\documentclass[11pt,a4paper,twoside]{article}
\usepackage{amsmath,amsthm,amssymb,mathtools}
\usepackage{newtxtext}
\usepackage[mathscr]{euscript}
\usepackage{microtype}
\usepackage[a4paper,textwidth=150mm,textheight=225mm,hmarginratio=1:1,headheight=14pt,headsep=7mm]{geometry}
\usepackage{xcolor,titlesec,titling,tocloft,fancyhdr,needspace,aliascnt}
\usepackage{tikz,float,booktabs,array}
\usepackage{enumitem}
\setlist[enumerate,1]{label=(\alph*)}
\setlist[enumerate,2]{label=(\roman*)}
\usepackage{pgfplots}\pgfplotsset{compat=1.16}
\definecolor{PaperInk}{HTML}{213747}
\definecolor{PaperAccent}{HTML}{2B6570}
\definecolor{PaperMuted}{HTML}{667882}
\definecolor{PaperRule}{HTML}{CBD7DC}
\usepackage[colorlinks=true,linkcolor=black!75!blue,citecolor=black!75!blue,urlcolor=black!75!blue,bookmarksnumbered=true,bookmarksopen=true,bookmarksopenlevel=1]{hyperref}
\hypersetup{pdftitle={Stokes Corners and the Global Structure of Steady Rotating Waves in Annular Domains},pdfsubject={Stagnation, Stokes corners, Weiss monotonicity and extreme waves in annular domains}}
\allowdisplaybreaks[2]
\numberwithin{equation}{section}
\titleformat{\section}{\normalfont\large\bfseries}{\thesection.}{.6em}{}
\titleformat{\subsection}{\normalfont\normalsize\bfseries}{\thesubsection.}{.5em}{}
\titleformat{\subsubsection}[runin]{\normalfont\normalsize\itshape}{\thesubsubsection.}{.4em}{}[.]
\titlespacing*{\section}{0pt}{1.6\baselineskip}{.7\baselineskip}
\titlespacing*{\subsection}{0pt}{1.1\baselineskip}{.5\baselineskip}
\titlespacing*{\subsubsection}{0pt}{.8\baselineskip}{.6em}
\AddToHook{cmd/subsection/before}{\Needspace{7\baselineskip}}
\theoremstyle{plain}
\newtheorem{theorem}{Theorem}[section]
\newaliascnt{lemma}{theorem}
\newtheorem{lemma}[lemma]{Lemma}
\aliascntresetthe{lemma}
\newaliascnt{proposition}{theorem}
\newtheorem{proposition}[proposition]{Proposition}
\aliascntresetthe{proposition}
\newaliascnt{corollary}{theorem}
\newtheorem{corollary}[corollary]{Corollary}
\aliascntresetthe{corollary}
\newaliascnt{conjecture}{theorem}
\newtheorem{conjecture}[conjecture]{Conjecture}
\aliascntresetthe{conjecture}

\AtBeginDocument{}
\theoremstyle{definition}
\newaliascnt{definition}{theorem}
\newtheorem{definition}[definition]{Definition}
\aliascntresetthe{definition}
\theoremstyle{remark}
\newaliascnt{remark}{theorem}
\newtheorem{remark}[remark]{Remark}
\aliascntresetthe{remark}

\AddToHook{env/theorem/before}{\Needspace{5\baselineskip}}
\AddToHook{env/lemma/before}{\Needspace{4\baselineskip}}
\AddToHook{env/proposition/before}{\Needspace{5\baselineskip}}
\AddToHook{env/definition/before}{\Needspace{4\baselineskip}}
\fancypagestyle{plain}{\fancyhf{}\fancyfoot[C]{\footnotesize\thepage}}

\renewenvironment{abstract}{\par\begin{quote}\small\noindent\textsc{Abstract.}\ \ignorespaces}{\end{quote}\par}
\newcommand{\cst}[1]{\ensuremath{\csname const@#1\endcsname}}
\newcommand{\registerconstant}[2]{\expandafter\gdef\csname const@#1\endcsname{#2}}
\AddToHook{cmd/appendix/before}{\clearpage}
\AddToHook{env/thebibliography/before}{%
  \phantomsection\addcontentsline{toc}{section}{References}}
\AddToHook{env/thebibliography/begin}{\small}

\makeatletter
\DeclareRobustCommand{\fint}{\mathpalette\paper@meanbar\relax\!\int}
\newcommand{\paper@meanbar}[2]{%
  \setbox\z@=\hbox{$\m@th#1\int$}%
  \vcenter{\hbox to\wd\z@{\hss$\m@th#1-$\hss}}%
  \kern-\wd\z@}
\makeatother
\NewDocumentCommand{\qty}{d() d[]}{%
\IfNoValueTF{#1}{%
\IfNoValueTF{#2}{\PackageError{qty}{Missing delimiter}{}}{\left[#2\right]}}%
{\left(#1\right)}}
\begin{document}
\thispagestyle{plain}
\begin{center}
{\Large\bfseries Stokes corners and the global structure\\[2pt] of steady rotating waves in annular domains\par}
\vspace{1.4em}
{\scshape Quan Wang}\par
\vspace{.6em}
{\small School of Mathematics, Sichuan University, Chengdu 610064, China\\
\texttt{xihujunzi@scu.edu.cn}\par}
\end{center}
\vspace{.6em}

\begin{abstract}
We study steady rotating water waves in a planar annulus, in the model of radial gravity and constant vorticity introduced by Li and the author. The stream function is harmonic, and the free surface carries the Bernoulli condition $p_0^2|\nabla\Psi|^2=Q_\alpha(S)$ with $Q_\alpha(R)=(R-r_-)(R-r_+)$.
We first settle the crest shape. For $\alpha\neq\alpha_\pm$, every stagnation point lies on the circle $R=r_-$, and the surface forms a symmetric $120^\circ$ corner there. At $\alpha=\alpha_\pm$ the corner is a right angle. The same holds for limits of smooth solutions and, in the lower region, for overhanging waves.
We then study the global bifurcation continuum. Its waves are monotone from crest to trough, and they always cross the laminar circle. Hence there is no secondary bifurcation, and the parameter stays above an explicit floor. The continuum ends in an extreme wave, overturns, or degenerates; the last case is excluded for $c_0>\gamma^2e^{2\gamma}/(1+\gamma)^2$. An annular Nekrasov equation and an invariant Keady--Norbury cone give a slope bound below $90^\circ$ whenever the Nekrasov weight is nonnegative on the waves. A sharpened angular flow-force estimate keeps the troughs where this weight is nonnegative, for every $\gamma\in(0,1)$ above the floor. Hence, whenever $c_0>\gamma^2e^{2\gamma}/(1+\gamma)^2$, the existence half of the Stokes conjecture holds: the continuum contains an extreme wave with a $120^\circ$ Stokes crest. In the remaining band $\sinh^2\gamma<c_0\le\gamma^2e^{2\gamma}/(1+\gamma)^2$ the result is conditional: the continuum ends in such a wave unless its parameter degenerates to the value at which the effective gravity vanishes on the stagnation circle, and we cannot exclude this degeneration. We also show that the bifurcation is supercritical for every $\gamma\in(0,1)$.
\end{abstract}
\begin{quote}\small
\end{quote}
{\small\setlength{\parskip}{0pt}\tableofcontents}
\clearpage

\section{Introduction}\label{sec_intro}

\subsection{The Stokes conjecture}
In 1880 Stokes argued that the steady periodic gravity wave of greatest height has a stagnation point at each crest, and that there the free surface forms a symmetric corner with an included angle of $120^\circ$ \cite{Stokes1880}. His argument was local. Near a stagnation point the Bernoulli condition $|\nabla\psi|^2=-2g\,y$ forces the stream function to be homogeneous of degree $3/2$. The only homogeneous harmonic function of that degree that vanishes on two rays and is positive between them lives in a sector of opening $2\pi/3$.

Turning this heuristic into a theorem took a century, and the problem has two parts. The first is the \emph{existence} of an extreme wave, as a limit of regular waves along a global branch. The second is the \emph{shape} of the crest of any such wave.

For flat irrotational waves both parts are now classical. Nekrasov reformulated the problem as a nonlinear integral equation for the angle $\vartheta$ between the free surface and the horizontal \cite{Nekrasov1921}. In his formulation $\vartheta$ is the image, under a positive integral operator, of an explicit nonlinear expression in $\vartheta$. Krasovskii used this structure to obtain large-amplitude waves with every maximal slope strictly between $0$ and $\pi/6$ \cite{Krasovskii1961}. Keady and Norbury worked in a cone of functions on which Nekrasov's operator is order preserving. Using degree theory there, they obtained a global branch of waves for every crest speed between the critical value and zero \cite{KeadyNorbury1978}. Toland then proved that a wave of greatest height exists, as a limit along this branch that solves Nekrasov's equation \cite{Toland1978}. Amick, Fraenkel and Toland \cite{AFT1982}, and independently Plotnikov \cite{Plotnikov1982}, proved that its crest is a corner of exactly $120^\circ$. The first proof analyses Nekrasov's equation near the crest, and the second uses a different reduction of the free-boundary problem. The extreme wave was studied further by McLeod \cite{McLeod1997} and by Fraenkel and Harwin \cite{FraenkelHarwin2010}. Toland's survey \cite{Toland1996} gives a self-contained account, including a maximum-principle proof that the Keady--Norbury cone is invariant; we transplant this proof to the annulus. The surveys \cite{Strauss2010,HaziotEtAl2022} and the monograph \cite{Constantin2011} describe the wider theory.

Near-extreme waves are subtle. Longuet-Higgins and Fox described the inner region of the almost-highest wave \cite{LonguetHigginsFox1977}. The maximal slope of a near-extreme irrotational wave may exceed $30^\circ$, but Amick proved that it always stays below $31.15^\circ$ \cite{Amick1987}. High-precision computations locate the complex singularities that approach the crest as the wave steepens \cite{DLK2016,Lushnikov2016}.

Varv\u{a}ruc\u{a} and Weiss gave a geometric proof of the crest-shape statement in the original physical variables \cite{VW2011}. It uses neither Nekrasov's equation nor any structural assumption on the wave. Its main tools are a monotonicity formula of Weiss type \cite{Weiss1999} for the Alt--Caffarelli functional \cite{AltCaffarelli1981} (see also \cite{CaffarelliSalsa2005}), and a non-linear frequency formula in the spirit of Garofalo and Lin \cite{GarofaloLin1987}. The blow-up limits at a stagnation point are classified, and the Stokes corner is singled out as the only possible limit under a growth assumption on the velocity. Singularities of general Bernoulli free boundaries were analysed earlier in \cite{Varvaruca2006}. For waves with vorticity the stagnation set is richer. At a stagnation point about which the profile is symmetric, the profile has a corner of $120^\circ$ or a horizontal tangent \cite{Varvaruca2009}. When the free surface is an injective curve, the possible asymptotics are the Stokes corner, a horizontal cusp, or a horizontally flat point \cite{VW2012}. The method also extends to axisymmetric flows \cite{VW2014}. For flat waves, the growth hypothesis of \cite{VW2011} has recently been removed for variational solutions \cite{DYZ2026}. That proof uses the Hopf-differential method of \cite{DHN2026}. Cusps have also been excluded under that hypothesis \cite{McCurdy2025,DPY2026}.

\subsection{Global bifurcation and its alternatives}
Existence of large waves is obtained by global bifurcation. This goes back to Rabinowitz's global theorem \cite{Rabinowitz1971} and the local theorem of Crandall and Rabinowitz \cite{CrandallRabinowitz1971}. The analytic version of the theory, in which the continuum is a locally analytic curve, is due to Dancer \cite{Dancer1973} and Buffoni and Toland \cite{BuffoniToland2003}; see also \cite{Kielhoefer2012}. For waves with vorticity, the modern framework was set up by Constantin and Strauss \cite{ConstantinStrauss2004}. They used the hodograph variables of Dubreil-Jacotin \cite{DubreilJacotin1934}, in which the free boundary becomes fixed. Along their global continuum the waves approach stagnation, unless other alternatives are excluded by a priori bounds. The regularity of the resulting waves and of their streamlines was established in \cite{ConstantinEscher2011,ConstantinVarvaruca2011,BDT2000}.

What happens at the end of a continuum depends on such a priori bounds. With critical layers the continuum can end in overhanging profiles instead of at an extreme wave; this was established for constant vorticity in \cite{CSV2016}, building on \cite{Wahlen2009}. Overhanging and touching waves are studied in \cite{HurWheeler2022}. Slope bounds are therefore a central issue. Strauss and Wheeler proved a bound of $45^\circ$ for a large class of waves with favourable vorticity \cite{StraussWheeler2016}, and So and Strauss extended it to small adverse vorticity \cite{SoStrauss2018}. For solitary waves, the bound $F>1$ on the Froude number \cite{AmickToland1981,Wheeler2015} and related qualitative results \cite{CWW2018} play the role that a parameter floor plays here. Integral invariants are another recurring tool. The flow force of Benjamin and Lighthill \cite{BenjaminLighthill1954}, which is conserved because of translation invariance, constrains the possible conjugate states of a flow; its rotational analogue is central in this paper.

\subsection{Curved geometry: the annular model}
Steady waves in curved geometries arise for flows on rotating planets, in planetary rings and in equatorial ocean currents. Henry and Martin constructed exact free-surface equatorial flows in spherical coordinates \cite{HenryMartin2019}. Moon and Wu constructed global curves of rotating capillary waves on a two-dimensional droplet \cite{MoonWu2024}. In that problem gravity is absent and surface tension drives the waves. Li and the author \cite{LiWang2026} introduced a model of gravity waves in a curved annulus. A fluid with constant vorticity occupies the region between a rigid inner circle and a free outer surface, under radially inward gravity. In \cite{LiWang2026} we proved local and global bifurcation of travelling waves, using a hodograph formulation and a modified Leray--Schauder degree in the spirit of Healey and Simpson \cite{HealeySimpson1998}. In each open set of the degree argument, the continuum is unbounded, returns to another laminar flow, or reaches the boundary of the set. The purpose of this paper is to determine what the waves at the end of the continuum look like.

\subsection{Main results}
The annular model differs from the flat one in three ways that matter for the Stokes problem. The effective Bernoulli function $Q_\alpha(R)=R^2-2\alpha R+2\lambda$ is quadratic, not linear, in the distance to the centre. The bottom is a circle, so the problem is not invariant under translations. Finally, rotation invariance replaces horizontal translation invariance. We show that each of these differences can be turned to advantage.

Our first result concerns the shape of the crest, and it holds at all parameters. Away from the two degenerate values $\alpha_\pm$, stagnation can occur only on one explicit circle $R=r_-(\alpha)$ and only at global maxima of the free surface, and every stagnation point is a symmetric corner of $120^\circ$ with the explicit amplitude $2\Delta^{1/4}/(3|p_0|)$. At $\alpha_\pm$ every stagnation point is a symmetric right-angle corner. The same conclusions hold for limits of regular waves and, in the lower region, for overhanging waves. Horizontally flat points and cusps are excluded by the radial-graph structure, not by assumption. For overhanging waves we adapt the Hopf-differential and isoperimetric methods of \cite{DHN2026,DYZ2026,DPY2026} to the radial weight $Q_\alpha(|x|)$.

Our second result concerns the global continuum of \cite{LiWang2026}. It consists of waves that are strictly monotone from crest to trough. Every wave crosses the laminar circle, so there is no secondary bifurcation, and the parameter stays between an explicit floor $\alpha_\circ$ and an explicit upper bound. The second-order coefficient of the parameter at the bifurcation point is computed in closed form and is positive for every $\gamma\in(0,1)$, so the bifurcation is supercritical.

Our main result is the existence half of the Stokes conjecture. Let $0<\gamma<1$ and $\gamma^2e^{2\gamma}/(1+\gamma)^2<c_0<e^{2\gamma}$. Then the continuum contains a symmetric extreme wave, monotone from crest to trough, whose single stagnation point per period is a $120^\circ$ corner with the Stokes asymptotics. For every $c_0>\sinh^2\gamma$, the continuum ends in such a wave unless its parameter degenerates to the value $\alpha_+$, at which the effective gravity vanishes on the stagnation circle. The proof rests on two new tools. The first is an annular Nekrasov equation for the flow angle. Its weight obeys $VF'(V)\le\kappa_\alpha<\frac12$ at every parameter with $\Delta>0$ and is monotone without any condition, so the annular Keady--Norbury cone is invariant along the continuum as long as the weight is nonnegative on the waves. The second is a lower bound for the trough, obtained from a sharpened angular flow force, which keeps the waves in the region where the weight is nonnegative. The final step reduces to the positivity of an explicit polynomial in two variables, which we prove by hand; no interval arithmetic or other computer-assisted step is used. As an independent check, the polynomial identities and inequalities of this step are also verified in Lean~4. This formalization covers only the algebra, not the analytic arguments.

We also treat the window $\Delta<0$ and report numerical continuation for $23$ parameter pairs, which agrees with these results. What remains open is the thin band $\sinh^2\gamma<c_0\le\gamma^2e^{2\gamma}/(1+\gamma)^2$, where we cannot exclude that the parameter degenerates to $\alpha_+$. Thus the paper settles the existence half of the conjecture in an explicit parameter region and reduces the rest of the range $c_0>\sinh^2\gamma$ to one a priori bound; it does not settle the conjecture for all parameters. \autoref{sec_open} separates the statements that are proved unconditionally, those proved under explicit hypotheses, and those supported only by numerics. Precise statements are given in \autoref{sec_problem}, and \autoref{fig_dependencies} shows how they depend on each other.

\subsection{Main ideas}
Six ideas are new compared with the flat theory.
\begin{itemize}
\item \emph{A weighted maximum principle.} The functions $R^k|\nabla\Psi|^2$ are subharmonic for every $k$, because they are squared moduli of holomorphic functions in the cylinder variables (\autoref{lem_subharmonic}). Compared with $R^kQ_\alpha(R)$, this gives a local Bernstein bound, and hence the growth $\Psi=O(r^{3/2})$ at stagnation points, without Nekrasov's equation. For $k=2$ it gives an exact maximum-velocity principle; compare Spielvogel's inequality for flat waves \cite{Spielvogel1970}.
\item \emph{Radial phase monotonicity.} The fluid region is a radial subgraph. After blow-up this becomes a vertical-line property, which excludes the cusp and flat alternatives of \cite{VW2012} and replaces the frequency formula of \cite{VW2011}. The non-linearity of $Q_\alpha$ enters only as an error term of order $r$ in a perturbed Weiss formula.
\item \emph{A rotational Jacobi field.} The angular derivative $\Psi_\theta$ satisfies a homogeneous Robin condition on the free surface. In the hodograph variables this becomes an oblique condition whose tangential coefficient vanishes on the symmetry axes. Serrin's corner lemma, in the form of Gidas, Ni and Nirenberg \cite{Serrin1971,GNN1979}, then shows that monotonicity is open and closed along the continuum. This replaces the nodal analysis of \cite{ConstantinStrauss2004} and the unilateral bifurcation theorem.
\item \emph{A crossing invariant.} The strict crest--trough comparison, together with the fact that the laminar radius is a root of the circular relation for every $\alpha$, forces every wave to straddle the laminar circle. This excludes secondary bifurcation and yields the parameter floor by a log-concavity argument.
\item \emph{An algebraic structure of the Nekrasov weight.} Written with the distances $a=r_--\sigma$ and $b=r_+-\sigma$ to the two stagnation radii, the quantities that control the cone argument become quadratic forms in $\sigma$ with negative discriminant, or sums of nonnegative terms. This gives the constant $\kappa_\alpha<\frac12$ for every parameter and the monotonicity of the weight without conditions (\autoref{lem_structural}, \autoref{lem_weight_monotone}).
\item \emph{The angular flow force.} Rotation invariance yields a conserved flux between the crest and trough sections, the annular analogue of the flow force of \cite{BenjaminLighthill1954}. Estimated from above on the trough ray by the maximum-velocity principle and from below on the crest ray by the Cauchy--Schwarz inequality, it bounds the parameter from above and keeps the trough away from the region where the Nekrasov weight is negative (\autoref{thm_alpha_bound}, \autoref{lem_trough_bound}).
\end{itemize}

\subsection{Relation to earlier work of Li and the author}
The present paper is a sequel to \cite{LiWang2026}, and the boundary between the two is as follows. From \cite{LiWang2026} we use the model \eqref{eq_model}, the parameter path \eqref{eq_parameter_path}, the hodograph formulation in the space $X=C^{3+s}_{\rm per}$, the open sets $\mathcal O_\delta$ defined by the two hodograph constraints, the local bifurcation at $\alpha_c$, and the global alternatives of \cite[Theorems~2.1 and~2.2]{LiWang2026}: in each $\mathcal O_\delta$ the continuum is unbounded, returns to another laminar flow, or reaches $\partial\mathcal O_\delta$. Nothing else is taken from \cite{LiWang2026}. In particular, \cite{LiWang2026} contains no analysis of stagnation points or of the crest, no monotonicity or nodal property of the waves, no bound on the parameter along the continuum, no slope bound, no limit $\delta\downarrow0$, and no extreme wave. \autoref{tab_relation} summarizes the division.
\begin{table}[ht]
\centering\small
\begin{tabular}{@{}p{0.44\linewidth}p{0.50\linewidth}@{}}
\toprule
Taken from \cite{LiWang2026} & New in this paper\\
\midrule
Model \eqref{eq_model}, path \eqref{eq_parameter_path}, hodograph formulation, space $X$, sets $\mathcal O_\delta$ & Localization of stagnation and the annular Stokes corner, including degenerate parameters and overhanging waves (Sections~\ref{sec_prelim}, \ref{sec_corner})\\
Local bifurcation at $\alpha_c$ & Direction of bifurcation, $\alpha_2>0$ (\autoref{sec_supercritical})\\
Global alternatives in each $\mathcal O_\delta$ & Monotone half-continuum, no laminar return, crossing invariant, parameter bounds (Sections~\ref{sec_continuum}, \ref{sec_final})\\
--- & Limit $\delta\downarrow0$, Nekrasov equation, invariant cone, slope bound, extreme wave (Sections~\ref{sec_corner}, \ref{sec_final}, \ref{sec_towards})\\
\bottomrule
\end{tabular}
\caption{Division of results between \cite{LiWang2026} and the present paper.}
\label{tab_relation}
\end{table}

\subsection{Organization}
\autoref{sec_problem} states the model, the standing assumptions and the main theorems. \autoref{sec_prelim} collects comparison principles, localizes stagnation points and proves the gradient bounds. \autoref{sec_corner} proves the Stokes corner theorem, its degenerate and overhanging versions, and the limit theorem. \autoref{sec_continuum} shows that the global continuum consists of monotone waves and bounds the parameter. \autoref{sec_final} proves the crossing invariant and the parameter floor, treats regularity, and proves the global theorem. \autoref{sec_towards} derives the annular Nekrasov equation, proves the slope bound in an invariant cone and the existence of the extreme Stokes wave under explicit conditions, and computes the direction of bifurcation. \autoref{sec_numerics} reports the numerical continuation, and \autoref{sec_open} lists the open problems.

\section{The model and the main results}\label{sec_problem}

\subsection{The annular model}\label{subsec_model}

\begin{figure}[!ht]
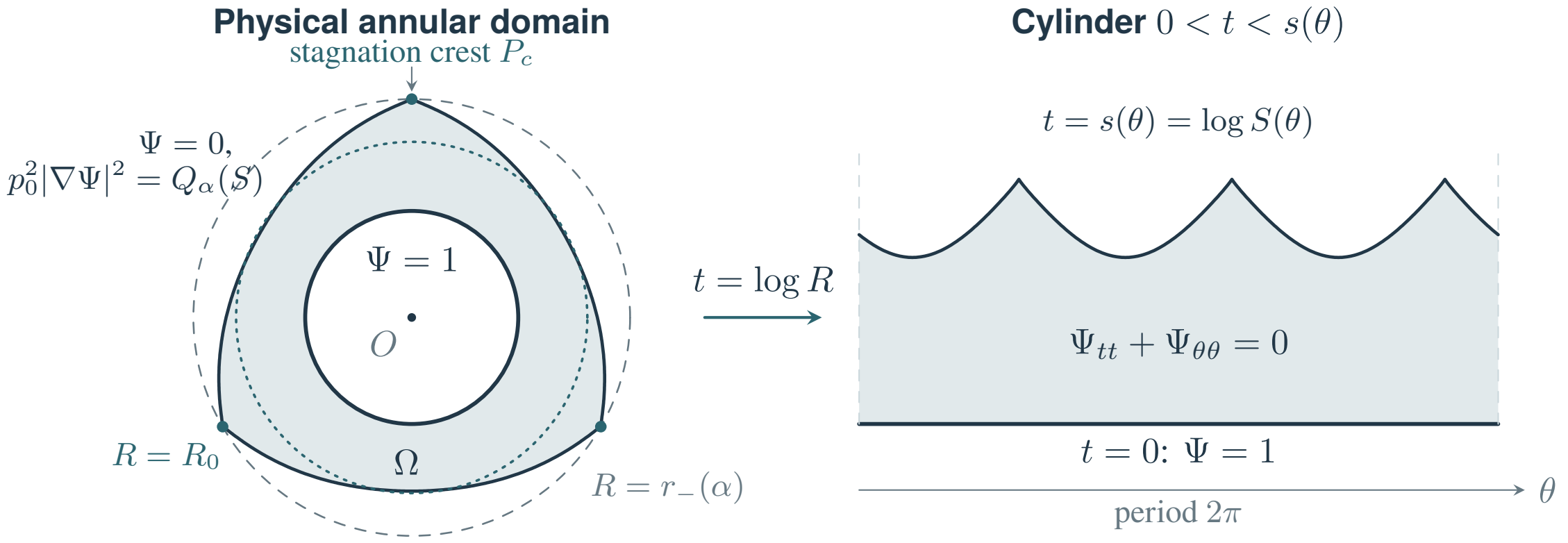

\centering
\resizebox{\linewidth}{!}{%
}
\caption{Schematic geometry (not to scale). The fluid occupies $1<R<S(\theta)$. An extreme wave touches the circle $R=r_-(\alpha)$ (dashed) exactly at its stagnation crests, where the surface has $120^\circ$ corners (\autoref{thm_main_stokes}). The dotted circle is the laminar radius $R_0$. Right: the conformal cylinder $t=\log R$, in which the inner boundary becomes flat.}
\label{fig_geometry}
\end{figure}

Let $(R,\theta)$ be polar coordinates in the plane. A steady wave is described by a $2\pi$-periodic free surface $R=S(\theta)>1$ bounding the fluid region
\[
\Omega=\{(R,\theta):1<R<S(\theta)\},
\]
together with a stream function $\Psi$ of the relative velocity. Following \cite{LiWang2026},
\begin{equation}\label{eq_model}
\left\{
\begin{aligned}
&\Delta\Psi=\Psi_{RR}+R^{-1}\Psi_R+R^{-2}\Psi_{\theta\theta}=0 &&\text{in }\Omega,\\
&\Psi=1 &&\text{on }R=1,\\
&\Psi=0,\qquad p_0^2|\nabla\Psi|^2=Q_\alpha(S) &&\text{on }R=S(\theta),
\end{aligned}\right.
\qquad Q_\alpha(R)=R^2-2\alpha R+2\lambda .
\end{equation}
Here $p_0<0$ is the relative mass flux, and the relative velocities are $V-R=-p_0\Psi_R$ and $U=p_0R^{-1}\Psi_\theta$. A \emph{stagnation point} is a point $P_c=(R_c,\theta_c)$ of the free surface at which $\nabla\Psi(P_c)=0$. \autoref{fig_geometry} illustrates the geometry and its cylindrical representation.

Write $R_0=e^\gamma$ and set $c_0=p_0^2/(R_0^2\gamma^2)$. As in \cite{LiWang2026}, we choose $\lambda$ along the path
\begin{equation}\label{eq_parameter_path}
2\lambda=c_0-R_0^2+2\alpha R_0 .
\end{equation}
On this path the circle $S\equiv R_0$ with $\Psi=1-\log R/\gamma$ solves \eqref{eq_model} for every $\alpha$, and a branch of non-circular solutions bifurcates at
\[
\alpha_c=R_0+\frac{2c_0R_0}{R_0^2-1}.
\]
Completing the square,
\begin{equation}\label{eq_completed_square}
Q_\alpha(R)=(R-\alpha)^2-\Delta(\alpha),\qquad \Delta(\alpha)=(\alpha-R_0)^2-c_0 .
\end{equation}
\emph{Standing assumptions.} The global bifurcation theorem of \cite[Theorems~2.1 and~2.2]{LiWang2026} is stated for
\begin{equation}\label{eq_standing}
0<\gamma<1,\qquad p_0^2<\gamma^2e^{4\gamma}\ \ \text{(equivalently } c_0<R_0^2\text{)},
\end{equation}
on the parameter interval $\alpha>\alpha_s:=\frac12\big(R_0-c_0/R_0\big)$. We assume \eqref{eq_standing} whenever we use the continuum of \cite{LiWang2026}, that is, in Sections~\ref{sec_continuum} and~\ref{sec_final}. Sections~\ref{sec_prelim} and~\ref{sec_corner} concern arbitrary classical solutions and use neither condition. The local computation of \autoref{sec_supercritical} is valid for every $\gamma>0$ (\autoref{rem_all_gamma}).

When $\Delta>0$ we write $r_\pm(\alpha)=\alpha\pm\sqrt{\Delta}$, so that $Q_\alpha=(R-r_-)(R-r_+)$, $2\alpha=r_-+r_+$ and $2\lambda=r_-r_+$. The two \emph{degenerate parameters} are $\alpha_\pm=R_0\pm\sqrt{c_0}$, where $\Delta=0$.

\begin{definition}[Classical solutions]\label{def_classical}
A pair $(\Psi,S)$ is a \emph{classical solution} if it satisfies the following conditions.
\begin{enumerate}
\item[(S1)] $S\in C(\mathbb T)$ and $S>1$.
\item[(S2)] $\Psi\in C^2(\Omega)\cap C^1(\overline\Omega)$ solves the first two lines of \eqref{eq_model}.
\item[(S3)] The Bernoulli condition holds at every point of the free surface.
\item[(S4)] $S$ is $C^1$ away from finitely many points.
\end{enumerate}
It belongs to the \emph{class} (L) if, in addition, $S\in C^1(\mathbb T)$ and $\Psi_R<0$ in $\overline\Omega$. The latter condition is equivalent to the hodograph condition $h_p>0$ of \cite{LiWang2026}.
\end{definition}

By the maximum principle, $0<\Psi<1$ in $\Omega$. In cylindrical variables $t=\log R$ the bottom is the line $t=0$. Odd reflection $\Psi(-t,\theta)=2-\Psi(t,\theta)$ extends $\Psi$ harmonically across it, and the Hopf lemma gives $\Psi_t(0,\theta)<0$.

\subsection{Principal conclusions}\label{subsec_results}

The classical Stokes conjecture asserts that the steady gravity wave of greatest height has a stagnation point at its crest, where the free surface forms a symmetric corner of $120^\circ$. For flat irrotational waves this was proved by Amick, Fraenkel and Toland \cite{AFT1982} and by Plotnikov \cite{Plotnikov1982}. Varvaruca and Weiss \cite{VW2011,VW2012} gave a geometric proof in the original variables. In the annular model the effective gravity is not constant and the bottom is curved, so none of these proofs applies directly. Our results are the following.

\begin{theorem}[Localization of stagnation]\label{thm_intro_localization}
Let $(\Psi,S)$ be a classical solution with a stagnation point $P_c$. Then $\Delta(\alpha)\ge0$. If $\Delta(\alpha)>0$, then
\[
R_c=r_-(\alpha)<\alpha,\qquad S(\theta)\le r_-(\alpha)\ \text{ for all }\theta,
\]
and the stagnation points are exactly the points where $S$ attains its maximum.
\end{theorem}

\noindent\emph{Source.} \autoref{thm_intro_localization} is \autoref{thm_localization}, proved in \autoref{sec_prelim}. The proof combines the two-region dichotomy of \autoref{lem_dichotomy} with the Hopf lemma.

\begin{theorem}[Annular Stokes corner theorem]\label{thm_intro_stokes}
Let $(\Psi,S)$ be a classical solution with a stagnation point $P_c$, and assume $\alpha\ne\alpha_\pm$. Then:
\begin{enumerate}
\item $P_c$ lies on the circle $R=r_-(\alpha)$ and is a global maximum of $S$;
\item the fluid has at $P_c$ a tangent cone which is the sector of opening $120^\circ$, symmetric about $-\mathbf e_R(\theta_c)$, with tangent rays
\[
\tau_\pm=-\tfrac12\mathbf e_R(\theta_c)\pm\tfrac{\sqrt3}{2}\mathbf e_\theta(\theta_c);
\]
\item in the frame $e_1=\mathbf e_\theta(\theta_c)$, $e_2=-\mathbf e_R(\theta_c)$, as $r\to0$,
\[
r^{-3/2}\Psi(P_c+ry)\to A_*\varrho^{3/2}\cos\!\Big(\tfrac32\big(\varphi-\tfrac\pi2\big)\Big)\chi_{\{\pi/6<\varphi<5\pi/6\}},
\qquad A_*=\frac{2\Delta(\alpha)^{1/4}}{3|p_0|},
\]
locally uniformly, with convergence of gradients in $L^2_{\rm loc}$;
\item stagnation points are isolated, and hence finite in number.
\end{enumerate}
\end{theorem}

\noindent\emph{Source.} \autoref{thm_intro_stokes} is \autoref{thm_main_stokes}, proved in \autoref{sec_corner}. Part~(a) follows from \autoref{thm_localization}. Parts (b)--(d) follow from the perturbed Weiss monotonicity formula (\autoref{lem_weiss}), the radial phase monotonicity (\autoref{lem_phase}), the angle bounds (\autoref{thm_angle_bounds}), and the compactness and classification of blow-ups (\autoref{prop_blowup}, \autoref{prop_classification}).

\begin{theorem}[Limit theorem]\label{thm_intro_limit}
Let $(\alpha_n,\Psi_n,S_n)$ be solutions of class (L) with $\max S_n<r_-(\alpha_n)$. Suppose $\alpha_n\to\alpha_\infty$ with $\Delta(\alpha_\infty)>0$ and $r_-(\alpha_\infty)>1$, and suppose $r_-(\alpha_n)-\max S_n\to0$. Then a subsequence converges to a non-circular variational solution at $\alpha_\infty$, in the sense of \autoref{def_variational}. This limit has at least one stagnation point, and each of its stagnation points satisfies conclusions (b)--(d) of \autoref{thm_intro_stokes}.
\end{theorem}

\noindent\emph{Source.} \autoref{thm_intro_limit} is \autoref{thm_limit}, proved in \autoref{sec_corner}. It follows from the uniform estimates of \autoref{prop_uniform}, the blow-up analysis of \autoref{prop_blowup}, and \autoref{thm_main_stokes} applied to the limit. The crest theory for variational solutions without the graph property (\autoref{thm_variational_crest}) extends it to overhanging limits.

\begin{theorem}[Global alternatives]\label{thm_intro_global}
Assume \eqref{eq_standing} and $c_0>\sinh^2\gamma$. The global continuum splits, under the half-period shift, into two copies of a connected set of symmetric waves. These waves are strictly monotone from crest to trough, and each of them crosses the laminar circle: $\min S<R_0<\max S$. The continuum never returns to the laminar branch, never meets any other circular flow, and its parameter obeys an explicit upper bound $\alpha<A(\gamma,p_0)$. Along these continua, as $\delta\downarrow0$, at least one of the following occurs:
\begin{enumerate}
\item[(I)] an extreme wave with a single symmetric $120^\circ$ Stokes crest per period exists;
\item[(II)] the waves begin to overturn;
\item[(III)] the parameter degenerates to $\alpha_+$.
\end{enumerate}
If (III) does not occur, then (I) holds; thus overturning can occur only together with degeneration. If moreover $c_0>\gamma^2e^{2\gamma}/(1+\gamma)^2$, for instance if $c_0\ge e^{2\gamma}/4$, then (III) is impossible, and the extreme wave with a $120^\circ$ Stokes crest exists.
\end{theorem}

\noindent\emph{Source.} \autoref{thm_intro_global} combines \autoref{thm_global} with \autoref{thm_existence} and \autoref{cor_thin}. The monotonicity of the waves and the absence of laminar returns are \autoref{thm_oriented}, which rests on \autoref{prop_nodal_open_closed}. The crossing $\min S<R_0<\max S$, which excludes secondary bifurcation and other circular flows, is \autoref{thm_crossing}. The upper bound is \autoref{thm_alpha_bound}. The alternatives (I)--(III) follow from \autoref{thm_final_reduction} and \autoref{cor_dichotomy}. The exclusion of (III) under $c_0>\gamma^2e^{2\gamma}/(1+\gamma)^2$ is \autoref{cor_floor}, and the implication from the failure of (III) to (I) is \autoref{thm_band}. The exclusion of (II) rests on four ingredients: the annular Nekrasov equation (\autoref{thm_nekrasov}); the structural inequality and the slope bound in the Keady--Norbury cone (\autoref{lem_structural}, \autoref{thm_cone_slope}); the monotonicity of the weight and the invariance of the cone (\autoref{lem_weight_monotone}, \autoref{thm_cone_invariance}); and the lower bound for the trough (\autoref{lem_trough_bound}). The resulting condition (C3) is verified in \autoref{cor_thin}.

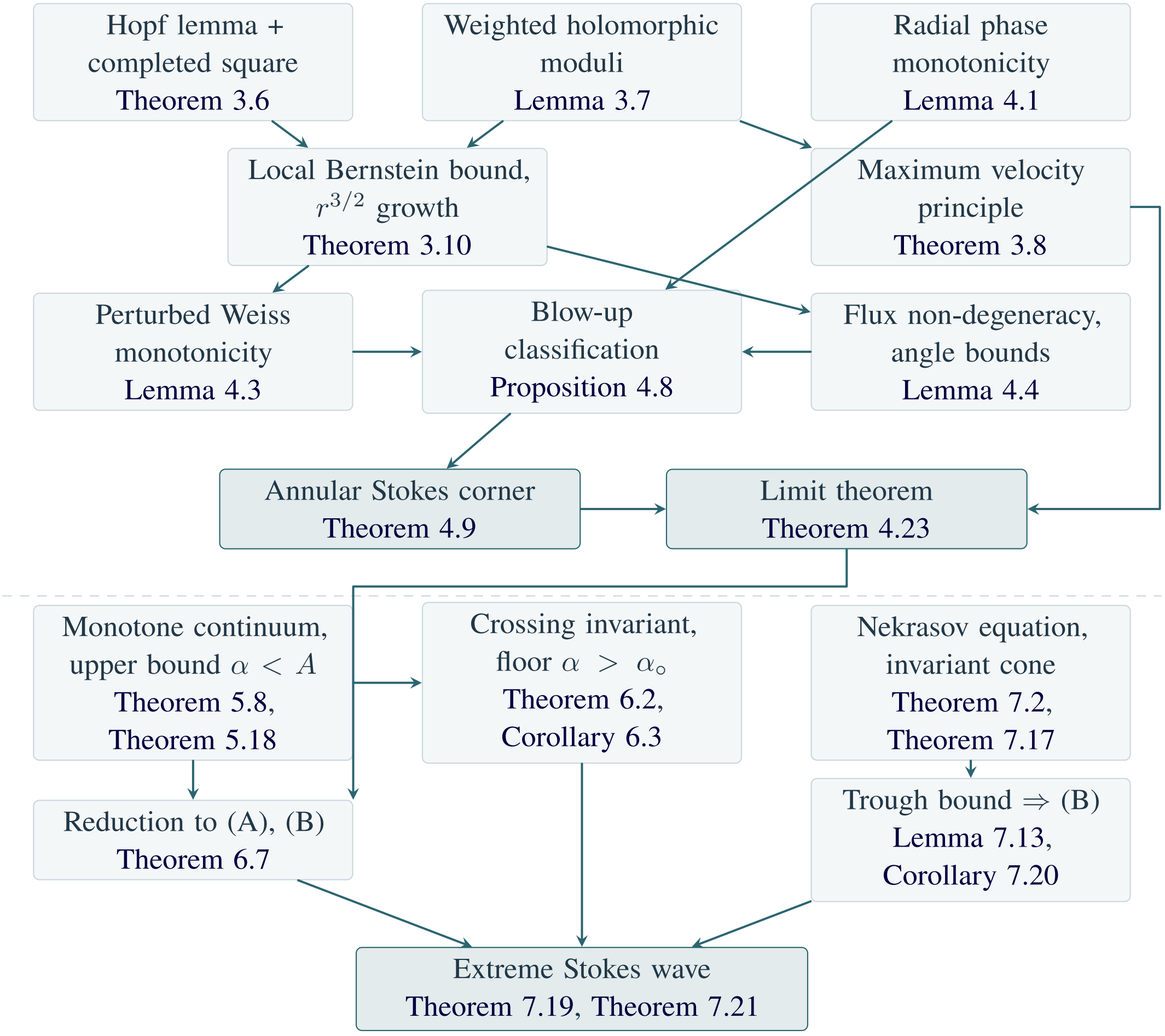
\begin{figure}[H]
\centering
\begin{tikzpicture}[
  >=stealth, font=\small, line cap=round, line join=round,
  paper step/.style={draw=PaperRule, fill=PaperAccent!5, text=PaperInk,
    rounded corners=2pt, align=center, text width=3.6cm, minimum height=8mm},
  result/.style={paper step, draw=PaperAccent, fill=PaperAccent!13},
  input/.style={->, draw=PaperAccent, line width=.75pt}]
\node[paper step] (hopf) at (-4.7,0) {Hopf lemma +\\completed square\\\autoref{thm_localization}};
\node[paper step] (sub) at (0,0) {Weighted holomorphic\\moduli\\\autoref{lem_subharmonic}};
\node[paper step] (graph) at (4.7,0) {Radial phase\\monotonicity\\\autoref{lem_phase}};
\node[paper step] (bern) at (-2.35,-1.75) {Local Bernstein bound,\\$r^{3/2}$ growth\\\autoref{thm_local_bernstein}};
\node[paper step] (vel) at (4.7,-1.75) {Maximum velocity\\principle\\\autoref{thm_max_velocity}};
\node[paper step] (weiss) at (-4.7,-3.5) {Perturbed Weiss\\monotonicity\\\autoref{lem_weiss}};
\node[paper step] (class) at (0,-3.5) {Blow-up\\classification\\\autoref{prop_classification}};
\node[paper step] (nondeg) at (4.7,-3.5) {Flux non-degeneracy,\\angle bounds\\\autoref{lem_nondeg}};
\node[result,text width=4.1cm] (stokes) at (-2.2,-5.4) {Annular Stokes corner\\\autoref{thm_main_stokes}};
\node[result,text width=4.1cm] (limit) at (3.2,-5.4) {Limit theorem\\\autoref{thm_limit}};
\draw[input] (hopf) -- (bern);
\draw[input] (sub) -- (bern);
\draw[input] (sub) -- (vel);
\draw[input] (bern) -- (weiss);
\draw[input] (bern) -- (nondeg);
\draw[input] (graph) -- (class);
\draw[input] (weiss) -- (class);
\draw[input] (nondeg) -- (class);
\draw[input] (class) -- (stokes);
\draw[input] (stokes) -- (limit);
\draw[input] (vel.east) -- ++(.35,0) |- (limit.east);
\draw[PaperRule,dashed] (-7,-6.45) -- (7,-6.45);
\node[paper step] (mono) at (-4.7,-7.5) {Monotone continuum,\\upper bound $\alpha<A$\\\autoref{thm_oriented}, \autoref{thm_alpha_bound}};
\node[paper step] (cross) at (0,-7.5) {Crossing invariant,\\floor $\alpha>\alpha_\circ$\\\autoref{thm_crossing}, \autoref{cor_floor}};
\node[paper step] (nek) at (4.7,-7.5) {Nekrasov equation,\\invariant cone\\\autoref{thm_nekrasov}, \autoref{thm_cone_invariance}};
\node[paper step] (red) at (-4.7,-9.4) {Reduction to (A), (B)\\\autoref{thm_final_reduction}};
\node[paper step] (trough) at (4.7,-9.4) {Trough bound $\Rightarrow$ (B)\\\autoref{lem_trough_bound}, \autoref{cor_thin}};
\node[result,text width=5.2cm] (exist) at (0,-11.2) {Extreme Stokes wave\\\autoref{thm_existence}, \autoref{thm_band}};
\draw[input] (mono) -- (red);
\draw[input] (mono) -- (cross);
\draw[input] (cross) -- (exist);
\draw[input] (nek) -- (trough);
\draw[input] (red) -- (exist);
\draw[input] (trough) -- (exist);
\draw[input] (limit.south) -- ++(0,-.45) -| (red.north east);
\end{tikzpicture}
\caption{Logical structure of the proofs. Above the dashed line: the crest shape. Below it: the existence of the extreme wave; the floor (A) is proved for $c_0>\gamma^2e^{2\gamma}/(1+\gamma)^2$, and (B) follows from (A).}
\label{fig_dependencies}
\end{figure}

\autoref{thm_intro_stokes} settles the crest-shape part of the Stokes conjecture for classical annular waves at non-degenerate parameters; at $\alpha=\alpha_\pm$ the corner is a right angle, with profile $\Psi\approx\frac{\sqrt2}{4|p_0|}(y_2^2-y_1^2)$ (\autoref{thm_degenerate}). The proof of \autoref{thm_intro_global} proceeds in three steps. First, the continuum consists of monotone waves and stays in a bounded parameter range (\autoref{sec_continuum}); this reduces the existence of the extreme wave to a parameter floor (A) and a slope bound (B) (\autoref{thm_final_reduction}). Second, the crossing invariant gives (A) under the floor condition (\autoref{cor_floor}). Third, the annular Nekrasov equation, the invariant cone and the lower bound for the trough give (B) whenever (A) holds (\autoref{thm_existence}, \autoref{thm_band}). Without the graph property, \autoref{cor_overhanging} shows that the crests of overhanging waves are still corners. \autoref{fig_dependencies} summarizes the logical structure.

\section{Comparison principles, stagnation points and gradient bounds}\label{sec_prelim}

This section collects the tools used throughout: comparison with circular flows and the flow-force identity, the localization of stagnation points, and the weighted gradient bounds that give the $3/2$ growth at stagnation.

\subsection{Comparison with circular flows}\label{sec_geometry}

\begin{lemma}[Two-region dichotomy]\label{lem_dichotomy}
Let $\Delta(\alpha)>0$ and let $(\Psi,S)$ be a classical solution. Then
\[
Q_\alpha(S(\theta))=p_0^2|\nabla\Psi|^2\ge0\quad\text{for all }\theta .
\]
Since $S(\mathbb T)$ is connected, either $\max S\le r_-$ or $\min S\ge r_+$.
\end{lemma}
\begin{proof}
$Q_\alpha<0$ exactly on $(r_-,r_+)$, and a connected set avoiding this interval lies entirely on one side of it.
\end{proof}

\begin{figure}[!ht]
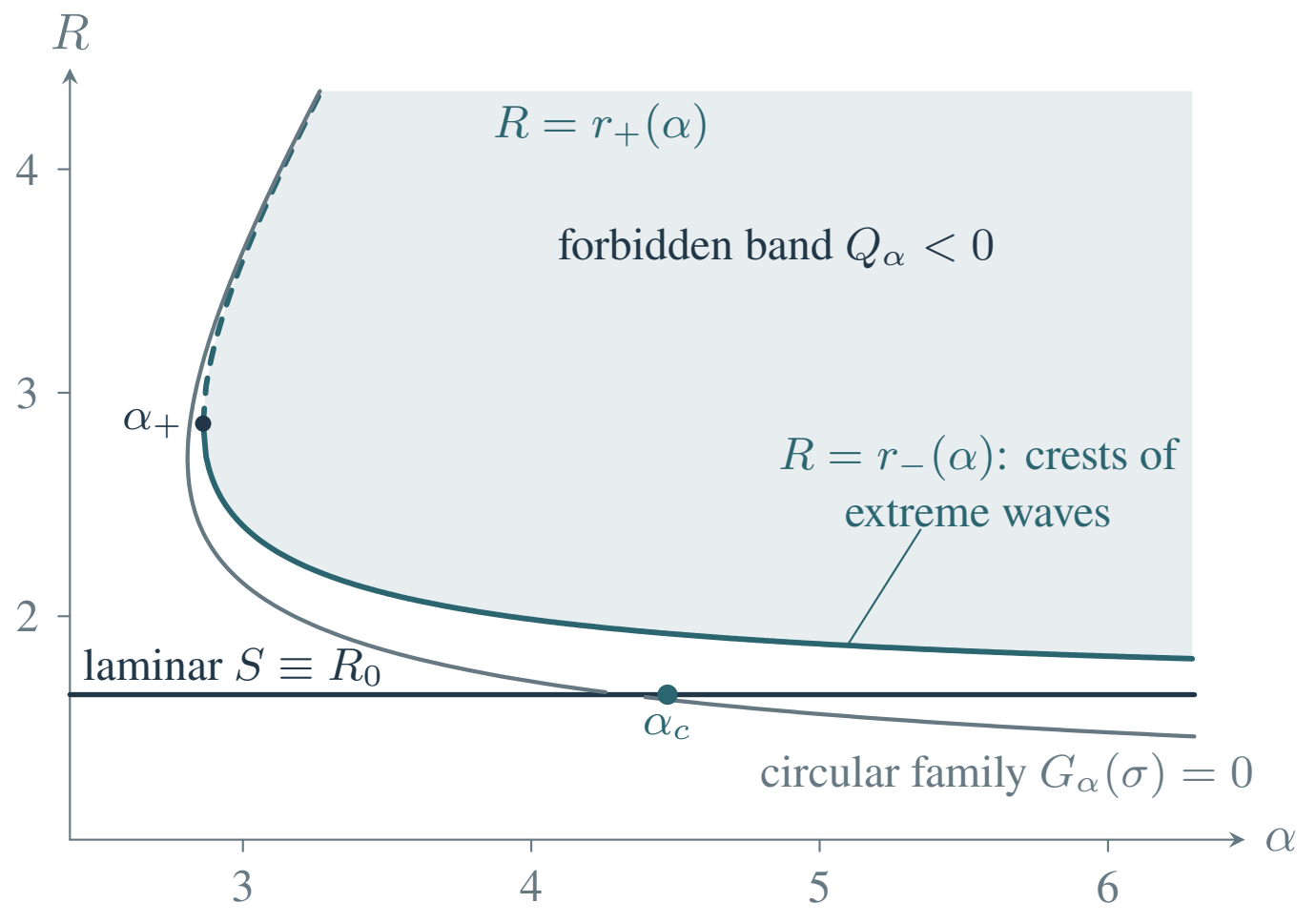

\centering

\caption{Parameter plane for $\gamma=1/2$, $p_0^2=1$; curves computed from the exact formulas. No free surface may enter the shaded band (\autoref{lem_dichotomy}). Stagnation can occur only on the solid curve $R=r_-(\alpha)$ (\autoref{thm_localization}); the lower region lies below it. The thin curve is the non-laminar circular family \eqref{eq_circular_relation}; the first non-circular branch bifurcates at $\alpha_c$.}
\label{fig_parameter_plane}
\end{figure}

We say that a solution lies in the \emph{lower region} if $\max S\le r_-(\alpha)$; this includes solutions with a stagnation crest. Circular solutions $S\equiv\sigma$, $\Psi=1-\log R/\log\sigma$ exist on the path \eqref{eq_parameter_path} exactly when
\begin{equation}\label{eq_circular_relation}
G_\alpha(\sigma):=\sigma^2\log^2\sigma\;Q_\alpha(\sigma)-p_0^2=0 .
\end{equation}
For $\sigma\ne R_0$ this relation can be solved explicitly for $\alpha$. Besides the laminar circle $\sigma=R_0$, it gives a second circular family, shown in \autoref{fig_parameter_plane}. The next proposition compares every solution with these circular flows.

\begin{proposition}[Crest--trough comparison]\label{prop_crest_trough}
Let $(\Psi,S)$ be a classical solution with $S\in C^1(\mathbb T)$. Then
\[
G_\alpha(\max S)\le0\le G_\alpha(\min S).
\]
If moreover $S\in C^2(\mathbb T)$ and the solution is not circular, both inequalities are strict.
\end{proposition}
\begin{proof}
Let $\sigma=\max S$ and $\Psi_\sigma=1-\log R/\log\sigma$. Then $\Omega$ is contained in the annulus $1<R<\sigma$, and $\Psi_\sigma\ge0=\Psi$ on the free surface, while $\Psi_\sigma=\Psi$ on $R=1$. By the maximum principle, $\Psi\le\Psi_\sigma$ in $\Omega$.

At a crest $P=(\sigma,\theta_*)$ both functions vanish. Moreover $S'(\theta_*)=0$, so $\nabla\Psi(P)$ is radial. Comparing difference quotients along the inward radial segment gives
\[
|\nabla\Psi(P)|=|\Psi_R(P)|\le\frac1{\sigma\log\sigma}.
\]
By the Bernoulli condition, $Q_\alpha(\sigma)\le p_0^2/(\sigma^2\log^2\sigma)$, which is $G_\alpha(\sigma)\le0$.

For the trough $\sigma'=\min S$, the annulus $1<R<\sigma'$ lies inside $\Omega$, and the reverse comparison $\Psi\ge1-\log R/\log\sigma'$ gives $G_\alpha(\sigma')\ge0$.

For strictness, suppose the solution is not circular. If $\Psi-\Psi_{\sigma'}$ vanished identically on the annulus $1<R<\sigma'$, it would vanish on $\Omega$ by analyticity. Then $\Psi_{\sigma'}=0$ on the free surface, so $S\equiv\sigma'$, which is excluded. Hence $\Psi-\Psi_{\sigma'}>0$ in the annulus. At a trough point the Hopf lemma, applied from inside the annulus, gives $|\Psi_R|>1/(\sigma'\log\sigma')$, and so $G_\alpha(\sigma')>0$. At a crest point, $\Psi_\sigma-\Psi>0$ in $\Omega$ for the same reason. Since $S\in C^2$, $\Omega$ satisfies an interior ball condition there, and the Hopf lemma gives $|\Psi_R|<1/(\sigma\log\sigma)$, so $G_\alpha(\sigma)<0$.
\end{proof}

\begin{proposition}[Flow-force identity]\label{prop_flow_force}
For every classical solution,
\begin{equation}\label{eq_flow_force}
\int_0^{2\pi}S(\theta)^2\,Q_\alpha\big(S(\theta)\big)\,d\theta=p_0^2\int_0^{2\pi}\Psi_R(1,\theta)^2\,d\theta .
\end{equation}
Consequently,
\begin{equation}\label{eq_flow_force_bounds}
\frac{2\pi p_0^2}{\log^2(\max S)}\le\int_0^{2\pi}S^2Q_\alpha(S)\,d\theta\le\frac{2\pi p_0^2}{\log^2(\min S)} .
\end{equation}
If moreover $\Delta(\alpha)>0$ and $\max S\le r_-(\alpha)$, then
\[
\frac1{2\pi}\int_0^{2\pi}\big(r_--S\big)\,d\theta\le\frac{p_0^2}{2\sqrt{\Delta}\,\log^2(\min S)} .
\]
\end{proposition}
\begin{proof}
In the cylinder variables, the field $\mathbf F=\Psi_t\nabla_c\Psi-\frac12|\nabla_c\Psi|^2e_t$ is divergence free for harmonic $\Psi$, because the problem is invariant under translations in $t$. Integrate $\operatorname{div}\mathbf F=0$ over the fluid, removing small balls around the finitely many non-smooth points.

On the free surface $\nabla_c\Psi=-|\nabla_c\Psi|\nu$, so $\mathbf F\cdot\nu=\frac12|\nabla_c\Psi|^2\,\nu_t$. Since the surface is the graph $t=s(\theta)$, we have $\nu_t\,ds=d\theta$. By the Bernoulli condition, $|\nabla_c\Psi|^2=S^2Q_\alpha(S)/p_0^2$ there. On the bottom $\nu=-e_t$ and $\Psi_\theta=0$, so $\mathbf F\cdot\nu=-\frac12\Psi_t^2$. This proves \eqref{eq_flow_force}.

For \eqref{eq_flow_force_bounds}, the comparisons of \autoref{prop_crest_trough} give
\[
1-\frac{\log R}{\log\min S}\le\Psi\le1-\frac{\log R}{\log\max S}
\]
near the bottom, with equality on $R=1$. Hence $1/\log\max S\le|\Psi_R(1,\theta)|\le1/\log\min S$.

For the last bound, use $S^2\ge1$ and $Q_\alpha(S)=(r_--S)(r_+-S)\ge2\sqrt\Delta\,(r_--S)$ on $[1,r_-]$.
\end{proof}

\begin{remark}
On circular solutions \eqref{eq_flow_force} reduces to the circular relation \eqref{eq_circular_relation}. It is an exact identity, not an estimate. On the computed near-extreme waves of \autoref{sec_numerics}, the two sides agree to a relative accuracy of $3\times10^{-7}$. The last inequality says that in the lower region, as $\alpha\to\infty$ (so $\sqrt\Delta\sim\alpha$), waves must flatten in $L^1$ onto the circle $R=r_-(\alpha)\to R_0$, unless their troughs approach the bottom.
\end{remark}

\begin{corollary}[Mean energy sandwich]\label{cor_energy_sandwich}
Let $(\Psi,S)$ be a classical solution with $S\in C^1(\mathbb T)$, and put $h_\alpha(\sigma)=\sigma^2Q_\alpha(\sigma)$, $S_c=\max S$, $S_t=\min S$. Then
\[
h_\alpha(S_c)\ \le\ \frac1{2\pi}\int_0^{2\pi}h_\alpha\big(S(\theta)\big)\,d\theta\ \le\ h_\alpha(S_t).
\]
In particular, if $h_\alpha$ is strictly increasing on $[S_t,S_c]$, the solution is circular. If $\Delta(\alpha)>0$, then $h_\alpha$ is either decreasing on $(1,r_-)$ or increasing and then decreasing, so for a non-circular solution in the lower region the crest lies beyond the maximum point of $h_\alpha$ on $(1,r_-)$, if there is one.
\end{corollary}
\begin{proof}
By \eqref{eq_flow_force_bounds}, the mean of $h_\alpha(S)$ lies between $p_0^2/\log^2S_c$ and $p_0^2/\log^2S_t$. By \autoref{prop_crest_trough}, $G_\alpha(S_c)\le0\le G_\alpha(S_t)$, i.e.\ $h_\alpha(S_c)\le p_0^2/\log^2S_c$ and $p_0^2/\log^2S_t\le h_\alpha(S_t)$. This proves the two inequalities. If $h_\alpha$ is strictly increasing on $[S_t,S_c]$ and $S_c>S_t$, then $h_\alpha(S_c)>h_\alpha(S_t)$, which contradicts them.

For the shape of $h_\alpha$, note that $h_\alpha'(\sigma)=2\sigma\big(2\sigma^2-3\alpha\sigma+2\lambda\big)$. The quadratic factor has at most two roots. At $\sigma=r_-$, using $2\lambda=r_-r_+$ and $r_+=2\alpha-r_-$, it equals $r_-(r_--\alpha)<0$. So on $(1,r_-)$ the factor is either negative throughout or positive and then negative. Finally, $h_\alpha(S_c)\le h_\alpha(S_t)$ with $S_c>S_t$ is impossible if $S_c$ lies at or before the maximum point.
\end{proof}

\subsection{Localization of stagnation points}\label{sec_localization}

\begin{theorem}[Localization]\label{thm_localization}
Let $(\Psi,S)$ be a classical solution with a stagnation point $P_c=(R_c,\theta_c)$. Then $\Delta(\alpha)\ge0$, and the following hold.
\begin{enumerate}
\item If $\Delta>0$, then $R_c=r_-(\alpha)<\alpha$ and $S\le r_-$. Each stagnation point is a global maximum of $S$, and $\{\theta:S(\theta)=r_-\}$ is exactly the set of stagnation points. Consequently $\Omega$ lies in the closed disc $\{R\le R_c\}$.
\item If $\Delta=0$, the stagnation points are the points where $S=\alpha$, and $\min S<\alpha$.
\end{enumerate}
\end{theorem}
\begin{proof}
At a stagnation point, $Q_\alpha(R_c)=0$, which by \eqref{eq_completed_square} forces $\Delta\ge0$.

(a) By \autoref{lem_dichotomy}, either $S\le r_-$ or $S\ge r_+$. Suppose $R_c=r_+$, so that $S\ge R_c$. Then $\{1<R<R_c\}\subset\Omega$, and for small $\varepsilon$ the disc $B=B_\varepsilon(P_c-\varepsilon\mathbf e_R)$ lies in $\Omega$ and touches $\partial\Omega$ at $P_c$. The Hopf lemma \cite{GilbargTrudinger2001} gives
\[
\liminf_{t\downarrow0}t^{-1}\Psi(P_c-t\mathbf e_R)>0 .
\]
On the other hand, $\nabla\Psi$ is continuous on $\overline\Omega$ and vanishes at $P_c$, so $\Psi(P_c-t\mathbf e_R)=o(t)$. This contradiction shows $R_c=r_-$ and $S\le r_-$. Finally, $Q_\alpha(S(\theta))=0$ if and only if $S(\theta)=r_-$.

(b) Here $Q_\alpha=(R-\alpha)^2$. If $\min S\ge\alpha$, the same Hopf argument applies.
\end{proof}

\subsection{Gradient bounds}\label{sec_gradient}

\begin{lemma}[Weighted subharmonicity]\label{lem_subharmonic}
Let $t=\log R$, and let $\Delta_c$ and $\nabla_c$ denote the flat operators in $(t,\theta)$.
\begin{enumerate}
\item For every $k\in\mathbb R$, the function $e^{(k-2)t}|\nabla_c\Psi|^2=R^k|\nabla\Psi|^2$ is subharmonic.
\item $\Delta_c\big(R^kQ_\alpha(R)\big)=R^kf_k(R)$, where
\[
f_k(R)=(k+2)^2R^2-2\alpha(k+1)^2R+2\lambda k^2 .
\]
\end{enumerate}
\end{lemma}
\begin{proof}
(a) The function $F=\Psi_t-i\Psi_\theta$ is holomorphic in $w=t+i\theta$, and
\[
e^{(k-2)t}|\nabla_c\Psi|^2=|e^{(k-2)w/2}F|^2 .
\]
(b) Use $\partial_t=R\partial_R$ and $(R\partial_R)^2R^j=j^2R^j$.
\end{proof}

\begin{theorem}[Maximum-velocity principle]\label{thm_max_velocity}
For every classical solution, $\sup_\Omega R^2|\nabla\Psi|^2=\max_\theta S^2Q_\alpha(S)/p_0^2$. For circular flows, $R^2|\nabla\Psi|^2$ is constant and equal to this value.
\end{theorem}
\begin{proof}
Fix $k>2$ and let $V_k=R^k|\nabla\Psi|^2$. On the bottom, $\Psi_\theta=\Psi_{tt}=0$, so
\[
\partial_tV_k=(k-2)\Psi_t^2>0 .
\]
Hence $V_k$ is not constant, and it cannot attain its maximum at the bottom: the Hopf lemma would force $\partial_tV_k<0$ there. By the strong maximum principle, the maximum is therefore attained on the free surface, where $V_k=S^kQ_\alpha(S)/p_0^2$. Letting $k\downarrow2$ gives the inequality; equality is attained on the free surface. For circular flows, $R^2|\nabla\Psi|^2=\log^{-2}S$.
\end{proof}

\begin{corollary}[First hodograph constraint]\label{cor_first_constraint}
For solutions of class (L),
\[
(h+1)h_p=\frac{|p_0|R}{|V-R|}\ge\frac{|p_0|}{\sqrt{\max_\theta S^2Q_\alpha(S)}} .
\]
\end{corollary}
\begin{proof}
Since $|V-R|=|p_0||\Psi_R|\le|p_0||\nabla\Psi|$ and $R\ge1$, this follows from \autoref{thm_max_velocity}.
\end{proof}

\begin{theorem}[Local Bernstein bound]\label{thm_local_bernstein}
Let $\Delta(\alpha)>0$, and let $(\Psi,S)$ be a classical solution with $\max S\le r_-$; this holds in particular when a stagnation point exists. Choose
\[
k>k_0=\tfrac12\big(r_-/\sqrt\Delta-1\big)
\]
and $\eta\in(0,r_--1)$ with $f_k<0$ on $[r_--\eta,r_-]$. Then, in $U=\Omega\cap\{R>r_--\eta\}$,
\[
R^k\big(p_0^2|\nabla\Psi|^2-Q_\alpha(R)\big)\le\frac M\eta\,(r_--R),
\qquad M=\max_{\overline\Omega}\big|R^k\big(p_0^2|\nabla\Psi|^2-Q_\alpha\big)\big| .
\]
Consequently $p_0^2|\nabla\Psi|^2\le C_1(r_--R)$ and $\Psi\le C_2(r_--R)^{3/2}$ in $U$.
\end{theorem}
\begin{proof}
From
\[
f_k(r_-)=r_-\big[(2k+3)r_--(2k+1)r_+\big]=r_-\big[2r_--(2k+1)2\sqrt\Delta\big],
\]
the choice $k>k_0$ makes $f_k(r_-)<0$, so an admissible $\eta$ exists.

Let $P=R^k(p_0^2|\nabla\Psi|^2-Q_\alpha)$ and $W=\frac M\eta(r_--R)$. By \autoref{lem_subharmonic}, $P$ is subharmonic in $U$. Since $\Delta|x|=1/|x|$, the function $W$ is superharmonic, so $P-W$ is subharmonic in $U$. On the free surface $P=0\le W$, because $S\le r_-$. On the inner arc $R=r_--\eta$ we have $P\le M=W$. The maximum principle gives $P\le W$ in $U$.

For the first consequence, use $Q_\alpha(R)\le(r_+-1)(r_--R)$ on $[1,r_-]$. For the second, integrate along the radial segment from $x$ to the free surface, on which $\Psi=0$:
\[
\Psi(x)\le\int_R^{r_-}\sqrt{C_1(r_--\rho)}\,d\rho/|p_0| .
\]
\end{proof}

\begin{lemma}[Phragm\'en--Lindel\"of at a boundary point]\label{lem_PL}
Let $D\subset\mathbb C$ be a bounded domain and $P\in\partial D$. Suppose that $\partial D$ contains a connected set $K\ni P$ with $\operatorname{diam}K\ge\rho_0>0$. Let $G$ be holomorphic in $D$ and continuous on $\overline D\setminus\{P\}$, with
\[
|G|\le M\ \text{ on }\partial D\setminus\{P\},\qquad |G(z)|\le C|z-P|^{-m}\ \text{ for }z\in D,
\]
for some constants $M\ge1$, $C$ and $m>0$. Then $|G|\le M$ in $D$.
\end{lemma}
\begin{proof}
The function $v=\log|G|-\log M$ is subharmonic in $D$. Fix $z\in D$ with $|z-P|<\rho_0/2$, and let $0<r<|z-P|/2$. Put $D_r=D\setminus\overline{B_r(P)}$ and $\Gamma_r=D\cap\partial B_r(P)$. On $\partial D_r\setminus\Gamma_r\subset\partial D\setminus\{P\}$ we have $v\le0$, and on $\Gamma_r$ we have $v\le\log^+(Cr^{-m}/M)$. By the maximum principle,
\[
v(z)\le\log^+\!\big(Cr^{-m}/M\big)\,\omega_r(z),
\]
where $\omega_r$ is the harmonic measure of $\Gamma_r$ in $D_r$.

The domain $D_r$ is contained in $\mathbb C\setminus(K\cup\overline{B_r(P)})$, and $K$ joins $P$ to the circle $\partial B_{\rho_0/2}(P)$. By the maximum principle and Beurling's estimate \cite{GarnettMarshall2005},
\[
\omega_r(z)\le C_0\Big(\frac r{|z-P|}\Big)^{1/2}
\]
with an absolute constant $C_0$. Hence $v(z)\le C_0\log^+(Cr^{-m}/M)\,(r/|z-P|)^{1/2}\to0$ as $r\to0$, and so $v\le0$ in $D\cap B_{\rho_0/2}(P)$.

On the rest of $D$, $v$ is subharmonic, bounded above and $\le0$ on the boundary of $D\setminus\overline{B_{\rho_0/2}(P)}$, so $v\le0$ there as well.
\end{proof}

\begin{proposition}[Interior Bernoulli inequality]\label{prop_interior_bernoulli}
Let $\Delta(\alpha)\ge0$, and let $(\Psi,S)$ be a classical solution in the lower region, $\max S\le r_-$. Put
\[
M=\max\Big\{1,\ \frac{p_0^2\max_\theta\Psi_R(1,\theta)^2}{Q_\alpha(1)}\Big\}.
\]
Then
\[
p_0^2|\nabla\Psi(x)|^2\le M\,Q_\alpha(|x|)\qquad\text{for all }x\in\Omega .
\]
In particular $p_0^2|\nabla\Psi|^2\le M(r_+-1)(r_--R)$ in $\Omega$, so the growth bound of \autoref{thm_local_bernstein} holds in all of $\Omega$, with an explicit constant and with no condition on $k$.
\end{proposition}
\begin{proof}
Let $V=\Psi_{x_1}-i\Psi_{x_2}$, which is holomorphic in $\Omega$, continuous on $\overline\Omega$, and satisfies $|V|=|\nabla\Psi|$. Let $\mathcal P(\zeta)=(r_--\zeta)(r_+-\zeta)$, so that $\mathcal P(R)=Q_\alpha(R)$ for real $R$.

Fix $\theta_0$ and set $G(z)=p_0^2V(z)^2/\mathcal P(ze^{-i\theta_0})$. The zeros of the denominator are $z=r_\pm e^{i\theta_0}$. Since $|z|\le S\le r_-$ on $\overline\Omega$, the only one that can lie in $\overline\Omega$ is $P=r_-e^{i\theta_0}$, and only when $S(\theta_0)=r_-$, in which case $P$ is a stagnation point on $\partial\Omega$. We bound $G$ on the boundary.
\begin{itemize}
\item On the free surface $z=Se^{i\theta}$, using $|r-Se^{i\varphi}|\ge|r-S|$ for real $r\ge S$,
\[
|\mathcal P(ze^{-i\theta_0})|\ge(r_--S)(r_+-S)=Q_\alpha(S)=p_0^2|V|^2,
\]
so $|G|\le1$.
\item On the bottom, $|V|=|\Psi_R|$ and $|\mathcal P(ze^{-i\theta_0})|\ge Q_\alpha(1)$, so $|G|\le M$.
\end{itemize}
Near $P$, $|G|\le C|z-P|^{-2}$ since $V$ is bounded. The free surface is a closed curve around the origin, so \autoref{lem_PL} applies with $\rho_0\ge2$ and gives $|G|\le M$ in $\Omega$. For $x=|x|e^{i\theta_0}$ this reads $p_0^2|\nabla\Psi(x)|^2\le M\,\mathcal P(|x|)=M\,Q_\alpha(|x|)$. Since $\theta_0$ was arbitrary, the inequality holds everywhere.
\end{proof}

\begin{remark}
The constant $M$ cannot be taken below $1$, since equality holds on the free surface. The inequality is also visible numerically. On the computed waves of \autoref{sec_numerics}, $\sup_\Omega p_0^2|\nabla\Psi|^2/Q_\alpha(|x|)$ equals $M$ to the accuracy of the computation. The supremum is attained on the free surface in runs~(a) and~(c), where $M=1$, and on the bottom in run~(b), where $M=1.469$.
\end{remark}

Since $R_c-R(x)=|P_c|-|x|\le|x-P_c|$, \autoref{thm_local_bernstein} implies, near a non-degenerate stagnation point,
\begin{equation}\label{eq_growth}
|\nabla\Psi(x)|^2\le C|x-P_c|,\qquad\Psi(x)\le C|x-P_c|^{3/2},\qquad\int_{B_r(P_c)\cap\Omega}|\nabla\Psi|^2\le Cr^3 .
\end{equation}

\section{The crest shape and limits of regular waves}\label{sec_corner}

This section proves the crest-shape half of the Stokes conjecture and the limit theorem, in five steps.
\begin{itemize}
\item A perturbed Weiss monotonicity formula and non-degeneracy (\autoref{sec_variational}).
\item The classification of blow-up limits, which gives the Stokes corner (\autoref{sec_blowup}).
\item The right-angle corner at the degenerate parameters (\autoref{subsec_degenerate}).
\item The extension of both corner theorems to profiles that are not radial graphs (\autoref{subsec_variational_crest}).
\item Compactness of families of regular waves (\autoref{sec_limit}).
\end{itemize}

\subsection{Weiss monotonicity and non-degeneracy}\label{sec_variational}

Fix a stagnation point $P_c$ with $\Delta>0$. Use the frame $e_1=\mathbf e_\theta(\theta_c)$, $e_2=-\mathbf e_R(\theta_c)$ and coordinates $x=P_c+y_1e_1+y_2e_2$, so that $y_2$ is the depth below the crest. Define
\[
d=R_c-R(x)=y_2-\frac{y_1^2}{2R_c}+O(|y|^3),\qquad
u=\frac{|p_0|\Psi}{(2\sqrt\Delta)^{1/2}},\qquad
q=d+\frac{d^2}{2\sqrt\Delta}.
\]
Then $Q_\alpha=2\sqrt\Delta\,d+d^2$, the Bernoulli condition becomes $|\nabla u|^2=q$ on $\partial\Omega$, and $q=y_2+O(|y|^2)$ is smooth. We extend $u$ by zero outside $\Omega$ and write $\chi=\chi_\Omega$.

\begin{lemma}[Radial phase monotonicity]\label{lem_phase}
Air is radially closed outward: if $x\notin\Omega$ with $R(x)\ge S(\theta(x))$, then $x+s\mathbf e_R\notin\Omega$ for all $s\ge0$. Fluid is radially closed inward: if $x\in\Omega$, then $x-s\mathbf e_R\in\Omega$ for $0\le s<R(x)-1$. Under the scaling $y\mapsto ry$, radial lines through points with $|y|\sim1$ converge to vertical lines at rate $O(r)$.
\end{lemma}
\begin{proof}
This follows immediately from $\Omega=\{1<R<S(\theta)\}$.
\end{proof}

\begin{lemma}[Domain variation]\label{lem_domain_variation}
For every $\phi\in C^1_c(B_{\rho_0};\mathbb R^2)$,
\begin{equation}\label{eq_domain_variation}
\int\Big[(|\nabla u|^2+q\chi)\operatorname{div}\phi-2\nabla u\cdot D\phi\,\nabla u+\chi\nabla q\cdot\phi\Big]=0 .
\end{equation}
\end{lemma}
\begin{proof}
For harmonic $u$,
\[
\operatorname{div}\big(|\nabla u|^2\phi-2(\nabla u\cdot\phi)\nabla u\big)=|\nabla u|^2\operatorname{div}\phi-2\nabla u\cdot D\phi\nabla u ,
\]
and $\operatorname{div}(q\phi)=q\operatorname{div}\phi+\nabla q\cdot\phi$. Integrate over $\Omega$ with small balls removed around the finitely many non-$C^1$ points. On the free surface $\nabla u=-|\nabla u|\nu$, so the boundary integrand is $(q-|\nabla u|^2)\phi\cdot\nu=0$. The small circles contribute $O(\varepsilon)$.
\end{proof}

\begin{lemma}[Perturbed Weiss monotonicity, cf.\ \cite{Weiss1999,VW2011}]\label{lem_weiss}
Let
\[
\Phi(r)=r^{-3}\int_{B_r}(|\nabla u|^2+q\chi)-\tfrac32r^{-4}\int_{\partial B_r}u^2 .
\]
Then, for almost every $r$,
\[
\Phi'(r)=2r^{-3}\int_{\partial B_r}\Big(\partial_\nu u-\frac{3u}{2r}\Big)^2+r^{-4}\int_{B_r}\chi\,(y\cdot\nabla q-q).
\]
The last term is bounded by a constant $C_W$, because $y\cdot\nabla q-q=O(|y|^2)$. Hence $\Phi(r)+C_Wr$ is non-decreasing and $\Phi(0+)$ exists.
\end{lemma}
\begin{proof}
Insert $\phi=y\,\eta_\kappa(|y|)$ into \eqref{eq_domain_variation} and let $\eta_\kappa\to\mathbf 1_{[0,r)}$. In two dimensions the volume terms in $|\nabla u|^2$ cancel, leaving
\[
\int_{B_r}\chi(2q+y\cdot\nabla q)=r\int_{\partial B_r}(|\nabla u|^2+q\chi)-2r\int_{\partial B_r}(\partial_\nu u)^2 .
\]
Write $I(r)=\int_{B_r}(|\nabla u|^2+q\chi)$ and $J(r)=\int_{\partial B_r}u^2$. Since $u$ is harmonic in $\Omega$ and vanishes on $\partial\Omega$, we have $\int_{B_r}|\nabla u|^2=\int_{\partial B_r}u\,\partial_\nu u$ for almost every $r$. Using the identity above,
\[
\frac{d}{dr}\big(r^{-3}I\big)=-3r^{-4}I+r^{-3}I'
=-3r^{-4}\!\int_{\partial B_r}u\,\partial_\nu u+2r^{-3}\!\int_{\partial B_r}(\partial_\nu u)^2+r^{-4}\!\int_{B_r}\chi\,(y\cdot\nabla q-q).
\]
In two dimensions $J(r)=r\int_{\partial B_1}u(rz)^2\,dz$, so
\[
\frac{d}{dr}\Big(-\frac32r^{-4}J\Big)=\frac92r^{-5}J-3r^{-4}\int_{\partial B_r}u\,\partial_\nu u .
\]
Adding the two lines and using $2r^{-3}\big(a-\frac{3u}{2r}\big)^2=2r^{-3}a^2-6r^{-4}au+\frac92r^{-5}u^2$ with $a=\partial_\nu u$ gives the formula for $\Phi'$. Since $q(0)=0$ and $\nabla q(0)=e_2$, the function $y\cdot\nabla q-q$ vanishes to second order at $0$, so $|y\cdot\nabla q-q|\le C|y|^2$ and the last term is at most $Cr^{-4}\cdot r^4=C_W$. Boundedness of $\Phi$ follows from \eqref{eq_growth}.
\end{proof}

\begin{lemma}[Bernoulli-flux non-degeneracy]\label{lem_nondeg}
Let $\sigma\in(0,1]$. There are $c_N(\sigma)>0$ and $\rho_2>0$ with the following property. If $z\in\partial\Omega$, $0<|z|<\rho_2$, $d(z)\ge\sigma|z|$ and $\rho\le\sigma|z|/4$, then
\[
\sup_{B_{4\rho}(z)}u\ge c_N\,\rho\,|z|^{1/2}.
\]
\end{lemma}
\begin{proof}
The free surface leaves $B_\rho(z)$, so it contains an arc of length at least $\rho$ starting at $z$. On this arc $d\ge\sigma|z|/2$, hence $|\nabla u|\ge(\sigma|z|/2)^{1/2}$.

For almost every $t\in(\rho,2\rho)$ the divergence theorem gives
\[
\int_{\partial\Omega\cap B_t(z)}|\nabla u|=\int_{\Omega\cap\partial B_t(z)}\partial_ru\le\int_{\Omega\cap\partial B_t(z)}|\nabla u| .
\]
Averaging over $t$ and applying the Cauchy--Schwarz inequality, then Caccioppoli's inequality $\int\eta^2|\nabla u|^2\le4\int u^2|\nabla\eta|^2$ for the subharmonic extension, we obtain
\[
\rho^2(\sigma|z|/2)^{1/2}\le C\rho\Big(\int_{B_{2\rho}(z)}|\nabla u|^2\Big)^{1/2}\le C'\rho\sup_{B_{4\rho}(z)}u .
\]
\end{proof}

\begin{theorem}[Angle bounds]\label{thm_angle_bounds}
Let $P_c$ be a stagnation point with $\Delta>0$.
\begin{enumerate}
\item No open sector with vertex $P_c$ and opening larger than $2\pi/3$ is contained, near $P_c$, in $\Omega$.
\item $\Omega$ is not contained, near $P_c$, in any closed sector of opening smaller than $2\pi/3$ that lies in the open lower half-plane $\{y_2>0\}$ apart from its vertex.
\end{enumerate}
\end{theorem}
\begin{proof}
(a) Suppose such a sector exists. Choose a slightly smaller sector of opening $\beta''>2\pi/3$, and let $h=r^{\pi/\beta''}\sin(\pi\phi/\beta'')$ in it. Since $\Psi\ge\mu_0>0$ on a compact arc inside $\Omega$, the maximum principle gives $\Psi\ge ch$. Along the bisector this means $\Psi\ge cr^{\pi/\beta''}$ with $\pi/\beta''<3/2$, contradicting \eqref{eq_growth}.

(b) Suppose such a sector exists. Since the extension of $u$ by zero is subharmonic, a harmonic barrier on a slightly larger sector gives $\sup_{B_{4r}}u\le Cr^{\pi/\beta'''}$ with $\pi/\beta'''>3/2$. On the other hand, one branch of the free surface crosses $B_r\setminus B_{r/2}$ at depth $d\ge\sigma r$. The flux argument of \autoref{lem_nondeg} then gives $\sup_{B_{4r}}u\ge cr^{3/2}$, a contradiction.
\end{proof}

\subsection{Blow-up classification and the Stokes corner}\label{sec_blowup}

Set
\[
u_r(y)=r^{-3/2}u(ry),\qquad\Omega_r=r^{-1}\Omega,\qquad\chi_r=\chi_{\Omega_r},\qquad q_r(y)=r^{-1}q(ry).
\]
Then $q_r\to y_2$ in $C^1_{\rm loc}$, and $\Phi(rs)$ equals the Weiss functional of $(u_r,q_r,\chi_r)$ evaluated at $s$.

\begin{proposition}[Compactness and homogeneity]\label{prop_blowup}
Every sequence $r_j\downarrow0$ has a subsequence along which
\[
u_{r_j}\to u_0\ \text{locally uniformly},\qquad \nabla u_{r_j}\to\nabla u_0\ \text{in }L^2_{\rm loc},\qquad \chi_{r_j}\overset{*}{\rightharpoonup}\chi_0 .
\]
The limits have the following properties.
\begin{enumerate}
\item $u_0\ge0$ is continuous and harmonic in $\{u_0>0\}$. Moreover $u_0=0=\chi_0$ on $\{y_2<0\}$, and $\chi_0=1$ on $\{u_0>0\}$.
\item $u_0$ is homogeneous of degree $3/2$.
\item The limiting domain-variation identity holds with $q$ replaced by $y_2$.
\item $\Phi(0+)=\int_{B_1}y_2\chi_0$.
\end{enumerate}
\end{proposition}
\begin{proof}
The bounds \eqref{eq_growth} make the $u_r$ locally uniformly Lipschitz, so a subsequence converges locally uniformly. Since $\Omega\subset\{d\ge0\}$, we have $\chi_r=0$ eventually on compact subsets of $\{y_2<0\}$.

For strong convergence of gradients, note that $\int\eta|\nabla u_r|^2=-\int u_r\nabla u_r\cdot\nabla\eta$, because $u_r$ is harmonic in $\Omega_r$ and vanishes on $\partial\Omega_r$. The right-hand side converges, and $u_0$ satisfies the same identity (test with $\eta(u_0-\varepsilon)^+$ and let $\varepsilon\to0$). Hence the $L^2$ norms converge, and weak convergence upgrades to strong convergence.

For (b), integrate \autoref{lem_weiss} between $rs_1$ and $rs_2$. The left-hand side tends to $\Phi(0+)-\Phi(0+)=0$, and the error term is $O(r)$. By strong convergence,
\[
\int_{B_{s_2}\setminus B_{s_1}}|y|^{-3}\big(\nabla u_0\cdot y/|y|-\tfrac32u_0/|y|\big)^2=0 .
\]
Parts (c) and (d) follow by passing to the limit, using $\int_{B_1}|\nabla u_0|^2=\frac32\int_{\partial B_1}u_0^2$ for (d).
\end{proof}

\begin{lemma}[Constant density and jump condition]\label{lem_jump}
Let $Z$ be a connected component of the interior of $\{u_0=0\}$ inside $\{y_2>0\}$. Then $\chi_0$ equals a constant $c_Z$ a.e.\ on $Z$. If $\Gamma$ is a ray of $\partial\{u_0>0\}$ separating the positivity set from $Z$, then
\[
|\nabla u_0|^2=(1-c_Z)\,y_2\qquad\text{on }\Gamma .
\]
If $\Gamma\subset\{y_2=0\}$, then $|\nabla u_0|=0$ on $\Gamma$.
\end{lemma}
\begin{proof}
For $\phi$ supported in $Z$, the limiting identity reduces to $\int\chi_0\operatorname{div}(y_2\phi)=0$. Since $y_2>0$ on $Z$, the field $y_2\phi$ ranges over all of $C^1_c(Z)$, so $\nabla\chi_0=0$.

For $\phi$ supported near a point of $\Gamma$, integrate by parts on each side. The positivity side contributes $\int_\Gamma(y_2-|\nabla u_0|^2)\phi\cdot\nu$, and the side $Z$ contributes $-c_Z\int_\Gamma y_2\phi\cdot\nu$. Since $\phi$ is arbitrary, the jump condition follows.
\end{proof}

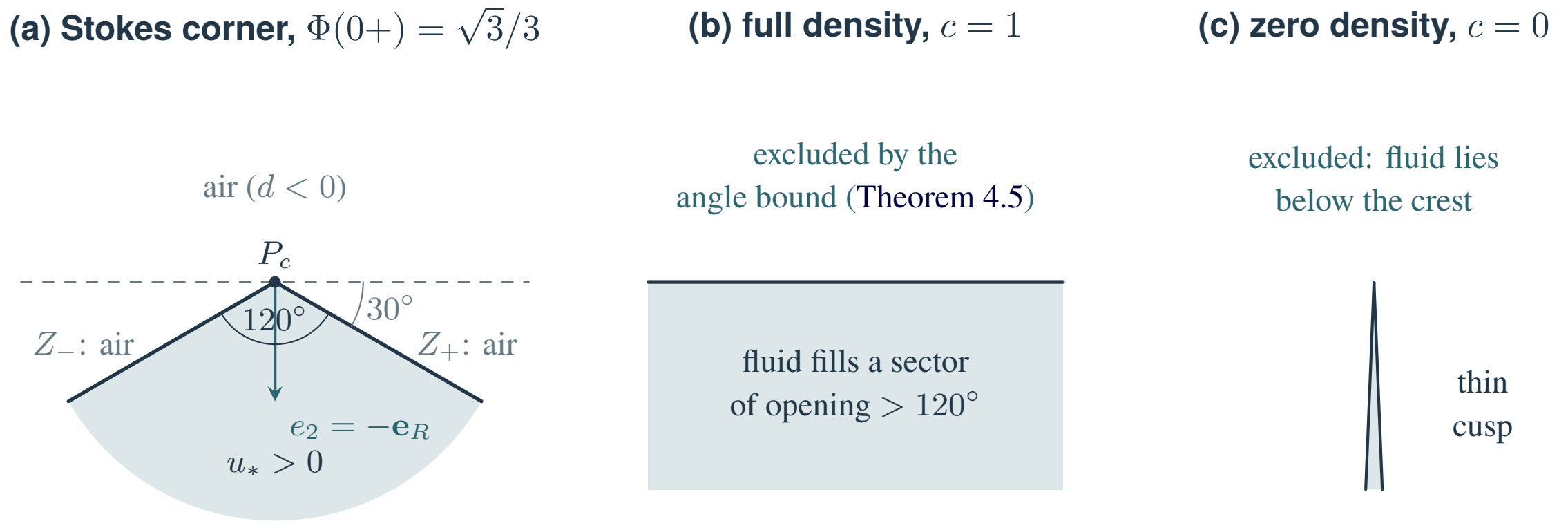
\begin{figure}[!ht]
\centering
\resizebox{\linewidth}{!}{%
\begin{tikzpicture}[x=1cm,y=1cm,>=stealth,font=\small,line cap=round,line join=round]
\begin{scope}
\node[font=\sffamily\small\bfseries,text=PaperInk] at (0,2.45) {(a) Stokes corner, $\Phi(0+)=\sqrt3/3$};
\fill[PaperAccent!16] (0,0) -- (-30:0) -- (210:2.3) arc (210:330:2.3) -- cycle;
\draw[PaperInk,line width=1pt] (210:2.3) -- (0,0) -- (330:2.3);
\draw[PaperMuted,dashed] (-2.45,0) -- (2.45,0);
\draw[PaperAccent,->,line width=.8pt] (0,0) -- (0,-1.15) node[below right,font=\footnotesize] {$e_2=-\mathbf e_R$};
\draw[PaperMuted] (330:.85) arc (330:360:.85);
\node[text=PaperMuted,font=\footnotesize] at (1.12,-.24) {$30^\circ$};
\draw[PaperInk] (210:.6) arc (210:330:.6);
\node[text=PaperInk,font=\footnotesize] at (0,-.35) {$120^\circ$};
\node[text=PaperInk,font=\footnotesize] at (0,-1.75) {$u_*>0$};
\node[text=PaperMuted,font=\footnotesize] at (1.85,-.62) {$Z_+$: air};
\node[text=PaperMuted,font=\footnotesize] at (-1.85,-.62) {$Z_-$: air};
\fill[PaperInk] (0,0) circle (1.6pt) node[above,font=\footnotesize] {$P_c$};
\node[text=PaperMuted,font=\footnotesize] at (0,.9) {air ($d<0$)};
\end{scope}
\begin{scope}[shift={(5.6,0)}]
\node[font=\sffamily\small\bfseries,text=PaperInk] at (0,2.45) {(b) full density, $c=1$};
\fill[PaperAccent!16] (-2,0) rectangle (2,-2.0);
\draw[PaperInk,line width=1pt] (-2,0)--(2,0);
\node[text=PaperInk,font=\footnotesize,align=center] at (0,-1.0) {fluid fills a sector\\of opening $>120^\circ$};
\node[text=PaperAccent,font=\footnotesize,align=center] at (0,1.0) {excluded by the\\angle bound (\autoref{thm_angle_bounds})};
\end{scope}
\begin{scope}[shift={(10.6,0)}]
\node[font=\sffamily\small\bfseries,text=PaperInk] at (0,2.45) {(c) zero density, $c=0$};
\fill[PaperAccent!16] (0,0) -- (-.08,-2.0) -- (.08,-2.0) -- cycle;
\draw[PaperInk,line width=.8pt] (-.08,-2.0) -- (0,0) -- (.08,-2.0);
\node[text=PaperInk,font=\footnotesize,align=center] at (1.05,-1.2) {thin\\cusp};
\node[text=PaperAccent,font=\footnotesize,align=center] at (0,1.0) {excluded: fluid lies\\below the crest};
\end{scope}
\end{tikzpicture}}
\caption{The three possible blow-up densities (\autoref{prop_classification}). Only the Stokes corner (a) survives (\autoref{thm_main_stokes}). In (a) the sector $S_*$ has its tangent rays $30^\circ$ below the horizontal and is bisected by the inward radial direction.}
\label{fig_blowups}
\end{figure}

\begin{proposition}[Classification]\label{prop_classification}
Let
\[
S_*=\{\pi/6<\varphi<5\pi/6\},\qquad u_*=\tfrac{\sqrt2}{3}\varrho^{3/2}\cos\big(\tfrac32(\varphi-\tfrac\pi2)\big)\chi_{S_*}.
\]
Every blow-up limit is of one of the following two types:
\begin{enumerate}
\item $u_0=u_*$ and $\chi_0=\chi_{S_*}$, so that $\Phi(0+)=\sqrt3/3$;
\item $u_0\equiv0$ and $\chi_0\equiv c\in\{0,1\}$ on $\{y_2>0\}$, so that $\Phi(0+)=2c/3$.
\end{enumerate}
\end{proposition}
\begin{proof}
\emph{Step 1 (shape of the positivity set).} By homogeneity, each angular component $(a,b)$ of $\{u_0>0\}$ carries $u_0=\varrho^{3/2}f(\varphi)$ with $f''+\frac94f=0$, $f>0$ on $(a,b)$ and $f(a)=f(b)=0$. Hence $b-a=2\pi/3$, and at most one such sector fits in $(0,\pi)$. So either $u_0\equiv0$, or
\[
u_0=A\varrho^{3/2}\cos\big(\tfrac32(\varphi-\varphi_r-\tfrac\pi3)\big)\quad\text{on }(\varphi_r,\varphi_r+2\pi/3),\quad A>0 .
\]
In the second case \autoref{lem_jump} excludes $\varphi_r=0$, so $\varphi_r\in(0,\pi/3)$.

\emph{Step 2 (the zero regions carry no fluid).} The three possible outcomes are sketched in \autoref{fig_blowups}. Let
\[
K_+=\{\tfrac12\le|y|\le5,\ \delta\le\varphi\le\varphi_r-\delta\}.
\]
Since $u_{r_j}\to0$ near $K_+$, \autoref{lem_nondeg} shows that $\partial\Omega_{r_j}$ does not meet $K_+$ for large $j$. As $K_+$ is connected, eventually either $K_+\subset\Omega_{r_j}$ or $K_+$ lies in the air.

Now consider $z_0=2(\cos\varphi_r,\sin\varphi_r)$. The function $u_0$ is not harmonic near $z_0$, so every $B_\varepsilon(z_0)$ contains a point that is not in $\Omega_{r_j}$. By \autoref{lem_phase}, the nearly vertical segment from that point towards $\{y_2=0\}$ lies in the air. Because $\cos\varphi_r>\frac12$, this segment crosses $K_+$. Hence $K_+$ is eventually in the air, and $c_{Z_+}=0$. The same argument on the left gives $c_{Z_-}=0$.

\emph{Step 3 (symmetry).} \autoref{lem_jump} now gives
\[
\tfrac94A^2=\sin\varphi_r=\sin(\varphi_r+2\pi/3).
\]
Hence $\varphi_r=\pi/6$ and $A=\sqrt2/3$.

\emph{Step 4 (the zero case).} Here $\{y_2>0\}$ is a single zero region, so $\chi_0\equiv c$. With $K=\{1\le|y|\le2,\ \delta\le\varphi\le\pi-\delta\}$, the argument of Step 2 shows that $K\subset\Omega_{r_j}$ or $K$ lies in the air, eventually. Weak-$*$ convergence forces this alternative to stabilize, so $c\in\{0,1\}$.
\end{proof}

\begin{theorem}[Annular Stokes corner]\label{thm_main_stokes}
Under the assumptions of \autoref{thm_intro_stokes}, conclusions (a)--(d) of that theorem hold.
\end{theorem}
\begin{proof}
Part (a) is \autoref{thm_localization}. By \autoref{prop_classification}, $\Phi(0+)\in\{\sqrt3/3,\,0,\,2/3\}$, and we exclude the last two values.

Suppose $\Phi(0+)=0$, so that every blow-up is of type (b) with $c=0$. By the stabilization argument of Step 4, the set $\{1\le|y|\le2,\ \pi/4\le\varphi\le3\pi/4\}$ lies in the air for all small $r$. This contradicts the fact that $P_c-\frac32r\,\mathbf e_R\in\Omega$.

Suppose $\Phi(0+)=2/3$. Then, in the same way, $\{1\le|y|\le2,\ \delta\le\varphi\le\pi-\delta\}\subset\Omega_r$ for all small $r$. The union of these sets over $r$ is an open sector of opening $\pi-2\delta>2\pi/3$ inside $\Omega$, contradicting \autoref{thm_angle_bounds}(a).

Hence $\Phi(0+)=\sqrt3/3$. This value differs from $0$ and $2/3$, so every blow-up equals $u_*$, and the whole family converges. Undoing the normalization gives $A_*=(2\sqrt\Delta)^{1/2}\frac{\sqrt2}{3}/|p_0|=2\Delta^{1/4}/(3|p_0|)$, which is (c).

For (b), fix $\delta>0$ and consider the annulus $1\le|y|\le2$. Points of $S_*$ at angular distance at least $\delta$ from its rays lie in $\Omega_r$, since $u_r\ge c(\delta)>0$ there. The regions $K_\pm$ of Step 2 lie in the air for all small $r$, and not only along a subsequence: otherwise a sequence $r_j\to0$ with fluid in $K_+$ would have a blow-up limit with $c_{Z_+}>0$, contradicting the uniqueness of the limit $u_*$ and Step 2. A fluid point within angle $\delta$ of the horizontal would, by \autoref{lem_phase}, force fluid into $K_\pm$. Hence $\Omega_r$ is squeezed between the $\delta$-core of $S_*$ and its $\delta$-neighbourhood, which is the tangent-cone statement. The direction $\varphi=\pi/6$ corresponds to $\tau_+$.

For (d), part (b) gives $d\ge c|y|>0$ on $\partial\Omega\setminus\{P_c\}$ near $P_c$, so $S<R_c$ there. By \autoref{thm_localization}, $P_c$ is therefore an isolated stagnation point, and compactness of $\mathbb T$ gives finiteness.
\end{proof}

\begin{remark}
In the normalization $|\nabla u|^2=y_2$, the profile $u_*$ is exactly the Stokes corner flow of \cite{VW2011}. The effective gravity at the crest is $g_{\rm eff}=\sqrt\Delta$, read off from $Q_\alpha\approx2\sqrt\Delta\,d$. The proof uses neither the frequency formula nor any assumption on the number of air components. Radial phase monotonicity rules out both horizontally flat and cusp-like stagnation. For flat waves with vorticity these alternatives were left open in \cite{VW2012}; cusps have since been excluded there \cite{DPY2026}. For the annular problem without the graph structure, see \autoref{subsec_variational_crest}.
\end{remark}

\subsection{Degenerate parameters: right-angle corners}\label{subsec_degenerate}

At $\alpha=\alpha_\pm$ we have $\Delta=0$ and $Q_\alpha(R)=(\alpha-R)^2$. The effective gravity then vanishes to first order at the stagnation circle, the natural homogeneity changes from $3/2$ to $2$, and the Stokes corner is replaced by a right angle. The local Bernstein bound of \autoref{thm_local_bernstein} is not available here: for $\Delta=0$ one computes $f_k(\alpha)=2\alpha^2>0$ for every $k$. The growth is instead supplied by \autoref{lem_PL}.

\begin{theorem}[Right-angle corners at degenerate parameters]\label{thm_degenerate}
Let $\alpha\in\{\alpha_+,\alpha_-\}$, and let $(\Psi,S)$ be a classical solution with a stagnation point $P_c=(R_c,\theta_c)$. Then $\alpha>1$, and:
\begin{enumerate}
\item $R_c=\alpha$, $S\le\alpha$ everywhere, and $P_c$ is a global maximum of $S$;
\item the fluid has at $P_c$ a tangent cone which is the sector of opening $90^\circ$, symmetric about $-\mathbf e_R(\theta_c)$;
\item in the frame $e_1=\mathbf e_\theta(\theta_c)$, $e_2=-\mathbf e_R(\theta_c)$, as $r\to0$,
\[
r^{-2}\Psi(P_c+ry)\to\frac{\sqrt2}{4|p_0|}\,\big(y_2^2-y_1^2\big)\,\chi_{\{|y_1|<y_2\}}
\]
locally uniformly, with convergence of gradients in $L^2_{\rm loc}$;
\item stagnation points are isolated, and hence finite in number.
\end{enumerate}
\end{theorem}

\begin{proof}
By \autoref{thm_localization}(b), $R_c=\alpha$. Moreover $\alpha>1$, since by the Hopf lemma $\Psi_R<0$ on the bottom. Set
\[
u=|p_0|\Psi,\qquad d=\alpha-R(x)=y_2-\frac{y_1^2}{2\alpha}+O(|y|^3),\qquad q=d^2=y_2^2+O(|y|^3).
\]
Then $|\nabla u|^2=q$ on $\partial\Omega$, and \autoref{lem_phase} and \autoref{lem_domain_variation} hold verbatim, since their proofs do not use the form of $q$. We follow the proof of \autoref{thm_main_stokes} and record only the changes.

\emph{Step 1: quadratic growth.} Let $V=\Psi_{x_1}-i\Psi_{x_2}$ and $G=|p_0|V/(\alpha-ze^{-i\theta_c})$. The denominator vanishes only at $z=P_c$. On the free surface, $|G|=|\alpha-S|/|\alpha-Se^{i(\theta-\theta_c)}|\le1$. On the bottom, $|G|\le|p_0|\max|\Psi_R(1,\cdot)|/(\alpha-1)$. Near $P_c$, $|G|\le C|z-P_c|^{-1}$. By \autoref{lem_PL}, $|G|\le M_1$ in $\Omega$, that is,
\[
|\nabla u(x)|\le M_1|x-P_c|\qquad(x\in\Omega).
\]
This uses no information on the sign of $S-\alpha$. Integrating along the segment from $x$ to its nearest boundary point gives $u\le 2M_1|x-P_c|^2$. Hence $\int_{B_r}|\nabla u|^2\le Cr^4$.

\emph{Step 2: Weiss formula of homogeneity $2$.} Let
\[
\Phi_2(r)=r^{-4}\int_{B_r}\big(|\nabla u|^2+q\chi\big)-2r^{-5}\int_{\partial B_r}u^2 .
\]
The computation in the proof of \autoref{lem_weiss}, with the exponent $\kappa=3/2$ replaced by $\kappa=2$, gives
\[
\Phi_2'(r)=2r^{-4}\int_{\partial B_r}\Big(\partial_\nu u-\frac{2u}r\Big)^2+r^{-5}\int_{B_r}\chi\,(y\cdot\nabla q-2q).
\]
Since $y\cdot\nabla q-2q=O(|y|^3)$, the last term is bounded. By Step 1, $\Phi_2$ is bounded.

\emph{Step 3: non-degeneracy and angle bounds.} The flux argument of \autoref{lem_nondeg}, with $|\nabla u|=|d|$ in place of $q^{1/2}$, gives
\[
\sup_{B_{4\rho}(z)}u\ge c\rho|z|
\]
for $z\in\partial\Omega$ with $d(z)\ge\sigma|z|$ and $\rho=\sigma|z|/4$. As in \autoref{thm_angle_bounds}, comparison with $\varrho^{\pi/\beta}\sin(\pi\varphi/\beta)$ shows two things. First, $\Omega$ contains no truncated sector at $P_c$ of opening larger than $\pi/2$, because $u\le Cr^2$. Second, $\Omega$ is not contained, near $P_c$, in a sector of opening smaller than $\pi/2$ inside $\{y_2\ge\sigma|y|\}$, by non-degeneracy.

\emph{Step 4: blow-ups.} The rescalings $u_r(y)=r^{-2}u(P_c+ry)$ satisfy $|\nabla u_r(y)|\le M_1|y|$. As in \autoref{prop_blowup}, every sequence $r_j\to0$ has a subsequence along which $u_{r_j}\to u_0$ locally uniformly, $\nabla u_{r_j}\to\nabla u_0$ in $L^2_{\rm loc}$ and $\chi_{r_j}\rightharpoonup\chi_0$ weakly-$*$. By Step 2, $u_0$ is homogeneous of degree $2$ and
\[
\Phi_2(0+)=\int_{B_1}y_2^2\,\chi_0 .
\]
By \autoref{lem_phase}, the positivity set of $u_0$ and the density $\chi_0$ are non-decreasing under $y\mapsto y+te_2$, $t\ge0$. The argument of \autoref{lem_jump}, with $q_0=y_2^2$, shows the following. First, $\chi_0$ is constant on each component of $\{u_0=0\}^\circ\setminus\{y_2=0\}$. Second, on a ray bounding the positivity set and adjacent to such a component of density $c$, $|\nabla u_0|^2=(1-c)y_2^2$. Third, on a ray with positivity on both sides, the two one-sided gradients have equal length.

\emph{Step 5: classification.} On each component of $\{u_0>0\}$, $u_0=A\varrho^2\sin(2(\varphi-a))$ on a sector $(a,a+\pi/2)$.

Suppose two such sectors share a ray $\Gamma$. If $\Gamma$ is not vertical, a point $z\in\Gamma$ with $|z|=1$ has, for small $r$, fluid of $\Omega_r$ on both sides and directly above it (in the $-e_2$ direction). By \autoref{lem_phase}, $\Omega_r$ then contains a neighbourhood of $z$, so $u_0$ is harmonic near $z$, non-negative, and zero at $z$, hence identically zero, which is absurd. If $\Gamma$ is vertical, the invariance under $+e_2$ forces the union of the sectors to be bounded by horizontal rays. On these rays $y_2=0$, while $|\nabla u_0|\ne0$, contradicting the jump condition.

So $\{u_0>0\}$ is either empty or a single sector $(a,a+\pi/2)$ which, by the $+e_2$ invariance, contains the direction $e_2$ in its closure, i.e.\ $a\in[0,\pi/2]$. The jump condition excludes $a\in\{0,\pi/2\}$.

For a single sector, the adjacent zero components inside $\{y_2>0\}$ are air, by the argument of Step 2 of the proof of \autoref{prop_classification}: a boundary point near a bounding ray produces an air point, whose upward segment crosses the region in question. Then
\[
4A^2\varrho^2=\varrho^2\sin^2a=\varrho^2\cos^2a
\]
on the two bounding rays, so $a=\pi/4$ and $A=\sqrt2/4$, i.e.\ $u_0=\frac{\sqrt2}4(y_2^2-y_1^2)$ on $\{|y_1|<y_2\}$. Its Weiss density is
\[
\int_{B_1}y_2^2\chi_{\{|y_1|<y_2\}}=\frac\pi{16}+\frac18 .
\]

If $u_0\equiv0$, then $\chi_0=c_\pm$ on $\{\pm y_2>0\}$ with $c_\pm\in\{0,1\}$, by the stabilization argument of Step 4 of that proof, so $\Phi_2(0+)\in\{0,\pi/8,\pi/4\}$. Each case is excluded.
\begin{itemize}
\item $c_+=c_-=0$ contradicts $P_c-s\mathbf e_R\in\Omega$ for small $s>0$.
\item $c_+=1$ makes $\Omega$ contain a truncated sector of opening $\pi-2\delta>\pi/2$, contradicting Step 3.
\item $c_+=0$, $c_-=1$ contradicts the inward radial closure of the fluid.
\end{itemize}
Since $\Phi_2(0+)$ does not depend on the sequence, every blow-up equals $u_*=\frac{\sqrt2}4(y_2^2-y_1^2)\chi_{\{|y_1|<y_2\}}$. This proves (c), since $\Psi=u/|p_0|$.

\emph{Step 6: conclusions.} The tangent-cone statement (b) follows from (c) exactly as in \autoref{thm_main_stokes}(b). In particular, $S<\alpha$ in a punctured neighbourhood of $\theta_c$, which gives (d).

For (a), suppose $S(\theta_1)>\alpha$ for some $\theta_1$. Since $\min S<\alpha$ by \autoref{thm_localization}(b), there is $\theta_\times$ with $S(\theta_\times)=\alpha$ and points arbitrarily close to $\theta_\times$ where $S>\alpha$. The point $(\alpha,\theta_\times)$ is a stagnation point, since $Q_\alpha(\alpha)=0$. This contradicts what was just shown at that point. Hence $S\le\alpha$, and every stagnation point is a global maximum.
\end{proof}

\begin{remark}
The amplitude $\sqrt2/(4|p_0|)$ does not involve $\Delta$, unlike $A_*=2\Delta^{1/4}/(3|p_0|)$ in \autoref{thm_intro_stokes}. As $\alpha\to\alpha_+$ the Stokes amplitude $A_*\to0$, consistent with the change of homogeneity from $3/2$ to $2$. Right-angle corners also appear in other Bernoulli problems whose gravity term vanishes on the free surface; here the degeneracy comes from the double root of $Q_\alpha$.
\end{remark}

\subsection{Crests without the graph structure}\label{subsec_variational_crest}

The proofs of \autoref{thm_main_stokes} and \autoref{thm_degenerate} use the radial-graph structure through \autoref{lem_phase}. This property fails if the waves overturn, which is alternative (II) of \autoref{thm_intro_global}. In this subsection we show that the crest analysis survives without it.

For classical solutions the growth bound needs no change, because the Phragm\'en--Lindel\"of argument of \autoref{prop_interior_bernoulli} never used the graph property. For the rest we adapt three recent ideas from the flat problem; the arguments there concern the weight $x_2$ in the plane, so we carry them out in full for the radial weight $Q_\alpha(|x|)$.
\begin{itemize}
\item The Hopf-differential method of De Philippis, Hirsch and Nahon \cite{DHN2026}, as used by Du, Yang and Zhang for flat gravity waves \cite{DYZ2026}, gives the growth bound for \emph{every} variational solution. Only the abstract approximation theorem of \cite{DHN2026}, which is independent of the domain and of the weight, is used as a black box.
\item The weighted isoperimetric argument of Du, Pu and Yang \cite{DPY2026} gives a positive lower bound for the Weiss density and so excludes cusps; for an earlier two-dimensional argument see \cite{McCurdy2025}.
\item The flux non-degeneracy of \autoref{lem_nondeg} excludes full density for regular interfaces.
\end{itemize}

Throughout, $\Delta(\alpha)\ge0$ and $r_*=r_-(\alpha)$, so that $r_*=\alpha$ when $\Delta=0$. We set
\[
\kappa=\tfrac32,\ \ g=2\sqrt\Delta\quad(\Delta>0),\qquad\qquad \kappa=2,\ \ g=1\quad(\Delta=0).
\]

\begin{definition}[Overhanging classical solutions]\label{def_general_classical}
A pair $(\Psi,\Omega)$ is a \emph{classical solution with general free surface} if the following hold:
\begin{itemize}
\item $\Omega$ is bounded by the unit circle and a Jordan curve $\Gamma$ surrounding it;
\item $\Gamma$ is $C^2$ away from $\Gamma\cap\{|x|=r_*\}$;
\item $\Psi\in C^2(\Omega)\cap C^1(\overline\Omega)$ is harmonic, with $\Psi=1$ on $|x|=1$, $\Psi=0$ on $\Gamma$, and $p_0^2|\nabla\Psi|^2=Q_\alpha(|x|)$ on $\Gamma$.
\end{itemize}
It lies in the \emph{lower region} if $\Gamma\subset\{|x|\le r_*\}$.
\end{definition}

Radial subgraphs $\Gamma=\{R=S(\theta)\}$ are a special case. As in \autoref{lem_dichotomy}, if $\Delta>0$ then $Q_\alpha(|x|)\ge0$ on the connected set $\Gamma$, so either $\Gamma\subset\{|x|\le r_-\}$ or $\Gamma\subset\{|x|\ge r_+\}$. In the lower region, every stagnation point lies on the circle $|x|=r_*$.

\begin{lemma}[Interior Bernoulli inequality without graph structure]\label{lem_general_bernoulli}
Let $(\Psi,\Omega)$ be a classical solution with general free surface in the lower region, and let
\[
M=\max\Big\{1,\ \frac{p_0^2\max_{|x|=1}|\nabla\Psi|^2}{Q_\alpha(1)}\Big\}.
\]
Then $p_0^2|\nabla\Psi(x)|^2\le M\,Q_\alpha(|x|)$ for all $x\in\Omega$.
\end{lemma}
\begin{proof}
The proof of \autoref{prop_interior_bernoulli} applies verbatim. It uses three facts. First, $|\mathcal P(ze^{-i\theta_0})|\ge(r_--|z|)(r_+-|z|)=Q_\alpha(|z|)$ for $|z|\le r_-$; this holds for every such $z$, because $|r-z|\ge r-|z|$ for real $r\ge|z|$. Second, the Bernoulli condition on $\Gamma$. Third, $\Gamma$ is a connected set of diameter at least $2$, as required in \autoref{lem_PL}. For $\Delta=0$ one takes $\mathcal P(\zeta)=(\alpha-\zeta)^2$. None of these facts uses the graph property.
\end{proof}

\begin{definition}[Variational solutions]\label{def_variational}
Let $D\subset\{|x|>1\}$ be open. A function $\Psi$ is a \emph{variational solution} of \eqref{eq_model} in $D$ at the parameter $\alpha$ if the following hold:
\begin{itemize}
\item $\Psi\in W^{1,2}(D)\cap C^0(D)$, $\Psi\ge0$, $\Psi\in C^2(D\cap\{\Psi>0\})$, and $\Psi=0$ on $D\cap\{|x|\ge r_*\}$;
\item for every $\phi\in C^1_c(D;\mathbb R^2)$, with $\chi=\chi_{\{\Psi>0\}}$,
\begin{equation}\label{eq_VS}
\int_D\Big[\big(p_0^2|\nabla\Psi|^2+Q_\alpha(|x|)\chi\big)\operatorname{div}\phi-2p_0^2\,\nabla\Psi\cdot D\phi\,\nabla\Psi+\chi\,\nabla\big(Q_\alpha(|x|)\big)\cdot\phi\Big]dx=0 .
\end{equation}
\end{itemize}
Put $\Gamma=\partial\{\Psi>0\}\cap D\cap\{|x|<r_*\}$. We say that the solution satisfies (G) if the first condition below holds, and that it is \emph{regular} if it satisfies all three:
\begin{enumerate}
\item[(G)] $p_0^2|\nabla\Psi|^2\le M\,Q_\alpha(|x|)$ a.e.\ in $D$, for some constant $M$;
\item[(H1)] every point of $\Gamma$ has a neighbourhood in which $\Gamma$ is a $C^2$ arc with $\{\Psi>0\}$ on one side, $\Psi$ is $C^1$ up to the arc from that side, and $p_0^2|\nabla\Psi|^2=Q_\alpha(|x|)$ on the arc;
\item[(H2)] no component of $\Gamma$ is a closed curve contained in $D$.
\end{enumerate}
\end{definition}

A closed curve in $\Gamma$ cannot bound fluid: $\Psi$ would be harmonic inside it and vanish on it, hence vanish identically. So (H2) only forbids air bubbles.

A classical solution with general free surface in the lower region is a regular variational solution in every small ball $D$ centred on the circle $|x|=r_*$. Identity \eqref{eq_VS} follows by the integration by parts in the proof of \autoref{lem_domain_variation}, which uses only harmonicity, the Bernoulli condition and the regularity of $\Gamma$ off the circle. Condition (G) is \autoref{lem_general_bernoulli}, and (H1) is part of \autoref{def_general_classical}. Condition (H2) holds because $\Gamma$ is a single Jordan curve surrounding the unit disc, which is not contained in a small ball.

\begin{proposition}[Growth for variational solutions]\label{prop_variational_growth}
Let $|x_0|=r_*$ and $0<\rho\le\min\{1,(r_*-1)/2\}$, and let $\Psi$ be a variational solution in $B_\rho(x_0)$ with $\int_{B_\rho(x_0)}|\nabla\Psi|^2\le E$. Then there is $C=C(\alpha,p_0,\rho,E)$ such that
\[
p_0^2|\nabla\Psi(x)|^2\le C\,Q_\alpha(|x|)\qquad\text{for }x\in B_{\rho/32}(x_0)\cap\{\Psi>0\}.
\]
In particular, $\Psi$ satisfies (G) in $B_{\rho/32}(x_0)$.
\end{proposition}

\begin{proof}
Put $u=|p_0|\Psi$, $q(x)=Q_\alpha(|x|)$ and $\chi^0=1-\chi$. Identify $\mathbb R^2$ with $\mathbb C$, with $\partial_z=\frac12(\partial_1-i\partial_2)$ and $\partial_{\bar z}=\frac12(\partial_1+i\partial_2)$. Then \eqref{eq_VS} reads $\operatorname{div}T=\chi\nabla q$ in $\mathcal D'$, where
\[
T=(|\nabla u|^2+q\chi)I-2\nabla u\otimes\nabla u .
\]

\emph{Step 1: the Hopf differential.} Let $\Phi=(u_{x_1}-iu_{x_2})^2=A-iB$, with $A=u_{x_1}^2-u_{x_2}^2$ and $B=2u_{x_1}u_{x_2}$. Then $|\Phi|=|\nabla u|^2$, and $\Phi=0$ a.e.\ on $\{u=0\}$. The two components of $\operatorname{div}T=\chi\nabla q$ read
\[
A_{x_1}+B_{x_2}=q\chi_{x_1},\qquad A_{x_2}-B_{x_1}=-q\chi_{x_2},
\]
by the product rule for the smooth function $q$. Hence
\begin{equation}\label{eq_hopf}
\partial_{\bar z}\Phi=\tfrac12(A_{x_1}+B_{x_2})+\tfrac i2(A_{x_2}-B_{x_1})=q\,\partial_z\chi=-q\,\partial_z\chi^0 .
\end{equation}
For the flat weight $q=x_2$ this is the equation of \cite[Lemma~3.1]{DYZ2026}.

\emph{Step 2: holomorphic approximation.} Let $\mathfrak Tg(z)=\frac1\pi\int g(\zeta)(z-\zeta)^{-1}d\zeta$ and $\mathfrak B=\partial_z\mathfrak T$ be the Cauchy and Beurling transforms. Then $\partial_{\bar z}\mathfrak Tg=g$ and $\partial_{\bar z}\mathfrak Bg=\partial_zg$. Moreover, $\mathfrak B$ is an isometry of $L^2(\mathbb R^2)$, and by Young's inequality $\|\mathfrak Tg\|_{L^2(B_s(b))}\le4s\|g\|_{L^2}$ when $\operatorname{supp}g\subset B_s(b)$ \cite{AIM2009}. For a disc $B=B_s(b)\subset B_\rho(x_0)$ put
\[
E_B=\mathfrak B\big(q\chi^0\mathbf 1_B\big)-\mathfrak T\big((\partial_zq)\chi^0\mathbf 1_B\big),\qquad h_B=\Phi+E_B .
\]
Then $\partial_{\bar z}E_B=\partial_z(q\chi^0)-(\partial_zq)\chi^0=q\,\partial_z\chi^0$ in $B$, so $\partial_{\bar z}h_B=0$ by \eqref{eq_hopf}. By Weyl's lemma, $h_B$ is holomorphic in $B$. Since $|\partial_zq|=\frac12|\nabla q|$,
\begin{equation}\label{eq_approx_error}
\|\Phi-h_B\|_{L^2(B)}\le\Lambda_B\,|B\cap\{u=0\}|^{1/2},\qquad \Lambda_B=\sup_B|q|+2s\sup_B|\nabla q| .
\end{equation}
Take $B=B_\rho(x_0)$. Then $\|h_B\|_{L^1(B)}\le p_0^2E+C$. The mean-value inequality for $|h_B|$ on discs of radius $\rho/8$, together with $\Phi=h_B-E_B$, gives
\begin{equation}\label{eq_Phi_L2}
\|\Phi\|_{L^2(B_{3\rho/4}(x_0))}\le C(\alpha,\rho)\,(1+p_0^2E).
\end{equation}

\emph{Step 3: decay at the stagnation circle.} Let $|a|=r_*$, $n=a/|a|$, and $0<r\le\rho/4$ with $B_r(a)\subset B_{3\rho/4}(x_0)$. On $B_r(a)$ we have $|r_*-|x||\le r$. Since $q=(r_--|x|)(r_+-|x|)$ and $|\nabla q|=2\big||x|-\alpha\big|$, we get
\[
\Lambda_{B_s(b)}\le C_\alpha r^{2\kappa-2}\qquad\text{for all }B_s(b)\subset B_r(a).
\]
Indeed, for $\Delta>0$ we have $|q|\le r(2\sqrt\Delta+r)$ and $|\nabla q|\le2(\sqrt\Delta+r)$, while for $\Delta=0$ we have $|q|\le r^2$ and $|\nabla q|\le2r$. Hence $\|\Phi-h_{B_r(a)}\|_{L^2(B_r(a))}\le Cr^{2\kappa-1}$.

The disc $D_r=B_{r/4}(a+\frac r2n)$ lies in $B_r(a)\cap\{|x|>r_*\}$, where $\Phi=0$. We use a unique-continuation inequality: for every $\varepsilon\in(0,1)$ there is $C_\varepsilon$ such that, for every unit vector $n$ and every holomorphic $h\in L^2(B_1)$,
\[
\|h\|_{L^2(B_{1/2})}\le\varepsilon\|h\|_{L^2(B_1)}+C_\varepsilon\|h\|_{L^2(B_{1/4}(n/2))}.
\]
If this failed, there would be holomorphic $h_j$ and unit vectors $n_j\to n$ with $\|h_j\|_{L^2(B_{1/2})}=1$, $\|h_j\|_{L^2(B_1)}\le\varepsilon^{-1}$ and $\|h_j\|_{L^2(B_{1/4}(n_j/2))}\to0$. By the mean-value inequality, a subsequence converges locally uniformly to a holomorphic $h$ on $B_1$. This $h$ vanishes on $B_{1/5}(n/2)$, hence identically, contradicting $\|h\|_{L^2(B_{1/2})}=1$.

Scale this to $B_r(a)$ and apply it to $h=h_{B_r(a)}$. Since $\|h\|_{L^2(D_r)}=\|h-\Phi\|_{L^2(D_r)}\le Cr^{2\kappa-1}$, we obtain
\[
\|\Phi\|_{L^2(B_{r/2}(a))}\le\varepsilon\|\Phi\|_{L^2(B_r(a))}+C_\varepsilon r^{2\kappa-1}.
\]
Take $\varepsilon=2^{-2\kappa}$. Then $N(r)=r^{1-2\kappa}\|\Phi\|_{L^2(B_r(a))}$ satisfies $N(r/2)\le\frac12N(r)+C$. Iterating from $r_0=\rho/4$ and using \eqref{eq_Phi_L2},
\begin{equation}\label{eq_decay}
\|\Phi\|_{L^2(B_r(a))}\le C(1+E)\,r^{2\kappa-1}\qquad(0<r\le\rho/4).
\end{equation}
For the flat weight this is the decay $r^2$ of \cite[Lemma~5.1]{DYZ2026}. For $\Delta=0$ it is one order better, because the weight vanishes quadratically.

\emph{Step 4: the pointwise bound.} The approximation theorem of De Philippis, Hirsch and Nahon \cite[Theorem~1]{DHN2026}, after scaling to $B_r(a)$ \cite[Theorem~5.4 and Appendix~A]{DYZ2026}, gives the following. Suppose every disc $B_s(b)\subset B_r(a)$ carries an approximant $h_{B_s(b)}$ with
\[
\|W-h_{B_s(b)}\|_{L^2(B_s(b))}^2\le\Lambda^2|B_s(b)\cap\{W=0\}| ,
\]
and with the interior gradient estimate valid for holomorphic functions. Then
\[
\|W\|_{L^\infty(B_{r/2}(a))}\le C\big(\Lambda+r^{-1}\|W\|_{L^2(B_r(a))}\big).
\]
Apply this to $W=\Phi$ with $\Lambda=Cr^{2\kappa-2}$, using \eqref{eq_approx_error}, $|B\cap\{u=0\}|\le|B\cap\{\Phi=0\}|$ and \eqref{eq_decay}. We get $\|\Phi\|_{L^\infty(B_{r/2}(a))}\le C(1+E)r^{2\kappa-2}$.

Now let $x\in B_{\rho/32}(x_0)\cap\{u>0\}$. Then $\delta=r_*-|x|>0$ and $\delta\le|x-x_0|<\rho/32$. Put $a=r_*x/|x|$ and $r=4\delta$. Then $x\in B_{r/2}(a)$ and $B_{\rho/4}(a)\subset B_{3\rho/4}(x_0)$, so
\[
|\nabla u(x)|^2\le C\delta^{2\kappa-2}\le C'\,Q_\alpha(|x|),
\]
because $Q_\alpha(|x|)\ge2\sqrt\Delta\,\delta$ if $\Delta>0$ and $Q_\alpha(|x|)=\delta^2$ if $\Delta=0$.
\end{proof}

\begin{theorem}[Crests of variational solutions]\label{thm_variational_crest}
Let $|x_0|=r_*$ and $0<\rho\le\min\{1,(r_*-1)/2\}$. Let $\Psi$ be a variational solution in $B_\rho(x_0)$ satisfying (G) with constant $M$; by \autoref{prop_variational_growth}, this holds on a smaller ball for every variational solution. Assume $x_0\in\partial\{\Psi>0\}$, and, in the frame $e_1=\mathbf e_\theta(\theta_0)$, $e_2=-\mathbf e_R(\theta_0)$, put $\Psi_r(y)=r^{-\kappa}\Psi(x_0+ry)$ and $\chi_r(y)=\chi(x_0+ry)$.
\begin{enumerate}
\item \emph{(No cusps.)} There are $c,r_1>0$, depending only on $\alpha$, $p_0$ and $M$, such that
\[
r^{-2\kappa}\int_{B_r(x_0)}Q_\alpha(|x|)\,\chi\,dx\ge c\qquad(0<r<r_1).
\]
\item \emph{(Blow-up alternative.)} As $r\to0$, the whole family $(\Psi_r,\chi_r)$ converges, locally uniformly and in $W^{1,2}_{\rm loc}\times L^1_{\rm loc}$, to exactly one of:
\begin{enumerate}
\item the corner: for $\Delta>0$ the Stokes corner of \autoref{thm_intro_stokes}(c), with amplitude $2\Delta^{1/4}/(3|p_0|)$ and phase $\chi_{\{\pi/6<\varphi<5\pi/6\}}$; for $\Delta=0$ the right-angle corner $\frac{\sqrt2}{4|p_0|}(y_2^2-y_1^2)$ with phase $\chi_{\{|y_1|<y_2\}}$;
\item full density: $(0,\chi_{\{y_2>0\}})$.
\end{enumerate}
\item \emph{(Exclusion of full density.)} Alternative (ii) does not occur, so the blow-up is the corner, under either of the following hypotheses:
\begin{enumerate}
\item[(R)] $\Psi$ is regular;
\item[(F)] $\{\Psi=0\}\cap B_\rho(x_0)$ has finitely many connected components, and the set of points of $\Gamma$ that are not regular in the sense of (H1) has $\mathcal H^1$-measure zero.
\end{enumerate}
\end{enumerate}
Hypothesis (F) is the annular counterpart of the hypotheses of \cite[Theorem~A]{VW2011} and \cite[Theorem~II(ii)]{DYZ2026}. In the flat case, full density is excluded under it by means of a frequency formula. Here the angle bound of \autoref{thm_angle_bounds} takes the place of the frequency formula.
\end{theorem}

\begin{proof}
Put $u=|p_0|\Psi$, $q(x)=Q_\alpha(|x|)$ and $d(x)=r_*-|x|$. Then (G) reads $|\nabla u|^2\le Mq$. On $\{d\ge0\}$ near $x_0$ we have $q=gd+d^2$ if $\Delta>0$ and $q=d^2$ if $\Delta=0$, so $q\le Cd^{2\kappa-2}\le C|x-x_0|^{2\kappa-2}$.

\emph{Step 1: Lipschitz bound, energy identity and the Laplacian.} By (G), $u$ is Lipschitz on $B=B_\rho(x_0)$, with $|\nabla u|\le C|x-x_0|^{\kappa-1}$, and $\nabla u=0$ a.e.\ on $\{u=0\}$. Hence $u(x)\le C|x-x_0|^\kappa$.

For $\zeta\in C^\infty_c(B)$, testing harmonicity with $(u-\varepsilon)^+\zeta$ and letting $\varepsilon\to0$ gives
\begin{equation}\label{eq_energy_identity}
\int\zeta|\nabla u|^2=-\int u\,\nabla u\cdot\nabla\zeta .
\end{equation}
Fix a smooth non-decreasing $\eta$ with $\eta=0$ on $(-\infty,1]$ and $\eta=1$ on $[2,\infty)$, and put $\eta_\varepsilon(s)=\eta(s/\varepsilon)$. Testing harmonicity with $\eta_\varepsilon(u)\zeta$, $\zeta\ge0$, gives
\[
\int\zeta\,d\mu_\varepsilon=-\int\eta_\varepsilon(u)\nabla u\cdot\nabla\zeta,\qquad\mu_\varepsilon=\eta_\varepsilon'(u)|\nabla u|^2dx\ge0 .
\]
Letting $\varepsilon\to0$ shows that $\mu=\Delta u\ge0$ is a Radon measure carried by $\partial\{u>0\}$, with $\mu_\varepsilon\rightharpoonup\mu$. Testing with radial cut-offs gives, for a.e.\ $t$,
\begin{equation}\label{eq_mu_flux}
\mu(B_t(x_0))\le\int_{\partial B_t(x_0)}|\nabla u|\,d\mathcal H^1\le\sqrt M\,V'(t),\qquad V(t)=\int_{B_t(x_0)}\sqrt{q_+}\,\chi\,dx .
\end{equation}
In particular $\mu(B_t(x_0))\le Ct^\kappa$.

\emph{Step 2: weighted perimeter.} Let $\phi\in C^1_c(B\cap\{q>0\};\mathbb R^2)$, and put $X=|\nabla u|^2\phi-2(\nabla u\cdot\phi)\nabla u$. In $\{u>0\}$, harmonicity gives $\operatorname{div}X=|\nabla u|^2\operatorname{div}\phi-2\nabla u\cdot D\phi\,\nabla u$. Integration by parts gives
\[
\int\eta_\varepsilon(u)\operatorname{div}X=-\int\eta_\varepsilon'(u)\nabla u\cdot X=\int(\nabla u\cdot\phi)\,d\mu_\varepsilon .
\]
As $\varepsilon\to0$, the left-hand side tends to $\int(|\nabla u|^2\operatorname{div}\phi-2\nabla u\cdot D\phi\,\nabla u)$. By \eqref{eq_VS}, this equals $-\int\chi\operatorname{div}(q\phi)$. Since $|\nabla u|\le\sqrt{Mq}$ on the support of $\mu_\varepsilon$, it follows that $\big|\int\chi\operatorname{div}(q\phi)\big|\le\sqrt M\int\sqrt q\,|\phi|\,d\mu$. With $\phi=\zeta/q$ we obtain $\chi\in BV_{\rm loc}(B\cap\{q>0\})$ and
\begin{equation}\label{eq_weighted_perimeter}
\sqrt q\,|D\chi|\le\sqrt M\,\mu\qquad\text{as measures in }B\cap\{q>0\}.
\end{equation}
For the flat weight these are \cite[Lemma~6.1]{DYZ2026} and \cite{DPY2026}.

\emph{Step 3: proof of (a).} We follow the isoperimetric argument of \cite{DPY2026}, which there concerns the half-plane with weight $\sqrt{x_2}$. Let $P(t)=\int_{B_t(x_0)\cap\{q>0\}}\sqrt q\,d|D\chi|$. Then $P(t)\le MV'(t)$ for a.e.\ $t$, by \eqref{eq_weighted_perimeter} and \eqref{eq_mu_flux}.

To straighten the stagnation circle, put
\[
G(x)=\big(r_*(\theta(x)-\theta_0),\ y_2(x)\big),\qquad y_2=q/g\ (\Delta>0),\qquad y_2=d\ (\Delta=0).
\]
This is a smooth diffeomorphism near $x_0$ with orthogonal differential at $x_0$. It maps $\{q>0\}$ into $\{y_2>0\}$, and $\sqrt q=\sqrt g\,y_2^{\kappa-1}$ exactly. On $B_t(x_0)$, $G$ and $G^{-1}$ are $(1+Ct)$-Lipschitz and $|\det DG|\le1+Ct$.

The isoperimetric inequality for the monomial weight $y_2^{\kappa-1}$ in the half-plane \cite[Theorem~1.4]{CabreRosOton2013}, \cite[Theorem~1.3]{CRS2016} has exponent $\kappa/(\kappa+1)$, since the weighted dimension is $\kappa+1$. It therefore transfers to sets $A\subset B_t(x_0)\cap\{q>0\}$ of finite perimeter, in the form
\[
\int_{\partial^*A}\sqrt q\,d\mathcal H^1\ge c_\kappa(1-Ct)\Big(\int_A\sqrt q\Big)^{\kappa/(\kappa+1)} .
\]
Apply this to $A_\varepsilon=B_t(x_0)\cap\{u>0\}\cap\{q>\varepsilon\}$. By the slicing properties of $BV$ functions \cite[Sections~5.3--5.5]{EvansGariepy2015}, its weighted perimeter is at most $P(t)+V'(t)+Ct\sqrt\varepsilon$ for a.e.\ $t$ and $\varepsilon$. Letting $\varepsilon\to0$, we get
\[
P(t)+V'(t)\ge c_\kappa(1-Ct)V(t)^{\kappa/(\kappa+1)},\qquad\text{hence}\qquad V'(t)\ge\frac{c_\kappa}{2(1+M)}V(t)^{\kappa/(\kappa+1)}
\]
for a.e.\ small $t$. Since $x_0\in\partial\{u>0\}$, we have $V(t)>0$ for $t>0$. Integrating $\frac d{dt}V^{1/(\kappa+1)}\ge c$ from $0$ gives $V(t)\ge ct^{\kappa+1}$. By the Cauchy--Schwarz inequality, $\int_{B_t(x_0)}q\chi\ge V(t)^2/|B_t|\ge ct^{2\kappa}$.

\emph{Step 4: blow-ups.} Put $u_r(y)=r^{-\kappa}u(x_0+ry)$ and $q_r(y)=r^{2-2\kappa}q(x_0+ry)$. Then $q_r\to g\,y_2^{2\kappa-2}$ in $C^1_{\rm loc}$. The $u_r$ are locally uniformly Lipschitz by Step~1. Rescaling \eqref{eq_weighted_perimeter} and $\mu(B_t)\le Ct^\kappa$ bounds the total variation of $\chi_r$ on compact subsets of $\{y_2>0\}$. Moreover $\chi_r=0$ on $\{q_r\le0\}$, a set that converges to $\{y_2\le0\}$. So every sequence $r_j\to0$ has a subsequence with the following properties: $u_{r_j}\to u_0$ locally uniformly; $\nabla u_{r_j}\to\nabla u_0$ in $L^2_{\rm loc}$, by \eqref{eq_energy_identity} as in \autoref{prop_blowup}; and $\chi_{r_j}\to\chi_0$ in $L^1_{\rm loc}$, where $\chi_0$ is a characteristic function vanishing on $\{y_2<0\}$. The pair $(u_0,\chi_0)$ satisfies the domain-variation identity with weight $g\,y_2^{2\kappa-2}$, and $\chi_0=1$ on $\{u_0>0\}$.

Let
\[
\Phi_\kappa(r)=r^{-2\kappa}\int_{B_r(x_0)}(|\nabla u|^2+q\chi)-\kappa\,r^{-2\kappa-1}\int_{\partial B_r(x_0)}u^2 .
\]
The computation in the proof of \autoref{lem_weiss}, with $\phi=(x-x_0)\eta_k(|x-x_0|)$ in \eqref{eq_VS} and with \eqref{eq_energy_identity}, gives
\[
\Phi_\kappa'(r)=2r^{-2\kappa}\int_{\partial B_r}\Big(\partial_\nu u-\frac{\kappa u}r\Big)^2+r^{-2\kappa-1}\int_{B_r}\chi\,\big(y\cdot\nabla q-(2\kappa-2)q\big),\qquad y=x-x_0 .
\]
Since $y\cdot\nabla q-(2\kappa-2)q=O(|y|^{2\kappa-1})$, the last term is bounded by a constant $C_W$, uniformly for centres on the circle. As in \autoref{prop_blowup}, $u_0$ is homogeneous of degree $\kappa$, and
\[
\Phi_\kappa(0+)=g\int_{B_1}y_2^{2\kappa-2}\chi_0 .
\]

\emph{Step 5: proof of (b).} Each component of $\{u_0>0\}$ is a sector of opening $\pi/\kappa$ on which $u_0=A\varrho^\kappa\cos(\kappa(\varphi-m))$, where $m$ is the bisecting angle, and $|\nabla u_0|=\kappa A\varrho^{\kappa-1}\ne0$ on its bounding rays. The argument of \autoref{lem_jump}, with weight $g\,y_2^{2\kappa-2}$, gives two facts. On a bounding ray in $\{y_2=0\}$, $|\nabla u_0|=0$. On a bounding ray adjacent to a zero component on which $\chi_0=c\in\{0,1\}$, $|\nabla u_0|^2=(1-c)g\,y_2^{2\kappa-2}$.

Two positive sectors would need both rays $\varphi=0$ and $\varphi=\pi$. This is impossible for $\kappa=3/2$ and contradicts the first fact for $\kappa=2$. So $\{u_0>0\}$ is empty, or a single sector $(a,a+\pi/\kappa)$ with $0<a<\pi-\pi/\kappa$. In the second case, $c=1$ on an adjacent zero component would give $|\nabla u_0|=0$ on a bounding ray, so $c=0$ on both. Then $\kappa^2A^2\varrho^{2\kappa-2}=g\,y_2^{2\kappa-2}$ on both rays. For $\kappa=3/2$ this means $a=\pi/6$ and $A^2=2g/9$; for $\kappa=2$ it means $a=\pi/4$ and $A^2=1/8$. Dividing by $|p_0|$ gives (i).

If $u_0\equiv0$, then $\chi_0\equiv c\in\{0,1\}$ on $\{y_2>0\}$. Passing to the limit in (a) gives $c=1$, which is (ii).

The values of $\Phi_\kappa(0+)/g$ in the two cases are $\sqrt3/3$ and $2/3$ for $\kappa=3/2$, and $\pi/16+1/8$ and $\pi/8$ for $\kappa=2$. They are distinct, so all subsequential limits coincide.

\emph{Step 6: non-degeneracy for regular solutions.} Let $\Psi$ be regular. Let $z\in\Gamma$ with $d(z)\ge\sigma|z-x_0|>0$, and let $\rho'\le\sigma|z-x_0|/4$ with $B_{4\rho'}(z)\subset B$. By (H1), $\Gamma\cap\overline{B_{2\rho'}(z)}$ is a finite union of $C^2$ arcs and closed curves. By (H2), the arc through $z$ leaves $B_{\rho'}(z)$, so its length in $B_{\rho'}(z)$ is at least $\rho'$. On it, $|\nabla u|=\sqrt q\ge c\,d(z)^{\kappa-1}$. The divergence theorem on $\{u>0\}\cap B_t(z)$, $t\in(\rho',2\rho')$, applies because $u\in C^1$ up to the $C^2$ arcs. Averaging in $t$ and using Caccioppoli's inequality for the subharmonic extension of $u$, as in \autoref{lem_nondeg}, we get
\begin{equation}\label{eq_nondeg_var}
\sup_{B_{4\rho'}(z)}u\ge c\,\rho'\,d(z)^{\kappa-1}.
\end{equation}
In rescaled form: if $z_r\in\Gamma_r=r^{-1}(\Gamma-x_0)$ lies in a fixed compact set $K\subset\{y_2>0\}\setminus\{0\}$, then $\sup u_r\ge c(K)>0$ on a fixed compact neighbourhood of $K$.

\emph{Step 7: proof of (c) under (R).} Suppose (ii) holds for a regular solution. Fix $\delta>0$ and $K=\{1\le|y|\le2,\ \delta\le\varphi\le\pi-\delta\}$. Since $u_r\to0$, Step~6 gives $r_1$ such that $\Gamma_r\cap K=\emptyset$ for $r<r_1$. So $K$ lies either in $\{u_r>0\}$ or in $\{u_r=0\}$. The sets $rK$ and $r'K$ overlap for $r'\in(r/2,r)$, so the alternative is the same for all $r<r_1$, and by $\chi_r\to\chi_{\{y_2>0\}}$ it is the first. Hence $\{u>0\}$ contains a truncated sector of opening $\pi-2\delta$ at $x_0$.

For $\kappa=3/2$ and $\pi-2\delta>2\pi/3$, the barrier argument of \autoref{thm_angle_bounds}(a) gives $u\ge cr^{\pi/\beta}$ along the bisector, with $\pi/\beta<3/2$. This contradicts $u\le Cr^{3/2}$. For $\kappa=2$ and $\pi-2\delta>\pi/2$, the same argument contradicts $u\le Cr^2$, as in Step~3 of the proof of \autoref{thm_degenerate}. We call this the \emph{sector contradiction}.

\emph{Step 8: proof of (c) under (F).} Suppose (ii) holds, and fix $\delta_0>0$ with $\pi-2\delta_0>\pi/\kappa$. Let $K=\{1\le|y|\le2,\ \delta_0\le\varphi\le\pi-\delta_0\}$, and put $h_0=\sin\delta_0$.

If there were $r_1>0$ with $K\subset\{u_r>0\}$ for all $r<r_1$, then $\{u>0\}$ would contain the truncated sector $\bigcup_{r<r_1}(x_0+rK)$, and the sector contradiction of Step~7 would apply. Hence there are $r_j\to0$ and points $p_j\in K$ with $u_{r_j}(p_j)=0$.

Let $\mathcal Z$ be the component of $\{u=0\}\cap B_\rho(x_0)$ containing $x_0$. By (F), the other components are relatively closed and finite in number, so a small ball $B_{\rho_1}(x_0)$ meets no other component. For large $j$, $x_0+r_jp_j\in B_{\rho_1}(x_0)$, so $x_0+r_jp_j\in\mathcal Z$. Hence the rescaled set $\mathcal Z_j=r_j^{-1}(\mathcal Z-x_0)$ is connected and contains both $0$ and $p_j$.

Fix $0<\delta<h_0/8$, so that $B_{2\delta}(p_j)\subset K':=\{y_2\ge h_0/2,\ |y|\le3\}$. Since $|p_j|\ge1>2\delta$, the connected set $\mathcal Z_j$ meets every circle $\partial B_t(p_j)$ with $\delta\le t\le2\delta$. No such circle lies entirely in $\{u_{r_j}=0\}$ for large $j$. Otherwise, by subharmonicity of the zero extension, $u_{r_j}=0$ on $B_t(p_j)$, whose area is at least $\pi\delta^2$; this is impossible, because $|K'\cap\{u_{r_j}=0\}|\to0$ by the phase convergence in (b). Each circle therefore contains points of both phases, and hence a point of $\Gamma_{r_j}=\partial\{u_{r_j}>0\}$.

The map $y\mapsto|y-p_j|$ is $1$-Lipschitz and sends $\Gamma_{r_j}\cap(B_{2\delta}(p_j)\setminus B_\delta(p_j))$ onto a set containing $[\delta,2\delta]$. So this part of $\Gamma_{r_j}$ has $\mathcal H^1$-measure at least $\delta$. By (F), the same holds for its regular part. On a regular arc, the distributional Laplacian of the zero extension is $|\nabla u_{r_j}|\,\mathcal H^1=\sqrt{q_{r_j}}\,\mathcal H^1$, by Green's formula and the one-sided Bernoulli trace. Since $\mu_{r_j}=\Delta u_{r_j}\ge0$ and $q_{r_j}\ge c(h_0)>0$ on $K'$ for large $j$, we get
\[
\mu_{r_j}(K')\ge c(h_0)\,\delta>0 .
\]
On the other hand, take $\zeta\in C^\infty_c(\{y_2>h_0/4\})$ with $0\le\zeta\le1$ and $\zeta=1$ on $K'$. Then
\[
\int\zeta\,d\mu_{r_j}=-\int\nabla u_{r_j}\cdot\nabla\zeta\to0,
\]
since $u_{r_j}\to0$ in $W^{1,2}_{\rm loc}$. This contradiction proves (c) under (F).
\end{proof}

\begin{corollary}[Isolation and tangent rays]\label{cor_no_bubbles}
Let $\Psi$ be a variational solution in $B_\rho(x_0)$ satisfying (G) and either (R) or (F). Then the stagnation points in $\partial\{\Psi>0\}\cap\{|x|=r_*\}\cap B_{\rho/2}(x_0)$ are isolated, and each is a corner as in \autoref{thm_variational_crest}(b)(i). Under (R), moreover, $r^{-1}(\Gamma-x_0)$ converges as $r\to0$, at each such point $x_0$, to the two bounding rays of the corner, in the Hausdorff sense on compact subsets of $\{y_2>0\}$.
\end{corollary}
\begin{proof}
\emph{Isolation.} By \autoref{thm_variational_crest}, each such point $x'$ is a corner with Weiss density $\Phi_{\kappa,x'}(0+)=g\,\omega_\kappa$, where $\omega_{3/2}=\sqrt3/3$ and $\omega_2=\pi/16+1/8$. Suppose such points $x_j\ne x_0$ converge to $x_0$, and set $r_j=|x_j-x_0|$. After a subsequence, $(x_j-x_0)/r_j\to e\in\{\pm e_1\}$, since the $x_j$ lie on the circle. Choose $s>0$ with $B_{2s}(e)$ disjoint from the closed corner sector at $x_0$. Since $\Phi_{\kappa,x_j}(t)+C_Wt$ is non-decreasing, with $C_W$ independent of $j$, and the boundary term in $\Phi_\kappa$ is non-positive,
\[
g\,\omega_\kappa\le\Phi_{\kappa,x_j}(r_js)+C_Wr_js\le s^{-2\kappa}\int_{B_s((x_j-x_0)/r_j)}\big(|\nabla u_{r_j}|^2+q_{r_j}\chi_{r_j}\big)+C_Wr_js .
\]
The right-hand side tends to $0$, because $u_{r_j}\to u_*$ in $W^{1,2}_{\rm loc}$ and $\chi_{r_j}\to\chi_{\rm corner}$ in $L^1_{\rm loc}$, and both limits vanish on $B_{2s}(e)$. This is a contradiction.

\emph{Hausdorff convergence} (under (R)). Let $z_r\in\Gamma_r$ with $z_r\to z\in\{y_2>0\}\setminus\{0\}$. If $z$ were inside the corner sector, $u_r(z_r)=0$ would contradict $u_r\to u_*>0$ near $z$. If $z$ were outside the closed sector, \eqref{eq_nondeg_var} would give $\sup u_r\ge c>0$ near $z$, where $u_*=0$. Conversely, every point of a bounding ray is a limit of points of $\Gamma_r$, since nearby there are points with $u_r>0$ and, by $L^1$ convergence, points with $\chi_r=0$.
\end{proof}

\begin{corollary}[Crests of overhanging waves]\label{cor_overhanging}
Let $(\Psi,\Omega)$ be a classical solution with general free surface in the lower region (\autoref{def_general_classical}). Then every stagnation point is a corner symmetric about the inward radial direction. If $\Delta>0$ it is a $120^\circ$ Stokes corner with amplitude $2\Delta^{1/4}/(3|p_0|)$; if $\Delta=0$ it is a right angle with profile $\frac{\sqrt2}{4|p_0|}(y_2^2-y_1^2)$. The stagnation points are finitely many, and at each of them the rescaled free surface converges to the two bounding rays, as in \autoref{cor_no_bubbles}.
\end{corollary}
\begin{proof}
By the remarks after \autoref{def_variational}, $\Psi$ is a regular variational solution in a small ball around each stagnation point $x_0$, with $x_0\in\partial\{\Psi>0\}$. Apply \autoref{thm_variational_crest} and \autoref{cor_no_bubbles}. Finiteness follows from isolation and the compactness of $\Gamma\cap\{|x|=r_*\}$.
\end{proof}

\autoref{cor_overhanging} assumes no graph property, symmetry or monotonicity. For radial subgraphs it recovers \autoref{thm_main_stokes}(c)--(d) and \autoref{thm_degenerate}(c)--(d); the full tangent-cone statement (b) of those theorems also uses \autoref{lem_phase}.

\begin{remark}[What remains open]\label{rem_full_density}
For variational solutions, \autoref{prop_variational_growth} and \autoref{thm_variational_crest} leave only two alternatives: the corner, and full density. Full density is excluded under (R) or (F). For flat waves, Varv\u{a}ruc\u{a} and Weiss exclude full density under (F) with a frequency formula \cite[Sections~7--10]{VW2011}, \cite[Proposition~6.3]{DYZ2026}. In Step~8 the frequency formula is not needed, because the angle bound already produces interface points at heights comparable to the scale. What remains open is a variational solution whose zero set has infinitely many components accumulating at the crest and whose interface is not regular. Excluding full density in that case seems to require a genuinely new idea, even for flat waves.
\end{remark}

\begin{remark}[Relevance to overturning]\label{rem_overturning_crest}
Suppose alternative (II) of \autoref{thm_intro_global} occurs, and the overturning waves remain classical solutions with general free surface in the lower region. By \autoref{cor_overhanging}, every stagnation point that forms is still a corner of the predicted angle. The crest-shape half of the Stokes conjecture is therefore insensitive to the loss of the graph property. We do not claim that the continuum of \cite{LiWang2026} can be followed beyond the onset of overturning, since its hodograph formulation requires the graph property.
\end{remark}

\subsection{Limits of regular waves}\label{sec_limit}

\begin{proposition}[Uniform estimates in the lower region]\label{prop_uniform}
Let $\mathcal A$ be a compact subset of $\{\alpha:\Delta(\alpha)>0,\ r_-(\alpha)>1\}$. There are constants depending only on $\mathcal A$ and $p_0$ such that every lower-region classical solution with $\alpha\in\mathcal A$ satisfies:
\begin{enumerate}
\item $|\nabla\Psi|\le L$ and $\min S\ge1+1/L$;
\item the conclusions of \autoref{thm_local_bernstein} with uniform $\eta$, $C_1$ and $C_2$;
\item the non-degeneracy $\sup_{B_{4\rho}(z)}\Psi\ge c(q_0)\rho$ at every $z\in\partial\Omega$ with $Q_\alpha(R(z))\ge q_0$, for $\rho\le\rho_*(q_0)$.
\end{enumerate}
\end{proposition}
\begin{proof}
(a) \autoref{thm_max_velocity} gives $|\nabla\Psi|^2\le r_-^2Q_\alpha(1)/p_0^2$. The fixed flux $1=\int_1^S|\Psi_R|\,dR\le(S-1)L$ then gives $\min S\ge1+1/L$.

(b) Repeat the proof of \autoref{thm_local_bernstein} with $W=\frac M\eta(r_--R)$. On the free surface $W\ge0$ because the solution lies in the lower region; $M$ is uniformly bounded by (a).

(c) This is the flux argument of \autoref{lem_nondeg}, whose constants depend only on $q_0$ and $p_0$.
\end{proof}

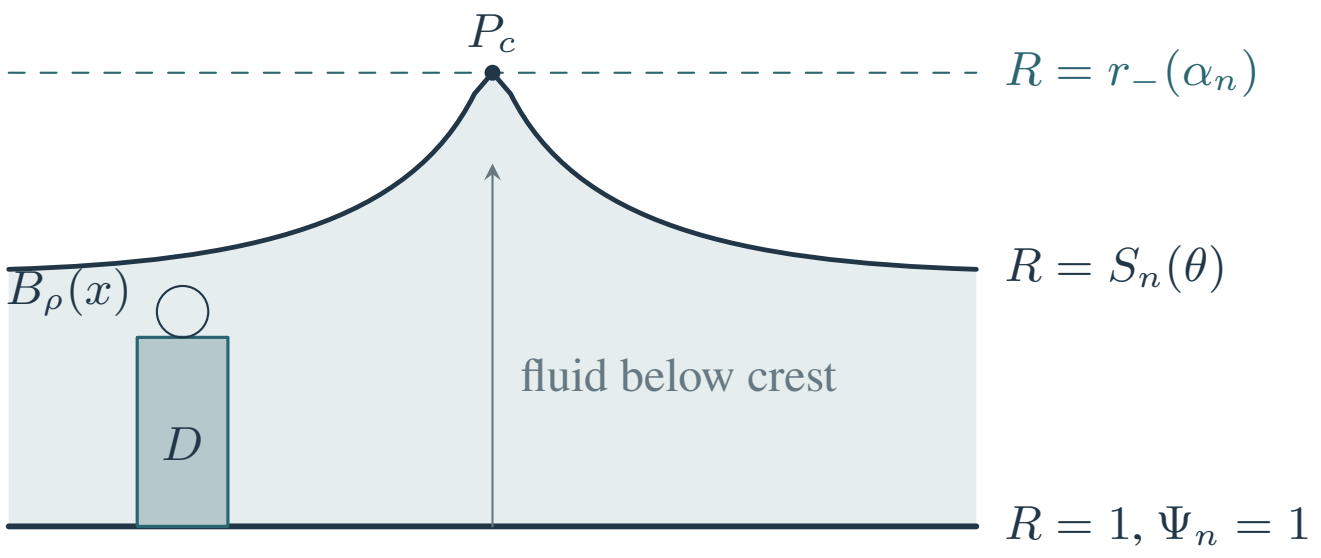
\begin{figure}[!ht]
\centering
\begin{tikzpicture}[x=1cm,y=1cm,>=stealth,font=\small,line cap=round,line join=round]
\fill[PaperAccent!12] (-3.2,0) -- (-3.2,1.7) .. controls (-1.6,1.75) and (-.55,2.0) .. (-.12,2.86) -- (0,3.0) -- (.12,2.86) .. controls (.55,2.0) and (1.6,1.75) .. (3.2,1.7) -- (3.2,0) -- cycle;
\draw[PaperInk,line width=.9pt] (-3.2,1.7) .. controls (-1.6,1.75) and (-.55,2.0) .. (-.12,2.86) -- (0,3.0) -- (.12,2.86) .. controls (.55,2.0) and (1.6,1.75) .. (3.2,1.7);
\draw[PaperInk,line width=1.1pt] (-3.2,0) -- (3.2,0);
\draw[PaperAccent,dashed] (-3.2,3.0) -- (3.2,3.0);
\node[anchor=west,text=PaperAccent,font=\footnotesize] at (3.25,3.0) {$R=r_-(\alpha_n)$};
\node[anchor=west,text=PaperInk,font=\footnotesize] at (3.25,1.7) {$R=S_n(\theta)$};
\node[anchor=west,text=PaperInk,font=\footnotesize] at (3.25,0) {$R=1$, $\Psi_n=1$};
\fill[PaperAccent!35] (-2.35,0) rectangle (-1.75,1.25);
\draw[PaperAccent,line width=.6pt] (-2.35,0) rectangle (-1.75,1.25);
\node[text=PaperInk,font=\footnotesize] at (-2.05,.55) {$D$};
\draw[PaperInk] (-2.05,1.42) circle (.17);
\node[text=PaperInk,font=\footnotesize,anchor=east] at (-2.25,1.55) {$B_\rho(x)$};
\draw[PaperMuted,->] (0,0) -- (0,2.4);
\node[text=PaperMuted,font=\footnotesize,anchor=west] at (.05,1.0) {fluid below crest};
\fill[PaperInk] (0,3.0) circle (1.5pt) node[above,font=\footnotesize] {$P_c$};
\end{tikzpicture}
\caption{The column argument in the limit theorem (schematic, radial coordinate vertical). If a ball $B_\rho(x)$ stays in $\Omega_n$, then by inward radial closure the fixed column $D$ beneath it lies in $\Omega_n$. The harmonic measure of the bottom then bounds $\Psi_n(x)$ from below, independently of $n$.}
\label{fig_column}
\end{figure}

\begin{theorem}[Limit theorem]\label{thm_limit}
Under the hypotheses of \autoref{thm_intro_limit}, a subsequence satisfies the following.
\begin{enumerate}
\item $\Psi_n\to\Psi_\infty$ uniformly, $\nabla\Psi_n\to\nabla\Psi_\infty$ in $L^2_{\rm loc}(\{R>1\})$, and $\chi_{\Omega_n}\to\chi_{\Omega_\infty}$ in $L^1_{\rm loc}$, where $\Omega_\infty=\{\Psi_\infty>0\}$.
\item $\Omega_\infty=\{1<R<S_\infty(\theta)\}$ with $S_\infty$ lower semicontinuous and $S_\infty\le r_-(\alpha_\infty)$; moreover $\Psi_\infty$ is Lipschitz, harmonic in $\Omega_\infty$ and non-increasing along rays.
\item $(\Psi_\infty,\Omega_\infty)$ satisfies the domain-variation identity at $\alpha_\infty$.
\item There is a stagnation point $P_c\in\partial\Omega_\infty$ on the circle $R=r_-(\alpha_\infty)$.
\item Every stagnation point is a Stokes corner as in \autoref{thm_intro_stokes}, and $\Omega_\infty$ is not an annulus.
\end{enumerate}
\end{theorem}
\begin{proof}
\emph{(a), (b).} By \autoref{prop_uniform}(a) and the Arzel\`a--Ascoli theorem, a subsequence converges uniformly. Monotonicity along rays passes to the limit, so $\Omega_\infty$ is a radial subgraph.

To prove $\chi_n\to\chi_\infty$, let $x$ be an interior point of $\{\Psi_\infty=0\}$ with $R(x)<r_--\varepsilon$, and choose $\rho$ with $B_{5\rho}(x)\subset\{\Psi_\infty=0\}$. If $B_\rho(x)$ met $\partial\Omega_n$, non-degeneracy would give $\sup_{B_{5\rho}(x)}\Psi_n\ge c\rho$, which is impossible for large $n$. Hence eventually either $B_\rho(x)\subset\Omega_n$ or $B_\rho(x)$ lies in the air.

In the first case (see \autoref{fig_column}), inward radial closure shows that $\Omega_n$ contains the fixed column
\[
D=\Big\{1<R<R(x)+\frac\rho2,\ |\theta-\theta(x)|<\frac{\rho}{4R(x)}\Big\},
\]
because every point of $D$ lies radially below a point of $B_\rho(x)$. The point $x$ is interior to $D$. On $\partial D$ we have $\Psi_n\ge0$, and $\Psi_n=1$ on the bottom arc of $D$. Comparison with the harmonic measure $\omega_D$ of that arc therefore gives $\Psi_n(x)\ge\omega_D(x)>0$, a bound independent of $n$, contradicting $\Psi_\infty(x)=0$. So $x$ is eventually in the air. Points of $\Omega_\infty$ are eventually in $\Omega_n$ by uniform convergence. It remains to see that $\partial\Omega_\infty\cap\{R<r_--\varepsilon\}$ is a null set, since the strip $\{r_--\varepsilon\le R\le r_-\}$ has measure $O(\varepsilon)$. Every $z\in\partial\Omega_\infty$ with $R(z)<r_--\varepsilon$ is a limit of points of $\partial\Omega_n$, by the column argument. Passing to the limit in the non-degeneracy of \autoref{prop_uniform}(c) gives $y\in B_{4\rho}(z)$ with $\Psi_\infty(y)\ge c\rho$. The Lipschitz bound $L$ then puts $B_{c\rho/L}(y)$ inside $\Omega_\infty$. So $\partial\Omega_\infty$ is porous there, hence of measure zero by the Lebesgue density theorem, and $\chi_n\to\chi_\infty$ in $L^1_{\rm loc}$. Strong convergence of gradients is proved as in \autoref{prop_blowup}.

\emph{(c).} Pass to the limit in the domain-variation identities of the $\Psi_n$.

\emph{(d).} Let $\theta_n$ be maximizers of $S_n$ with $\theta_n\to\theta_*$, and put $P_c=(r_-,\theta_*)$. Fix $h<r_-$. Applied at $(h,\theta_n)\in\Omega_n$, the same column/non-degeneracy dichotomy gives
\begin{equation}\label{eq_fluid_below}
\Omega_\infty\cap B_{6\rho}\big((h,\theta_*)\big)\ne\emptyset
\end{equation}
for all $h<r_-$ and all small $\rho$. Hence $P_c\in\partial\Omega_\infty$. The uniform bound $|\nabla\Psi_\infty|^2\le C(r_--R)$ shows that $P_c$ is a stagnation point.

\emph{(e).} The proof of \autoref{thm_main_stokes} uses only the following properties of the solution:
\begin{itemize}
\item the bound $S\le r_-$;
\item the growth bound \eqref{eq_growth};
\item the domain-variation identity, together with $\int\nabla\Psi\cdot\nabla(\eta\Psi)=0$;
\item non-degeneracy (every boundary point of $\Omega_\infty$ is a limit of boundary points of $\Omega_n$, by the column argument);
\item radial phase monotonicity;
\item positivity and harmonicity of $\Psi_\infty$ in $\Omega_\infty$.
\end{itemize}
All of these hold for the limit. The one step that used the segment below $P_c$ is replaced by \eqref{eq_fluid_below} with $h=r_--\frac32r$ and $\rho=\delta r$. Finally, an annulus $\{1<R<r_-\}$ is excluded, since its outer boundary carries the speed $1/(r_-\log r_-)\ne0$.
\end{proof}

Thus any family of waves of \cite{LiWang2026} in the lower region whose crests approach the stagnation circle produces an extreme wave with Stokes crests. We now compare this with the three alternatives of the global theorem of \cite{LiWang2026} in $\mathcal O_\delta$.
\begin{itemize}
\item \emph{Boundary contact.} By \autoref{cor_first_constraint}, in the lower region the first constraint cannot degenerate below $|p_0|/(r_-\sqrt{Q_\alpha(1)})$ while $\alpha$ stays bounded. Boundary contact therefore happens through the second constraint $\min SQ_\alpha(S)=\delta$, which is exactly the hypothesis of \autoref{thm_limit} as $\delta\downarrow0$.
\item \emph{Unboundedness.} This requires an upper bound for $\alpha$ and H\"older estimates away from stagnation.
\item \emph{Return to the laminar branch.} This is excluded for the oriented continuum in \autoref{sec_monotone}.
\end{itemize}
The continuum may also leave the lower region through the window $\Delta\le0$, or approach $\alpha_\pm$.

\begin{theorem}[Degenerate limits]\label{thm_degenerate_limit}
Let $(\Psi_n,S_n)$ be classical solutions in the lower region with parameters $\alpha_n$, and suppose that
\[
\alpha_n\to\alpha_*\in\{\alpha_+,\alpha_-\},\qquad \alpha_*>1,\qquad \max S_n-r_-(\alpha_n)\to0 .
\]
Then a subsequence converges, in the sense of \autoref{thm_limit}(a)--(c), to a limit $(\Psi_\infty,\Omega_\infty)$ with $S_\infty\le\alpha_*$. The limit has a stagnation point on the circle $R=\alpha_*$, and at every stagnation point the free surface forms the symmetric right-angle corner of \autoref{thm_degenerate}(b)--(c). If moreover each $S_n$ is even and non-increasing on $(0,\pi)$, the limit has exactly one stagnation point per period, at $\theta=0$.
\end{theorem}
\begin{proof}
\emph{Uniform bounds.} \autoref{thm_max_velocity} gives $|\nabla\Psi_n|\le L$ with $L$ independent of $n$, because $r_-(\alpha_n)$ and $Q_{\alpha_n}(1)$ converge. As in \autoref{prop_uniform}(a), it follows that $\min S_n\ge1+1/L$. By \eqref{eq_flow_force_bounds}, $|\Psi_{n,R}(1,\cdot)|\le1/\log(1+1/L)$. Since $Q_{\alpha_n}(1)\to(\alpha_*-1)^2>0$, the constants $M_n$ of \autoref{prop_interior_bernoulli} are bounded by some $M$. Hence
\[
p_0^2|\nabla\Psi_n(x)|^2\le M\,Q_{\alpha_n}(|x|)\qquad\text{in }\Omega_n .
\]
The non-degeneracy of \autoref{prop_uniform}(c) holds uniformly at boundary points with $Q_{\alpha_n}(R)\ge q_0$, because its proof uses only $q_0$ and $p_0$.

\emph{Convergence.} The proofs of parts (a)--(d) of \autoref{thm_limit} use only these three facts, together with \autoref{lem_phase}. The one exception is the bound $|\nabla\Psi_\infty|^2\le C(r_--R)$ at the end of part (d), which here is replaced by the limit of the display above:
\[
p_0^2|\nabla\Psi_\infty(x)|^2\le M(\alpha_*-|x|)^2\le M|x-P_c|^2 .
\]
This gives (a)--(c), the bound $S_\infty\le\alpha_*$, and a stagnation point $P_c$ with $|P_c|=\alpha_*$.

\emph{Blow-up.} Steps 2--6 of the proof of \autoref{thm_degenerate} use only four facts: the growth $|\nabla u|\le C|x-P_c|$, the domain-variation identity, non-degeneracy, and \autoref{lem_phase}. All four hold for the limit, the second by \autoref{thm_limit}(c) and the third by passing to the limit in the uniform non-degeneracy. Step 6 is not needed, since $S_\infty\le\alpha_*$ is already known. This gives the right-angle corner at every stagnation point.

For monotone $S_n$, the argument of \autoref{cor_extreme_shape} applies verbatim.
\end{proof}

\section{Monotone waves and an upper bound for the parameter}\label{sec_continuum}

From here on we work on the global continuum of \cite{LiWang2026} under the standing assumptions \eqref{eq_standing}. We first show that the continuum splits into two copies of a connected set of monotone waves, and then that the parameter is bounded above on that set.

\subsection{Monotone waves along the oriented continuum}\label{sec_monotone}

Throughout this subsection, solutions are even in $\theta$ and belong to the class (L). We also assume $Q_\alpha(S)>0$ on the free surface, which holds in every open set $\mathcal O_\delta$ of \cite{LiWang2026}. By the regularity theory for non-degenerate Bernoulli free boundaries \cite{KNS1978}, such solutions are then analytic up to the free surface. Write $\Omega_+=\Omega\cap\{0<\theta<\pi\}$, and call $P_0=(S(0),0)$ the \emph{crest} and $P_\pi=(S(\pi),\pi)$ the \emph{trough}.

\begin{lemma}[The rotational Jacobi field]\label{lem_rotation_field}
The function $w=\Psi_\theta$ has the following properties:
\begin{itemize}
\item it is harmonic in $\Omega$, odd in $\theta$, and vanishes on $R=1$ and on the axes $\theta\in\{0,\pi\}$;
\item on the free surface, $w=|\Psi_R|\,S'(\theta)$;
\item on the free surface it satisfies a homogeneous Robin condition $\partial_\nu w=c\,w$, where $c$ is a continuous function.
\end{itemize}
\end{lemma}
\begin{proof}
Polar Laplacian coefficients do not depend on $\theta$, so $w$ is harmonic. Oddness follows from evenness of the solution, and $w=0$ on the bottom because $\Psi\equiv1$ there.

Differentiating $\Psi(S(\theta),\theta)=0$ gives $w=-\Psi_RS'=|\Psi_R|S'$.

For the Robin condition, write $w=x^\perp\cdot\nabla\Psi$, where $x^\perp=R\,\mathbf e_\theta$. On the surface $\nabla\Psi=-|\nabla\Psi|\nu$, so the first-order term of $\nabla w$ is orthogonal to $\nu$, and
\[
\partial_\nu w=x^\perp\cdot D^2\Psi\,\nu .
\]
Decompose $x^\perp=a\tau+b\nu$; then $b=-w/|\nabla\Psi|$. Next, $\tau\cdot D^2\Psi\nu=-\partial_\tau|\nabla\Psi|^2/(2|\nabla\Psi|)$, and by the Bernoulli condition $\partial_\tau|\nabla\Psi|^2=Q_\alpha'(R)\,(\mathbf e_R\cdot\tau)/p_0^2$. In two dimensions $\mathbf e_R\cdot\tau=\pm(x^\perp\cdot\nu)/R=\pm b/R$. Hence both contributions to $\partial_\nu w$ are multiples of $b$, and therefore of $w$, with coefficients that are continuous because $|\nabla\Psi|\ge\mu_0>0$ on the free surface.
\end{proof}

\begin{figure}[!ht]
\centering
\begin{tikzpicture}[x=1.9cm,y=1.9cm,>=stealth,font=\small,line cap=round,line join=round]
\fill[PaperAccent!14] plot[domain=0:180,samples=90] ({(1.75+.22*cos(\x)-.05*cos(2*\x))*cos(\x)},{(1.75+.22*cos(\x)-.05*cos(2*\x))*sin(\x)}) -- (-1,0) arc (180:0:1) -- cycle;
\draw[PaperInk,line width=.9pt] plot[domain=0:180,samples=90] ({(1.75+.22*cos(\x)-.05*cos(2*\x))*cos(\x)},{(1.75+.22*cos(\x)-.05*cos(2*\x))*sin(\x)});
\draw[PaperInk,line width=1.1pt] (1,0) arc (0:180:1);
\draw[PaperMuted,dashed] (-2.0,0) -- (2.35,0);
\draw[PaperAccent,line width=1.2pt] (1,0) -- (1.92,0);
\draw[PaperAccent,line width=1.2pt] (-1,0) -- (-1.48,0);
\fill[PaperInk] (1.92,0) circle (1.4pt) node[below,font=\footnotesize] {crest $P_0$};
\fill[PaperInk] (-1.48,0) circle (1.4pt) node[below,font=\footnotesize] {trough $P_\pi$};
\node[text=PaperAccent,font=\footnotesize] at (1.46,.11) {$w=0$};
\node[text=PaperAccent,font=\footnotesize] at (-1.24,.11) {$w=0$};
\node[text=PaperInk] at (0,1.36) {$w=\Psi_\theta<0$ in $\Omega_+$};
\node[text=PaperInk,font=\footnotesize] at (0,.6) {$w=0$ on $R=1$};
\node[anchor=west,text=PaperInk,font=\footnotesize,align=left] at (1.55,1.62) {free surface:\\ $\partial_\nu w=c\,w$,\ \ $w=|\Psi_R|S'<0$};
\draw[->,PaperMuted] (1.6,1.55) -- (1.37,1.38);
\node[anchor=west,text=PaperMuted,font=\footnotesize] at (2.38,0) {axis $\theta=0$};
\node[anchor=east,text=PaperMuted,font=\footnotesize] at (-2.03,0) {axis $\theta=\pi$};
\end{tikzpicture}
\caption{The rotational Jacobi field $w=\Psi_\theta$ on the half domain $\Omega_+$ (schematic). It vanishes on the bottom and on both axes, and satisfies a homogeneous Robin condition on the free surface. Negativity of $w$ is preserved along $\mathcal C^+$ by the Hopf lemma, the oblique boundary condition on the free surface, and Serrin's corner lemma at all four corners (\autoref{prop_nodal_strict}).}
\label{fig_nodal}
\end{figure}
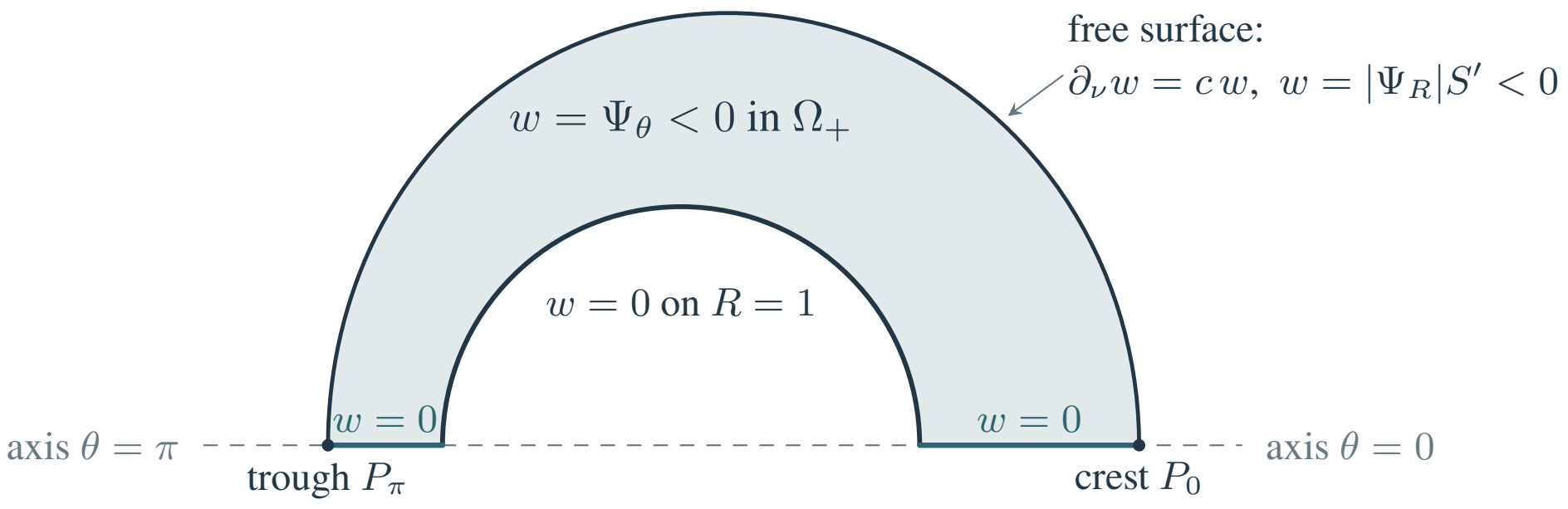

The field $w$ and its boundary conditions are shown in \autoref{fig_nodal}.

\subsubsection{Hodograph form of the Jacobi field}

Openness and closedness of the monotone class are most transparent in the hodograph variables of \cite{LiWang2026}, because the domain there is fixed. Set $q=\theta$ and $p=-\Psi$, so that the fluid occupies $D=\mathbb T\times(-1,0)$, the bottom is $p=-1$ and the free surface is $p=0$. As in \cite{LiWang2026}, we write $R=1+h(q,p)$. We use the logarithmic height
\[
H(q,p)=\log R=\log\big(1+h(q,p)\big),\qquad S(\theta)=e^{H(\theta,0)} .
\]
The map $h\mapsto H$ is a smooth diffeomorphism of the relevant open set of $X=C^{3+s}_{\rm per}(\overline D)$ onto its image, and $H=0$ on the bottom. In $\mathcal O_\delta$ the first constraint $(h+1)h_p>\delta$ gives $H_p=h_p/(1+h)>0$.

\begin{lemma}[Hodograph equations]\label{lem_hodograph}
Since $\Psi$ is harmonic in $(t,\theta)=(\log R,\theta)$, the function $H$ satisfies
\begin{equation}\label{eq_hod}
(1+H_q^2)H_{pp}-2H_qH_pH_{pq}+H_p^2H_{qq}=0\quad\text{in }D,
\end{equation}
together with $H=0$ on $p=-1$ and the Bernoulli condition
\begin{equation}\label{eq_hod_bc}
p_0^2\,(1+H_q^2)=H_p^2\,G(H)\quad\text{on }p=0,\qquad G(H):=e^{2H}Q_\alpha(e^H)=S^2Q_\alpha(S).
\end{equation}
The function $v:=H_q$ is odd in $q$, vanishes on $p=-1$ and on $q\in\{0,\pi\}$, and satisfies
\begin{align}
&\mathcal Lv:=H_p^2v_{qq}-2H_qH_pv_{qp}+(1+H_q^2)v_{pp}+b^1v_q+b^2v_p=0\quad\text{in }D, \label{eq_v}\\
&\mathcal Bv:=2p_0^2H_q\,v_q-2H_pG(H)\,v_p-H_p^2G'(H)\,v=0\quad\text{on }p=0, \label{eq_v_bc}
\end{align}
with $b^1=2H_qH_{pp}-2H_pH_{pq}$ and $b^2=2H_pH_{qq}-2H_qH_{pq}$. The operator $\mathcal L$ is uniformly elliptic, with determinant $H_p^2$, and has no zeroth-order term. In $\mathcal O_\delta$ the coefficient of $v_p$ in $\mathcal B$ is strictly negative, so $\mathcal B$ is oblique.
\end{lemma}
\begin{proof}
Equation \eqref{eq_hod} is the Dubreil-Jacotin equation \cite{DubreilJacotin1934,ConstantinStrauss2004} for a harmonic stream function in the plane $(t,\theta)$; it is invariant under $p\mapsto-p$. The Bernoulli condition becomes \eqref{eq_hod_bc} because $|\nabla\Psi|^2=e^{-2t}(\Psi_t^2+\Psi_\theta^2)=e^{-2H}(1+H_q^2)/H_p^2$. The coefficients of \eqref{eq_hod} and \eqref{eq_hod_bc} do not depend on $q$ explicitly, so differentiating in $q$ gives \eqref{eq_v} and \eqref{eq_v_bc}. Evenness of $H$ gives oddness of $v$, and $H=0$ on the bottom gives $v=0$ there. The determinant is $H_p^2(1+H_q^2)-H_q^2H_p^2=H_p^2$. Finally $H_p>0$ and $G(H)=S\cdot SQ_\alpha(S)>0$ in $\mathcal O_\delta$.
\end{proof}

Since $\Psi_\theta=H_q/H_p$ and $S'(\theta)=S\,H_q(\theta,0)$, the sign of $v$ is the sign of the rotational field $w=\Psi_\theta$ of \autoref{lem_rotation_field} and of $S'$.

\begin{definition}[Monotone waves]\label{def_nodal}
Let $\mathcal N$ be the set of even solutions in $\mathcal O_\delta$ whose field $v=H_q$ satisfies:
\begin{enumerate}
\item[(N1)] $v<0$ in $(0,\pi)\times(-1,0)$; equivalently $\Psi_\theta<0$ in $\Omega_+$;
\item[(N2)] $v(q,0)<0$ for $0<q<\pi$; equivalently $S'(\theta)<0$ for $0<\theta<\pi$;
\item[(N3)] $v_p(q,-1)<0$ for $0<q<\pi$;
\item[(N4)] $v_q(0,p)<0$ and $v_q(\pi,p)>0$ for $-1<p\le0$;
\item[(N5)] $v_{qp}(0,-1)<0$ and $v_{qp}(\pi,-1)>0$.
\end{enumerate}
Since $v_q(0,0)=S''(0)/S(0)$ at the crest, (N4) contains $S''(0)<0<S''(\pi)$. Waves in $\mathcal N$ have a single crest and a single trough per period and are strictly monotone between them.
\end{definition}

The five conditions combine into one. For $(q,p)\in(0,\pi)\times(-1,0]$ set
\begin{equation}\label{eq_omega}
\omega(q,p)=\frac{v(q,p)}{\sin q\,(p+1)} .
\end{equation}
Because $v$ vanishes on $q=0$ and on $p=-1$,
\[
v(q,p)=q\,(p+1)\int_0^1\!\!\int_0^1 v_{qp}\big(\sigma q,\,-1+\tau(p+1)\big)\,d\sigma\,d\tau ,
\]
and the analogous formula holds with base point $q=\pi$. Hence $\omega$ extends continuously to the closed rectangle $K=[0,\pi]\times[-1,0]$. Its boundary values are $v_p/\sin q$ on the bottom, $v_q/(p+1)$ on $q=0$, $-v_q/(p+1)$ on $q=\pi$, and $\pm v_{qp}$ at the two bottom corners. Comparing with (N1)--(N5) gives the following.

\begin{lemma}[A single sign condition]\label{lem_omega}
An even solution in $\mathcal O_\delta$ belongs to $\mathcal N$ if and only if $\omega<0$ on the compact set $K$. Moreover,
\[
\|\omega_1-\omega_2\|_{C^0(K)}\le C\,\|H_1-H_2\|_{C^3(\overline D)},
\]
with $C$ an absolute constant.
\end{lemma}
\begin{proof}
The equivalence is the list of boundary values above. For the estimate, use the integral formula with base point $0$ on $[0,2\pi/3]\times[-1,0]$, where $q/\sin q\le 2\pi/(3\sin(2\pi/3))$, and the formula with base point $\pi$ on $[\pi/3,\pi]\times[-1,0]$. In both cases $\omega$ is a bounded multiple of an average of $v_{qp}=H_{qqp}$.
\end{proof}

\begin{proposition}[Strict nodal property]\label{prop_nodal_strict}
Let $H$ be an even, non-circular solution in $\mathcal O_\delta$, and suppose that $v=H_q\le0$ in $(0,\pi)\times(-1,0)$. Then the solution belongs to $\mathcal N$.
\end{proposition}
\begin{proof}
All statements are applied to $v\in C^{2+s}(\overline D)$, which solves the elliptic equation \eqref{eq_v} without zeroth-order term.

\emph{(N1).} If $v\equiv0$, then $H$ depends only on $p$ and the solution is circular. Otherwise the strong maximum principle gives $v<0$ in the open half-rectangle.

\emph{(N3) and the interior of (N4).} On the bottom and on the axes $v=0$ is the maximum of $v$ over the closed half-rectangle, attained at points where the boundary is flat. The Hopf lemma gives $v_p<0$ on the bottom, $v_q<0$ on $q=0$ and $v_q>0$ on $q=\pi$, for $-1<p<0$.

\emph{(N2).} Suppose $v(q_1,0)=0$ for some $q_1\in(0,\pi)$. Then $v$ attains its maximum $0$ at the boundary point $(q_1,0)$, so $v_q(q_1,0)=0$, and the Hopf lemma gives $v_p(q_1,0)>0$. But \eqref{eq_v_bc} at this point reduces to $-2H_pG(H)\,v_p=0$, so $v_p(q_1,0)=0$. This is a contradiction.

\emph{Corners.} We use Serrin's corner lemma \cite{Serrin1971}, in the form given by Gidas, Ni and Nirenberg \cite[Lemma~S]{GNN1979}. Let the domain near a corner $Q$ be $\{\rho<0,\ \sigma<0\}$ for two $C^2$ functions with transversal zero sets. Let $w\in C^2$ near $Q$ satisfy $a^{ij}w_{ij}+b^iw_i+cw\ge0$, with bounded coefficients and uniform ellipticity, and suppose $w<0$ inside and $w(Q)=0$. Assume $a^{ij}\rho_i\sigma_j\ge0$ at $Q$; when equality holds, assume also $a^{ij}\in C^2$ near $Q$. In two dimensions no further condition is needed, since the intersection of the two zero sets is the point $Q$. Then for every direction $s$ entering the domain transversally to both edges, $\partial_sw<0$ or $\partial_s^2w<0$ at $Q$.

At each of our four corners the edges are coordinate lines, so $a^{ij}\rho_i\sigma_j=\pm(-H_qH_p)$, which vanishes at $Q$ because $H_q=0$ on the axes and on the bottom. The coefficients of $\mathcal L$ are smooth functions of $\nabla H$ and $D^2H$, hence of class $C^{1+s}$; the leading coefficients $a^{ij}$ depend on $\nabla H$ only and are of class $C^{2+s}$. So the lemma applies to $w=v$.

\emph{(N5).} At $(0,-1)$, $v$ vanishes on both edges. Hence $\nabla v=0$, $v_{qq}=0$ and $v_{pp}=0$ there. With $s=(1,1)/\sqrt2$, Serrin's lemma gives $\partial_s^2v=v_{qp}<0$. The corner $(\pi,-1)$ is symmetric.

\emph{(N4) at the crest and the trough.} At $(0,0)$ we have $v_q\le0$, because $v<0$ for $q>0$. Suppose $v_q(0,0)=0$. Since $v\equiv0$ on the axis, also $v_p=v_{pp}=0$ there, and $\nabla v(0,0)=0$. Moreover $v_{qq}(0,0)=H_{qqq}(0,0)=0$ because $H_{qqq}$ is odd. Differentiate \eqref{eq_v_bc} in $q$ along $p=0$ and evaluate at $(0,0)$. Every term except $-2H_pG(H)\,v_{qp}$ contains one of the factors $v$, $v_q$, $v_p$ or $H_q$, and these all vanish at $(0,0)$. Hence $v_{qp}(0,0)=0$, so all first and second derivatives of $v$ vanish at the corner. This contradicts Serrin's lemma. Therefore $v_q(0,0)<0$. The trough $(\pi,0)$ is handled in the same way and gives $v_q(\pi,0)>0$.
\end{proof}

\begin{remark}
In physical variables, the corner step at the crest is the statement $S''(0)<0$. The argument uses only the oddness of $v$, the obliqueness of \eqref{eq_v_bc}, and the vanishing of the coefficient $2p_0^2H_q$ of $v_q$ in \eqref{eq_v_bc} on the axis. No reflection across the free surface is needed.
\end{remark}

\begin{proposition}[Openness, closedness and returns]\label{prop_nodal_open_closed}
In the topology of the space $X=\{h\in C^{3+s}_{\rm per}(\overline D):h=0\text{ on the bottom}\}$, whose elements are even and $2\pi$-periodic in $q$:
\begin{enumerate}
\item $\mathcal N$ is open;
\item every non-circular limit in $\mathcal O_\delta$ of elements of $\mathcal N$ lies in $\mathcal N$;
\item the local branch at $\alpha_c$ lies in $\mathcal N$ for one sign of the amplitude, namely $S=R_0(1+\varepsilon\cos\theta)+o(\varepsilon)$ with $\varepsilon>0$;
\item if elements of $\mathcal N$ converge to a circular solution $(\alpha_*,S\equiv\sigma)$, then the linearization at that circle is singular in the mode $\cos\theta$. For the laminar circle $\sigma=R_0$ this means $\alpha_*=\alpha_c$.
\end{enumerate}
\end{proposition}
\begin{proof}
(a) Let $H\in\mathcal N$. By \autoref{lem_omega}, $\omega\le-\mu<0$ on the compact set $K$ for some $\mu>0$. The map $h\mapsto H$ is continuous into $C^3$, so the same lemma shows that every $\tilde h$ with $\|\tilde h-h\|_X$ small enough satisfies $\tilde\omega\le-\mu/2$ on $K$. Such $\tilde h$ lie in $\mathcal N$.

(b) If $H_n\to H$ with $H_n\in\mathcal N$, then $v_n\to v$ uniformly, so $v\le0$. Now apply \autoref{prop_nodal_strict}.

(c) The laminar flow is $H_0(p)=\gamma(1+p)$. Linearizing $\Psi(e^H,q)=-p$ in the mode $\cos q$ gives
\[
H=H_0+\varepsilon V(p)\cos q+o(\varepsilon)\ \text{in }X,\qquad V(p)=\frac{R-R^{-1}}{R_0-R_0^{-1}},\quad R=e^{\gamma(1+p)} .
\]
Here $V(-1)=0$, $V>0$ for $p>-1$, and $V'(-1)=2\gamma/(R_0-R_0^{-1})>0$. Hence $\varepsilon^{-1}\omega\to-V(p)/(p+1)$ uniformly on $K$, and this limit is negative on $K$. So $\omega<0$ on $K$ for small $\varepsilon>0$, and \autoref{lem_omega} gives (c).

(d) Let $H_n\in\mathcal N$ converge in $C^{3+s}$ to the circular solution $H_*(p)=(1+p)\log\sigma$ at $\alpha_*$, and let $v_n=\partial_qH_n\not\equiv0$. Each $v_n$ solves the linear problem \eqref{eq_v}--\eqref{eq_v_bc}, with $v_n=0$ on the bottom. The coefficients of this problem converge in $C^{1+s}$ to those at $H_*$, which are independent of $q$. The domain $D$ has two boundary components and no corners, the operator is uniformly elliptic, and the top condition is uniformly oblique. The Schauder estimate for oblique problems \cite[Theorem~6.30]{GilbargTrudinger2001} therefore gives
\[
\|v\|_{C^{2+s}(\overline D)}\le C\|v\|_{C^0(\overline D)}
\]
uniformly in $n$. The normalized fields $\hat v_n=v_n/\|v_n\|_{C^0}$ have a subsequence converging in $C^2$ to some $\hat v$ with $\|\hat v\|_{C^0}=1$. The limit solves the problem at $H_*$, vanishes on the bottom, is odd in $q$, and satisfies $\hat v\le0$ on $(0,\pi)\times[-1,0]$.

At $H_*$ the problem separates. Write $\hat v=\sum_{k\ge1}V_k(p)\sin kq$. Each $V_k$ solves a second-order ODE on $(-1,0)$ with $V_k(-1)=0$ and a homogeneous Robin condition at $p=0$. Since $q$-differentiation maps the mode $\cos kq$ of the linearization of \eqref{eq_hod}--\eqref{eq_hod_bc} at $H_*$ onto the mode $\sin kq$ of \eqref{eq_v}--\eqref{eq_v_bc}, $V_k\not\equiv0$ exactly when the linearization at the circle is singular in the mode $\cos kq$. In cylinder variables the mode-$k$ kernel is $\sinh(kt)\cos k\theta$, and the order-one computation of \autoref{thm_alpha2}, with $(R_0,c_0)$ replaced by $(\sigma,Q_{\alpha_*}(\sigma))$, shows that it is singular if and only if
\[
\alpha_*-\sigma=\frac{Q_{\alpha_*}(\sigma)}{\sigma}\big(1+k\coth(k\log\sigma)\big).
\]
The right-hand side is strictly increasing in $k$, so at most one mode $k\ge1$ is singular, and $\hat v=bV_k(p)\sin kq$ with $b\neq0$. Moreover $V_k$ is unique up to a factor, since it solves a second-order ODE with $V_k(-1)=0$. Linearizing $\Psi(e^{H},q)=-p$ at the circle shows that $V_k$ is a non-zero multiple of $\sinh\big(k(1+p)\log\sigma\big)$, which does not vanish on $(-1,0]$. Since $\hat v\le0$ on $(0,\pi)$ and $\sin kq$ changes sign there when $k\ge2$, we get $k=1$. If no mode $k\ge1$ were singular, $\hat v$ would vanish, contradicting $\|\hat v\|_{C^0}=1$. At $\sigma=R_0$ the mode-$1$ condition reads $\alpha_*=R_0+c_0(1+\coth\gamma)/R_0=R_0+c_0/\sinh\gamma=\alpha_c$.
\end{proof}

\begin{theorem}[The oriented continuum consists of monotone waves]\label{thm_oriented}
Assume \eqref{eq_standing}. Fix $\delta\in(0,\gamma)$, and let $\mathcal C_\delta$ be the continuum of \cite[Theorem~2.2]{LiWang2026}: the connected component, containing $(\alpha_c,R_0)$, of the closure of the non-trivial solutions in $\mathcal O_\delta$. Let $T$ denote the half-period shift $\theta\mapsto\theta+\pi$. Then exactly one of the following holds:
\begin{enumerate}
\item $\mathcal C_\delta$ contains a non-laminar circular solution at which the mode-$1$ linearization is singular (a secondary bifurcation from the second circular family);
\item $\mathcal C_\delta\subset\mathcal N\cup T(\mathcal N)\cup\{(\alpha_c,R_0)\}$. In this case the set $\mathcal C^+_\delta:=(\mathcal C_\delta\cap\mathcal N)\cup\{(\alpha_c,R_0)\}$ is connected, satisfies $\mathcal C_\delta=\mathcal C^+_\delta\cup T(\mathcal C^+_\delta)$ and $\mathcal C^+_\delta\cap T(\mathcal C^+_\delta)=\{(\alpha_c,R_0)\}$, contains the $\varepsilon>0$ half of the local branch, and is unbounded in $(\alpha_s,\infty)\times X$ or meets $\partial\mathcal O_\delta$.
\end{enumerate}
By \autoref{thm_crossing}, case (a) never occurs. No unilateral bifurcation theorem is needed; the argument uses only \cite[Theorem~2.2]{LiWang2026} and the symmetry $T$.
\end{theorem}
\begin{proof}
\emph{Symmetry.} $T$ preserves the equations, the class of even functions, both constraints defining $\mathcal O_\delta$, and the laminar line. It fixes $(\alpha_c,R_0)$ and maps the kernel element $\cos\theta$ to $-\cos\theta$. Hence $T(\mathcal C_\delta)=\mathcal C_\delta$, and near the bifurcation point $T$ exchanges the halves $\varepsilon>0$ and $\varepsilon<0$ of the local Crandall--Rabinowitz curve.

\emph{Clopen decomposition.} Suppose (a) fails, and write $b=(\alpha_c,R_0)$ and $M=\mathcal N\cup T(\mathcal N)$. By \autoref{prop_nodal_open_closed}(a), $M\cap\mathcal C_\delta$ is open in $\mathcal C_\delta$. Let $x\in\mathcal C_\delta\setminus\{b\}$ be a limit point of $M\cap\mathcal C_\delta$. If $x$ is non-circular, then $x\in M$ by \autoref{prop_nodal_open_closed}(b) applied to $\mathcal N$ or to $T(\mathcal N)$. If $x$ is circular, \autoref{prop_nodal_open_closed}(d) forces it to be laminar with $\alpha=\alpha_c$, hence $x=b$, which we have excluded, or a non-laminar circle with a singular mode-$1$ linearization, which (a) excludes. Thus $E_1:=(M\cap\mathcal C_\delta)\cup\{b\}$ is closed in $\mathcal C_\delta$.

Near $b$, local uniqueness in the Crandall--Rabinowitz theorem \cite{CrandallRabinowitz1971} shows that $\mathcal C_\delta$ consists only of the local curve, which lies in $M\cup\{b\}$ by \autoref{prop_nodal_open_closed}(c) and the symmetry. Hence $E_2:=\mathcal C_\delta\setminus E_1$ stays a positive distance from $b$, so it is closed. Since $\mathcal C_\delta=E_1\sqcup E_2$ is connected and $E_1\ne\emptyset$, we get $E_2=\emptyset$, i.e.\ $\mathcal C_\delta\subset M\cup\{b\}$.

\emph{The half $\mathcal C^+_\delta$.} Since $\mathcal N\cap T(\mathcal N)=\emptyset$, we have $\mathcal C_\delta=\mathcal C^+_\delta\cup T(\mathcal C^+_\delta)$ with intersection $\{b\}$, and both sets are closed. Suppose $\mathcal C^+_\delta=A\sqcup B$ with $A$, $B$ closed, $b\in A$ and $B\ne\emptyset$. Then $B$ is closed in $\mathcal C_\delta$ and disjoint from $A\cup T(\mathcal C^+_\delta)$, contradicting the connectedness of $\mathcal C_\delta$. So $\mathcal C^+_\delta$ is connected.

\emph{Alternatives.} \cite[Theorem~2.2]{LiWang2026} offers three alternatives for $\mathcal C_\delta$. The second, a trivial point $(\alpha_k,H)$ with $\alpha_k\neq\alpha_c$, is impossible because $\mathcal C_\delta\subset M\cup\{b\}$ and trivial points are circular. The other two are invariant under the isometry $T$, so whichever holds for $\mathcal C_\delta$ also holds for $\mathcal C^+_\delta$.
\end{proof}

\begin{corollary}[Shape of the extreme wave]\label{cor_extreme_shape}
Suppose a sequence in $\mathcal C^+$ satisfies the hypotheses of \autoref{thm_intro_limit}. Then the limiting extreme wave is symmetric, and $S_\infty$ is non-increasing on $(0,\pi)$. It has exactly one stagnation point per period, located at the crest $\theta=0$, where the free surface forms the symmetric $120^\circ$ Stokes corner.
\end{corollary}
\begin{proof}
Monotonicity of the $S_n$ and the radial-subgraph structure pass to the limit. By \autoref{thm_limit}, the stagnation points are exactly the points where $S_\infty$ attains its maximum $r_-$, and they are isolated. A non-increasing function on $[0,\pi)$ whose maximum set is a finite set containing $0$ can attain its maximum only at $0$.
\end{proof}

\begin{theorem}[Refined reduction of existence]\label{thm_refined_reduction}
Let $\mathcal C^+=\mathcal C^+_\delta$ be as in \autoref{thm_oriented}, and suppose that the following hold along $\mathcal C^+$:
\begin{enumerate}
\item[(R1)] $\alpha$ stays bounded and $\Delta(\alpha)$ stays bounded away from $0$;
\item[(R2)] solutions remain in the lower region;
\item[(R3)] norms in $X=C^{3+s}_{\rm per}$ stay bounded on every part of $\mathcal C^+$ where $\min_\theta Q_\alpha(S)$ is bounded below;
\end{enumerate}
Condition (R1) is required uniformly in $\delta$: there are $\bar\alpha$ and $\bar\Delta>0$, independent of $\delta$, such that $\alpha\le\bar\alpha$ and $\Delta(\alpha)\ge\bar\Delta$ on $\mathcal C^+_\delta$ for all small $\delta$. Conditions (R2) and (R3) are required for each fixed small $\delta$, and the bound in (R3) may depend on $\delta$.
Then, for every $\delta<\delta_1:=\inf_{\mathcal C^+}|p_0|/\big(r_-\sqrt{Q_\alpha(1)}\big)$, the component $\mathcal C^+$ meets $\partial\mathcal O_\delta$ through the second constraint. Letting $\delta\downarrow0$ and applying \autoref{thm_limit} and \autoref{cor_extreme_shape}, a symmetric monotone extreme wave exists, with exactly one stagnation point per period, at which the free surface forms a symmetric $120^\circ$ Stokes corner.
\end{theorem}
\begin{proof}
By \autoref{thm_crossing}(a), case (b) of \autoref{thm_oriented} holds. In $\mathcal O_\delta$ we have $Q_\alpha(S)>\delta/S\ge\delta/r_-$. Together with (R1)--(R3), this bounds $\alpha$, the free-surface radius and the norm in $X$ on $\mathcal C^+_\delta$, so $\mathcal C^+_\delta$ is bounded and must meet $\partial\mathcal O_\delta$. In the lower region the first constraint stays above $\delta_1$ by \autoref{cor_first_constraint}, so the contact happens through the second constraint: $\min_\theta S\,Q_\alpha(S)=\delta$ at some point of $\mathcal C^+_\delta$.

Since (R1) is uniform, $\delta_1>0$ does not depend on $\delta$. Take $\delta_k\downarrow0$ and corresponding contact solutions. The bound of (R3) is not used in this limit; only the $\delta$-independent bounds of \autoref{prop_uniform} are. By (R1), after passing to a subsequence, $\alpha_k\to\alpha_\infty$ with $\Delta(\alpha_\infty)>0$. Since $Q_{\alpha_k}(S_k)\to0$ at some point, $\max S_k-r_-(\alpha_k)\to0$, and $r_-(\alpha_\infty)>1$ by \autoref{prop_uniform}(a). Thus the hypotheses of \autoref{thm_intro_limit} hold, and \autoref{cor_extreme_shape} applies.
\end{proof}

\subsection{An upper bound for the parameter}\label{sec_abound}

This subsection proves condition (R1) of \autoref{thm_refined_reduction} in the form of an explicit upper bound for $\alpha$ on monotone waves in the lower region. The argument rests on three exact structures:
\begin{itemize}
\item the Robin coefficient of the rotational field, which yields a ground-state inequality;
\item an angular flow-force invariant coming from rotational symmetry;
\item the unimodality of the circular comparison function.
\end{itemize}
Throughout, $(\Psi,S)\in\mathcal N$ lies in the lower region, so that $\Delta(\alpha)>0$ and $S<r_-(\alpha)$, and $Q_\alpha(S)>0$ on the free surface. We write $S_c=S(0)=\max S$ for the crest radius and $S_t=S(\pi)=\min S$ for the trough radius. Since $r_-$ decreases on $[\alpha_+,\infty)$ from $\alpha_+$ to $R_0$, we have $R_0<r_-\le\alpha_+$ there.

\subsubsection{The Robin coefficient and a ground-state inequality}

\begin{lemma}[Explicit Robin coefficient]\label{lem_robin_coefficient}
The coefficient in \autoref{lem_rotation_field} is
\[
c=\frac{a}{S}\,\frac{\alpha-S}{Q_\alpha(S)}-\kappa,\qquad \frac aS=\frac{S}{\sqrt{S^2+S'^2}} ,
\]
where $\kappa$ is the curvature of the free surface, positive for a circle.
\end{lemma}
\begin{proof}
In the computation of \autoref{lem_rotation_field}, $\partial_\nu w=a\,\tau\cdot D^2\Psi\nu+b\,\Psi_{\nu\nu}$. On the level set $\Psi_\nu=-|\nabla\Psi|$, and since $\Psi$ is harmonic, $\Psi_{\nu\nu}=-\tau\cdot D^2\Psi\tau=\kappa|\nabla\Psi|$. Moreover $\tau\cdot D^2\Psi\nu=-\partial_\tau|\nabla\Psi|$, which by the Bernoulli condition equals $-Q_\alpha'(R)(\mathbf e_R\cdot\tau)/(2p_0^2|\nabla\Psi|)$. Finally $\mathbf e_R\cdot\tau=-b/R$, $b=-w/|\nabla\Psi|$ and $p_0^2|\nabla\Psi|^2=Q_\alpha(S)$. Substituting gives the formula.
\end{proof}

\begin{remark}
At the laminar flow, $c=(\alpha-R_0)/c_0-1/R_0$. For the rotational field of the bifurcating mode, $w\propto(R-R^{-1})\sin\theta$, the Robin condition reads $c=(1+R_0^{-2})/(R_0-R_0^{-1})$. These two expressions agree exactly when $\alpha=\alpha_c$. The Robin formula therefore recovers the bifurcation condition of \cite{LiWang2026}; we also checked this numerically for several values of $(\gamma,p_0)$.
\end{remark}

\begin{lemma}[Ground-state inequality]\label{lem_picone}
Let $\phi\in C^1(\overline{\Omega_+})$ vanish near the bottom and near the axes, and let $\phi_{\rm top}(\theta)=\phi(S(\theta),\theta)$. Then
\[
\int_0^\pi\frac{S(\alpha-S)}{Q_\alpha(S)}\,\phi_{\rm top}^2\,d\theta\le\int_{\Omega_+}|\nabla\phi|^2+\int\kappa\,\phi_{\rm top}^2\,ds .
\]
Moreover,
\[
\Big|\int\kappa\,\phi_{\rm top}^2\,ds\Big|\le\int_0^\pi\phi_{\rm top}^2\,d\theta+\frac\pi2\,\mathrm{TV}_{(0,\pi)}\big(\phi_{\rm top}^2\big).
\]
\end{lemma}
\begin{proof}
By (N1) and (N2), $w<0$ on $\Omega_+$ and on the open top arc, so $\psi=\phi/w$ is $C^1$ on the support of $\phi$. Since $w$ is harmonic and satisfies $\partial_\nu w=cw$, integration by parts gives Picone's identity
\[
\int_{\Omega_+}|\nabla\phi|^2=\int_{\Omega_+}w^2|\nabla\psi|^2+\int_{\rm top}c\,\phi^2\,ds\ \ge\ \int_{\rm top}c\,\phi^2\,ds .
\]
By \autoref{lem_robin_coefficient} and the relation $(a/S)\,ds=S\,d\theta$, this is the first inequality.

For the curvature term, write $\kappa\,ds=d\vartheta$, where $\vartheta=\theta+\frac\pi2+\psi_0$ is the tangent angle and $\psi_0=-\arctan(S'/S)\in(-\frac\pi2,\frac\pi2)$. Integrating by parts against the compactly supported function $\phi_{\rm top}^2$ gives
\[
\int\phi_{\rm top}^2\,d\vartheta=\int\phi_{\rm top}^2\,d\theta-\int\psi_0\,(\phi_{\rm top}^2)'\,d\theta ,
\]
which yields the bound.
\end{proof}

\begin{proposition}[Definite trough depth]\label{prop_definite_depth}
Fix radii $1<s_2<s_1<R_0$ and angles $\delta_0,\delta>0$ with $\delta_0+\delta<\pi/2$, and set $\ell=\pi-2\delta_0$. Put
\[
C_*=\frac2\delta\log\alpha_++\frac{\ell\,(s_1+s_2)}{2(s_1-s_2)}+\ell+\frac{3\pi}2,\qquad
d_*:=\min\Big\{R_0-s_1,\ \frac{\ell-2\delta}{2C_*}\Big\}.
\]
Every monotone wave in the lower region with $\alpha\ge\alpha_+$ satisfies $r_--S_t\ge d_*$.
\end{proposition}
\begin{proof}
Let $I=[\delta_0+\delta,\ \pi-\delta_0-\delta]$ and $\theta_b=\pi-\delta_0-\delta$. Since $S$ is decreasing on $(0,\pi)$ and $r_-\ge R_0$ for $\alpha\ge\alpha_+$, we have $r_--S_t\ge r_--S(\theta_b)$.

If $S<s_1$ somewhere on $I$, then $S(\theta_b)<s_1$, so $r_--S_t\ge R_0-s_1$.

Otherwise $S\ge s_1$ on $I$. Take $\phi=\chi(\theta)\zeta(R)$ with the following choices:
\begin{itemize}
\item $\chi=1$ on $I$, $\chi=0$ outside $(\delta_0,\pi-\delta_0)$, and $\chi$ linear on the two ramps of width $\delta$, so that $\int\chi'^2\,d\theta=2/\delta$;
\item $\zeta=0$ for $R\le s_2$, $\zeta=1$ for $R\ge s_1$, and $\zeta$ linear in between.
\end{itemize}
Then $\phi$ is Lipschitz and vanishes in a neighbourhood of the bottom and of the axes, where $w=0$, so \autoref{lem_picone} applies. Its Dirichlet energy in $\Omega_+$ is bounded as follows.
\begin{itemize}
\item The angular part is $\iint\chi'^2\zeta^2R^{-1}\,dR\,d\theta\le\frac2\delta\log\frac{r_-}{s_2}\le\frac2\delta\log\alpha_+$, using $r_-\le\alpha_+$.
\item The radial part is $\iint\chi^2\zeta'^2R\,dR\,d\theta\le\ell\,\frac{s_1^2-s_2^2}{2(s_1-s_2)^2}=\frac{\ell(s_1+s_2)}{2(s_1-s_2)}$.
\end{itemize}
For the curvature term, write $\kappa\,ds=d\psi$ with $\psi$ the angle of the tangent. On a radial graph $\psi-\theta-\pi/2$ stays in $(-\pi/2,\pi/2)$. Integrating by parts against this function, and using that $\phi_{\rm top}$ vanishes at the ends of its support,
\[
\Big|\int\kappa\phi^2\,ds\Big|\le\int\phi_{\rm top}^2\,d\theta+\frac\pi2\,\mathrm{TV}(\phi_{\rm top}^2)\le\ell+\frac{3\pi}2 .
\]
Here we used $\mathrm{TV}(\phi_{\rm top}^2)\le\mathrm{TV}(\chi^2)+\mathrm{TV}(\zeta(S)^2)\le2+1$, which holds because $\zeta(S)^2$ is monotone in $\theta$.

In the lower region $(\alpha-S)/Q_\alpha(S)\ge1/(2(r_--S))$, because $r_+-S\le2(\alpha-S)$ when $S\le r_-$. The integrand on the left of \autoref{lem_picone} is non-negative, and $\phi_{\rm top}=1$ on $I$. Hence
\[
\frac{|I|}{2(r_--S(\theta_b))}\le\int_I\frac{d\theta}{2(r_--S)}\le C_* ,
\]
by monotonicity, with $|I|=\ell-2\delta$. Therefore $r_--S_t\ge r_--S(\theta_b)\ge(\ell-2\delta)/(2C_*)$.
\end{proof}

\subsubsection{An angular flow-force invariant}

\begin{proposition}[Angular flow force]\label{prop_angular_flow_force}
For every even classical solution that is smooth up to the free surface, and for any antiderivative $H$ of $\sigma\mapsto\sigma Q_\alpha(\sigma)$,
\begin{equation}\label{eq_angular_flow_force}
p_0^2\Big(\int_1^{S_t}R\,\Psi_R(R,\pi)^2\,dR-\int_1^{S_c}R\,\Psi_R(R,0)^2\,dR\Big)=H(S_c)-H(S_t)=\int_{S_t}^{S_c}\sigma\,Q_\alpha(\sigma)\,d\sigma .
\end{equation}
\end{proposition}
\begin{proof}
Rotations preserve both the Laplacian and $Q_\alpha(|x|)$. The corresponding Noether field
\[
\mathbf J=(x^\perp\cdot\nabla\Psi)\nabla\Psi-\tfrac12|\nabla\Psi|^2x^\perp
\]
is divergence free for harmonic $\Psi$, because $\operatorname{div}x^\perp=0$ and $Dx^\perp$ is antisymmetric. Integrate $\operatorname{div}\mathbf J=0$ over $\Omega_+$.
\begin{itemize}
\item On the bottom, $\Psi_\theta=0$ and $x^\perp\cdot\nu=0$, so the flux vanishes.
\item On the axes, $\Psi_\theta=0$, so $\mathbf J\cdot\mathbf e_\theta=-\frac R2\Psi_R^2$.
\item On the free surface, $\mathbf J\cdot\nu=\frac12|\nabla\Psi|^2x^\perp\cdot\nu$ and $x^\perp\cdot\nu\,ds=-SS'\,d\theta$, so by the Bernoulli condition the flux is $-\frac1{2p_0^2}\int_0^\pi Q_\alpha(S)SS'\,d\theta$.
\end{itemize}
Adding the three contributions gives \eqref{eq_angular_flow_force}. We also verified the underlying flux identity, before the Bernoulli substitution, numerically on non-circular harmonic configurations, to relative accuracy $10^{-10}$.
\end{proof}

\begin{remark}
Identity \eqref{eq_angular_flow_force} is a conjugate-flow condition between the crest and trough sections, the annular analogue of the flow-force constant of Benjamin and Lighthill \cite{BenjaminLighthill1954}. It rules out the scenario that obstructs a direct bound on $\alpha$: a wave that is nearly laminar near the crest and nearly the thin circular state of the second family near the trough, which would be a periodic \emph{bore}.
\end{remark}

\subsubsection{The crest--trough sandwich and the bound}

\begin{lemma}[Unimodality]\label{lem_unimodal}
For $\Delta(\alpha)>0$, the function $k(\sigma)=\sigma\log\sigma\sqrt{Q_\alpha(\sigma)}$ is strictly log-concave on $(1,r_-)$. Hence the set $\{G_\alpha\ge0\}\cap(1,r_-)=\{k\ge|p_0|\}$ is a closed interval $[\sigma_{\rm th},\sigma_{\rm dp}]$, possibly empty. For every non-circular solution in the lower region,
\[
\sigma_{\rm th}\le S_t,\qquad \sigma_{\rm dp}\le S_c<r_- .
\]
\end{lemma}
\begin{proof}
We have
\[
(\log k)'=\frac1\sigma+\frac1{\sigma\log\sigma}-\frac1{2(r_--\sigma)}-\frac1{2(r_+-\sigma)},
\]
and each of the four terms is strictly decreasing on $(1,r_-)$. By \autoref{prop_crest_trough}, $G_\alpha(S_t)\ge0$, so $S_t\in[\sigma_{\rm th},\sigma_{\rm dp}]$. Also $G_\alpha(S_c)\le0$, so either $S_c\le\sigma_{\rm th}$ or $S_c\ge\sigma_{\rm dp}$. The first option would give $S_c\le S_t$, which is impossible for a non-circular wave.
\end{proof}

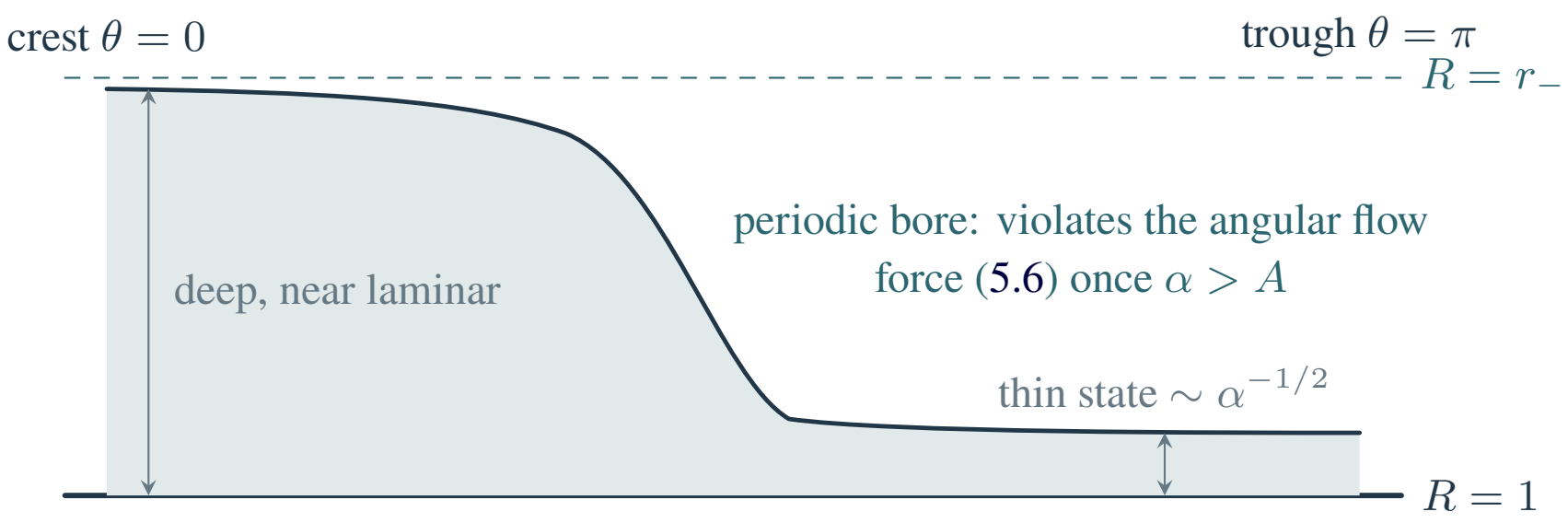
\begin{figure}[!ht]
\centering
\begin{tikzpicture}[x=1cm,y=1cm,>=stealth,font=\small,line cap=round,line join=round]
\draw[PaperAccent,dashed] (-.3,3.0) -- (9.3,3.0) node[right,font=\footnotesize] {$R=r_-$};
\draw[PaperInk,line width=1.1pt] (-.3,0) -- (9.3,0) node[right,font=\footnotesize] {$R=1$};
\fill[PaperAccent!14] (0,0) -- (0,2.92) .. controls (1.5,2.9) and (2.6,2.85) .. (3.3,2.6) .. controls (4.0,2.3) and (4.4,.85) .. (4.9,.55) .. controls (5.6,.45) and (8,.45) .. (9,.45) -- (9,0) -- cycle;
\draw[PaperInk,line width=.9pt] (0,2.92) .. controls (1.5,2.9) and (2.6,2.85) .. (3.3,2.6) .. controls (4.0,2.3) and (4.4,.85) .. (4.9,.55) .. controls (5.6,.45) and (8,.45) .. (9,.45);
\draw[PaperMuted,<->] (.3,0) -- (.3,2.92); \node[anchor=west,font=\footnotesize,text=PaperMuted] at (.35,1.45) {deep, near laminar};
\draw[PaperMuted,<->] (7.6,0) -- (7.6,.45); \node[anchor=south,font=\footnotesize,text=PaperMuted] at (7.6,.5) {thin state $\sim\alpha^{-1/2}$};
\node[font=\footnotesize,text=PaperInk] at (0,3.3) {crest $\theta=0$};
\node[font=\footnotesize,text=PaperInk] at (9,3.3) {trough $\theta=\pi$};
\node[text=PaperAccent,font=\footnotesize,align=center] at (7.0,1.75) {periodic bore: violates the angular flow\\force \eqref{eq_angular_flow_force} once $\alpha>A$};
\end{tikzpicture}
\caption{The large-$\alpha$ scenario excluded by \autoref{thm_alpha_bound} (schematic; angle horizontal, radius vertical). By \autoref{lem_unimodal} the crest stays within $O(\Delta^{-1/2})$ of $r_-$, and by \autoref{prop_definite_depth} the trough lies a definite distance $d_*$ below it. The angular flow force then forces $\sqrt\Delta\,d_*^{3/2}\lesssim|p_0|r_-^2\sqrt{r_+}$, which fails for large $\alpha$.}
\label{fig_bore}
\end{figure}

\begin{theorem}[Upper bound for $\alpha$]\label{thm_alpha_bound}
There is an explicit $A=A(\gamma,p_0)$ such that every monotone wave in $\mathcal N$ lying in the lower region with $\alpha\ge\alpha_+$ satisfies $\alpha<A$. One may take $A$ to be the smallest $\hat\alpha\ge\alpha_+$ such that the following three conditions hold for every $\alpha\ge\hat\alpha$:
\begin{equation}\label{eq_A_conditions}
\sqrt\Delta\ge\frac{p_0^2}{s_1^2\log^2s_1\,(R_0-1)},\qquad
\frac{p_0^2}{2\sqrt\Delta\,s_1^2\log^2s_1}\le\frac{d_*}2,\qquad
d_*^{3/2}>\frac{4|p_0|\,r_-^2\sqrt{r_+-1}}{3\sqrt\Delta} .
\end{equation}
\end{theorem}
\begin{proof}
Suppose $\alpha\ge A$, and let $D=r_--S_t$ and $d_c=r_--S_c$. By \autoref{prop_definite_depth}, $D\ge d_*$.

\emph{Crest depth.} The scenario to be excluded is sketched in \autoref{fig_bore}. We have $k(s_1)^2=s_1^2\log^2s_1\,Q_\alpha(s_1)\ge s_1^2\log^2s_1\cdot2\sqrt\Delta\,(R_0-1)/2\ge p_0^2$ by the first condition in \eqref{eq_A_conditions}, so $s_1\in[\sigma_{\rm th},\sigma_{\rm dp}]$ and $\sigma_{\rm dp}\ge s_1$. At $\sigma_{\rm dp}$,
\[
p_0^2=\sigma^2\log^2\sigma\,(r_--\sigma)(r_+-\sigma)\ge s_1^2\log^2s_1\,(r_--\sigma_{\rm dp})\,2\sqrt\Delta .
\]
By \autoref{lem_unimodal} and the second condition, $d_c\le r_--\sigma_{\rm dp}\le d_*/2\le D/2$.

\emph{Flow-force bound.} Drop the negative term in \eqref{eq_angular_flow_force}. Using $\int_1^{S_t}|\Psi_R|\,dR=1$ and the maximum-velocity principle (\autoref{thm_max_velocity}), which gives $|\nabla\Psi|\le r_-\sqrt{Q_\alpha(S_t)}/|p_0|$ because $Q_\alpha$ decreases on $[1,r_-]$, we obtain
\[
\int_{S_t}^{S_c}\sigma Q_\alpha(\sigma)\,d\sigma\le|p_0|\,r_-^2\sqrt{Q_\alpha(S_t)}\le|p_0|\,r_-^2\sqrt{(r_+-1)D} .
\]
On the other hand $\sigma Q_\alpha(\sigma)\ge2\sqrt\Delta\,(r_--\sigma)$ on $[1,r_-]$, so
\[
\int_{S_t}^{S_c}\sigma Q_\alpha(\sigma)\,d\sigma\ge\sqrt\Delta\,(D^2-d_c^2)\ge\tfrac34\sqrt\Delta\,D^2 .
\]
Combining, $D^{3/2}\le4|p_0|r_-^2\sqrt{r_+-1}/(3\sqrt\Delta)$. Together with $D\ge d_*$, this contradicts the third condition in \eqref{eq_A_conditions}.
\end{proof}

\begin{remark}[Size of $A$]
Each condition in \eqref{eq_A_conditions} holds for large $\alpha$, since $\sqrt\Delta\sim\alpha$ and $r_+\sim2\alpha$; the right-hand side of the third is $O(\alpha^{-1/2})$. Since that condition reads $d_*^{3/2}>O(\alpha^{-1/2})$, the bound scales like $A\sim d_*^{-3}$, and it pays to optimize the free parameters of \autoref{prop_definite_depth}. We take $\delta_0=0.02$, $s_2=1+\frac15(s_1-1)$, and maximize over $\delta$ and $s_1$. The script \texttt{verify\_A.py} then gives:
\begin{itemize}
\item $d_*=0.054$ and $A\approx9.1\cdot10^4$ for $\gamma=\frac12$, $p_0^2=0.554$ (where $\alpha_c=3.21$);
\item $d_*=0.035$ and $A\approx8.2\cdot10^3$ for $\gamma=0.2$, $p_0^2=0.0445$;
\item $d_*=0.055$ and $A\approx3.4\cdot10^7$ for $\gamma=0.95$, $p_0^2=36.3$.
\end{itemize}
The fixed choice $s_1=\frac{1+R_0}2$, $s_2=\frac{1+s_1}2$ with a window $[\pi/4,\pi/2]$ gives bounds $25$ to $50$ times larger. The bound is still far above the values $\alpha\le10$ observed numerically (\autoref{sec_numerics}). This is expected: it rests on a simple test function and on a flow-force estimate that discards the crest term.
\end{remark}

\subsubsection{Consequence for the existence problem}

\begin{corollary}[Updated reduction]\label{cor_updated_reduction}
Assume \eqref{eq_standing} and $\alpha_c>\alpha_+$; the latter holds if and only if $c_0>\sinh^2\gamma$. Conditions (R1) and (R2) of \autoref{thm_refined_reduction} follow from the single condition
\[
\text{(R1$^*$)}\qquad \inf_{\mathcal C^+_\delta}\alpha\ \ge\ \alpha_++\eta\quad\text{for some }\eta>0\text{ independent of }\delta .
\]
Consequently, (R1$^*$) and (R3) together imply the existence of a symmetric monotone extreme wave whose single stagnation crest per period is a $120^\circ$ Stokes corner.
\end{corollary}
\begin{proof}
At $\alpha_c$ the laminar radius satisfies $R_0<r_-(\alpha_c)$, so $\mathcal C^+$ starts in the lower region. Within $\mathcal O_\delta$, a solution can pass between the lower region and the region $\min S>r_+$ only through $\Delta\le0$, since $S$ cannot meet $[r_-,r_+]$. Under (R1$^*$), $\mathcal C^+$ therefore stays in the lower region, which is (R2). By \autoref{thm_alpha_bound}, $\alpha<A$ on $\mathcal C^+$, and $\Delta\ge\Delta(\alpha_++\eta)>0$, which is (R1). The final statement is \autoref{thm_refined_reduction}.
\end{proof}

\section{The crossing invariant, regularity and the global theorem}\label{sec_final}

This section proves that the continuum always crosses the laminar circle. This excludes secondary bifurcation, and it gives a parameter floor. It also reduces (R3) to a bound on the slope of the free surface.

\subsection{The laminar crossing invariant}\label{subsec_crossing}

\begin{lemma}[Radial crossing]\label{lem_radial_crossing}
The second circular family $\{G_\alpha(\sigma)=0,\ \sigma\neq R_0\}$ meets the laminar line exactly at
\[
\alpha_\circ=R_0+\frac{c_0(1+\gamma)}{\gamma R_0},
\]
which is where the radial mode of the linearization at the laminar flow becomes singular. Moreover $\alpha_\circ<\alpha_c$ for every $\gamma>0$. For $\alpha>\alpha_\circ$ with $\Delta(\alpha)>0$, the laminar radius $R_0$ is the larger root $\sigma_{\rm dp}$ of \autoref{lem_unimodal}, and the second-family circle lying in the lower region is the thin root $\sigma_{\rm th}<R_0$.
\end{lemma}
\begin{proof}
Circular solutions on the path \eqref{eq_parameter_path} are the roots of $G_\alpha$, equivalently of $k(\sigma)=|p_0|$ with $k$ as in \autoref{lem_unimodal}. Linearizing \eqref{eq_model} at the laminar flow in the class of radial functions gives a one-dimensional problem, which is singular exactly when $\sigma\mapsto G_\alpha(\sigma)$ has a double root at $R_0$, i.e.\ when $k'(R_0)=0$. At the same parameter the second root of $k=|p_0|$ merges with $R_0$. Explicitly, solving $G_\alpha(\sigma)=0$ for $\alpha$ gives
\[
\alpha(\sigma)=\frac{\sigma+R_0}2+\frac{c_0-P(\sigma)}{2(\sigma-R_0)},\qquad P(\sigma)=\frac{p_0^2}{\sigma^2\log^2\sigma},
\]
and $\alpha(\sigma)\to R_0-\frac12P'(R_0)=\alpha_\circ$ as $\sigma\to R_0$, since $P(R_0)=c_0$ and $P'(R_0)=-2c_0(1+\gamma)/(\gamma R_0)$.

We have $k(R_0)=|p_0|$ for every $\alpha$. Differentiating,
\[
(\log k)'(R_0)=\frac1{R_0}\Big(1+\frac1\gamma\Big)-\frac{\alpha-R_0}{c_0},
\]
which vanishes exactly at $\alpha_\circ$ and is negative for $\alpha>\alpha_\circ$. By unimodality, $R_0$ then lies on the decreasing side of $k$, so $R_0=\sigma_{\rm dp}$.

To compare with $\alpha_c$, note $\alpha_c-R_0=c_0/\sinh\gamma$ and $\alpha_\circ-R_0=c_0(1+\gamma)e^{-\gamma}/\gamma$. The inequality $\gamma e^\gamma/(1+\gamma)>\sinh\gamma$ is equivalent to
\[
\varphi(\gamma):=(\gamma-1)e^{2\gamma}+1+\gamma>0 .
\]
Here $\varphi(0)=\varphi'(0)=0$ and $\varphi''(\gamma)=4\gamma e^{2\gamma}>0$, so $\varphi>0$ for $\gamma>0$.
\end{proof}

\begin{theorem}[The continuum crosses the laminar circle]\label{thm_crossing}
Assume \eqref{eq_standing}, fix $\delta\in(0,\gamma)$, and let $b=(\alpha_c,R_0)$. Let $\mathcal N'$ be the set of waves in $\mathcal N$ that cross the laminar circle,
\[
\mathcal N'=\big\{(\Psi,S)\in\mathcal N:\ \min S<R_0<\max S\big\}.
\]
\begin{enumerate}
\item $\mathcal C_\delta\subset\mathcal N'\cup T(\mathcal N')\cup\{b\}$. In particular, case (b) of \autoref{thm_oriented} always holds, and $\mathcal C_\delta$ contains no circular solution other than $b$: there is no secondary bifurcation.
\item Every element of $\mathcal C^+_\delta$ in the lower region satisfies $\alpha>\alpha_\circ$.
\item If $\alpha_\circ>\alpha_+$, equivalently
\begin{equation}\label{eq_floor_condition}
c_0>\Big(\frac{\gamma R_0}{1+\gamma}\Big)^2=\frac{\gamma^2e^{2\gamma}}{(1+\gamma)^2},
\end{equation}
then all of $\mathcal C^+_\delta$ lies in the lower region, and $\alpha>\alpha_\circ$ on $\mathcal C^+_\delta$.
\end{enumerate}
\end{theorem}

\begin{figure}[!ht]
\centering
\begin{tikzpicture}
\begin{axis}[width=11cm,height=6.8cm,xmin=5.4,xmax=9.6,ymin=1.38,ymax=1.95,
  xlabel={$\alpha$},ylabel={radius},tick label style={font=\footnotesize},label style={font=\small},
  legend style={font=\footnotesize,at={(0.98,0.02)},anchor=south east,draw=PaperRule},legend cell align=left]
\addplot[PaperInk,line width=.6pt,domain=5.4:9.6] {1.64872};
\addlegendentry{laminar circle $R_0$}
\addplot[PaperAccent,line width=1.1pt] coordinates {(6.6044,1.6487) (6.6165,1.6536) (6.6305,1.6556) (6.6491,1.6576) (6.6738,1.6596) (6.7065,1.6618) (6.7498,1.6643) (6.8063,1.6671) (6.8661,1.6699) (6.9259,1.6725) (6.9858,1.6750) (7.0456,1.6775) (7.1054,1.6799) (7.1653,1.6823) (7.2251,1.6846) (7.2850,1.6870) (7.3448,1.6893) (7.4047,1.6916) (7.4646,1.6939) (7.5244,1.6962) (7.5843,1.6985) (7.6442,1.7007) (7.7041,1.7030) (7.7640,1.7053) (7.8239,1.7076) (7.8838,1.7098) (7.9437,1.7121) (8.0036,1.7144) (8.0635,1.7167) (8.1234,1.7190) (8.1833,1.7213) (8.2432,1.7237) (8.3031,1.7260) (8.3631,1.7284) (8.4230,1.7308) (8.4829,1.7333) (8.5428,1.7358) (8.6028,1.7383) (8.6627,1.7409) (8.7226,1.7436) (8.7826,1.7463) (8.8725,1.7506) (8.9324,1.7536) (8.9923,1.7568) (9.0523,1.7602) (9.1122,1.7639) (9.1721,1.7679) (9.2620,1.7752) (9.3220,1.7817)};
\addlegendentry{crest $\max S$ on $\mathcal C^+$}
\addplot[PaperAccent,line width=1.1pt,dashed] coordinates {(6.6044,1.6487) (6.6165,1.6431) (6.6305,1.6402) (6.6491,1.6371) (6.6738,1.6336) (6.7065,1.6296) (6.7498,1.6248) (6.8063,1.6191) (6.8661,1.6135) (6.9259,1.6083) (6.9858,1.6033) (7.0456,1.5986) (7.1054,1.5940) (7.1653,1.5896) (7.2251,1.5854) (7.2850,1.5813) (7.3448,1.5773) (7.4047,1.5734) (7.4646,1.5696) (7.5244,1.5658) (7.5843,1.5622) (7.6442,1.5587) (7.7041,1.5552) (7.7640,1.5518) (7.8239,1.5485) (7.8838,1.5452) (7.9437,1.5420) (8.0036,1.5389) (8.0635,1.5358) (8.1234,1.5328) (8.1833,1.5298) (8.2432,1.5269) (8.3031,1.5240) (8.3631,1.5212) (8.4230,1.5184) (8.4829,1.5157) (8.5428,1.5130) (8.6028,1.5103) (8.6627,1.5077) (8.7226,1.5052) (8.7826,1.5026) (8.8725,1.4989) (8.9324,1.4965) (8.9923,1.4941) (9.0523,1.4917) (9.1122,1.4894) (9.1721,1.4871) (9.2620,1.4837) (9.3220,1.4814)};
\addlegendentry{trough $\min S$ on $\mathcal C^+$}
\addplot[black!45,line width=1.3pt] coordinates {(6.3025,1.6527) (6.2733,1.6553) (6.2445,1.6580) (6.2161,1.6606) (6.1880,1.6632) (6.1601,1.6658) (6.1327,1.6685) (6.1055,1.6711) (6.0786,1.6737) (6.0521,1.6763) (6.0258,1.6790) (5.9999,1.6816) (5.9742,1.6842) (5.9488,1.6868) (5.9237,1.6895) (5.8989,1.6921) (5.8744,1.6947) (5.8501,1.6973) (5.8261,1.6999) (5.8024,1.7026) (5.7789,1.7052) (5.7557,1.7078) (5.7327,1.7104) (5.7100,1.7131) (5.6875,1.7157) (5.6653,1.7183) (5.6433,1.7209) (5.6216,1.7236) (5.6001,1.7262) (5.5788,1.7288) (5.5577,1.7314) (5.5369,1.7341) (5.5162,1.7367) (5.4958,1.7393) (5.4756,1.7419) (5.4557,1.7446) (5.4359,1.7472) (5.4163,1.7498)};
\addlegendentry{second circular family}
\addplot[black!45,line width=1.3pt,forget plot] coordinates {(9.5920,1.4717) (9.4892,1.4753) (9.3886,1.4789) (9.2902,1.4825) (9.1940,1.4861) (9.0997,1.4897) (9.0075,1.4934) (8.9172,1.4970) (8.8289,1.5006) (8.7423,1.5042) (8.6575,1.5078) (8.5745,1.5114) (8.4932,1.5150) (8.4135,1.5186) (8.3354,1.5222) (8.2588,1.5258) (8.1838,1.5294) (8.1102,1.5330) (8.0381,1.5366) (7.9674,1.5402) (7.8981,1.5438) (7.8300,1.5474) (7.7633,1.5510) (7.6978,1.5546) (7.6336,1.5582) (7.5706,1.5618) (7.5087,1.5654) (7.4480,1.5690) (7.3884,1.5726) (7.3298,1.5762) (7.2724,1.5798) (7.2159,1.5835) (7.1605,1.5871) (7.1061,1.5907) (7.0526,1.5943) (7.0000,1.5979) (6.9484,1.6015) (6.8977,1.6051) (6.8479,1.6087) (6.7989,1.6123) (6.7507,1.6159) (6.7034,1.6195) (6.6569,1.6231) (6.6111,1.6267) (6.5661,1.6303) (6.5219,1.6339) (6.4784,1.6375) (6.4356,1.6411) (6.3935,1.6447)};
\addplot[PaperInk,dotted,line width=.8pt] coordinates {(5.5000,2.0000) (5.5678,1.9933) (5.6356,1.9869) (5.7034,1.9808) (5.7712,1.9748) (5.8390,1.9691) (5.9068,1.9636) (5.9746,1.9583) (6.0424,1.9531) (6.1102,1.9482) (6.1780,1.9434) (6.2458,1.9387) (6.3136,1.9343) (6.3814,1.9299) (6.4492,1.9257) (6.5169,1.9216) (6.5847,1.9176) (6.6525,1.9138) (6.7203,1.9100) (6.7881,1.9064) (6.8559,1.9029) (6.9237,1.8995) (6.9915,1.8961) (7.0593,1.8929) (7.1271,1.8897) (7.1949,1.8866) (7.2627,1.8836) (7.3305,1.8807) (7.3983,1.8779) (7.4661,1.8751) (7.5339,1.8724) (7.6017,1.8697) (7.6695,1.8671) (7.7373,1.8646) (7.8051,1.8622) (7.8729,1.8597) (7.9407,1.8574) (8.0085,1.8551) (8.0763,1.8528) (8.1441,1.8506) (8.2119,1.8485) (8.2797,1.8464) (8.3475,1.8443) (8.4153,1.8423) (8.4831,1.8403) (8.5508,1.8384) (8.6186,1.8365) (8.6864,1.8346) (8.7542,1.8328) (8.8220,1.8310) (8.8898,1.8293) (8.9576,1.8276) (9.0254,1.8259) (9.0932,1.8242) (9.1610,1.8226) (9.2288,1.8210) (9.2966,1.8195) (9.3644,1.8179) (9.4322,1.8164) (9.5000,1.8149)};
\addlegendentry{stagnation circle $r_-(\alpha)$}
\addplot[PaperInk,dashed,line width=.5pt,forget plot] coordinates {(6.3476,1.38) (6.3476,1.95)};
\node[font=\footnotesize,anchor=south east] at (axis cs:6.3476,1.39) {$\alpha_\circ$};
\addplot[only marks,mark=*,mark size=2pt,color=PaperInk,forget plot] coordinates {(6.0955,1.6721)};
\node[font=\footnotesize,anchor=north,align=center] at (axis cs:5.95,1.655) {singular\\circle};
\addplot[only marks,mark=*,mark size=2pt,color=PaperAccent,forget plot] coordinates {(6.6044,1.6487)};
\node[font=\footnotesize,anchor=north west] at (axis cs:6.6044,1.6487) {};
\node[font=\footnotesize,text=PaperAccent] at (axis cs:6.45,1.69) {$\alpha_c$};
\end{axis}
\end{tikzpicture}
\caption{The crossing invariant (\autoref{thm_crossing}) for $\gamma=0.5$, $c_0=0.95\,e^{2\gamma}$, computed. Along the primary continuum, which starts at $(\alpha_c,R_0)$, the crest stays above and the trough below the laminar circle. The second circular family has a mode-$1$ singular circle on the deep side at $\alpha_*=6.096<\alpha_\circ$ (\autoref{rem_secondary_points}). Waves branching from it lie above $R_0$, so they cannot be connected to the primary continuum. By \autoref{thm_crossing}(b), the primary continuum also never enters $\alpha\le\alpha_\circ$ while it stays in the lower region, i.e.\ below the stagnation circle $r_-(\alpha)$.}
\label{fig_crossing}
\end{figure}
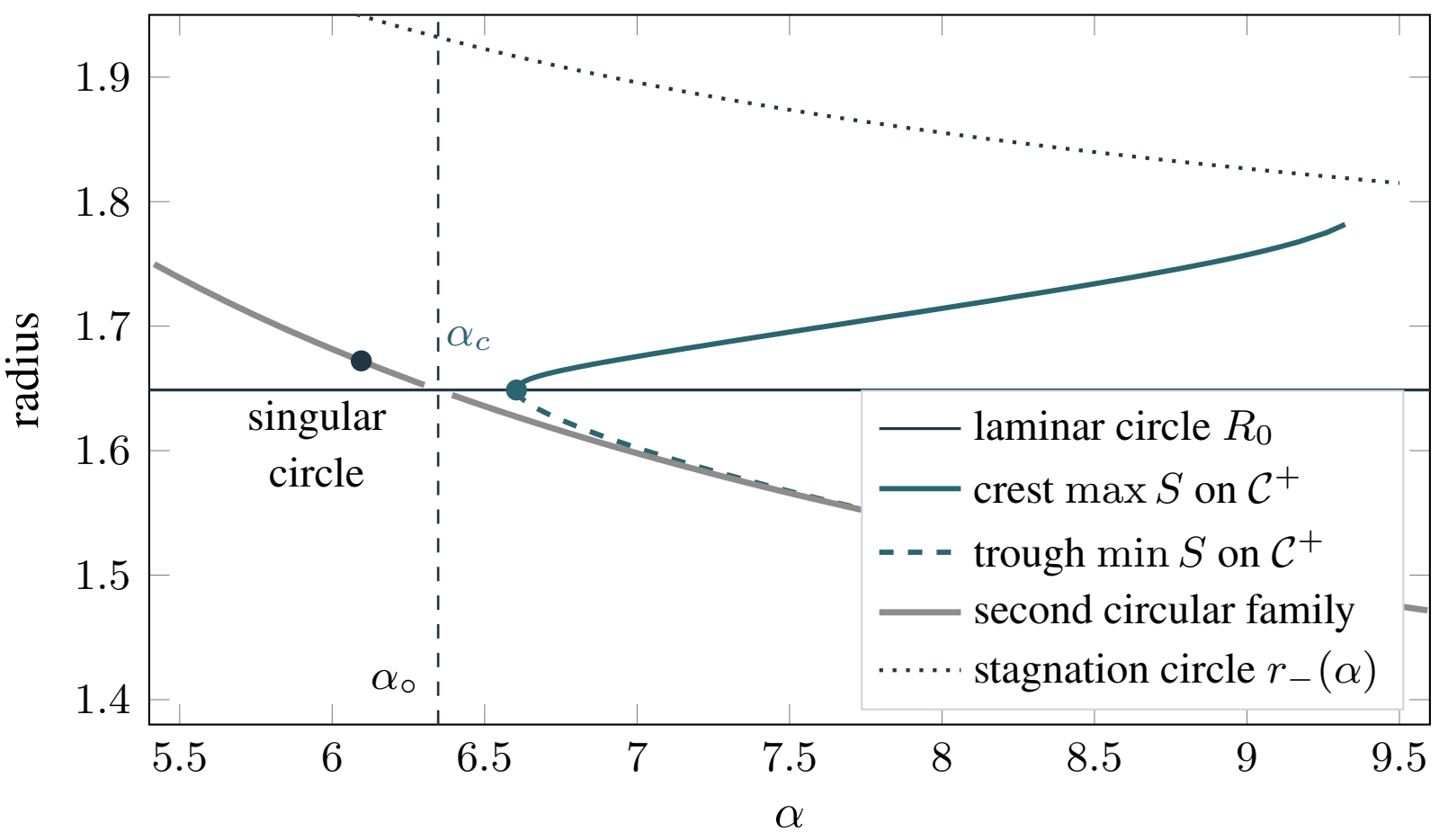

The invariant has a simple meaning. Waves on secondary branches, which bifurcate from circles $S\equiv\sigma\ne R_0$, lie entirely on one side of the laminar circle $R=R_0$. Waves on the primary continuum always straddle it. The singular circles of \autoref{rem_secondary_points} therefore exist, but the primary continuum cannot reach them. \autoref{fig_crossing} shows a computed example.

\begin{proof}
(a) Put $M'=\mathcal N'\cup T(\mathcal N')$. Since $\min S$ and $\max S$ depend continuously on the solution in $X$, and $\mathcal N$ is open (\autoref{prop_nodal_open_closed}(a)), $M'$ is open. Near $b$ the continuum consists of the local curve $S=R_0(1+\varepsilon\cos\theta)+o(\varepsilon)$, which lies in $M'\cup\{b\}$. Let $E_1=(M'\cap\mathcal C_\delta)\cup\{b\}$. We claim that $E_1$ is closed in $\mathcal C_\delta$. Let $x\in\mathcal C_\delta\setminus\{b\}$ be a limit of elements of $M'\cap\mathcal C_\delta$. Then $\min S\le R_0\le\max S$ at $x$.
\begin{itemize}
\item If $x$ is non-circular, then $x\in\mathcal N\cup T(\mathcal N)$ by \autoref{prop_nodal_open_closed}(b). It is analytic up to the free surface, and $G_\alpha(R_0)=0$ on the path \eqref{eq_parameter_path}. The strict comparison of \autoref{prop_crest_trough} gives $G_\alpha(\min S)>0>G_\alpha(\max S)$, so $\min S\ne R_0\ne\max S$, and hence $x\in M'$.
\item If $x$ is circular, $S\equiv\sigma$, then $\min S=\max S=\sigma$ forces $\sigma=R_0$. By \autoref{prop_nodal_open_closed}(d), applied to $\mathcal N$ or $T(\mathcal N)$, this gives $\alpha=\alpha_c$, so $x=b$, which was excluded.
\end{itemize}
Thus $E_1$ is closed. As in the proof of \autoref{thm_oriented}, $\mathcal C_\delta\setminus E_1$ stays away from $b$ and is closed. By connectedness, $\mathcal C_\delta=E_1$.

(b) Let $x\in\mathcal C^+_\delta\setminus\{b\}$ lie in the lower region, so that $\Delta(\alpha)>0$ and $\max S<r_-$. By (a), $S_t=\min S<R_0<\max S<r_-$, so $(1,R_0]\subset(1,r_-)$. Moreover $G_\alpha(S_t)>0$, i.e.\ $k(S_t)>|p_0|$, with $k$ as in \autoref{lem_unimodal}, and $k(R_0)=|p_0|$ because $G_\alpha(R_0)=0$ on the path \eqref{eq_parameter_path}. Suppose $\alpha\le\alpha_\circ$. By the computation in the proof of \autoref{lem_radial_crossing},
\[
(\log k)'(R_0)=\frac1{R_0}\Big(1+\frac1\gamma\Big)-\frac{\alpha-R_0}{c_0}\ge0 .
\]
Since $(\log k)'$ is strictly decreasing on $(1,r_-)$, by \autoref{lem_unimodal}, we get $(\log k)'>0$ on $(1,R_0)$. Hence $k<k(R_0)=|p_0|$ on $(1,R_0)$, contradicting $k(S_t)>|p_0|$. So $\alpha>\alpha_\circ$. At $b$ itself, $\alpha_c>\alpha_\circ$ by \autoref{lem_radial_crossing}.

(c) Let $W$ be the set of $x\in\mathcal C^+_\delta$ with $\alpha(x)>\alpha_\circ$ and $\max S<r_-(\alpha(x))$. It contains $b$, since $R_0<r_-(\alpha_c)$. It is open in $\mathcal C^+_\delta$. To see that it is closed, let $x$ be a limit of elements of $W$. Then $\alpha(x)\ge\alpha_\circ>\alpha_+$, so $\Delta>0$, and $\max S\le r_-$. Since $x\in\mathcal O_\delta$, $Q_\alpha(S)>0$ on the free surface, so in fact $\max S<r_-$. By (b), $\alpha(x)>\alpha_\circ$. So $x\in W$. Since $\mathcal C^+_\delta$ is connected, $W=\mathcal C^+_\delta$. Finally, $\alpha_\circ>\alpha_+$ means $c_0(1+\gamma)/(\gamma R_0)>\sqrt{c_0}$, which is \eqref{eq_floor_condition}.
\end{proof}

\begin{corollary}[An unconditional parameter floor]\label{cor_floor}
Assume \eqref{eq_standing} and \eqref{eq_floor_condition}. Then, for every $\delta\in(0,\gamma)$, every element of $\mathcal C^+_\delta$ lies in the lower region and satisfies
\[
\alpha_\circ<\alpha<A(\gamma,p_0),\qquad \Delta(\alpha)>\Delta(\alpha_\circ)>0 .
\]
The bounds are independent of $\delta$.
\end{corollary}
\begin{proof}
Combine \autoref{thm_crossing}(c) with \autoref{thm_alpha_bound}.
\end{proof}

Condition \eqref{eq_floor_condition} is stronger than $c_0>\sinh^2\gamma$, because $\gamma e^\gamma/(1+\gamma)>\sinh\gamma$ (\autoref{lem_radial_crossing}). In terms of $f=c_0e^{-2\gamma}$, it reads $f>(\gamma/(1+\gamma))^2$, which is less than $1/4$ in the range $\gamma<1$ of the standing assumptions. For example, it holds whenever $f\ge1/4$. Of the $20$ parameter pairs in the scan of \autoref{subsec_scan}, $19$ satisfy it.

\begin{remark}[Singular circles that the continuum cannot reach]\label{rem_secondary_points}
At a circle $S\equiv\sigma$ of the second family, the mode-$1$ linearization is singular exactly when
\[
\alpha-\sigma=\frac{2\sigma\,Q_\alpha(\sigma)}{\sigma^2-1} .
\]
This is the laminar condition $\alpha_c-R_0=2c_0R_0/(R_0^2-1)$ with $(R_0,c_0)$ replaced by $(\sigma,Q_\alpha(\sigma))$, and it reduces to it at $\sigma=R_0$. Solving it numerically along the second family produces genuine secondary bifurcation points $(\sigma_*,\alpha_*)$ on the deep side, in the lower region, with $\alpha_+<\alpha_*<\alpha_\circ$:
\begin{center}
\begin{tabular}{cccccccc}
\hline
$\gamma$ & $p_0^2$ & $R_0$ & $\sigma_*$ & $\alpha_*$ & $\alpha_+$ & $\alpha_\circ$ & $\alpha_c$\\
\hline
0.3 & 0.0896 & 1.3499 & 1.3554 & 3.0637 & 2.0890 & 3.1038 & 3.1441\\
0.5 & 0.5542 & 1.6487 & 1.6750 & 3.0527 & 2.5518 & 3.1326 & 3.2137\\
0.5 & 1.7549 & 1.6487 & 1.6721 & 6.0955 & 3.2557 & 6.3476 & 6.6044\\
0.7 & 2.4174 & 2.0138 & 2.0920 & 3.3484 & 3.1167 & 3.4809 & 3.6175\\
0.9 & 17.787 & 2.4596 & 2.6040 & 5.1842 & 4.3648 & 5.5751 & 5.9957\\
\hline
\end{tabular}
\end{center}
Secondary branches of monotone waves therefore issue from the second circular family, at parameters below $\alpha_\circ$. Their waves lie on the deep side, above the laminar circle. By \autoref{thm_crossing}, the primary continuum never meets them. The existence of these circles also shows that no inequality holding separately at each $\alpha$ can give a parameter floor: the floor of \autoref{thm_crossing} uses the crossing invariant, which is a property of the primary continuum.
\end{remark}

\subsection{The flow angle and the regularity condition}

\begin{lemma}[Flow angle]\label{lem_flow_angle}
For $(\Psi,S)\in\mathcal N$, let $F=\Psi_t-i\Psi_\theta$ in the cylinder and set $\vartheta=\pi-\arg F$. Then:
\begin{itemize}
\item $\vartheta$ is harmonic in $\Omega_+$;
\item $\vartheta=0$ on the bottom and on both axes;
\item $0<\vartheta<\pi/2$ in $\Omega_+$ and on the open top arc;
\item on the free surface $\vartheta=\arctan|s'(\theta)|$, where $t=s(\theta)$ is the surface in the cylinder; that is, $\vartheta$ is the slope angle of the free surface.
\end{itemize}
Consequently $\sup_{\Omega_+}\vartheta=\sup_{(0,\pi)}\arctan|s'|$.
\end{lemma}
\begin{proof}
$F$ is holomorphic in $t+i\theta$ and does not vanish, since $\Psi_t<0$. In $\Omega_+$, (N1) and class (L) give $\operatorname{Re}F=\Psi_t<0$ and $\operatorname{Im}F=-\Psi_\theta>0$, so $\arg F\in(\pi/2,\pi)$. The function $\arg F$ is harmonic, so $\vartheta$ is harmonic as well.

On the bottom and on the axes $\Psi_\theta=0$, so $\vartheta=0$. Differentiating $\Psi(s(\theta),\theta)=0$ gives $s'=-\Psi_\theta/\Psi_t$, hence $\tan\vartheta=|\Psi_\theta|/|\Psi_t|=|s'|$ on the top.
\end{proof}

\begin{proposition}[(R3) follows from a slope bound]\label{prop_slope_regularity}
Fix $\delta>0$, a compact parameter interval, and $\eta>0$. Suppose that solutions of $\mathcal N$ in $\mathcal O_\delta$, in the lower region and with $\alpha$ in that interval, satisfy the slope bound $\vartheta\le\pi/2-\eta$. Then their hodograph heights $h$ are bounded in $X=C^{3+s}_{\rm per}(\overline D)$, by a constant depending only on $\delta$, $\eta$ and the parameter interval.
\end{proposition}
\begin{proof}
Recall $h(q,p)=R-1$ with $p=-\Psi$, so that $h_p=1/|\Psi_R|=R/|\Psi_t|$ and $|h_q|=R\tan\vartheta$ on $\overline{\Omega_+}$, and by evenness on the whole period.

\emph{Lower bound for $|\Psi_t|$.} The function $\Psi_t$ is harmonic and negative, so $|\Psi_t|$ attains its infimum on the boundary. On the bottom, $|\Psi_t|\ge1/\log S_c\ge1/\log r_-$ by the comparison of \autoref{prop_crest_trough}. On the top,
\[
|\Psi_t|=|\nabla_c\Psi|\cos\vartheta=\frac{S\sqrt{Q_\alpha(S)}}{|p_0|}\cos\vartheta\ge\frac{\sqrt\delta\,\sin\eta}{|p_0|},
\]
using $SQ_\alpha(S)>\delta$ in $\mathcal O_\delta$.

\emph{Bounds on $H$.} In the variables of \autoref{lem_hodograph}, $H_p=h_p/R=1/|\Psi_t|$, so $H_p\le C$. The first constraint gives $H_p=(h+1)h_p/R^2\ge\delta/r_-^2$. Finally $|H_q|=|h_q|/R=\tan\vartheta$, so $|H_q|\le\cot\eta$. Thus $H$ is uniformly Lipschitz, with $H_p$ bounded above and below.

\emph{Regularity.} With these bounds, \eqref{eq_hod} is a quasilinear equation with smooth coefficients in $\nabla H$. It is uniformly elliptic, since its determinant is $H_p^2$ and its trace is $1+H_q^2+H_p^2$. The top condition \eqref{eq_hod_bc} has the form $F(H,\nabla H)=0$ with $\partial F/\partial H_p=-2H_pG(H)\le-2\delta^2/r_-^2<0$, using $G(H)=S\cdot SQ_\alpha(S)\ge\delta$. It is therefore uniformly oblique on the Lipschitz bounds just obtained. The bottom condition is $H=0$, and the problem is periodic in $q$.

The Lieberman--Trudinger theory for nonlinear oblique problems \cite{LiebermanTrudinger1986} gives a uniform $C^{1,\beta}(\overline D)$ bound for some $\beta\in(0,1)$; this step uses only the Lipschitz bounds. The equation is linear in $D^2H$ with $C^\beta$ coefficients. The linearized boundary operator has $C^{\beta}$ coefficients and is oblique. Schauder theory for oblique problems \cite[Theorem~6.30 and Lemma~6.29]{GilbargTrudinger2001} (see also \cite{ADN1959}) gives a uniform $C^{2,\beta}$ bound.

Next, $v=H_q$ solves the linear oblique problem \eqref{eq_v}--\eqref{eq_v_bc}, whose coefficients are now bounded in $C^{\beta}$, with $v=0$ on the bottom. The same Schauder estimate bounds $v$ in $C^{2,\beta}$. The equation then expresses $H_{pp}$ through $H_{qq}$ and $H_{qp}$, giving a $C^{3,\beta}$ bound for $H$. One more iteration, with coefficients in $C^{1,\beta}$, gives bounds in $C^{3,s}$ for any $s<1$. Since $h=e^H-1$, the same bounds hold for $h$ in $X$.
\end{proof}

\subsection{The final reduction}

By \autoref{thm_crossing}, secondary bifurcation never occurs, and $\alpha>\alpha_\circ$ holds automatically in the lower region. What remains of the parameter floor is the distance from the degenerate value $\alpha_+$.

\begin{theorem}[Existence of the extreme wave modulo two a priori bounds]\label{thm_final_reduction}
Assume \eqref{eq_standing} and $c_0>\sinh^2\gamma$. Suppose there is $\eta>0$ such that, for all small $\delta\in(0,\gamma)$, every element of $\mathcal C^+_\delta$ satisfies:
\begin{enumerate}
\item[(A)] \emph{(parameter floor)} $\alpha\ge\alpha_++\eta$;
\item[(B)] \emph{(no vertical tangents)} $\vartheta\le\pi/2-\eta$, i.e.\ the free-surface slope in the cylinder satisfies $|s'|\le\cot\eta$.
\end{enumerate}
Then there exists a symmetric extreme wave, monotone from crest to trough, with exactly one stagnation point per period. At that point the free surface forms a corner of $120^\circ$, symmetric about the radial axis, with the explicit Stokes asymptotics of \autoref{thm_intro_stokes}. Condition (A) holds automatically, with $\eta=\alpha_\circ-\alpha_+$, if $c_0$ satisfies \eqref{eq_floor_condition}.
\end{theorem}
\begin{proof}
\emph{Step 1: lower region.} By (A), $\Delta(\alpha)\ge\Delta(\alpha_++\eta)>0$ on $\mathcal C^+_\delta$. At the bifurcation point $R_0<r_-(\alpha_c)$. Within $\mathcal O_\delta$ and $\{\Delta>0\}$ the lower region and the region $\{\min S>r_+\}$ are separated, because $S$ cannot meet $[r_-,r_+]$. Since $\mathcal C^+_\delta$ is connected, it lies entirely in the lower region.

\emph{Step 2: conditions (R1)--(R3).} By \autoref{thm_alpha_bound}, $\alpha\in[\alpha_++\eta,A]$ on $\mathcal C^+_\delta$, which gives (R1). Step 1 is (R2). Condition (B) and \autoref{prop_slope_regularity} give (R3). Now apply \autoref{thm_refined_reduction} and \autoref{cor_extreme_shape}. The last statement is \autoref{cor_floor}.
\end{proof}

\begin{corollary}[Extreme wave or overturning]\label{cor_dichotomy}
Assume \eqref{eq_standing}, and either \eqref{eq_floor_condition} or, more generally, $c_0>\sinh^2\gamma$ together with condition (A). Then at least one of the following holds:
\begin{enumerate}
\item there exists a symmetric monotone extreme wave with a single $120^\circ$ Stokes crest per period, as in \autoref{thm_final_reduction};
\item there are $\delta_k\downarrow0$ and monotone waves in $\mathcal C^+_{\delta_k}$ whose free-surface slope angle tends to $\pi/2$, i.e.\ the onset of overturning.
\end{enumerate}
\end{corollary}
\begin{proof}
If (b) fails, there are $\eta>0$ and $\delta_0>0$ such that $\vartheta\le\pi/2-\eta$ on $\mathcal C^+_\delta$ for every $\delta<\delta_0$. This is condition (B), and \autoref{thm_final_reduction} gives (a).
\end{proof}

\begin{remark}
Alternative (b) is not an artefact of the method. For flat waves with constant vorticity, global continua can end in overhanging profiles rather than at a stagnation point \cite{CSV2016}, and the flows of \cite{LiWang2026} carry constant vorticity in the absolute frame. Deciding between (a) and (b) in the annular model therefore requires an annular Nekrasov-type analysis of the flow angle, or a parameter restriction under which overturning is excluded.
\end{remark}

\subsection{An unconditional global theorem}

The conditional statements above can be packaged into an unconditional description of the global continuum. The parameter floor is no longer assumed; failures of each condition appear instead as explicit alternatives.

\begin{theorem}[Global alternatives for the monotone continuum]\label{thm_global}
Assume \eqref{eq_standing} and $c_0>\sinh^2\gamma$, equivalently $\alpha_c>\alpha_+$. For each $\delta\in(0,\gamma)$ let $\mathcal C_\delta$ be the continuum of \cite[Theorem~2.2]{LiWang2026}. Then at least one of the following holds.
\begin{enumerate}
\item[(I)] \emph{Extreme wave.} There exists a symmetric extreme wave, monotone from crest to trough, with exactly one stagnation point per period. The crest lies on the circle $R=r_-(\alpha_\infty)$ and is a corner of $120^\circ$, symmetric about the radial axis, with the Stokes asymptotics of \autoref{thm_intro_stokes}.
\item[(II)] \emph{Overturning.} There are $\delta_k\downarrow0$ and monotone waves in $\mathcal C^+_{\delta_k}$ whose maximal free-surface slope angle tends to $\pi/2$.
\item[(III)] \emph{Degeneration of the effective gravity.} $\liminf_{\delta\downarrow0}\inf_{\mathcal C^+_\delta}\alpha\le\alpha_+$.
\end{enumerate}
Secondary bifurcation does not occur (\autoref{thm_crossing}). If $c_0$ satisfies \eqref{eq_floor_condition}, then (III) does not occur either, so (I) or (II) holds. In every case $\alpha>\alpha_\circ$ on the part of $\mathcal C^+_\delta$ in the lower region, and $\alpha<A(\gamma,p_0)$ on it (\autoref{thm_alpha_bound}).
\end{theorem}
\begin{proof}
If (III) fails, there are $\eta>0$ and $\delta_0>0$ such that $\alpha\ge\alpha_++\eta$ on $\mathcal C^+_\delta$ for every $\delta<\delta_0$. This is (A), and \autoref{cor_dichotomy} gives (I) or (II). Under \eqref{eq_floor_condition}, \autoref{cor_floor} gives $\alpha>\alpha_\circ>\alpha_+$ on $\mathcal C^+_\delta$ uniformly in $\delta$, so (III) fails.
\end{proof}

\begin{figure}[!ht]
\centering
\begin{tikzpicture}[x=1cm,y=1cm,>=stealth,font=\small,line cap=round,line join=round]
\tikzset{lab/.style={font=\footnotesize,fill=white,inner sep=1.2pt}}
\draw[->,PaperMuted] (0,0) -- (11.2,0) node[right] {$\alpha$};
\draw[->,PaperMuted] (0,0) -- (0,4.9) node[above] {amplitude};
\draw[PaperInk,line width=1pt] (0,0) -- (11,0);
\node[anchor=north,lab,text=PaperInk] at (9.3,-.08) {laminar flows};
\draw[PaperAccent,dashed] (1.2,0) -- (1.2,3.9) node[above,font=\footnotesize] {$\alpha_+$};
\draw[PaperMuted,dashed] (3.9,0) -- (3.9,3.9) node[above,font=\footnotesize] {$\alpha_\circ$};
\draw[PaperMuted,dashed] (0,4.5) -- (10.8,4.5);
\node[anchor=south west,lab,text=PaperMuted] at (.1,4.52) {upper bound $\alpha<A$ (off scale)};
\draw[PaperMuted,line width=.8pt] (.4,.32) .. controls (1.6,.5) and (3.0,.55) .. (3.9,0);
\node[anchor=north,lab,text=PaperMuted] at (2.0,-.08) {second circular family (never reached)};
\fill[PaperAccent] (5.0,0) circle (2pt) node[below=2pt,font=\footnotesize,text=PaperAccent] {$\alpha_c$};
\draw[PaperAccent,line width=1.1pt] (5.0,0) .. controls (5.6,1.0) and (6.6,1.8) .. (7.4,2.6);
\draw[PaperAccent,line width=1.1pt,dotted] (7.4,2.6) .. controls (8.0,3.1) and (8.6,3.4) .. (9.3,3.7);
\fill[PaperInk] (9.3,3.7) circle (2pt);
\node[anchor=west,lab,text=PaperInk,align=left] at (9.45,3.7) {(I) extreme wave,\\$120^\circ$ Stokes crest};
\draw[PaperAccent,line width=1.1pt,dotted] (7.4,2.6) .. controls (7.6,3.2) and (7.4,3.7) .. (7.0,4.0);
\node[anchor=south,lab,text=PaperInk] at (7.0,4.02) {(II) overturning};
\draw[PaperAccent,line width=1.1pt,dotted] (6.6,1.8) .. controls (5.0,2.6) and (3.0,2.85) .. (1.3,2.75);
\node[anchor=south west,lab,text=PaperInk] at (1.3,2.85) {(III) $\alpha\downarrow\alpha_+$, only if $\alpha_\circ\le\alpha_+$};
\end{tikzpicture}
\caption{The alternatives of \autoref{thm_global} (schematic). The monotone continuum $\mathcal C^+$ leaves the laminar line at $\alpha_c$. It ends in an extreme wave with a Stokes crest (I), overturns (II), or its parameter degenerates to $\alpha_+$ (III). Alternative (III) is possible only when $\alpha_\circ\le\alpha_+$; the schematic shows this case. Every wave of $\mathcal C^+$ crosses the laminar circle, so the continuum never meets the second circular family (\autoref{thm_crossing}). Return to the laminar line is impossible (\autoref{thm_oriented}), and $\alpha$ is bounded above (\autoref{thm_alpha_bound}).}
\label{fig_global}
\end{figure}
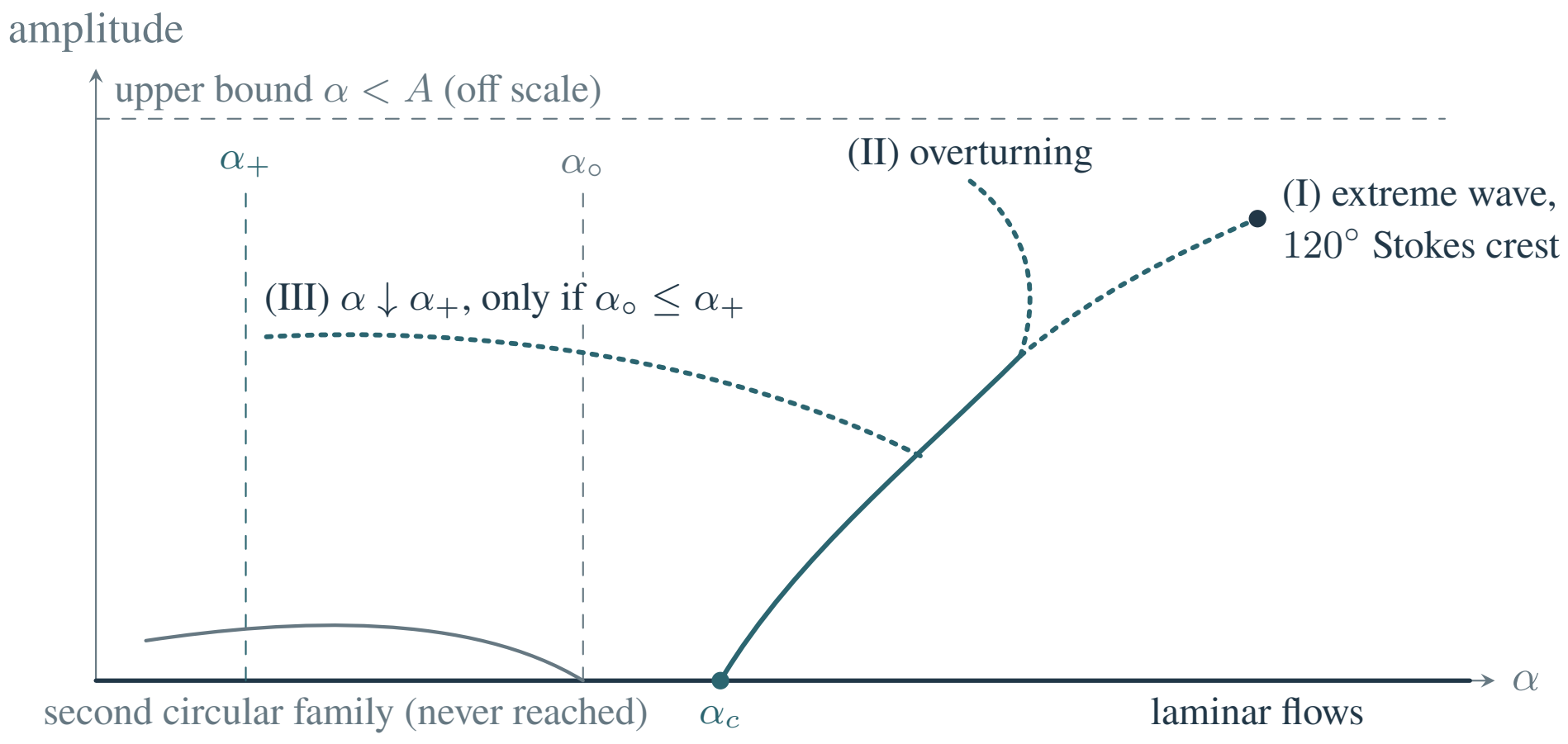

The alternatives are illustrated in \autoref{fig_global}.

\begin{corollary}[The degenerate alternative]\label{cor_IV}
In \autoref{thm_global}, suppose that (III) holds through monotone waves of $\mathcal C^+_{\delta_k}$ in the lower region whose crests approach stagnation, i.e.\ $\alpha_k\to\alpha_+$ and $\max S_k-r_-(\alpha_k)\to0$. Then there exists a symmetric extreme wave at $\alpha=\alpha_+$, monotone from crest to trough, whose single crest per period is a symmetric right-angle corner.
\end{corollary}
\begin{proof}
Apply \autoref{thm_degenerate_limit} with $\alpha_*=\alpha_+>R_0>1$.
\end{proof}

\begin{theorem}[The critical case $c_0=\sinh^2\gamma$]\label{thm_critical}
Assume \eqref{eq_standing} and $c_0=\sinh^2\gamma$, so that $\alpha_c=\alpha_+$. Then the bifurcation is supercritical: $\alpha=\alpha_++\alpha_2\varepsilon^2+O(\varepsilon^4)$ with $\alpha_2>0$. Moreover, the conclusion of \autoref{thm_global} holds with (III) replaced by the following.
\begin{enumerate}
\item[(III$'$)] There are $\varepsilon_0>0$, $\delta_k\downarrow0$ and $x_k\in\mathcal C^+_{\delta_k}$ such that $\|x_k-(\alpha_c,R_0)\|_{\mathbb R\times X}\ge\varepsilon_0$ and $\limsup_k\alpha(x_k)\le\alpha_+$.
\end{enumerate}
\end{theorem}
\begin{proof}
\autoref{thm_supercritical} gives $\alpha_2>0$; by \autoref{rem_all_gamma}, \autoref{thm_alpha2} applies at $c_0=\sinh^2\gamma$.

Suppose (II) and (III$'$) both fail. Since (III$'$) fails, for every $\varepsilon_0>0$ there are $\eta>0$ and $\delta_0>0$ such that every $x\in\mathcal C^+_\delta$, $\delta<\delta_0$, at distance at least $\varepsilon_0$ from $b=(\alpha_c,R_0)$ satisfies $\alpha(x)\ge\alpha_++\eta$. Choose $\varepsilon_0$ so small that, within distance $\varepsilon_0$ of $b$, the continuum consists of the local curve. On that curve $\alpha\ge\alpha_+$, $S$ is close to $R_0<\alpha_+=r_-(\alpha_+)$, and $Q_\alpha(S)$ is close to $c_0$.

In particular $\alpha\ge\alpha_+$ on all of $\mathcal C^+_\delta$, and the continuum lies in the lower region. The rest of the proof of \autoref{thm_global} goes through with three remarks.
\begin{itemize}
\item \autoref{prop_slope_regularity} and \autoref{thm_alpha_bound} need $\alpha\ge\alpha_+$ but not $\Delta>0$.
\item Contact with $\partial\mathcal O_\delta$ through the second constraint, for small $\delta$, can only occur at distance at least $\varepsilon_0$ from $b$, because $Q_\alpha(S)$ is close to $c_0$ near $b$. There $\Delta\ge\Delta(\alpha_++\eta)>0$.
\item The limit theorem therefore applies exactly as before and gives (I).
\end{itemize}
\end{proof}

\begin{remark}[How generic is each alternative?]\label{rem_alternatives}
The alternatives have different status.
\begin{itemize}
\item Secondary bifurcation is excluded by \autoref{thm_crossing}, although singular circles exist: on a grid of $2140$ parameter pairs with $0.05\le\gamma\le0.99$, $p_0^2=f\gamma^2e^{4\gamma}$ with $0.02\le f\le0.99$ (so that \eqref{eq_standing} holds) and $c_0>\sinh^2\gamma$, a singular circle on the deep side, in the lower region, with parameter in $(\alpha_+,\alpha_\circ)$ exists for $1885$ of them (\autoref{rem_secondary_points}). The branches issuing from these circles never meet the primary continuum.
\item Alternative (II) has a flat counterpart: the global theory of \cite{CSV2016} allows overhanging profiles, and overhanging and touching waves with constant vorticity are studied in \cite{HurWheeler2022}. In the annular model, \autoref{thm_band} below shows that (II) can occur only together with (III).
\item Alternative (III) is the annular analogue of a critical Froude number. It is excluded under \eqref{eq_floor_condition}, which holds, for instance, whenever $c_0\ge e^{2\gamma}/4$.
\end{itemize}
\autoref{thm_global} is the annular counterpart of the global theorems for flat waves: extreme waves, overhanging waves, or loss of the structure that drives the continuation.
\end{remark}

\begin{remark}[Comparison with the flat Stokes problem]
In the classical theory for flat irrotational waves, the two analogous facts were the hardest a priori estimates on the way to the extreme wave:
\begin{itemize}
\item a bound for the Froude number along the branch;
\item a bound on the slope of the free surface, below $\pi/2$; Amick proved that it is in fact below $31.15^\circ$ \cite{Amick1987}, using Nekrasov's integral equation \cite{Nekrasov1921,AFT1982}. For waves with vorticity, slope bounds are known for favourable and for small adverse vorticity \cite{StraussWheeler2016,SoStrauss2018}.
\end{itemize}
In the annular model, \autoref{lem_flow_angle} provides the analogue of the Nekrasov unknown: a harmonic angle that vanishes on the bottom and on the axes. Its boundary condition on the free surface follows from the Bernoulli relation $\log|F|=\frac12\log\big(S^2Q_\alpha(S)/p_0^2\big)$ through the Cauchy--Riemann equations. This structure leads to the annular Nekrasov equation of \autoref{sec_nekrasov}, through which (B) is proved under explicit conditions (\autoref{thm_cone_invariance}).
\end{remark}

\subsection{The case \texorpdfstring{$c_0<\sinh^2\gamma$}{c0 < sinh^2 gamma}}\label{subsec_window}

If $c_0<\sinh^2\gamma$, then $\alpha_-<\alpha_c<\alpha_+$, so the bifurcation takes place in the window where $\Delta<0$. There $Q_\alpha\ge-\Delta(\alpha)>0$, and stagnation is impossible. The following theorem describes how the continuum can leave the window.

\begin{theorem}[Continuum in the window $\Delta<0$]\label{thm_window}
Assume \eqref{eq_standing} and $c_0<\sinh^2\gamma$. Let $J=[a,b]\subset(\alpha_-,\alpha_+)$ with $a<\alpha_c<b$, and let $\eta\in(0,\pi/2)$. There is $\delta_J>0$ such that for every $\delta\in(0,\delta_J)$ at least one of the following holds:
\begin{enumerate}
\item $\mathcal C^+_\delta$ contains a wave with $\alpha\in\{a,b\}$;
\item $\mathcal C^+_\delta$ contains a wave with $\alpha\in J$ whose maximal slope satisfies $\vartheta>\pi/2-\eta$.
\end{enumerate}
Moreover, no solution with $\alpha\in(\alpha_-,\alpha_+)$ has a stagnation point. Consequently a stagnation point can develop along $\mathcal C^+_\delta$ only after $\alpha$ has reached $\alpha_+$ or $\alpha_-$, where the corner, if it forms, is a right angle (\autoref{thm_degenerate}).
\end{theorem}
\begin{proof}
Put $\Delta_J=\min_J(-\Delta)>0$. For $\alpha\in J$ we have $Q_\alpha\ge\Delta_J$, so there is no stagnation point.

\emph{Bound on the crest.} \autoref{prop_crest_trough} gives $S_c^2\log^2S_c\,Q_\alpha(S_c)\le p_0^2$, hence $S_c\le\bar\sigma_J$, where $\bar\sigma_J^2\log^2\bar\sigma_J=p_0^2/\Delta_J$.

\emph{The constraints.} The second constraint satisfies $SQ_\alpha(S)\ge\Delta_J$. For the first, \autoref{cor_first_constraint} gives $(h+1)h_p\ge\delta_{1,J}:=|p_0|/\big(\bar\sigma_J\max_{[1,\bar\sigma_J]\times J}Q^{1/2}\big)$. Set $\delta_J=\min\{\Delta_J,\delta_{1,J}\}/2$. Then no solution with $\alpha\in J$ lies on $\partial\mathcal O_\delta$ for $\delta<\delta_J$.

\emph{The alternatives.} Suppose (a) and (b) both fail. By \autoref{thm_crossing}(a), case (b) of \autoref{thm_oriented} holds. Since $\alpha$ is continuous on the connected set $\mathcal C^+_\delta$, takes the value $\alpha_c\in(a,b)$, and avoids $a$ and $b$, every element has $\alpha\in(a,b)$. Moreover every element has slope at most $\pi/2-\eta$. The proof of \autoref{prop_slope_regularity} then applies with $r_-$ replaced by $\bar\sigma_J$ and $\delta$ by $\Delta_J$. On the top,
\[
|\Psi_t|\ge\sqrt{\Delta_J}\sin\eta/|p_0| .
\]
Also $H_p\ge\delta_{1,J}/\bar\sigma_J^2$ and $G(H)\ge\Delta_J$. This bounds $\mathcal C^+_\delta$ in $(\alpha_s,\infty)\times X$; note that $J\subset(\alpha_s,\infty)$, because $\alpha_--\alpha_s=(R_0-\sqrt{c_0})^2/(2R_0)>0$. But $\mathcal C^+_\delta$ does not meet $\partial\mathcal O_\delta$, which contradicts \autoref{thm_oriented}(b).
\end{proof}

\section{The slope bound and the direction of bifurcation}\label{sec_towards}

\autoref{thm_final_reduction} reduces the existence of the extreme wave to a parameter floor (A) and a slope bound (B). This section treats both. For (B), we derive an annular Nekrasov equation for the flow angle and prove a slope bound in an invariant Keady--Norbury cone: $\vartheta<\kappa_\alpha\pi$ with $\kappa_\alpha<\frac12$ in general, and $\vartheta<\pi/3$ under an algebraic condition. This gives the existence theorem, \autoref{thm_existence}. For (A), we compute the direction of bifurcation in closed form.

\subsection{An annular Nekrasov equation}\label{sec_nekrasov}

Condition (B) of \autoref{thm_final_reduction} is a uniform bound on the slope of the free surface. For flat irrotational waves the corresponding bound comes from Nekrasov's integral equation for the flow angle \cite{Nekrasov1921,Toland1978,AFT1982,Amick1987}. In this subsection we derive the annular counterpart of that equation. We show that it reproduces the bifurcation condition, prove two exact consequences, and locate the point of maximal slope. We do not prove (B) here; the last subsection explains what remains.

\subsubsection{Conformal parametrization}
Let $(\Psi,S)$ be a solution of class (L) with $S\in C^{1,\beta}(\mathbb T)$. The fluid region is doubly connected, so there are $\varrho>1$ and a conformal map $f$ from the annulus $A_\varrho=\{1<|\zeta|<\varrho\}$ onto $\Omega$. This map takes $|\zeta|=1$ to the bottom and $|\zeta|=\varrho$ to the free surface, and by the Kellogg--Warschawski theorem it is $C^{1,\beta}$ up to the boundary. The function $\Psi\circ f$ is harmonic in $A_\varrho$, equals $1$ on $|\zeta|=1$ and $0$ on $|\zeta|=\varrho$, so
\[
\Psi\circ f(\zeta)=1-\frac{\log|\zeta|}{L},\qquad L:=\log\varrho .
\]
For even solutions we normalize $f(\varrho)=S(0)$, so that $f$ commutes with complex conjugation and the crest corresponds to $\zeta=\varrho$.

We use the strip variable $\tau=\log\zeta=s+i\phi$, with $0<s<L$ and $\phi\in\mathbb T$, and the cylinder variable $w=\log z=t+i\theta$. Write $w=W(\tau)$, so that $W(\tau+2\pi i)=W(\tau)+2\pi i$ and $W'\neq0$.

\begin{lemma}[Flow angle in conformal variables]\label{lem_conformal_angle}
Let $F=\Psi_t-i\Psi_\theta$ as in \autoref{lem_flow_angle}. Then
\[
F\big(W(\tau)\big)=-\frac{1}{L\,W'(\tau)},\qquad \vartheta=\arg W'(\tau).
\]
The function $\Lambda=\log W'=U+i\vartheta$ is holomorphic in the strip, with $\vartheta=0$ on $s=0$. On the free surface $s=L$,
\[
|W'|=\frac{|p_0|}{L\sqrt{h_\alpha(S)}},\qquad \frac{d}{d\phi}\log S=-|W'|\sin\vartheta ,
\]
where $h_\alpha(\sigma)=\sigma^2Q_\alpha(\sigma)$ as in \autoref{cor_energy_sandwich}.
\end{lemma}
\begin{proof}
The function $\Phi(\tau)=1-\tau/L$ is holomorphic and $\operatorname{Re}\Phi=\Psi\circ f$. Hence $F=d\Phi/dw=\Phi'(\tau)/W'(\tau)=-1/(LW')$. Since $\vartheta=\pi-\arg F$ and $\arg(-1/W')=\pi-\arg W'$, we get $\vartheta=\arg W'$.

On $s=0$ we have $\operatorname{Re}W=0$, so $\partial_\phi\operatorname{Re}W=-\operatorname{Im}W'=0$. Thus $W'$ is real there, and positive by orientation, so $\vartheta=0$.

On $s=L$, the Bernoulli condition gives $|F|^2=|\nabla_c\Psi|^2=S^2Q_\alpha(S)/p_0^2$, which is the formula for $|W'|$. Finally $\log S(\phi)=\operatorname{Re}W(L+i\phi)$, and $\partial_\phi\operatorname{Re}W=\operatorname{Re}(iW')=-|W'|\sin\vartheta$.
\end{proof}

\subsubsection{The equation}
For $\vartheta=\sum_{k\ge1}a_k\sin k\phi$ define
\[
\mathcal K_L\vartheta=\sum_{k\ge1}k\coth(kL)\,a_k\sin k\phi,\qquad \mathcal C_L\vartheta=\sum_{k\ge1}\coth(kL)\,a_k\cos k\phi .
\]
These are the Dirichlet-to-Neumann operator and the conjugation operator of the strip with Dirichlet data on $s=0$. Since $\frac{d}{d\phi}\mathcal C_L\vartheta=-\mathcal K_L\vartheta$, for $L=\infty$ they reduce to the operators $|D|$ and the Hilbert transform that appear in Nekrasov's equation.

\begin{theorem}[Annular Nekrasov equation]\label{thm_nekrasov}
Let $(\Psi,S)\in\mathcal N$, and let $\vartheta(\phi)$ be the flow angle on the free surface in the conformal parametrization. Then $\vartheta$ is odd, $0<\vartheta<\pi/2$ on $(0,\pi)$, and
\begin{equation}\label{eq_nekrasov}
\mathcal K_L\vartheta=\mathfrak g_\alpha(S)\,\sin\vartheta,\qquad \mathfrak g_\alpha(\sigma)=-\frac{|p_0|\,\sigma\,h_\alpha'(\sigma)}{2L\,h_\alpha(\sigma)^{3/2}} .
\end{equation}
Here $S=S(\phi)$ solves
\begin{equation}\label{eq_nekrasov_S}
\frac{d}{d\phi}\log S=-\frac{|p_0|\,\sin\vartheta}{L\sqrt{h_\alpha(S)}},\qquad S(0)=S_c .
\end{equation}
\end{theorem}
\begin{proof}
Oddness and the bounds on $\vartheta$ come from \autoref{lem_flow_angle}. The function $\vartheta(s,\phi)$ is harmonic in the strip, vanishes on $s=0$ and is odd in $\phi$, so
\[
\vartheta(s,\phi)=\sum_k a_k\frac{\sinh ks}{\sinh kL}\sin k\phi ,\qquad \partial_s\vartheta(L,\cdot)=\mathcal K_L\vartheta .
\]
By the Cauchy--Riemann equations for $\Lambda=U+i\vartheta$, $\partial_s\vartheta=-\partial_\phi U$. On $s=L$, \autoref{lem_conformal_angle} gives
\[
U=\log|p_0|-\log L-\tfrac12\log h_\alpha(S),
\]
hence
\[
-\partial_\phi U=\frac{h_\alpha'(S)}{2h_\alpha(S)}\,S\,\partial_\phi\log S=-\frac{S\,h_\alpha'(S)}{2h_\alpha(S)}\,|W'|\sin\vartheta=\mathfrak g_\alpha(S)\sin\vartheta .
\]
Equation \eqref{eq_nekrasov_S} is the last identity of \autoref{lem_conformal_angle}.
\end{proof}

\begin{remark}[Consistency checks]\label{rem_nekrasov_checks}
\emph{Bifurcation.} At the laminar flow $L=\gamma$ and $S\equiv R_0$. Linearizing \eqref{eq_nekrasov} in the mode $\sin\phi$ gives $\coth\gamma=\mathfrak g_\alpha(R_0)=\big(R_0(\alpha-R_0)-c_0\big)/c_0$, that is, $\alpha-R_0=c_0/\sinh\gamma$. This is exactly $\alpha=\alpha_c$.

\emph{Numerics.} On the computed wave of run~(b) with $N=16384$, the residual of \eqref{eq_nekrasov} is $4\times10^{-5}$ relative to $\max|\mathfrak g_\alpha(S)\sin\vartheta|$. The equation is not imposed by the scheme.
\end{remark}

\subsubsection{Two exact consequences}

\begin{proposition}[Total weight]\label{prop_total_weight}
For every $(\Psi,S)\in\mathcal N$,
\[
\int_0^\pi\mathfrak g_\alpha\big(S(\phi)\big)\sin\vartheta(\phi)\,d\phi=\frac12\log\frac{h_\alpha(S_t)}{h_\alpha(S_c)} .
\]
In particular the total weight is positive, by \autoref{cor_energy_sandwich}. It tends to $+\infty$ exactly when $h_\alpha(S_c)\to0$, that is, when the crest approaches the stagnation circle.
\end{proposition}
\begin{proof}
Integrate \eqref{eq_nekrasov} over $(0,\pi)$. Since $\mathcal K_L\vartheta=-\frac{d}{d\phi}\mathcal C_L\vartheta$, the left-hand side equals $(\mathcal C_L\vartheta)(0)-(\mathcal C_L\vartheta)(\pi)$. The function $\Lambda_0(\tau)=\sum_ka_k\cosh(k\tau)/\sinh(kL)$ is holomorphic, $2\pi i$-periodic, and has imaginary part $\vartheta$; its real part on $s=L$ is $\mathcal C_L\vartheta$. The difference $\Lambda-\Lambda_0$ is holomorphic and periodic with vanishing imaginary part, hence a real constant. A term of the form $-ic\tau$ is excluded because it would add $-cs$, which is not odd in $\phi$. So the integral equals $U(L,0)-U(L,\pi)$, which is $\frac12\log h_\alpha(S_t)-\frac12\log h_\alpha(S_c)$ by \autoref{lem_conformal_angle}. On the computed wave of run~(b) the two sides are $2.0888$ and $2.0894$; the discrepancy is concentrated at the crest.
\end{proof}

\begin{proposition}[Location of the maximal slope]\label{prop_max_slope}
Let $(\Psi,S)\in\mathcal N$. The supremum of $\vartheta$ over the closed half-strip is attained only on the free surface, at points $\phi_*\in(0,\pi)$ with
\[
\mathfrak g_\alpha\big(S(\phi_*)\big)>0,\qquad\text{equivalently}\qquad h_\alpha'\big(S(\phi_*)\big)<0 .
\]
Since $h_\alpha$ increases and then decreases on $(1,r_-)$ (\autoref{cor_energy_sandwich}), this means that $S(\phi_*)$ exceeds the maximum point of $h_\alpha$, if there is one: the steepest point lies between the crest and the point where $S$ crosses that maximum point.
\end{proposition}
\begin{proof}
$\vartheta$ is harmonic, positive in the open half-strip, and zero on the bottom and on the axes $\phi\in\{0,\pi\}$. By the strong maximum principle its maximum is attained only on the top. At such a point the Hopf lemma gives $\partial_s\vartheta>0$, and \eqref{eq_nekrasov} gives $\mathfrak g_\alpha(S)\sin\vartheta>0$. Since $0<\vartheta<\pi/2$, $\mathfrak g_\alpha(S)>0$. By the definition of $\mathfrak g_\alpha$, this is $h_\alpha'(S)<0$.
\end{proof}

\begin{lemma}[Positive kernel]\label{lem_kernel}
The inverse of $\mathcal K_L$ on odd functions has the representation
\[
\mathcal K_L^{-1}f(\phi)=\int_0^\pi K_L(\phi,\phi')f(\phi')\,d\phi',\qquad K_L(\phi,\phi')=\frac2\pi\sum_{k\ge1}\frac{\tanh kL}{k}\sin k\phi\,\sin k\phi',
\]
and $K_L\ge0$ on $(0,\pi)^2$, with a logarithmic singularity on the diagonal. Consequently every $(\Psi,S)\in\mathcal N$ satisfies
\[
\vartheta(\phi)=\int_0^\pi K_L(\phi,\phi')\,\mathfrak g_\alpha\big(S(\phi')\big)\sin\vartheta(\phi')\,d\phi' .
\]
\end{lemma}
\begin{proof}
The series formula follows from the definition of $\mathcal K_L$, which is invertible on odd functions since $k\coth kL\ge1/L$. For positivity, let $f\ge0$ be odd, continuous and not identically zero on $(0,\pi)$. Let $u$ be harmonic in the half-strip $(0,L)\times(0,\pi)$, with $u=0$ on $s=0$ and on $\phi\in\{0,\pi\}$, and $\partial_su=f$ on $s=L$. Then $u(L,\cdot)=\mathcal K_L^{-1}f$. If $u$ had a negative minimum, it would be attained on the open top, where the Hopf lemma gives $\partial_su<0$, contradicting $\partial_su=f\ge0$. So $u\ge0$, and $u>0$ in the interior and on the open top by the strong maximum principle and the Hopf lemma. Since this holds for every such $f$, the kernel $K_L$ is non-negative almost everywhere, and hence everywhere off the diagonal by continuity there. The integral form of the equation follows by applying $\mathcal K_L^{-1}$ to \eqref{eq_nekrasov}.
\end{proof}

\subsubsection{A slope bound from the Keady--Norbury cone}\label{subsec_cone}
\autoref{prop_max_slope} and \autoref{lem_kernel} reduce condition (B) to a uniform estimate for the integral equation of \autoref{lem_kernel}. That equation has the structure of Nekrasov's equation in finite depth: a positive kernel, and a weight $\mathfrak g_\alpha(S)$ that is positive near the crest and blows up there as the wave approaches stagnation. For flat waves, Keady and Norbury obtained the bound $\pi/3$ by working in an invariant cone whose decisive condition is that $\Phi(t)/\sin(t/2)$ is nonincreasing \cite{KeadyNorbury1978,Toland1996}. We show that the same mechanism works in the annulus. The decisive hypothesis is the cone condition (C2) on the flow angle. In the annulus it gives a slope bound below $90^\circ$ for every parameter with $\Delta>0$, and the bound $60^\circ$ under an explicit algebraic condition (C1).

\emph{The Nekrasov variables.} For $(\Psi,S)\in\mathcal N$ put
\[
V(\phi)=\frac{|p_0|}{L}\int_0^\phi\sin\vartheta,\qquad F(V)=\tfrac12\log h_\alpha\big(S(\phi)\big).
\]
By \eqref{eq_nekrasov_S}, $dV=-\frac{\sqrt{h_\alpha(S)}}{S}\,dS$. Hence $V$ is the increasing function $V=\int_S^{S_c}\sqrt{h_\alpha(\sigma)}\,\sigma^{-1}d\sigma$ of $S$, and
\[
F'(V)=-\frac{S\,h_\alpha'(S)}{2h_\alpha(S)^{3/2}}=\frac L{|p_0|}\,\mathfrak g_\alpha(S),\qquad
\frac{d}{d\phi}F(V(\phi))=F'(V)\,V'=\mathfrak g_\alpha(S)\sin\vartheta .
\]
\autoref{lem_kernel} therefore reads $\vartheta=K_L\big[(F\circ V)'\big]$. In flat deep water, $h$ is linear in the height and $F(V)=\frac13\log(a+bV)$, which is Nekrasov's equation. Here $F$ is determined by $h_\alpha$. Near a stagnation point $h_\alpha\sim2\sigma^2\sqrt\Delta\,(r_--\sigma)$, so again $VF'(V)\to\frac13$; this is the source of the $30^\circ$ slope at the Stokes corner.

\begin{lemma}[Structural inequality]\label{lem_structural}
Let $\alpha>R_0$ with $\Delta(\alpha)>0$, and put
\[
\delta_\alpha=\frac{\Delta(\alpha)}{(\alpha-1)^2},\qquad \kappa_\alpha=\frac1{2+\delta_\alpha}<\frac12 .
\]
Then every $(\Psi,S)\in\mathcal N$ in the lower region at the parameter $\alpha$ satisfies
\[
(F\circ V)'(\phi)\le\kappa_\alpha\,\frac{V'(\phi)}{V(\phi)}\qquad(0<\phi<\pi).
\]
If moreover
\begin{enumerate}
\item[(C1)] $(\sigma h_\alpha'(\sigma))'\le0$ for $\sigma\in[s_h,r_-]$,
\end{enumerate}
where $s_h$ is the maximum point of $h_\alpha$ on $(1,r_-)$ ($s_h=1$ if $h_\alpha$ is decreasing there), then the same holds with $\kappa_\alpha$ replaced by $\frac13$. Since $(\sigma h_\alpha')'=2\sigma\,(8\sigma^2-9\alpha\sigma+4\lambda)$, condition (C1) holds if and only if $8\sigma^2-9\alpha\sigma+4\lambda\le0$ at $\sigma=r_-$ and at $\sigma=s_h$. The first condition is equivalent to $r_-\le\frac56\alpha$, i.e.\ to $\sqrt{\Delta}\ge\alpha/6$. The second holds automatically if $s_h>1$.
\end{lemma}
\begin{proof}
Where $h_\alpha'(S)\ge0$, the left-hand side is $\le0$ and there is nothing to prove. The set where $h_\alpha'(S(\phi))<0$ is an interval $(0,\phi_h)$, on which $S\in(s_h,S_c]$. There, with $p=-h_\alpha'>0$, the claim with a constant $\kappa$ reads $V\le2\kappa h_\alpha^{3/2}/(Sp)$. Put
\[
\Xi_\kappa(S)=\frac{2\kappa h_\alpha(S)^{3/2}}{Sp(S)}-\int_S^{S_c}\frac{\sqrt{h_\alpha(\sigma)}}\sigma\,d\sigma .
\]
Then $\Xi_\kappa(S_c)>0$, and a direct computation gives
\[
\Xi_\kappa'(S)=\frac{\sqrt{h_\alpha}}S\,\big(1-3\kappa+\kappa J\big),\qquad J:=\frac{2h_\alpha\,(\sigma h_\alpha')'}{\sigma h_\alpha'^2}\Big|_{\sigma=S}.
\]
Hence $\Xi_\kappa>0$ on $(s_h,S_c]$, which is the claim, as soon as $J\le3-1/\kappa$ on $(s_h,S_c)$. Under (C1) we have $J\le0=3-1/\frac13$.

In general we show that $1-J\ge\delta_\alpha$ wherever $1<\sigma<r_-$ and $h_\alpha'(\sigma)<0$. Since $3-1/\kappa_\alpha=1-\delta_\alpha$, this proves the lemma. Write $h_\alpha=\sigma^2Q_\alpha$, $a=r_--\sigma>0$ and $b=r_+-\sigma>a$. Then $Q_\alpha=ab$, $Q_\alpha'=-(a+b)$ and $Q_\alpha''=2$. Moreover
\[
h_\alpha'=-\sigma D,\qquad D:=\sigma(a+b)-2ab,\qquad (\sigma h_\alpha')'=\sigma\big(4Q_\alpha+5\sigma Q_\alpha'+2\sigma^2\big).
\]
So $h_\alpha'<0$ means $D>0$, and a short computation gives
\[
1-J=\frac{8a^2b^2+6ab\,D+\sigma^2(b-a)^2}{D^2}\ \ge\ \frac{\sigma^2(b-a)^2}{D^2}\ \ge\ \frac{(b-a)^2}{(a+b)^2},
\]
because $0<D<\sigma(a+b)$. Finally $b-a=2\sqrt\Delta$ and $a+b=2(\alpha-\sigma)<2(\alpha-1)$, so $1-J\ge\delta_\alpha$.

For the formula in (C1), note that $h_\alpha=\sigma^4-2\alpha\sigma^3+2\lambda\sigma^2$, so $\sigma h_\alpha'=4\sigma^4-6\alpha\sigma^3+4\lambda\sigma^2$. The quadratic $8\sigma^2-9\alpha\sigma+4\lambda$ is convex, so it is $\le0$ on an interval if and only if it is $\le0$ at both ends. At $r_-$, $4\lambda=4\alpha r_--2r_-^2$ gives the value $r_-(6r_--5\alpha)$. At $s_h>1$, $h_\alpha'(s_h)=0$ gives $2s_h^2=3\alpha s_h-2\lambda$, so the value is $3\alpha s_h-4\lambda$. Since $s_h=\frac14\big(3\alpha-\sqrt{9\alpha^2-16\lambda}\big)$, the inequality $3\alpha s_h\le4\lambda$ is equivalent to $\sqrt{9\alpha^2-16\lambda}\le3\alpha$, that is, to $\lambda\ge0$. This holds on the path \eqref{eq_parameter_path} for $\alpha>\alpha_s$.
\end{proof}

\begin{remark}[Uniformity as $\Delta\to0$]\label{rem_kappa_crest}
The constant $\kappa_\alpha$ tends to $\frac12$ as $\Delta\to0$. The proof also gives a bound that stays uniform as $\Delta\to0$, as long as the crest stays away from stagnation. On the range $\sigma\le S_c$ we have $a\ge d_c:=r_--S_c$ and $b\ge a$, so $ab/(a+b)\ge a/2$. Since $D<\sigma(a+b)$,
\[
1-J\ \ge\ \frac{8a^2b^2}{D^2}\ \ge\ \frac{8a^2b^2}{\sigma^2(a+b)^2}\ \ge\ \frac{2d_c^2}{r_-^2}.
\]
Hence \autoref{lem_structural} holds with $\kappa=1/\big(2+\max\{\delta_\alpha,\,2d_c^2/r_-^2\}\big)$. Along a family of waves whose parameter approaches $\alpha_+$, the slope bound of \autoref{thm_cone_slope} can therefore degenerate only if the crests approach stagnation.
\end{remark}

Note that $\delta_\alpha$ is increasing in $\alpha>R_0$. Indeed,
\[
\frac{d}{d\alpha}\frac{(\alpha-R_0)^2-c_0}{(\alpha-1)^2}=\frac{2\big((\alpha-R_0)(R_0-1)+c_0\big)}{(\alpha-1)^3}>0 .
\]

\begin{lemma}[Two kernel estimates]\label{lem_kernel_estimates}
Let $K_\infty(\phi,\phi')=\frac2\pi\sum_{k\ge1}k^{-1}\sin k\phi\sin k\phi'$ be the kernel of \autoref{lem_kernel} with $L=\infty$.
\begin{enumerate}
\item For every $L>0$ and every $g\ge0$, $K_Lg\le K_\infty g$.
\item $K_\infty\big[\tfrac12\cot(\cdot/4)\big](s)\le\pi-\frac s2$ for $0<s<\pi$.
\end{enumerate}
\end{lemma}
\begin{proof}
(a) Let $u_\infty$ be the bounded harmonic function in $(-\infty,L)\times(0,\pi)$ that vanishes on the sides and satisfies $\partial_su_\infty=g$ on $s=L$. Let $u_L$ be the corresponding function in $(0,L)\times(0,\pi)$ that vanishes on $s=0$. Their traces on $s=L$ are $K_\infty g$ and $K_Lg$, and $u_\infty\ge0$ by the argument of \autoref{lem_kernel}. The difference $w=u_\infty-u_L$ is harmonic in the finite strip, vanishes on the sides, satisfies $\partial_sw=0$ on the top and $w=u_\infty\ge0$ on the bottom. A negative minimum would lie on the top, where the Hopf lemma would give $\partial_sw<0$. Hence $w\ge0$.

(b) From $\cot(t/4)=\cot(t/2)+1/\sin(t/2)$ and $1/\sin(t/2)=\cot(t/2)+\tan(t/4)$, we get $\frac12\cot(t/4)=\cot(t/2)+\frac12\tan(t/4)\le\cot(t/2)+\frac12\tan(t/2)$ on $(0,\pi)$. The sine coefficients of $\cot(t/2)$ on $(0,\pi)$ all equal $2$, so $K_\infty[\cot(\cdot/2)](s)=2\sum_kk^{-1}\sin ks=\pi-s$. Moreover $K_\infty[\tan(\cdot/2)](s)=s$; this is Lemma~5 of \cite{Toland1996}, whose kernel is half of ours. Since $K_\infty$ is positive, adding the two parts gives (b).
\end{proof}

\begin{theorem}[Slope bound in the Keady--Norbury cone]\label{thm_cone_slope}
Let $(\Psi,S)\in\mathcal N$ lie in the lower region at a parameter $\alpha>R_0$ with $\Delta(\alpha)>0$, and suppose that
\begin{enumerate}
\item[(C2)] $\vartheta(\phi)/\sin(\phi/2)$ is nonincreasing on $(0,\pi)$.
\end{enumerate}
Then $\vartheta(\phi)<\kappa_\alpha\big(\pi-\frac\phi2\big)<\kappa_\alpha\pi$ for $0<\phi<\pi$, with $\kappa_\alpha<\frac12$ as in \autoref{lem_structural}. In particular (B) holds with $\eta=\frac\pi2-\kappa_\alpha\pi=\frac{\pi\delta_\alpha}{2(2+\delta_\alpha)}$. If (C1) also holds, then $\vartheta(\phi)<\frac13\big(\pi-\frac\phi2\big)<\frac\pi3$.
\end{theorem}
\begin{proof}
Fix $0<\tau<\phi<\pi$ and put $\lambda=\sin(\tau/2)/\sin(\phi/2)\in(0,1)$. By (C2), $\vartheta(\tau)\ge\lambda\vartheta(\phi)$. Since $0<\vartheta<\pi/2$ (\autoref{thm_nekrasov}) and $\sin$ is concave and increasing on $[0,\pi/2]$, $\sin\vartheta(\tau)\ge\sin(\lambda\vartheta(\phi))\ge\lambda\sin\vartheta(\phi)$. Integrating in $\tau$,
\[
\int_0^\phi\sin\vartheta\ge\frac{\sin\vartheta(\phi)}{\sin(\phi/2)}\int_0^\phi\sin\frac\tau2\,d\tau=2\tan\frac\phi4\,\sin\vartheta(\phi),
\qquad\text{so}\qquad\frac{V'(\phi)}{V(\phi)}\le\frac12\cot\frac\phi4 .
\]
By \autoref{lem_kernel} and the positivity of $K_L$, then \autoref{lem_structural} and \autoref{lem_kernel_estimates},
\[
\vartheta=K_L\big[(F\circ V)'\big]\le K_L\big[(F\circ V)'_+\big]\le\kappa\,K_L\big[\tfrac12\cot(\cdot/4)\big]\le\kappa\,K_\infty\big[\tfrac12\cot(\cdot/4)\big]\le\kappa\Big(\pi-\frac\phi2\Big),
\]
with $\kappa=\kappa_\alpha$, or $\kappa=\frac13$ under (C1). The inequality is strict, because $\sin(\lambda x)>\lambda\sin x$ for $\lambda\in(0,1)$ and $x\in(0,\pi/2]$.
\end{proof}

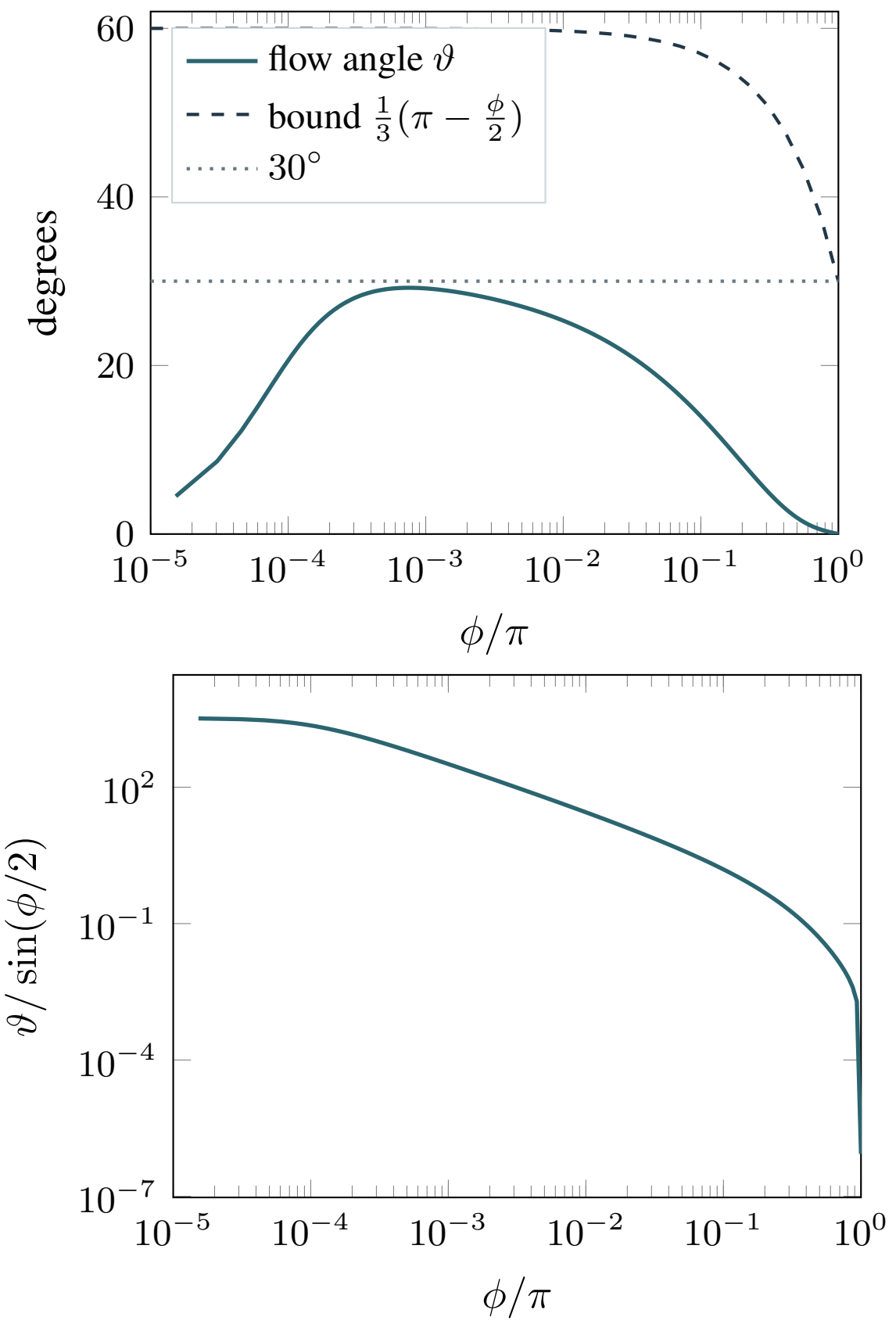
\begin{figure}[!ht]
\centering
\begin{tikzpicture}
\begin{axis}[width=7.4cm,height=6cm,xmode=log,xmin=1e-5,xmax=1,ymin=0,ymax=62,xlabel={$\phi/\pi$},ylabel={degrees},
  tick label style={font=\footnotesize},label style={font=\small},legend style={font=\footnotesize,at={(0.03,0.97)},anchor=north west,draw=PaperRule},legend cell align=left]
\addplot[PaperAccent,line width=1pt] coordinates {(1.52588e-05,4.482) (3.05176e-05,8.645) (4.57764e-05,12.271) (6.10352e-05,15.285) (7.62939e-05,17.730) (9.15527e-05,19.693) (0.000106812,21.266) (0.00012207,22.529) (0.000137329,23.552) (0.000152588,24.392) (0.000167847,25.084) (0.000183105,25.657) (0.000198364,26.139) (0.000213623,26.547) (0.000228882,26.895) (0.000244141,27.190) (0.000259399,27.445) (0.000289917,27.859) (0.000305176,28.025) (0.000335693,28.300) (0.000350952,28.413) (0.00038147,28.600) (0.000411987,28.748) (0.000427246,28.808) (0.000473022,28.954) (0.00050354,29.025) (0.000534058,29.081) (0.000579834,29.142) (0.000610352,29.170) (0.000656128,29.199) (0.000701904,29.214) (0.000762939,29.220) (0.000808716,29.215) (0.000869751,29.200) (0.000930786,29.178) (0.00100708,29.143) (0.00106812,29.111) (0.00115967,29.058) (0.00123596,29.011) (0.00132751,28.952) (0.00141907,28.892) (0.00152588,28.821) (0.00163269,28.751) (0.00175476,28.670) (0.00187683,28.591) (0.00201416,28.503) (0.00215149,28.417) (0.00231934,28.315) (0.00248718,28.215) (0.00265503,28.119) (0.00285339,28.008) (0.00306702,27.893) (0.00328064,27.782) (0.00352478,27.659) (0.00376892,27.541) (0.00404358,27.412) (0.0043335,27.282) (0.00465393,27.142) (0.00498962,27.002) (0.00534058,26.861) (0.0057373,26.707) (0.00614929,26.553) (0.0065918,26.393) (0.00706482,26.229) (0.00758362,26.056) (0.00811768,25.884) (0.00871277,25.699) (0.00933838,25.512) (0.0100098,25.319) (0.0107422,25.115) (0.0115204,24.907) (0.0123444,24.695) (0.0132446,24.471) (0.0141907,24.244) (0.0152283,24.005) (0.0163269,23.761) (0.0175018,23.509) (0.0187683,23.248) (0.0201263,22.978) (0.0215759,22.700) (0.0231323,22.413) (0.0247955,22.117) (0.0265961,21.810) (0.0285187,21.493) (0.0305786,21.167) (0.0327911,20.831) (0.0351562,20.484) (0.0376892,20.128) (0.0404205,19.759) (0.043335,19.381) (0.046463,18.991) (0.0498199,18.590) (0.053421,18.178) (0.0572815,17.754) (0.0614166,17.319) (0.0658569,16.872) (0.0706177,16.414) (0.0757141,15.945) (0.0811768,15.466) (0.0870514,14.974) (0.093338,14.472) (0.100082,13.960) (0.107315,13.437) (0.115067,12.905) (0.123367,12.365) (0.132278,11.816) (0.141846,11.259) (0.152084,10.696) (0.163071,10.128) (0.17485,9.556) (0.187485,8.982) (0.201035,8.406) (0.215561,7.832) (0.231125,7.261) (0.247818,6.696) (0.265732,6.139) (0.284927,5.593) (0.305511,5.061) (0.327576,4.545) (0.351242,4.050) (0.376617,3.577) (0.403824,3.129) (0.432999,2.709) (0.464279,2.320) (0.497818,1.962) (0.533768,1.638) (0.572327,1.347) (0.613678,1.088) (0.658005,0.862) (0.705551,0.666) (0.756516,0.497) (0.811172,0.352) (0.869766,0.225) (0.932602,0.111) (0.999969,0.000)};
\addlegendentry{flow angle $\vartheta$}
\addplot[PaperInk,dashed,line width=.8pt,domain=1e-5:1,samples=40] {(180-90*x)/3};
\addlegendentry{bound $\frac13(\pi-\frac\phi2)$}
\addplot[PaperMuted,dotted,line width=.8pt,domain=1e-5:1,samples=2] {30};
\addlegendentry{$30^\circ$}
\end{axis}
\end{tikzpicture}\hfill
\begin{tikzpicture}
\begin{axis}[width=7.4cm,height=6cm,xmode=log,ymode=log,xmin=1e-5,xmax=1,xlabel={$\phi/\pi$},ylabel={$\vartheta/\sin(\phi/2)$},
  tick label style={font=\footnotesize},label style={font=\small}]
\addplot[PaperAccent,line width=1pt] coordinates {(1.52588e-05,3263.4) (3.05176e-05,3147.5) (4.57764e-05,2978.4) (6.10352e-05,2782.6) (7.62939e-05,2582.1) (9.15527e-05,2390) (0.000106812,2212.2) (0.00012207,2050.6) (0.000137329,1905.6) (0.000152588,1776.1) (0.000167847,1660.5) (0.000183105,1556.9) (0.000198364,1464.1) (0.000213623,1380.8) (0.000228882,1305.6) (0.000244141,1237.5) (0.000259399,1175.6) (0.000289917,1067.7) (0.000305176,1020.4) (0.000335693,936.7) (0.000350952,899.57) (0.00038147,833.04) (0.000411987,775.32) (0.000427246,749.2) (0.000473022,680.12) (0.00050354,640.46) (0.000534058,605.04) (0.000579834,558.43) (0.000610352,531.02) (0.000656128,494.47) (0.000701904,462.46) (0.000762939,425.54) (0.000808716,401.38) (0.000869751,373.03) (0.000930786,348.3) (0.00100708,321.54) (0.00106812,302.83) (0.00115967,278.41) (0.00123596,260.8) (0.00132751,242.33) (0.00141907,226.22) (0.00152588,209.87) (0.00163269,195.66) (0.00175476,181.54) (0.00187683,169.26) (0.00201416,157.24) (0.00215149,146.76) (0.00231934,135.65) (0.00248718,126.05) (0.00265503,117.68) (0.00285339,109.06) (0.00306702,101.05) (0.00328064,94.094) (0.00352478,87.19) (0.00376892,81.193) (0.00404358,75.325) (0.0043335,69.951) (0.00465393,64.802) (0.00498962,60.13) (0.00534058,55.884) (0.0057373,51.722) (0.00614929,47.978) (0.0065918,44.489) (0.00706482,41.252) (0.00758362,38.176) (0.00811768,35.429) (0.00871277,32.774) (0.00933838,30.356) (0.0100098,28.105) (0.0107422,25.979) (0.0115204,24.024) (0.0123444,22.229) (0.0132446,20.531) (0.0141907,18.985) (0.0152283,17.517) (0.0163269,16.172) (0.0175018,14.927) (0.0187683,13.765) (0.0201263,12.687) (0.0215759,11.692) (0.0231323,10.768) (0.0247955,9.9136) (0.0265961,9.1141) (0.0285187,8.3768) (0.0305786,7.6944) (0.0327911,7.0614) (0.0351562,6.4774) (0.0376892,5.9374) (0.0404205,5.4352) (0.043335,4.9731) (0.046463,4.5456) (0.0498199,4.1503) (0.053421,3.7852) (0.0572815,3.4484) (0.0614166,3.1381) (0.0658569,2.8517) (0.0706177,2.5879) (0.0757141,2.3455) (0.0811768,2.1226) (0.0870514,1.9172) (0.093338,1.729) (0.100082,1.5562) (0.107315,1.3978) (0.115067,1.253) (0.123367,1.1207) (0.132278,0.99969) (0.141846,0.88927) (0.152084,0.78895) (0.163071,0.69771) (0.17485,0.61497) (0.187485,0.54006) (0.201035,0.47241) (0.215561,0.41151) (0.231125,0.35687) (0.247818,0.30795) (0.265732,0.2643) (0.284927,0.22556) (0.305511,0.19131) (0.327576,0.1612) (0.351242,0.13484) (0.376617,0.11193) (0.403824,0.092149) (0.432999,0.075189) (0.464279,0.060766) (0.497818,0.048604) (0.533768,0.038443) (0.572327,0.030027) (0.613678,0.023125) (0.658005,0.017516) (0.705551,0.01299) (0.756516,0.0093541) (0.811172,0.0064216) (0.869766,0.0040135) (0.932602,0.0019434) (0.999969,8.586e-07)};
\end{axis}
\end{tikzpicture}
\caption{The Keady--Norbury cone on the near-extreme wave of \autoref{subsec_hp} ($\gamma=0.9$, $q_{\rm cr}/q_{\rm cr}^{\,0}=0.091$, $N=32768$), on logarithmic axes in $\phi$. Left: the flow angle on the free surface rises within $\phi/\pi<10^{-3}$ of the crest to $29.2^\circ$, and then decays towards the trough. It stays below the bound $\frac13(\pi-\frac\phi2)$ of \autoref{thm_cone_slope}. Right: the ratio $\vartheta/\sin(\phi/2)$ is strictly decreasing, which is condition (C2). The parameters lie above the floor, so by \autoref{cor_thin} the cone condition holds along the whole continuum, as observed here.}
\label{fig_cone}
\end{figure}

\autoref{fig_cone} shows the bound and the cone condition on a computed near-extreme wave.

\begin{corollary}[Existence modulo the cone condition]\label{cor_cone_existence}
Assume \eqref{eq_standing} and \eqref{eq_floor_condition}. If, for all small $\delta$, every wave of $\mathcal C^+_\delta$ satisfies (C2), then there exists a symmetric monotone extreme wave whose single crest per period is a $120^\circ$ Stokes corner.
\end{corollary}
\begin{proof}
By \autoref{cor_floor}, $\mathcal C^+_\delta$ lies in the lower region with $\alpha\in(\alpha_\circ,A)$, so $\delta_\alpha>\delta_\circ:=\Delta(\alpha_\circ)/(\alpha_\circ-1)^2>0$. \autoref{thm_cone_slope} gives (B) with $\eta=\pi\delta_\circ/(2(2+\delta_\circ))$, and \autoref{cor_floor} gives (A). Now apply \autoref{thm_final_reduction}.
\end{proof}

\subsubsection{Invariance of the cone along the continuum}\label{subsec_cone_invariance}
We now show that (C2) propagates along $\mathcal C^+_\delta$. We need one explicit condition on the parameter:
\begin{enumerate}
\item[(C3)] either $s_h(\alpha)\le\sigma_{\rm th}(\alpha)$, or $s_h(\alpha)>1$ and
\[
\Phi_\alpha:=H_\alpha(r_-)-H_\alpha(s_h)+\frac{p_0^2}{\log r_-}-|p_0|\sqrt{h_\alpha(s_h)}\ \ge0 ,
\]
\end{enumerate}
where $s_h(\alpha)$ is the maximum point of $h_\alpha$ on $(1,r_-)$ ($s_h=1$ if $h_\alpha$ is decreasing there) and $H_\alpha(\sigma)=\frac14\sigma^4-\frac23\alpha\sigma^3+\lambda\sigma^2$ is an antiderivative of $\sigma Q_\alpha(\sigma)$. By \autoref{lem_trough_bound} below, (C3) keeps the trough of every wave beyond $s_h$, where $h_\alpha$ decreases. The monotonicity of the weight, which is the second ingredient of the flat argument, holds without any condition.

\begin{lemma}[The weight is monotone]\label{lem_weight_monotone}
Let $\Delta(\alpha)>0$. Then $\sigma\mapsto-\sigma h_\alpha'(\sigma)/h_\alpha(\sigma)^{3/2}$ is strictly increasing on $(0,r_-)$.
\end{lemma}
\begin{proof}
Put $P=-2(2\sigma^2-3\alpha\sigma+2\lambda)$, so that $h_\alpha'=-\sigma P$ and $-\sigma h_\alpha'/h_\alpha^{3/2}=P/(\sigma Q_\alpha^{3/2})$. A direct computation gives
\[
\frac{d}{d\sigma}\,\frac{P}{\sigma Q_\alpha^{3/2}}=\frac{N}{\sigma^2Q_\alpha^{5/2}},\qquad N=4Q_\alpha(\lambda-\sigma^2)+3\sigma(\alpha-\sigma)P .
\]
With $a=r_--\sigma$ and $b=r_+-\sigma$ we have $Q_\alpha=ab$, $2\lambda=r_-r_+$, $2\alpha=r_-+r_+$ and $P=\sigma(a+b)-2ab$. Substituting,
\[
N=\Big(\tfrac32(a+b)^2-2ab\Big)\sigma^2-ab(a+b)\,\sigma+2a^2b^2 .
\]
For $a,b>0$ the leading coefficient is positive, and the discriminant is $a^2b^2\big(16ab-11(a+b)^2\big)<0$. Hence $N>0$.
\end{proof}

\begin{lemma}[A lower bound for the trough]\label{lem_trough_bound}
Let $\alpha>\alpha_\circ$ with $\Delta(\alpha)>0$, and let $(\Psi,S)\in\mathcal N$ be a monotone wave in the lower region at the parameter $\alpha$. Then
\begin{equation}\label{eq_trough_bound}
H_\alpha(r_-)-H_\alpha(S_t)+\frac{p_0^2}{\log r_-}<|p_0|\max_{[S_t,r_-]}\sqrt{h_\alpha} .
\end{equation}
Consequently, if (C3) holds at $\alpha$, then $h_\alpha'\le0$ on $[S_t,r_-]$, and $h_\alpha'<0$ on $(S_t,r_-]$.
\end{lemma}
\begin{proof}
By the angular flow force \eqref{eq_angular_flow_force},
\[
H_\alpha(S_c)-H_\alpha(S_t)=p_0^2\Big(\int_1^{S_t}R\,\Psi_R(R,\pi)^2\,dR-\int_1^{S_c}R\,\Psi_R(R,0)^2\,dR\Big).
\]
Along both rays $\Psi$ decreases from $1$ to $0$, so $\int|\Psi_R|\,dR=1$ on each. On the trough ray, \autoref{thm_max_velocity} gives $R|\Psi_R|\le\max_\theta\sqrt{h_\alpha(S(\theta))}/|p_0|$, and the range of $S$ is $[S_t,S_c]\subset[S_t,r_-]$. Hence the first integral is at most $\max_{[S_t,r_-]}\sqrt{h_\alpha}/|p_0|$. On the crest ray, the Cauchy--Schwarz inequality gives
\[
1\le\Big(\int_1^{S_c}|\Psi_R|\,dR\Big)^2\le\int_1^{S_c}R\,\Psi_R^2\,dR\cdot\log S_c ,
\]
so the second integral is at least $1/\log S_c$. Hence $H_\alpha(S_c)+p_0^2/\log S_c<H_\alpha(S_t)+|p_0|\max_{[S_t,r_-]}\sqrt{h_\alpha}$. Finally, $S_c\in[\sigma_{\rm dp},r_-)$ with $\sigma_{\rm dp}=R_0$, by \autoref{lem_unimodal} and \autoref{lem_radial_crossing}. On this interval
\[
\frac{d}{dS}\Big(H_\alpha(S)+\frac{p_0^2}{\log S}\Big)=SQ_\alpha(S)-\frac{p_0^2}{S\log^2S}=\frac{G_\alpha(S)}{S\log^2S}<0,
\]
because $G_\alpha<0$ on $(\sigma_{\rm dp},r_-)$. So $H_\alpha(S_c)+p_0^2/\log S_c>H_\alpha(r_-)+p_0^2/\log r_-$, which proves \eqref{eq_trough_bound}.

For the consequence, recall that $h_\alpha$ increases on $[1,s_h]$ and decreases on $[s_h,r_-]$ (\autoref{cor_energy_sandwich}). If $s_h\le\sigma_{\rm th}$, then $S_t\ge\sigma_{\rm th}\ge s_h$ by \autoref{lem_unimodal}. Otherwise $\Phi_\alpha\ge0$. Suppose $S_t\le s_h$. Then $\max_{[S_t,r_-]}h_\alpha=h_\alpha(s_h)$ and $H_\alpha(S_t)\le H_\alpha(s_h)$. Hence \eqref{eq_trough_bound} gives $\Phi_\alpha<0$, a contradiction. So $S_t>s_h$ in this case.
\end{proof}

Under (C3), by \autoref{lem_trough_bound}, $h_\alpha'\le0$ on the range of every wave in the lower region at the parameter $\alpha$. By \autoref{lem_weight_monotone} and the monotonicity of $S$, the weight $\mathfrak g_\alpha(S(\phi))$ is nonincreasing in $\phi$; under (C3) it is also nonnegative, as in the flat case. The linear part of the argument is Toland's maximum-principle proof of the invariance of the Keady--Norbury cone \cite[Theorem~3]{Toland1996}, transplanted from the half-disc to the strip.

\begin{remark}[Structure of the invariance argument]\label{rem_cone_steps}
The proof that (C2) propagates has five steps, and each delicate point is handled in exactly one of them.
\begin{enumerate}
\item[(1)] \emph{Operator and axis.} The quotient $m=w/a^3$ satisfies an elliptic equation whose first-order coefficient $2\cot(\phi/2)$ is singular on the axis $\phi=0$; near the axis the operator is the five-dimensional Laplacian, so the strong maximum principle and the Hopf lemma hold up to $\phi=0$ (\autoref{lem_cone_equation}).
\item[(2)] \emph{Boundary signs.} On the bottom $m=0$; on the trough axis $m=2\vartheta_\phi\le0$; on the top $\partial_sm\le0$, which uses (C2), (C3), \autoref{lem_trough_bound} and \autoref{lem_weight_monotone} (first part of the proof of \autoref{lem_cone_strict}).
\item[(3)] \emph{Maximum principle.} The zeroth-order coefficient is negative, and every boundary point, including the four corners of $\overline Z$, is covered by one of the sign conditions of step (2) or by the Hopf lemma on the top (second part of that proof).
\item[(4)] \emph{Strictness.} The strong maximum principle, the Hopf lemma on the bottom and on the trough axis, (N4), and Serrin's corner lemma at $(0,\pi)$ give strict signs on all of $\overline Z$ (third part of that proof).
\item[(5)] \emph{Propagation.} The strict condition is open in $X$, the closed condition (C2) passes to limits within $\mathcal C^+_\delta$, and step (4) upgrades it again; connectedness of $\mathcal C^+_\delta$ finishes the argument (\autoref{thm_cone_invariance}).
\end{enumerate}
The argument is carried out for each fixed $\delta$. The slope constant it produces, through \autoref{thm_cone_slope}, depends only on a lower bound for $\delta_\alpha$, and hence not on $\delta$ (\autoref{rem_uniform_delta}).
\end{remark}

\begin{lemma}[A regular equation for the cone]\label{lem_cone_equation}
Let $\vartheta$ be harmonic in $Z=(0,L)\times(0,\pi)$ and smooth up to $\overline Z$, with $\vartheta=0$ on $s=0$ and on $\phi\in\{0,\pi\}$. Put $a=\sin(\phi/2)$ and
\[
w=a\,\vartheta_\phi-a'\vartheta=a^2\,\partial_\phi\Big(\frac\vartheta a\Big),\qquad m=\frac w{a^3}.
\]
Then $m$ extends smoothly to $\phi=0$, and
\begin{equation}\label{eq_cone_m}
\Delta m+2\cot(\phi/2)\,m_\phi-m=0\quad\text{in }Z,\qquad m=0\ \text{on }s=0,\qquad m=2\vartheta_\phi\ \text{on }\phi=\pi .
\end{equation}
Near $\phi=0$, the operator $\Delta+2\cot(\phi/2)\partial_\phi$ is the Laplacian of $\mathbb R^5$, acting on functions of $s$ and of the distance $\phi$ to a three-dimensional axis, up to a smooth lower-order term. Hence the strong maximum principle and the Hopf lemma hold for \eqref{eq_cone_m} up to $\phi=0$.
\end{lemma}
\begin{proof}
Since $a''=-a/4$ and $a'''=-a'/4$, a direct computation gives $\Delta w=\cot(\phi/2)\,w_\phi+\frac14w$; this is the strip version of the equation in \cite[Theorem~3]{Toland1996}. Substituting $w=a^3m$ and using $(a^3)'=3a^2a'$ and $(a^3)''=6a(a')^2-\frac34a^3$, the coefficient of $m_\phi$ becomes $2\cdot3a'/a-2a'/a=4a'/a=2\cot(\phi/2)$, and that of $m$ becomes $-1$. Since $\vartheta$ is odd in $\phi$, $w=O(\phi^3)$ at $\phi=0$, and $m$ is even and smooth there. The boundary values follow from $\vartheta=0$ on $s=0$ and from $a(\pi)=1$, $a'(\pi)=0$, $\vartheta(\cdot,\pi)=0$. Finally $2\cot(\phi/2)=4/\phi+O(\phi)$, and $\partial_\phi^2+\frac4\phi\partial_\phi$ is the radial Laplacian in $\mathbb R^5$.
\end{proof}

\begin{lemma}[Closed cone implies strict cone]\label{lem_cone_strict}
Let $(\Psi,S)\in\mathcal N$ be a monotone wave in the lower region at a parameter $\alpha>\alpha_\circ$ with $\Delta(\alpha)>0$ satisfying (C3), and suppose that $\vartheta$ satisfies (C2). Then $m<0$ on $(0,L]\times[0,\pi]$ and $\partial_sm<0$ on $\{0\}\times[0,\pi]$. In particular, $\vartheta(\phi)/\sin(\phi/2)$ is strictly decreasing on $(0,\pi)$, with $\partial_\phi^2\big(\vartheta/\sin(\phi/2)\big)(0)<0$.
\end{lemma}
\begin{proof}
\emph{Step (2): boundary signs.} Let $g=\mathfrak g_\alpha(S)\sin\vartheta$ be the Neumann data of $\vartheta$ on $s=L$ (\autoref{thm_nekrasov}). The function $\sin\vartheta/a$ is nonincreasing. Indeed,
\[
\Big(\frac{\sin\vartheta}a\Big)'=\frac{a\cos\vartheta\,\vartheta'-a'\sin\vartheta}{a^2},
\]
and by (C2) $a\vartheta'\le a'\vartheta$. If $\vartheta'\ge0$, then $a\cos\vartheta\,\vartheta'\le a'\vartheta\cos\vartheta\le a'\sin\vartheta$, because $0<\vartheta<\pi/2$ on the top for $0<\phi<\pi$: positivity follows from (N2), and $\vartheta<\pi/2$ because $\tan\vartheta=|s'|$ is finite (\autoref{lem_flow_angle}). If $\vartheta'<0$, the numerator is negative. By (C3) and \autoref{lem_trough_bound}, $h_\alpha'\le0$ on the range of $S$, so the factor $\mathfrak g_\alpha(S(\phi))$ is nonnegative; by \autoref{lem_weight_monotone} and the monotonicity of $S$, it is nonincreasing. Hence $g/a$ is nonincreasing, and on the top
\[
\partial_sm=\frac{a\,g'-a'g}{a^3}=\frac1a\Big(\frac ga\Big)'\le0 .
\]

\emph{Step (3): maximum principle.} The zeroth-order coefficient in \eqref{eq_cone_m} is negative. Hence $m$ has no positive interior maximum, by \autoref{lem_cone_equation} also on the axis $\phi=0$. On $s=0$, including the corners $(0,0)$ and $(0,\pi)$, we have $m=0$. On $\phi=\pi$, including the corner $(L,\pi)$, we have $m=2\vartheta_\phi\le0$, because $\vartheta>0$ in $Z$ and vanishes on $\phi=\pi$. It remains to consider a positive maximum at a point $(L,\phi_0)$ with $0\le\phi_0<\pi$. In the five-dimensional picture of \autoref{lem_cone_equation} the axis is an interior line and the top is a flat boundary piece, so the corner $(L,0)$ is an ordinary boundary point of the top, and the Hopf lemma gives $\partial_sm>0$ at $(L,\phi_0)$ in all cases. This contradicts $\partial_sm\le0$. So $m\le0$.

\emph{Step (4): strictness.} If $m\equiv0$, then $\vartheta=c(s)\sin(\phi/2)$, which contradicts $\vartheta(\cdot,\pi)=0$ unless $\vartheta\equiv0$. By the strong maximum principle, $m<0$ in $Z$, including the axis. At a point of the top where $m=0$, the Hopf lemma would give $\partial_sm>0$, so $m<0$ on the top. On $\phi=\pi$ and $0<s\le L$, $m=2\vartheta_\phi<0$: by the Hopf lemma for $\vartheta$ on the side, and at the corner $(L,\pi)$ by (N4) of \autoref{def_nodal}. On $s=0$, where $m=0$, the Hopf lemma gives $\partial_sm<0$ for $0\le\phi<\pi$. At the corner $(0,\pi)$, $\partial_sm=2\vartheta_{s\phi}(0,\pi)<0$, by Serrin's corner lemma applied to $\vartheta$ \cite{Serrin1971,GNN1979}.
\end{proof}

\begin{theorem}[The cone is invariant along the continuum]\label{thm_cone_invariance}
Assume \eqref{eq_standing} and $c_0>\sinh^2\gamma$, and fix $\delta\in(0,\gamma)$. Suppose that every element of $\mathcal C^+_\delta$ lies in the lower region and has $\Delta(\alpha)>0$, and that (C3) holds at every such $\alpha$. Then every wave of $\mathcal C^+_\delta$ other than the bifurcation point satisfies the strict cone condition of \autoref{lem_cone_strict}, in particular (C2).
\end{theorem}
\begin{proof}
This is step (5) of \autoref{rem_cone_steps}. Let $E\subset\mathcal C^+_\delta$ consist of $b=(\alpha_c,R_0)$ and the waves satisfying the conclusion of \autoref{lem_cone_strict}. By assumption and \autoref{thm_crossing}(b), every wave lies in the lower region with $\alpha>\alpha_\circ$ and $\Delta(\alpha)>0$, so \autoref{lem_cone_strict} applies to it.

\emph{$E$ contains a neighbourhood of $b$.} On the local branch, $\vartheta=\varepsilon A\,\frac{\sinh s}{\sinh\gamma}\sin\phi+O(\varepsilon^2)$ with $A>0$, in every $C^k$. For $\vartheta_0=\sin\phi\,\sinh s$, a direct computation gives $w=-a^3\sinh s$, i.e.\ $m=-\sinh s$, with $\partial_sm=-1$ on $s=0$. The strict inequalities therefore persist for small $\varepsilon>0$.

\emph{$E$ is open.} Waves in $\mathcal O_\delta$ are analytic. By elliptic regularity for the hodograph problem \cite{ADN1959,LiebermanTrudinger1986}, convergence in $X$ implies convergence in $C^k(\overline D)$ for every $k$. Hence $m$ and $\partial_sm$ converge uniformly on $\overline Z$. The conditions $m<0$ on $(0,L]\times[0,\pi]$ and $\partial_sm<0$ on $\{0\}\times[0,\pi]$ are together open under $C^1$ convergence of $m$ on the compact set $\overline Z$.

\emph{$E$ is closed.} Let $x$ be a limit of elements of $E$. If $x$ is circular, $x=b$ by \autoref{thm_crossing}(a). Otherwise, $x\in\mathcal N$, and its flow angle satisfies the closed condition (C2) as a limit. By \autoref{lem_cone_strict}, $x\in E$.

Since $\mathcal C^+_\delta$ is connected, $E=\mathcal C^+_\delta$.
\end{proof}

\begin{remark}[Uniformity in $\delta$]\label{rem_uniform_delta}
In \autoref{thm_cone_invariance} the parameter $\delta$ is fixed, and the closedness step takes limits within the closed set $\mathcal C^+_\delta\subset\mathcal O_\delta$, in the topology of $X$. No limit $\delta\downarrow0$ of the cone condition is needed. What is used as $\delta\downarrow0$ is only the slope bound of \autoref{thm_cone_slope}, $\vartheta<\kappa_\alpha\pi$, whose constant depends on $\alpha$ alone. Under (A) it is bounded by $\kappa_{\alpha_++\eta}<\frac12$, and under \eqref{eq_floor_condition} by $\kappa_{\alpha_\circ}$, uniformly in $\delta$. Strictness at the limit wave is not needed: the limit $\delta\downarrow0$ is taken in \autoref{thm_limit}, which does not use the cone.
\end{remark}

\begin{theorem}[Existence of the extreme Stokes wave]\label{thm_existence}
Assume \eqref{eq_standing} and \eqref{eq_floor_condition}, and that (C3) holds for every $\alpha\in(\alpha_\circ,A(\gamma,p_0))$. Then the global continuum contains, in the limit $\delta\downarrow0$, a symmetric extreme wave. It is monotone from crest to trough, and its single stagnation point per period lies on the circle $R=r_-(\alpha_\infty)$, where the free surface forms a corner of exactly $120^\circ$, symmetric about the radial axis, with the Stokes asymptotics of \autoref{thm_intro_stokes}. Along the continuum, every wave satisfies $\alpha_\circ<\alpha<A$ and has slope below $90^\circ-\eta_\circ$, where $\eta_\circ=\frac{\pi\delta_\circ}{2(2+\delta_\circ)}$ and $\delta_\circ=\Delta(\alpha_\circ)/(\alpha_\circ-1)^2>0$. If (C1) also holds for these $\alpha$, the slope stays below $60^\circ$.
\end{theorem}
\begin{proof}
By \autoref{cor_floor}, every wave of $\mathcal C^+_\delta$ lies in the lower region and has $\alpha\in(\alpha_\circ,A)$ and $\Delta(\alpha)>0$. By \autoref{cor_thin}, (C3) holds at each such $\alpha$, so \autoref{thm_cone_invariance} shows that every wave satisfies (C2). Moreover $\delta_\alpha>\delta_\circ$ for $\alpha>\alpha_\circ$, because $\delta_\alpha$ increases in $\alpha$, and $\delta_\circ>0$ because $\alpha_\circ>\alpha_+$. The constant $\eta_\circ$ depends only on $\gamma$ and $c_0$, not on $\delta$. \autoref{thm_cone_slope} then gives (B) with $\eta=\eta_\circ$, and \autoref{cor_floor} gives (A) with $\eta=\alpha_\circ-\alpha_+$. Now apply \autoref{thm_final_reduction}.
\end{proof}

\begin{corollary}[Every annulus above the floor]\label{cor_thin}
Let $0<\gamma<1$ and $0<c_0<e^{2\gamma}$. Then (C3) holds for every $\alpha>\max\{\alpha_\circ,\alpha_+\}$. In particular, if the floor condition \eqref{eq_floor_condition} holds, then the extreme Stokes wave of \autoref{thm_existence} exists: the existence half of the Stokes conjecture holds whenever \eqref{eq_floor_condition} holds.
\end{corollary}
\begin{proof}
If $s_h\le1$ there is nothing to prove, so let $s_h>1$. We show that $\Phi_\alpha>0$.

\emph{Scaled variables.} Write $f=c_0e^{-2\gamma}<1$, $u=f(1+\gamma)/\gamma=\alpha_\circ/R_0-1$, $\tau=R_0/\alpha$ and $\beta=1+\frac12(f-1)\tau$. Put $x=s_h/R_0$ and $\rho=r_-/R_0$. Then $x$ is the smaller zero of $2\tau z^2-3z+2\beta$, because $h_\alpha'=2\sigma(2\sigma^2-3\alpha\sigma+2\lambda)$. Similarly, $\rho$ is the smaller zero of $\tau z^2-2z+2\beta$, because $\tau Q_\alpha(R_0z)/R_0^2=\tau z^2-2z+2\beta$. Hence $4\tau x\le3$, $\tau\rho<1$, and $\rho>1$ because $Q_\alpha(R_0)=c_0>0$. Moreover $\beta\le1$, since $f\le1$. Finally, $\alpha>\alpha_\circ$ and $\gamma<1$ give $\tau<1/(1+u)<1/(1+2f)$, that is, $3\tau+4\beta<5$.

\emph{The quantity to estimate.} Put $D=\rho^2-2x^2$ and $Y=x^{3/2}\sqrt{1-\tau x}$. Subtracting the two quadratic relations gives $\tau D=2\rho-3x$. With $K(z)=\frac14\tau z^4-\frac23z^3+\beta z^2=\tau H_\alpha(R_0z)/R_0^4$, $p_0^2=\gamma^2R_0^4f$ and $\log r_-=\gamma+\log\rho$, a direct computation gives
\[
\frac{\tau\,\Phi_\alpha}{R_0^4}=I+\frac{\gamma^2f\tau}{L}-\gamma\sqrt{f\tau}\,Y,\qquad I=K(\rho)-K(x),\quad L=\gamma+\log\rho .
\]
Completing the square in $\gamma\sqrt{f\tau}$, and using $L<1+\log\rho\le1+\frac{\rho^2-1}{2\rho}$ (since $\gamma<1$ and $\rho\ge1$), we obtain
\[
\frac{\tau\,\Phi_\alpha}{R_0^4}\ \ge\ I-\frac L4Y^2\ >\ I-\frac{\rho^2+2\rho-1}{8\rho}\,Y^2 .
\]
The two quadratic relations also give the factorizations
\[
12D\,I=(\rho-x)^3\big(7x^2+3x\rho-2\rho^2\big),\qquad D\,(1-\tau x)=(\rho-x)^2 .
\]
With $z=\rho/x$ this yields
\[
24\rho D\Big(I-\frac{\rho^2+2\rho-1}{8\rho}Y^2\Big)=(\rho-x)^2x^3\,F(z,\rho),\qquad F(z,\rho)=2\rho\,g(z)-3(\rho^2+2\rho-1),
\]
where $g(z)=(z-1)(7+3z-2z^2)$. It therefore suffices to show $D>0$ and $F(z,\rho)>0$.

\emph{The domain.} Since $\tau\rho^2-2\rho+2\beta=0$ and $4\tau x\le3$,
\[
0=9\big(2\tau x^2-3x+2\beta\big)=-\tau\rho^2+(3x-2\rho)(6\tau x+4\tau\rho-9).
\]
If $3x\ge2\rho$, the right-hand side would be negative. Hence $3x<2\rho$, so $\tau D=2\rho-3x>0$, $D>0$ and $z>\frac32$. Next, the identities
\[
\begin{gathered}
D(1-\tau\rho)=-x^2(z-1)(z-2),\qquad 2D(1-\beta)=x^2\big(2(z^2-2)+(4-3z)\rho\big),\\
zD(5-3\tau-4\beta)=x\,T(z,\rho),
\end{gathered}
\]
with $T(z,\rho)=(8-6z)\rho^2+5(z^2-2)\rho-3z(2z-3)$, turn the constraints $\tau\rho<1$, $\beta\le1$ and $3\tau+4\beta<5$ into
\[
z<2,\qquad 2(z^2-2)+(4-3z)\rho\ge0,\qquad T(z,\rho)>0 .
\]

\emph{Positivity of $F$.} Note that $F$ is a concave quadratic in $\rho$. On $[\frac32,2]$ we have $2g(z)-7=(2z-3)\big(3+2(2-z)(z+1)\big)\ge0$, so $F(z,1)=2g(z)-6\ge1$.
\begin{itemize}
\item \emph{Case $\frac32<z\le\frac{31}{20}$.} Put $\rho_f=2(z^2-2)/(3z-4)$; the second constraint says $\rho\le\rho_f$. A computation gives $(3z-4)^2F(z,\rho_f)=P_f(z)$ with
\[
P_f\Big(\frac32+\frac w{20}\Big)=\frac{2\,000\,000-700\,000w-410\,000w^2+106\,000w^3-5\,800w^4-310w^5-3w^6}{8\cdot10^6}.
\]
For $w\in[0,1]$ this is at least $(2\,000\,000-1\,116\,113)/(8\cdot10^6)>0$. Since $F(z,\cdot)$ is concave and positive at $\rho=1$ and at $\rho=\rho_f$, it is positive on $[1,\rho_f]$.
\item \emph{Case $\frac{31}{20}\le z<2$.} Here we use the third constraint through the multiplier $5$:
\[
F-5T=(30z-43)\rho^2+(-4z^3-15z^2+8z+30)\rho+30z^2-45z+3 .
\]
Its leading coefficient is at least $\frac72$, and its discriminant is
\[
\mathrm{disc}(z)=16z^6+120z^5+161z^4-4080z^3+9724z^2-7620z+1416 .
\]
With $W=20z-31\in[0,9]$, $9^6\cdot4\cdot10^6\,\mathrm{disc}(z)$ equals
\[
\begin{aligned}
&-6\,951\,644\cdot9^6-291\,503\,196\cdot9^5W-869\,335\,335\cdot9^4W^2+361\,890\,180\cdot9^3W^3\\
&\quad+273\,528\,090\cdot9^2W^4+19\,840\,464\cdot9\,W^5+531\,441\,W^6 .
\end{aligned}
\]
Using $W^k\le9^{k-2}W^2$ for $k\ge2$, the sum of the terms of degree at least $2$ is at most $9^4W^2(-869\,335\,335+361\,890\,180+273\,528\,090+19\,840\,464+531\,441)<0$. Hence $\mathrm{disc}<0$, so $F-5T>0$ for all $\rho$, and $F>5T>0$.
\end{itemize}
In both cases $F(z,\rho)>0$, hence $\Phi_\alpha>0$ and (C3) holds.
\end{proof}

The proof of \autoref{cor_thin} is elementary and complete as written. It uses no interval arithmetic, no grid search and no computer-generated certificate; every inequality is reduced to the explicit polynomial estimates displayed above. As an independent check, its algebraic steps have also been verified in Lean~4, in the file \texttt{AnnularStokes.lean} (theorem \texttt{direct\_C3} and the lemmas before it): the scaling identity, the two factorizations, the three constraint identities, the two polynomial estimates and the final positivity. This is a partial formal verification. The inequality $\log\rho\le(\rho^2-1)/(2\rho)$ and the definitions of $x$ and $\rho$ as roots enter the Lean statements as hypotheses, and the analytic results on which \autoref{cor_thin} is used are not formalized. The margins are small: at $\tau=0$, i.e.\ in the limit $\alpha\to\infty$, $F(\frac32,1)=1$, and the smallest value of $F$ on the domain is about $0.22$. It is attained where $f=1$ and $\tau(1+2f)=1$ meet, at $z\approx1.558$, $\rho\approx1.268$. On the two near-extreme waves computed in \autoref{sec_numerics} with $\gamma=0.5$ and $\gamma=0.9$, the three estimates in the proof of \autoref{lem_trough_bound} hold with the trough term within $0.3\%$ of its bound, and $S_t$ exceeds $s_h$ by $0.29$ and $0.31$ respectively.

\begin{theorem}[No overturning without degeneration]\label{thm_band}
Assume \eqref{eq_standing} and $c_0>\sinh^2\gamma$, and suppose that (A) holds: there is $\eta>0$ with $\alpha\ge\alpha_++\eta$ on $\mathcal C^+_\delta$ for all small $\delta$. Then there exists a symmetric extreme wave, monotone from crest to trough, with a single $120^\circ$ Stokes crest per period, as in alternative (I) of \autoref{thm_global}. Consequently, for every $c_0>\sinh^2\gamma$, at least one of the alternatives (I) and (III) of \autoref{thm_global} holds; overturning (II) can occur only together with degeneration (III).
\end{theorem}
\begin{proof}
By Step~1 of the proof of \autoref{thm_final_reduction}, $\mathcal C^+_\delta$ lies in the lower region, and $\Delta(\alpha)\ge\Delta(\alpha_++\eta)>0$ on it. By \autoref{thm_crossing}(b), $\alpha>\alpha_\circ$ on it. Hence $\alpha>\max\{\alpha_\circ,\alpha_+\}$, and (C3) holds at every such $\alpha$ by \autoref{cor_thin}. By \autoref{thm_cone_invariance}, every wave satisfies (C2). Since $\delta_\alpha$ increases in $\alpha$, $\delta_\alpha\ge\delta_{\alpha_++\eta}>0$ along the continuum, and \autoref{thm_cone_slope} gives (B) with a uniform constant. \autoref{thm_final_reduction} now gives the extreme wave. For the consequence, if (III) fails, then (A) holds, as in the proof of \autoref{thm_global}.
\end{proof}

\begin{figure}[!ht]
\centering
\begin{tikzpicture}
\begin{axis}[width=9.2cm,height=8.2cm,xmin=0,xmax=1,ymin=0,ymax=1,
  xlabel={$\gamma$},ylabel={$f=c_0e^{-2\gamma}$},tick label style={font=\footnotesize},label style={font=\small},
  legend style={font=\footnotesize,at={(1.02,1)},anchor=north west,draw=PaperRule},legend cell align=left]
\addplot[draw=none,fill=PaperMuted!25,domain=0:1,samples=60,area legend] {((1-exp(-2*x))/2)^2} \closedcycle;
\addlegendentry{window $c_0\le\sinh^2\gamma$}
\addplot[PaperInk,dashed,line width=.8pt,domain=0:1,samples=60] {(x/(1+x))^2};
\addlegendentry{$\alpha_\circ=\alpha_+$}
\addplot[fill=PaperAccent!45,draw=PaperAccent,line width=.8pt,area legend] coordinates {(0.0000,0.00000) (0.0125,0.00015) (0.0250,0.00059) (0.0375,0.00131) (0.0500,0.00227) (0.0625,0.00346) (0.0750,0.00487) (0.0875,0.00647) (0.1000,0.00826) (0.1125,0.01023) (0.1250,0.01235) (0.1375,0.01461) (0.1500,0.01701) (0.1625,0.01954) (0.1750,0.02218) (0.1875,0.02493) (0.2000,0.02778) (0.2125,0.03072) (0.2250,0.03374) (0.2375,0.03683) (0.2500,0.04000) (0.2625,0.04323) (0.2750,0.04652) (0.2875,0.04986) (0.3000,0.05325) (0.3125,0.05669) (0.3250,0.06016) (0.3375,0.06367) (0.3500,0.06722) (0.3625,0.07079) (0.3750,0.07438) (0.3875,0.07800) (0.4000,0.08163) (0.4125,0.08528) (0.4250,0.08895) (0.4375,0.09263) (0.4500,0.09631) (0.4625,0.10001) (0.4750,0.10371) (0.4875,0.10741) (0.5000,0.11111) (0.5125,0.11481) (0.5250,0.11852) (0.5375,0.12222) (0.5500,0.12591) (0.5625,0.12960) (0.5750,0.13328) (0.5875,0.13696) (0.6000,0.14062) (0.6125,0.14428) (0.6250,0.14793) (0.6375,0.15156) (0.6500,0.15519) (0.6625,0.15880) (0.6750,0.16240) (0.6875,0.16598) (0.7000,0.16955) (0.7125,0.17310) (0.7250,0.17664) (0.7375,0.18017) (0.7500,0.18367) (0.7625,0.18716) (0.7750,0.19064) (0.7875,0.19409) (0.8000,0.19753) (0.8125,0.20095) (0.8250,0.20435) (0.8375,0.20774) (0.8500,0.21110) (0.8625,0.21445) (0.8750,0.21778) (0.8875,0.22109) (0.9000,0.22438) (0.9125,0.22765) (0.9250,0.23090) (0.9375,0.23413) (0.9500,0.23734) (0.9625,0.24054) (0.9750,0.24371) (0.9875,0.24687) (1.0000,0.25000) (1,1) (0,1) (0,0)};
\addlegendentry{extreme Stokes wave exists}
\addplot[PaperInk,dotted,line width=1pt] coordinates {(0,0) (0.0082,0.00166) (0.0163,0.00337) (0.0245,0.00514) (0.0327,0.00695) (0.0408,0.00881) (0.0490,0.01073) (0.0572,0.01268) (0.0653,0.01469) (0.0735,0.01674) (0.0817,0.01883) (0.0898,0.02096) (0.0980,0.02313) (0.1062,0.02534) (0.1143,0.02759) (0.1225,0.02988) (0.1307,0.03219) (0.1388,0.03454) (0.1470,0.03693) (0.1552,0.03934) (0.1633,0.04178) (0.1715,0.04425) (0.1797,0.04674) (0.1878,0.04926) (0.1960,0.05180) (0.2042,0.05437) (0.2123,0.05695) (0.2205,0.05955) (0.2287,0.06218) (0.2368,0.06482) (0.2450,0.06747) (0.2532,0.07015) (0.2613,0.07283) (0.2695,0.07553) (0.2777,0.07824) (0.2858,0.08097) (0.2940,0.08370) (0.3022,0.08644) (0.3103,0.08919) (0.3185,0.09195) (0.3267,0.09472) (0.3348,0.09749) (0.3430,0.10027) (0.3512,0.10306) (0.3593,0.10584) (0.3675,0.10863) (0.3757,0.11143) (0.3838,0.11423) (0.3920,0.11702) (0.4002,0.11982) (0.4083,0.12262) (0.4165,0.12542) (0.4247,0.12822) (0.4328,0.13102) (0.4410,0.13382) (0.4492,0.13662) (0.4573,0.13941) (0.4655,0.14220) (0.4737,0.14499) (0.4818,0.14777) (0.4900,0.15055)};
\addlegendentry{$f=f_1(\gamma)$: above, $\vartheta<60^\circ$}
\addplot[only marks,mark=o,mark size=2pt,color=PaperInk,line width=.6pt] coordinates {(0.2,0.2)(0.2,0.5)(0.2,0.8)(0.2,0.95)(0.4,0.2)(0.4,0.5)(0.4,0.8)(0.4,0.95)(0.6,0.2)(0.6,0.5)(0.6,0.8)(0.6,0.95)(0.8,0.2)(0.8,0.5)(0.8,0.8)(0.8,0.95)(0.95,0.2)(0.95,0.5)(0.95,0.8)(0.95,0.95)};
\addlegendentry{scan of \autoref{subsec_scan}}
\end{axis}
\end{tikzpicture}
\caption{The parameter plane of the standing assumptions, $0<\gamma<1$ and $f=c_0e^{-2\gamma}\in(0,1)$. Below the grey curve, $c_0\le\sinh^2\gamma$ and the bifurcation takes place in the window $\Delta<0$. Accent region (above the dashed curve): \eqref{eq_floor_condition} holds, the floor (A) is proved (\autoref{cor_floor}), and the extreme Stokes wave exists (\autoref{thm_existence}, \autoref{cor_thin}); \emph{proved}. Above the dotted curve $f=f_1(\gamma)$, drawn for $\gamma\le0.49$, all slopes along the continuum stay below $60^\circ$ (\autoref{rem_C1_region}); \emph{proved}. White band between the grey and the dashed curves: (A) implies the extreme wave (\autoref{thm_band}), but (A) itself is \emph{open} (\autoref{rem_cone_status}). Circles mark the $20$ computed branches of the scan (\autoref{subsec_scan}); that they stay supercritical and turn back once in $\alpha$ is \emph{numerical evidence} only.}
\label{fig_region}
\end{figure}
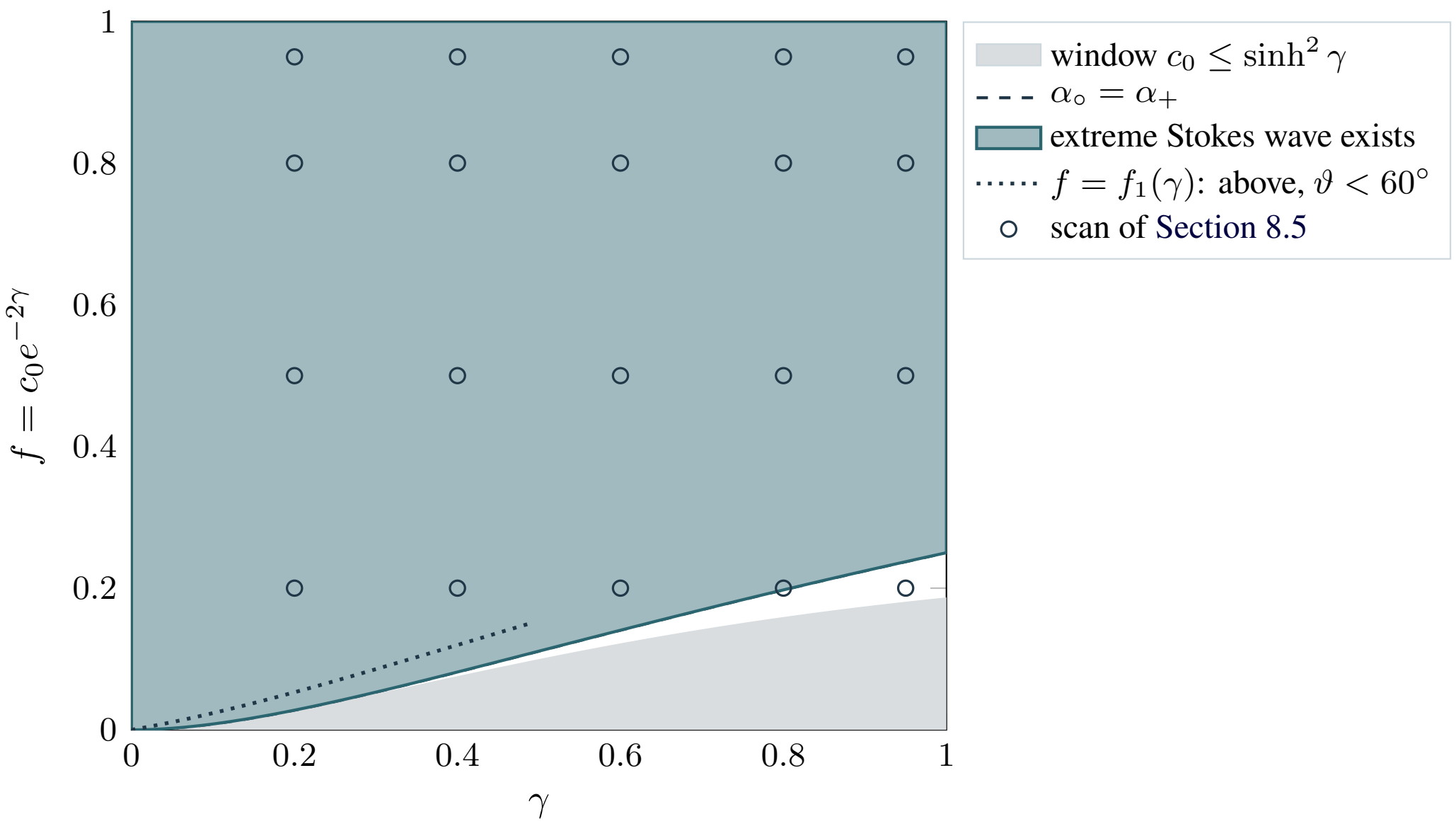

\begin{remark}[The sharper slope bound]\label{rem_C1_region}
Condition (C1) improves the slope bound of \autoref{thm_existence} from $90^\circ-\eta_\circ$ to $60^\circ$. Put $k=1+1/\gamma$ and $v(f)=\frac6{35}\big(6+\sqrt{1+35f}\big)-1$. Since $\lambda>0$, (C1) is $\sqrt\Delta\ge\alpha/6$ by \autoref{lem_structural}. For $\alpha>R_0$ this is the quadratic inequality $\frac{35}{36}\alpha^2-2R_0\alpha+R_0^2-c_0\ge0$, whose smaller root lies below $R_0$. So (C1) holds for all $\alpha\ge\alpha_\circ$ if and only if $u\ge v(f)$. Since $u-v(f)$ increases in $f$, this is equivalent to
\[
f\ge f_1(\gamma):=\frac{s^2-1}{35},\qquad s=\frac{3+\sqrt{9+k(1+k)}}{k}.
\]
For instance, $f_1(0.1)<0.024$, $f_1(0.2)<0.054$ and $f_1(0.49)<0.151$. The extreme waves of \autoref{cor_thin} with $f\ge f_1(\gamma)$ therefore have all slopes below $60^\circ$ along the continuum.
\end{remark}

\begin{remark}[The parameter region]\label{rem_cone_status}
\autoref{cor_thin} proves (C3) for every $\alpha>\max\{\alpha_\circ,\alpha_+\}$, at all $0<\gamma<1$ and $0<c_0<e^{2\gamma}$; see \autoref{fig_region}. In the floor region, where (A) is also proved, the extreme wave therefore exists. In the band $\sinh^2\gamma<c_0\le\gamma^2e^{2\gamma}/(1+\gamma)^2$, \autoref{thm_band} shows that (B) follows from (A), so only the floor (A) is open there: the parameter may a priori approach $\alpha_+$.
\end{remark}

\subsubsection{The existence chain at a glance}\label{subsec_audit}
\autoref{tab_chain} lists the steps from the continuum of \cite{LiWang2026} to the extreme wave, with the hypotheses of each step and the dependence of its constants on $\delta$. Only $\delta$-independent bounds are used in the limit $\delta\downarrow0$.

We comment on the passage to the limit, which is where the chain is most delicate.
\begin{itemize}
\item \emph{Function spaces.} For fixed $\delta$ the continuum lives in $(\alpha_s,\infty)\times X$, $X=C^{3+s}_{\rm per}(\overline D)$, and all arguments at fixed $\delta$ (steps 1, 6, 8, 9) use the topology of $X$. The limit $\delta\downarrow0$ is taken in weaker topologies: $\Psi_k\to\Psi_\infty$ uniformly, $\nabla\Psi_k\to\nabla\Psi_\infty$ in $L^2_{\rm loc}$, $\chi_{\Omega_k}\to\chi_{\Omega_\infty}$ in $L^1_{\rm loc}$ (\autoref{thm_limit}(a)). These need only the bounds of step 10, which depend on $\alpha$ through a compact set with $\Delta>0$, supplied by (R1) uniformly in $\delta$.
\item \emph{The free-boundary condition.} The Bernoulli condition is inherited in the weak form of the domain-variation identity \eqref{eq_VS}, whose terms converge under the convergences above together with $\alpha_k\to\alpha_\infty$ (\autoref{thm_limit}(c)). Thus the limit is a variational solution in the sense of \autoref{def_variational}; the crest analysis of \autoref{sec_corner} is formulated for this class.
\item \emph{Symmetry and monotonicity.} Evenness in $\theta$ and the monotonicity of $S$ on $(0,\pi)$ are closed conditions under uniform convergence of $\Psi$, so they pass to the limit (\autoref{cor_extreme_shape}). The strict conditions (N1)--(N5) are not claimed for the limit.
\item \emph{The branch reaches stagnation.} The local classification of \autoref{thm_main_stokes} describes any stagnation point; that the continuum actually produces one is a separate fact. It comes from step 9: the contact solutions have $\min S_kQ_{\alpha_k}(S_k)=\delta_k\to0$ with $\Delta(\alpha_k)\ge\bar\Delta$, so $r_-(\alpha_k)-\max S_k\to0$, and \autoref{thm_limit}(d) gives a stagnation point of the limit on the circle $R=r_-(\alpha_\infty)$.
\item \emph{No laminar limit.} The limit is not circular, by \autoref{thm_limit}(e). A circular limit with a stagnation point on $R=r_-(\alpha_\infty)$ would be the annulus $\{1<R<r_-\}$ with $\Psi=1-\log R/\log r_-$, whose outer boundary carries the speed $1/(r_-\log r_-)\ne0$; this contradicts the uniform bound $|\nabla\Psi_\infty|^2\le C(r_--R)$ used in the proof of \autoref{thm_limit}(d).
\end{itemize}

\subsection{The direction of bifurcation}\label{sec_supercritical}

The parameter floor (A) is the first of the two open bounds. This subsection decides it near the bifurcation point. We compute the second-order coefficient of $\alpha$ along the local branch in closed form. It is positive for every $\gamma\in(0,1)$ and every $c_0>0$.

\subsubsection{The coefficient $\alpha_2$}

Work in the cylinder variables $t=\log R$, and write $T(\theta)=\log S(\theta)$. Then $\Psi$ is harmonic in $\{0<t<T(\theta)\}$, $\Psi=1$ on $t=0$, and on $t=T$
\begin{equation}\label{eq_cyl_bc}
\Psi=0,\qquad \Psi_t^2+\Psi_\theta^2=g_\alpha(T),\qquad g_\alpha(T):=\frac{e^{2T}Q_\alpha(e^T)}{p_0^2} .
\end{equation}
Throughout this subsection
\[
C=\coth\gamma,\qquad \lambda=\coth\gamma-\frac1\gamma\in(0,1),\qquad m=\frac{R_0^2}{c_0} .
\]
Write $\hat g_n=\gamma^2\,\partial_T^ng_{\alpha_c}(\gamma)$. Since $p_0^2=c_0R_0^2\gamma^2$ and $\alpha_c=R_0+c_0/\sinh\gamma$, a direct computation gives
\[
\gamma^2g_{\alpha_c}(\gamma)=1,\quad \hat g_1=-2C,\quad \hat g_2=2m-6-10C,\quad \hat g_3=18m-30-38C,
\]
and $\gamma^2\partial_\alpha\partial_Tg_\alpha(\gamma)=-2m/R_0$.
Moreover $\partial_\alpha g_\alpha(\gamma)=0$, since $Q_\alpha(R_0)=c_0$ for every $\alpha$.

\begin{theorem}[Second-order coefficient]\label{thm_alpha2}
Let $0<\gamma<1$ and $c_0>0$. Normalize the local Crandall--Rabinowitz curve at $(\alpha_c,R_0)$ by
\[
\log S=\gamma+\varepsilon\cos\theta+O(\varepsilon^2),
\]
where the $O(\varepsilon^2)$ term has no $\cos\theta$ component. Then $\alpha=\alpha_c+\alpha_2\varepsilon^2+O(\varepsilon^4)$, with
\begin{equation}\label{eq_alpha2}
\alpha_2=\frac{R_0}{2m}\Big[\Big(\frac1{2\lambda}-\frac C4\Big)w^2+\Big(2C+\frac94\Big)w+\frac34(C+1)^2(4C+1)\Big],\qquad w=m-(C+1)(3C+2).
\end{equation}
\end{theorem}
\begin{proof}
The problem is analytic in $(\varepsilon,\alpha)$. By the Crandall--Rabinowitz theorem \cite{CrandallRabinowitz1971}, in its analytic form \cite[Chapter~8]{BuffoniToland2003}, the local curve is analytic, and its Taylor coefficients are those of the formal expansion below. The half-period shift $\theta\mapsto\theta+\pi$ maps $\varepsilon$ to $-\varepsilon$ and fixes $\alpha$, so $\alpha(\varepsilon)$ is even. Write
\[
\Psi=1-\frac t\gamma+\sum_{n\ge1}\varepsilon^n\psi_n,\qquad T=\gamma+\sum_{n\ge1}\varepsilon^n\eta_n,\qquad \eta_1=\cos\theta ,
\]
and transfer \eqref{eq_cyl_bc} to $t=\gamma$ by Taylor expansion. All $\psi_n$ are harmonic and vanish at $t=0$. All functions below are evaluated at $t=\gamma$.

\emph{Order $1$.} The Dirichlet condition gives $\psi_1=\eta_1/\gamma$. Hence $\psi_1=\sinh t\cos\theta/(\gamma\sinh\gamma)$ and $\partial_t\psi_1=(C/\gamma)\cos\theta$. The Bernoulli condition reads $-(2/\gamma)\partial_t\psi_1=g_1\eta_1$, i.e.\ $\hat g_1=-2C$. Using the formula for $\hat g_1$, this is $\alpha=R_0+c_0/\sinh\gamma=\alpha_c$.

\emph{Order $2$.} Write $\eta_2=b_0+b_2\cos2\theta$ and $\psi_2=A_0t+A_2\sinh2t\cos2\theta$. The two boundary conditions are
\[
\psi_2=\frac{\eta_2}\gamma-\eta_1\partial_t\psi_1,\qquad -\frac2\gamma\big(\partial_t\psi_2+\eta_1\partial_t^2\psi_1\big)+|\nabla\psi_1|^2=g_1\eta_2+\tfrac12g_2\eta_1^2 .
\]
Here $|\nabla\psi_1|^2=\gamma^{-2}\big(\tfrac{C^2+1}2+\tfrac{C^2-1}2\cos2\theta\big)$ and $\eta_1\partial_t^2\psi_1=\gamma^{-1}\cos^2\theta$. The modes $0$ and $2$ give two linear systems of size two. In mode $2$ we use $2C-4\coth2\gamma=-2/C$. The solutions are
\[
b_0=\frac{\hat g_2/4+(1-C^2)/2-C/\gamma}{2\lambda},\qquad b_2=-\frac C4\Big(\frac{\hat g_2}2+1-3C^2\Big),
\]
together with $\gamma^2A_0=b_0-C/2$ and $\gamma A_2\sinh2\gamma=b_2-C/2$.

\emph{Order $3$.} Only the $\cos\theta$ component is needed, and $\eta_3$ has none. The Dirichlet condition determines the $\cos\theta$ component of $\psi_3$ from
\[
\psi_3=\frac{\eta_3}{\gamma}-\eta_1\partial_t\psi_2-\eta_2\partial_t\psi_1-\tfrac12\eta_1^2\partial_t^2\psi_1 ,
\]
and then $\partial_t\psi_3=C\psi_3$ in that mode. The Bernoulli condition at order $3$ is
\[
-\frac2\gamma u_3+2\big(\partial_t\psi_1\,u_2+\partial_\theta\psi_1\,v_2\big)=g_1\eta_3+g_2\eta_1\eta_2+\tfrac16g_3\eta_1^3+\alpha_2\,\partial_\alpha g_1\,\eta_1 ,
\]
where
\[
u_2=\partial_t\psi_2+\eta_1\partial_t^2\psi_1,\qquad v_2=\partial_\theta\psi_2+\eta_1\partial_t\partial_\theta\psi_1,
\]
\[
u_3=\partial_t\psi_3+\eta_1\partial_t^2\psi_2+\eta_2\partial_t^2\psi_1+\tfrac12\eta_1^2\partial_t^3\psi_1 .
\]
Project on $\cos\theta$, multiply by $\gamma^2$, and insert the order-$2$ solution. Using $4C\coth2\gamma=2(1+C^2)$, the result is
\[
\frac{4C}{\gamma}\Big(b_0-\frac C2\Big)+2C^2\Big(b_2-\frac C2\Big)+2(C^2-1)\sigma+2C-\hat g_2\,\sigma-\frac{\hat g_3}8=-\frac{2m}{R_0}\,\alpha_2,\qquad \sigma=b_0+\frac{b_2}2 .
\]

\emph{Simplification.} Put $u=\hat g_2/4=(m-3-5C)/2$ and $\kappa=(1-C^2)/2-C/\gamma$. Move everything to the right-hand side. The terms containing $b_0$ combine to $8\lambda b_0^2$, because their coefficient is $4u+2-2C^2-4C/\gamma=8\lambda b_0$. The terms containing $b_2$ combine to $-\frac C4(2u+1-3C^2)^2$, because their coefficient is $2u+1-3C^2$ and $b_2=-\frac C4(2u+1-3C^2)$. Collecting terms,
\[
\frac{2m}{R_0}\alpha_2=\frac{2}{\lambda}(u+\kappa)^2-\frac C4(2u+1-3C^2)^2+\frac92u+3+\frac92C+C^3+\frac{2C^2}{\gamma} .
\]
Now substitute $C/\gamma=C^2-C\lambda$. Then $2(u+\kappa)=w+2C\lambda$ and $2u+1-3C^2=w$, where $w=m-(C+1)(3C+2)$. The terms containing $\lambda$ outside the first square cancel, and the constant term factors as $3C^3+\frac{27}4C^2+\frac92C+\frac34=\frac34(C+1)^2(4C+1)$. This gives \eqref{eq_alpha2}.
\end{proof}

\begin{remark}[The Crandall--Rabinowitz hypotheses]\label{rem_all_gamma}
The hypotheses of the Crandall--Rabinowitz theorem are verified in \cite{LiWang2026} for $0<\gamma<1$ and $c_0>\sinh^2\gamma$. The following argument shows that they hold for every $c_0>0$. By the order-one computation, the mode $\cos k\theta$, $k\ge1$, is singular at the laminar flow exactly when
\[
\alpha=R_0+\frac{c_0}{R_0}\big(1+k\coth k\gamma\big).
\]
These values are distinct, since $k\coth k\gamma$ is strictly increasing in $k$. The radial mode is singular only at $\alpha_\circ<\alpha_c$ (\autoref{lem_radial_crossing}). Hence the kernel at $\alpha_c$ is spanned by the mode $k=1$. The transversality condition reduces to $\partial_\alpha\hat g_1=-2R_0/c_0\neq0$, which is the coefficient of $\alpha_2$ in the order-three equation. Thus \autoref{thm_alpha2} holds for every $\gamma\in(0,1)$ and every $c_0>0$.
\end{remark}

\begin{remark}[Comparison with \cite{LiWang2026}]\label{rem_LW_direction}
In \cite{LiWang2026} the bifurcation is shown to be of pitchfork type. Its direction is determined by the sign of an explicit quantity $\mathcal O$ \cite[Theorem~2.1(3) and (4.31)]{LiWang2026}, a polynomial of degree $18$ in $e^\gamma$ divided by a positive factor. Formula \eqref{eq_alpha2} gives a factorized form in which the sign can be analysed. We compared the two signs on a grid of $1600$ parameter pairs with $0.05\le\gamma\le0.99$ and $1<m<(C+1)^2$, which is the overlap of \eqref{eq_standing} with $c_0>\sinh^2\gamma$. The signs agree at every point; by \autoref{thm_supercritical} both are positive there. Outside this range we use only \eqref{eq_alpha2}, which we also checked against the spectral solver; for example, at $\gamma=0.3$ and $m=3(C+1)^2$ both give $\alpha_2=3.88677$.
\end{remark}

\begin{remark}[Independent check]\label{rem_alpha2_check}
We compared \eqref{eq_alpha2} with $\alpha_2$ extracted from the spectral solver of \autoref{sec_numerics}, by fitting $\alpha-\alpha_c$ at small amplitude. The comparison includes values of $\gamma$ in $[0.05,1)$ and $c_0$ from $1.001\sinh^2\gamma$ to $100\sinh^2\gamma$. The two values agree to between three and seven significant digits. The agreement is weakest for the thinnest annulus, where the fit is least accurate.
\end{remark}

\subsubsection{Sign of $\alpha_2$}

In the admissible range $c_0>\sinh^2\gamma$ we have $m<R_0^2/\sinh^2\gamma=(C+1)^2$. Write $w=-(C+1)y$. As $c_0$ runs over $(\sinh^2\gamma,\infty)$, $y$ runs over $(2C+1,3C+2)$, and
\begin{equation}\label{eq_alpha2_y}
\begin{gathered}
\alpha_2=\frac{R_0(C+1)}{2m}\,\big(\mathsf a\,y^2-\mathsf b\,y+\mathsf c\big),\\
\mathsf a=(C+1)\Big(\frac1{2\lambda}-\frac C4\Big),\qquad \mathsf b=2C+\frac94,\qquad \mathsf c=\frac34(C+1)(4C+1).
\end{gathered}
\end{equation}

\begin{lemma}\label{lem_langevin}
For every $\gamma>0$: $\lambda<\gamma/3$, $C\lambda<1$, $C>1/\gamma$, and $\mathsf a>(C+1)/(4\lambda)>0$.
\end{lemma}
\begin{proof}
The bound $\lambda<\gamma/3$ is equivalent to $\gamma\cosh\gamma<(1+\gamma^2/3)\sinh\gamma$. Comparing the coefficients of $\gamma^{2n+1}$ in the two power series, we need $\frac1{(2n)!}\le\frac1{(2n+1)!}+\frac1{3(2n-1)!}$ for $n\ge1$, with equality for $n=0$. Multiplying by $(2n+1)!$, this is $2n\le\frac{2n(2n+1)}3$, i.e.\ $n\ge1$. Equality holds for $n=1$ and strict inequality for $n\ge2$, so the inequality between the series is strict for $\gamma>0$. The inequality $C\lambda<1$ is equivalent to $\coth^2\gamma-1<\coth\gamma/\gamma$, i.e.\ to $\gamma<\sinh\gamma\cosh\gamma$. The inequality $C>1/\gamma$ is $\lambda>0$. Finally $1/(2\lambda)-C/4=(2-C\lambda)/(4\lambda)>1/(4\lambda)$.
\end{proof}

\begin{theorem}[Supercriticality]\label{thm_supercritical}
Let $0<\gamma<1$ and $c_0>0$. Then $\alpha_2>0$. Hence the bifurcation is supercritical: near $(\alpha_c,R_0)$ the continuum $\mathcal C^+_\delta$ lies in $\{\alpha>\alpha_c\}$.
\end{theorem}
\begin{proof}
By \autoref{lem_langevin}, $\mathsf a>(C+1)/(4\lambda)>3(C+1)/(4\gamma)$, so
\[
4\mathsf a\mathsf c>\frac{9C}{\gamma}(C+1)^2>\frac9{\gamma^2}(C+1)^2 .
\]
Also $\mathsf b\le\frac94(C+1)$, so $\mathsf b^2\le\frac{81}{16}(C+1)^2$. Since $\gamma<1<4/3$, we have $9/\gamma^2>81/16$, hence $4\mathsf a\mathsf c>\mathsf b^2$. The quadratic $\mathsf ay^2-\mathsf by+\mathsf c$ is then positive for every real $y$, so $\alpha_2>0$ for every $c_0>0$ by \eqref{eq_alpha2_y}.

Finally, $\alpha=\alpha_c+\alpha_2\varepsilon^2+O(\varepsilon^4)$ with $\alpha_2>0$ gives $\alpha>\alpha_c$ for small $\varepsilon\ne0$. By local uniqueness, this describes all of $\mathcal C^+_\delta$ near the bifurcation point.
\end{proof}

\section{Numerical evidence}\label{sec_numerics}

Above the floor \eqref{eq_floor_condition} the theorems give the extreme wave; in the band $\sinh^2\gamma<c_0\le\gamma^2e^{2\gamma}/(1+\gamma)^2$ they leave alternatives (I) and (III). We computed the monotone continuum numerically for three parameter pairs above the floor, and for $20$ further pairs in \autoref{subsec_scan}, one of which lies in the band. The computations illustrate the continuum, test the exact identities of the paper on computed waves, and show how far the a priori bounds are from sharp. They are evidence only: no statement of Sections~\ref{sec_prelim}--\ref{sec_towards} depends on them.

\subsection{Method}
Since $\Psi$ is harmonic with $\Psi=1$ on $R=1$ and $\Psi=0$ on $R=S(\theta)$, the fluid domain is conformally equivalent to an annulus $\{1<|\zeta|<\varrho\}$, with unknown modulus $\varrho>1$. In that model domain $\Psi=1-\log|\zeta|/\log\varrho$. We write the conformal map as
\[
z=f(\zeta)=\zeta\,e^{g(\zeta)},\qquad g(\zeta)=\sum_{k=1}^{N}b_k\Big(\big(\tfrac{\zeta}{\varrho}\big)^{k}-(\varrho\,\zeta)^{-k}\Big),\qquad b_k\in\mathbb R .
\]
Then $|f|=1$ on $|\zeta|=1$, and the waves are even in $\theta$. On $|\zeta|=\varrho$ the Bernoulli condition becomes
\[
\frac{p_0^2}{\varrho^2\log^2\varrho\;|f'(\zeta)|^2}=Q_\alpha\big(|f(\zeta)|\big).
\]
We impose it at the $N+1$ points $\zeta_j=\varrho e^{i\pi j/N}$, $0\le j\le N$. The unknowns are $(b_1,\dots,b_N,\varrho,\alpha)$. Their number exceeds the number of equations by one, and pseudo-arclength continuation supplies the missing equation. The Jacobian with respect to $b$ is exact; the derivative in $\varrho$ uses a finite difference.

We start from the laminar flow at $\alpha=\alpha_c$, with $\varrho=R_0$ and $b=0$, in the direction of the first mode. Whenever the coefficients in the top eighth of the spectrum exceed $10^{-7}$ of the largest one, we double $N$, from $128$ up to $2048$. A run stops at the first of these events: the spectral tail exceeds $2\times10^{-6}$ at $N=2048$; monotonicity from crest to trough fails on a grid of $6000$ points; or $|\nabla\Psi|$ at the crest drops below $0.08$.

As an a posteriori check, we evaluated the Bernoulli residual on a grid eight times finer than the collocation grid at the last computed wave. The maximum is $3.9\times10^{-3}$ for run~(a) and $4.2\times10^{-3}$ for run~(b). In both runs, the residual is largest next to the crest.

The computed waves also satisfy the exact identities of the paper, which were not imposed by the scheme. We evaluated them at the last computed wave with $32768$ boundary points and adaptive quadrature along the axes. The results are as follows.
\begin{itemize}
\item The flow-force identity \eqref{eq_flow_force} holds to a relative accuracy of $1.5\times10^{-7}$ in run~(a) and $3.0\times10^{-7}$ in run~(b).
\item The angular flow-force identity \eqref{eq_angular_flow_force} holds to $1.7\times10^{-6}$ and $1.8\times10^{-6}$.
\item In the maximum-velocity principle (\autoref{thm_max_velocity}), $\sup R^2|\nabla\Psi|^2$ over $60$ interior circles and $\max S^2Q_\alpha(S)/p_0^2$ agree to eight digits.
\item The crest--trough comparison $G_\alpha(S_c)\le0\le G_\alpha(S_t)$ and the bound $S\le r_-$ (\autoref{prop_crest_trough}, \autoref{thm_localization}) hold at all $351$ computed branch points of the three runs.
\end{itemize}
The script \texttt{verify\_identities.py} accompanies the paper.

To test the resolution limit, we repeated run~(b) with $N$ allowed up to $4096$, halving the arclength step whenever the spectral tail exceeded $2\times10^{-6}$. The branch then continues to $q_{\rm cr}/q_{\rm cr}^{\,0}=0.13$ and $r_--S_c=7.5\times10^{-3}$. Along the way the maximal slope increases to $28.5^\circ$, $\alpha$ stays above $\alpha_c$, and the waves remain monotone. The fine-grid residual at the last wave is $8.6\times10^{-3}$. The trend towards a $30^\circ$ slope at the crest, and hence towards the $120^\circ$ corner, continues, but resolving the corner itself needs a crest-adapted method.

\subsection{Convergence and residuals}\label{subsec_convergence}
For waves away from stagnation the free surface is analytic, and the scheme converges geometrically in $N$. We tested this on two waves of run~(a) ($\gamma=0.5$, $f=0.6$), produced by continuation and then re-solved at fixed $\alpha$ for several $N$. The column ``residual'' is the maximum of the Bernoulli residual on a grid eight times finer than the collocation grid, and ``tail'' is the size of the top eighth of the coefficients relative to the largest.
\begin{center}
\small
\begin{tabular}{@{}rccccc@{}}
\toprule
$N$ & $\varrho$ & $S_c$ & max slope & tail & residual\\
\midrule
\multicolumn{6}{@{}l}{\emph{Wave 1:} $\alpha=5.75698$, $q_{\rm cr}/q_{\rm cr}^{\,0}=0.798$}\\
$32$ & $1.5637093838$ & $1.721470355$ & $7.756^\circ$ & $6.1\times10^{-5}$ & $5.2\times10^{-4}$\\
$64$ & $1.5637094531$ & $1.721466257$ & $7.748^\circ$ & $5.7\times10^{-8}$ & $5.5\times10^{-7}$\\
$128$ & $1.5637094531$ & $1.721466256$ & $7.748^\circ$ & $2.1\times10^{-13}$ & $8.2\times10^{-13}$\\
$256$ & $1.5637094531$ & $1.721466256$ & $7.748^\circ$ & $4.5\times10^{-18}$ & $9.3\times10^{-15}$\\
\midrule
\multicolumn{6}{@{}l}{\emph{Wave 2:} $\alpha=6.46483$, $q_{\rm cr}/q_{\rm cr}^{\,0}=0.489$}\\
$128$ & $1.5115458653$ & $1.779663813$ & $19.504^\circ$ & $2.3\times10^{-4}$ & $2.4\times10^{-2}$\\
$256$ & $1.5115598336$ & $1.779314022$ & $18.883^\circ$ & $2.9\times10^{-6}$ & $6.2\times10^{-4}$\\
$512$ & $1.5115598418$ & $1.779312875$ & $18.890^\circ$ & $2.9\times10^{-9}$ & $6.8\times10^{-7}$\\
$1024$ & $1.5115598418$ & $1.779312874$ & $18.890^\circ$ & $1.2\times10^{-14}$ & $1.2\times10^{-12}$\\
$2048$ & $1.5115598418$ & $1.779312874$ & $18.890^\circ$ & $6.6\times10^{-19}$ & $2.3\times10^{-14}$\\
\bottomrule
\end{tabular}
\end{center}
For wave 2, Newton's method at fixed $\alpha$ does not converge for $N\le64$. The rate of convergence deteriorates as the crest approaches stagnation, because the singularity of the analytic continuation of the free surface approaches it; with $N=2048$ the dense scheme stops near $q_{\rm cr}/q_{\rm cr}^{\,0}\approx0.2$. For the last computed waves we have no convergence study in $N$. We report instead the a~posteriori Bernoulli residual ($3.9\times10^{-3}$ and $4.2\times10^{-3}$ for the dense runs, $1.9\times10^{-6}$ for the high-resolution run of \autoref{subsec_hp}) and the errors in the exact identities listed above, which are not imposed by the scheme and serve as conservation checks. Values reported at those waves, such as the maximal slope, should be read together with these residuals. The script \texttt{convergence\_test.py} reproduces the table.

\subsection{Results}
In the table, $q_{\rm cr}$ is the fluid speed at the crest along the branch and $q_{\rm cr}^{\,0}$ is its laminar value at $\alpha_c$. The gap $r_--S_c$ is measured at the last computed wave.
\begin{center}
\small
\begin{tabular}{@{}lccccccccc@{}}
\toprule
run & $\gamma$ & $c_0$ & $\alpha_+$ & $\alpha_\circ$ & $\alpha_c$ & $\max\alpha$ & $q_{\rm cr}/q_{\rm cr}^{\,0}$ & $r_--S_c$ & max slope\\
\midrule
(a) & $0.50$ & $1.631$ & $2.926$ & $4.616$ & $4.779$ & $6.522$ & $0.222$ & $8.4\times10^{-3}$ & $26.7^\circ$\\
(b) & $0.90$ & $1.815$ & $3.807$ & $4.017$ & $4.228$ & $4.930$ & $0.181$ & $1.4\times10^{-2}$ & $27.3^\circ$\\
(c) & $0.25$ & $1.319$ & $2.432$ & $6.420$ & $6.505$ & $9.899$ & $0.345$ & $9.1\times10^{-3}$ & $24.2^\circ$\\
\bottomrule
\end{tabular}
\end{center}

\begin{figure}[!ht]
\centering
\begin{tikzpicture}
\begin{axis}[width=7.2cm,height=6.2cm,xlabel={$\alpha$},ylabel={$q_{\rm cr}/q_{\rm cr}^{\,0}$},
  ymin=0,ymax=1.05,xmin=4.0,xmax=10.2,tick label style={font=\footnotesize},label style={font=\small},
  legend style={font=\footnotesize,at={(0.03,0.03)},anchor=south west,draw=PaperRule},legend cell align=left]
\addplot[PaperAccent,line width=1pt] coordinates {(4.7846,1.0000) (4.7974,0.9938) (4.8160,0.9881) (4.8430,0.9818) (4.8805,0.9745) (4.9201,0.9676) (4.9598,0.9610) (4.9994,0.9544) (5.0391,0.9479) (5.0788,0.9413) (5.1186,0.9346) (5.1583,0.9278) (5.2179,0.9174) (5.2577,0.9103) (5.2975,0.9031) (5.3373,0.8957) (5.3771,0.8881) (5.4169,0.8803) (5.4567,0.8723) (5.4965,0.8642) (5.5363,0.8558) (5.5761,0.8472) (5.6160,0.8383) (5.6558,0.8292) (5.7156,0.8151) (5.7554,0.8052) (5.7953,0.7951) (5.8351,0.7845) (5.8750,0.7736) (5.9148,0.7623) (5.9547,0.7505) (5.9946,0.7382) (6.0345,0.7252) (6.0743,0.7117) (6.1142,0.6973) (6.1741,0.6742) (6.2339,0.6484) (6.2738,0.6295) (6.3137,0.6086) (6.3536,0.5852) (6.4134,0.5429) (6.4533,0.5062) (6.4932,0.4531) (6.5210,0.3605) (6.5216,0.3476) (6.5216,0.3390) (6.5215,0.3376) (6.5215,0.3359) (6.5214,0.3326) (6.5213,0.3299) (6.5211,0.3268) (6.5209,0.3230) (6.5206,0.3184) (6.5201,0.3128) (6.5194,0.3061) (6.5183,0.2978) (6.5169,0.2876) (6.5148,0.2747) (6.5117,0.2583) (6.5044,0.2223)};
\addlegendentry{(a) $\gamma=0.5$}
\addplot[PaperInk,line width=1pt,dashed] coordinates {(4.2290,1.0000) (4.2325,0.9902) (4.2377,0.9813) (4.2444,0.9727) (4.2529,0.9642) (4.2635,0.9554) (4.2764,0.9462) (4.2921,0.9363) (4.3099,0.9261) (4.3277,0.9166) (4.3457,0.9075) (4.3638,0.8986) (4.3820,0.8898) (4.4002,0.8811) (4.4184,0.8724) (4.4368,0.8637) (4.4551,0.8548) (4.4735,0.8459) (4.4920,0.8369) (4.5104,0.8276) (4.5290,0.8182) (4.5475,0.8086) (4.5661,0.7987) (4.5847,0.7885) (4.6033,0.7780) (4.6220,0.7672) (4.6407,0.7560) (4.6594,0.7443) (4.6781,0.7322) (4.6969,0.7196) (4.7156,0.7063) (4.7344,0.6923) (4.7532,0.6775) (4.7720,0.6618) (4.7909,0.6449) (4.8097,0.6266) (4.8285,0.6066) (4.8662,0.5586) (4.8849,0.5283) (4.9036,0.4894) (4.9266,0.4017) (4.9291,0.3753) (4.9300,0.3507) (4.9300,0.3420) (4.9298,0.3302) (4.9295,0.3198) (4.9290,0.3108) (4.9285,0.3020) (4.9278,0.2928) (4.9270,0.2831) (4.9252,0.2628) (4.9240,0.2517) (4.9227,0.2389) (4.9210,0.2237) (4.9191,0.2051) (4.9168,0.1805)};
\addlegendentry{(b) $\gamma=0.9$}
\addplot[PaperMuted,line width=1pt,dotted] coordinates {(6.5125,1.0000) (6.5256,0.9978) (6.5570,0.9945) (6.6097,0.9901) (6.6696,0.9855) (6.7295,0.9809) (6.7895,0.9762) (6.8494,0.9715) (6.9094,0.9667) (6.9493,0.9635) (7.0093,0.9585) (7.0692,0.9535) (7.1292,0.9484) (7.1892,0.9432) (7.2491,0.9379) (7.3091,0.9325) (7.3691,0.9271) (7.4090,0.9233) (7.4690,0.9177) (7.5290,0.9119) (7.5890,0.9060) (7.6489,0.9001) (7.7089,0.8940) (7.7689,0.8877) (7.8289,0.8814) (7.8888,0.8749) (7.9288,0.8705) (7.9888,0.8638) (8.0488,0.8569) (8.1088,0.8499) (8.1688,0.8427) (8.2287,0.8354) (8.2887,0.8279) (8.3487,0.8202) (8.3887,0.8150) (8.4487,0.8069) (8.5087,0.7987) (8.5687,0.7902) (8.6486,0.7785) (8.7086,0.7694) (8.7686,0.7600) (8.8286,0.7504) (8.8886,0.7403) (8.9286,0.7334) (8.9886,0.7228) (9.0486,0.7117) (9.1086,0.7001) (9.1685,0.6880) (9.2285,0.6752) (9.3085,0.6571) (9.3685,0.6426) (9.4085,0.6323) (9.4685,0.6160) (9.5285,0.5982) (9.5885,0.5787) (9.6485,0.5568) (9.7285,0.5221) (9.7884,0.4892) (9.8484,0.4431) (9.8992,0.3451)};
\addlegendentry{(c) $\gamma=0.25$}
\end{axis}
\end{tikzpicture}\hfill
\begin{tikzpicture}
\begin{axis}[width=6.6cm,height=6.6cm,axis equal image,hide axis,clip=false,xmin=-3.5,xmax=3.5,ymin=-3.5,ymax=3.5]
\draw[PaperMuted,fill=PaperMuted!15] (axis cs:0,0) circle[radius=1];
\draw[PaperAccent,dashed] (axis cs:0,0) circle[radius=2.8618];
\addplot[PaperInk,line width=1pt,fill=PaperAccent!10,fill opacity=.6] coordinates {(-0.0000,-2.0895) (-0.1089,-2.0866) (-0.2166,-2.0783) (-0.3228,-2.0646) (-0.4274,-2.0457) (-0.5299,-2.0218) (-0.6303,-1.9931) (-0.7284,-1.9598) (-0.8238,-1.9221) (-0.9165,-1.8801) (-1.0063,-1.8342) (-1.0930,-1.7845) (-1.1766,-1.7313) (-1.2569,-1.6748) (-1.3338,-1.6152) (-1.4073,-1.5527) (-1.4773,-1.4876) (-1.5436,-1.4201) (-1.6064,-1.3503) (-1.6655,-1.2785) (-1.7209,-1.2050) (-1.7727,-1.1299) (-1.8208,-1.0534) (-1.8653,-0.9757) (-1.9062,-0.8970) (-1.9434,-0.8174) (-1.9772,-0.7373) (-2.0074,-0.6566) (-2.0343,-0.5756) (-2.0577,-0.4945) (-2.0779,-0.4133) (-2.0949,-0.3323) (-2.1087,-0.2514) (-2.1195,-0.1710) (-2.1273,-0.0910) (-2.1323,-0.0116) (-2.1344,0.0671) (-2.1339,0.1450) (-2.1308,0.2221) (-2.1253,0.2982) (-2.1173,0.3733) (-2.1071,0.4474) (-2.0947,0.5203) (-2.0802,0.5921) (-2.0637,0.6626) (-2.0454,0.7318) (-2.0252,0.7998) (-2.0034,0.8664) (-1.9800,0.9317) (-1.9551,0.9956) (-1.9288,1.0582) (-1.9012,1.1193) (-1.8723,1.1791) (-1.8423,1.2374) (-1.8113,1.2944) (-1.7793,1.3499) (-1.7463,1.4041) (-1.7126,1.4569) (-1.6781,1.5083) (-1.6429,1.5584) (-1.6071,1.6071) (-1.5708,1.6545) (-1.5339,1.7006) (-1.4967,1.7455) (-1.4591,1.7891) (-1.4212,1.8314) (-1.3831,1.8726) (-1.3448,1.9125) (-1.3063,1.9514) (-1.2678,1.9890) (-1.2292,2.0256) (-1.1906,2.0611) (-1.1521,2.0956) (-1.1136,2.1290) (-1.0753,2.1615) (-1.0372,2.1929) (-0.9993,2.2235) (-0.9616,2.2531) (-0.9242,2.2818) (-0.8872,2.3096) (-0.8505,2.3366) (-0.8142,2.3628) (-0.7783,2.3882) (-0.7429,2.4128) (-0.7079,2.4367) (-0.6735,2.4598) (-0.6396,2.4822) (-0.6063,2.5039) (-0.5736,2.5249) (-0.5416,2.5452) (-0.5102,2.5649) (-0.4795,2.5840) (-0.4495,2.6024) (-0.4203,2.6203) (-0.3918,2.6375) (-0.3642,2.6541) (-0.3373,2.6701) (-0.3113,2.6855) (-0.2862,2.7004) (-0.2620,2.7146) (-0.2387,2.7283) (-0.2163,2.7413) (-0.1950,2.7538) (-0.1746,2.7656) (-0.1552,2.7768) (-0.1369,2.7873) (-0.1197,2.7972) (-0.1035,2.8063) (-0.0885,2.8145) (-0.0745,2.8220) (-0.0617,2.8286) (-0.0501,2.8342) (-0.0396,2.8387) (-0.0304,2.8421) (-0.0223,2.8446) (-0.0155,2.8462) (-0.0099,2.8472) (-0.0056,2.8476) (-0.0025,2.8478) (-0.0006,2.8478) (-0.0000,2.8478) (0.0006,2.8478) (0.0025,2.8478) (0.0056,2.8476) (0.0099,2.8472) (0.0155,2.8462) (0.0223,2.8446) (0.0304,2.8421) (0.0396,2.8387) (0.0501,2.8342) (0.0617,2.8286) (0.0745,2.8220) (0.0885,2.8145) (0.1035,2.8063) (0.1197,2.7972) (0.1369,2.7873) (0.1552,2.7768) (0.1746,2.7656) (0.1950,2.7538) (0.2163,2.7413) (0.2387,2.7283) (0.2620,2.7146) (0.2862,2.7004) (0.3113,2.6855) (0.3373,2.6701) (0.3642,2.6541) (0.3918,2.6375) (0.4203,2.6203) (0.4495,2.6024) (0.4795,2.5840) (0.5102,2.5649) (0.5416,2.5452) (0.5736,2.5249) (0.6063,2.5039) (0.6396,2.4822) (0.6735,2.4598) (0.7079,2.4367) (0.7429,2.4128) (0.7783,2.3882) (0.8142,2.3628) (0.8505,2.3366) (0.8872,2.3096) (0.9242,2.2818) (0.9616,2.2531) (0.9993,2.2235) (1.0372,2.1929) (1.0753,2.1615) (1.1136,2.1290) (1.1521,2.0956) (1.1906,2.0611) (1.2292,2.0256) (1.2678,1.9890) (1.3063,1.9514) (1.3448,1.9125) (1.3831,1.8726) (1.4212,1.8314) (1.4591,1.7891) (1.4967,1.7455) (1.5339,1.7006) (1.5708,1.6545) (1.6071,1.6071) (1.6429,1.5584) (1.6781,1.5083) (1.7126,1.4569) (1.7463,1.4041) (1.7793,1.3499) (1.8113,1.2944) (1.8423,1.2374) (1.8723,1.1791) (1.9012,1.1193) (1.9288,1.0582) (1.9551,0.9956) (1.9800,0.9317) (2.0034,0.8664) (2.0252,0.7998) (2.0454,0.7318) (2.0637,0.6626) (2.0802,0.5921) (2.0947,0.5203) (2.1071,0.4474) (2.1173,0.3733) (2.1253,0.2982) (2.1308,0.2221) (2.1339,0.1450) (2.1344,0.0671) (2.1323,-0.0116) (2.1273,-0.0910) (2.1195,-0.1710) (2.1087,-0.2514) (2.0949,-0.3323) (2.0779,-0.4133) (2.0577,-0.4945) (2.0343,-0.5756) (2.0074,-0.6566) (1.9772,-0.7373) (1.9434,-0.8174) (1.9062,-0.8970) (1.8653,-0.9757) (1.8208,-1.0534) (1.7727,-1.1299) (1.7209,-1.2050) (1.6655,-1.2785) (1.6064,-1.3503) (1.5436,-1.4201) (1.4773,-1.4876) (1.4073,-1.5527) (1.3338,-1.6152) (1.2569,-1.6748) (1.1766,-1.7313) (1.0930,-1.7845) (1.0063,-1.8342) (0.9165,-1.8801) (0.8238,-1.9221) (0.7284,-1.9598) (0.6303,-1.9931) (0.5299,-2.0218) (0.4274,-2.0457) (0.3228,-2.0646) (0.2166,-2.0783) (0.1089,-2.0866) (0.0000,-2.0895)} -- cycle;
\draw[PaperMuted,fill=PaperMuted!15] (axis cs:0,0) circle[radius=1];
\node[font=\footnotesize,text=PaperAccent,anchor=west] at (axis cs:2.15,2.0) {$r_-(\alpha)$};
\node[font=\footnotesize,text=PaperInk] at (axis cs:0,0) {$R=1$};
\node[font=\footnotesize,text=PaperInk] at (axis cs:0,-1.75) {fluid};
\end{axis}
\end{tikzpicture}
\caption{Numerical continuation of the monotone continuum. Left: normalized crest speed against $\alpha$ along runs (a)--(c). Each branch starts at $(\alpha_c,1)$; branches (a) and (b) turn back in $\alpha$ before the computation stops. Right: the last computed wave of run (b) ($\gamma=0.9$, $\alpha=4.917$), with the circle $R=r_-(\alpha)$ on which the Stokes crest must form (dashed). The crest lies $0.014$ below that circle. All three runs satisfy the standing assumptions \eqref{eq_standing} and $c_0>\sinh^2\gamma$.}
\label{fig_numerics}
\end{figure}
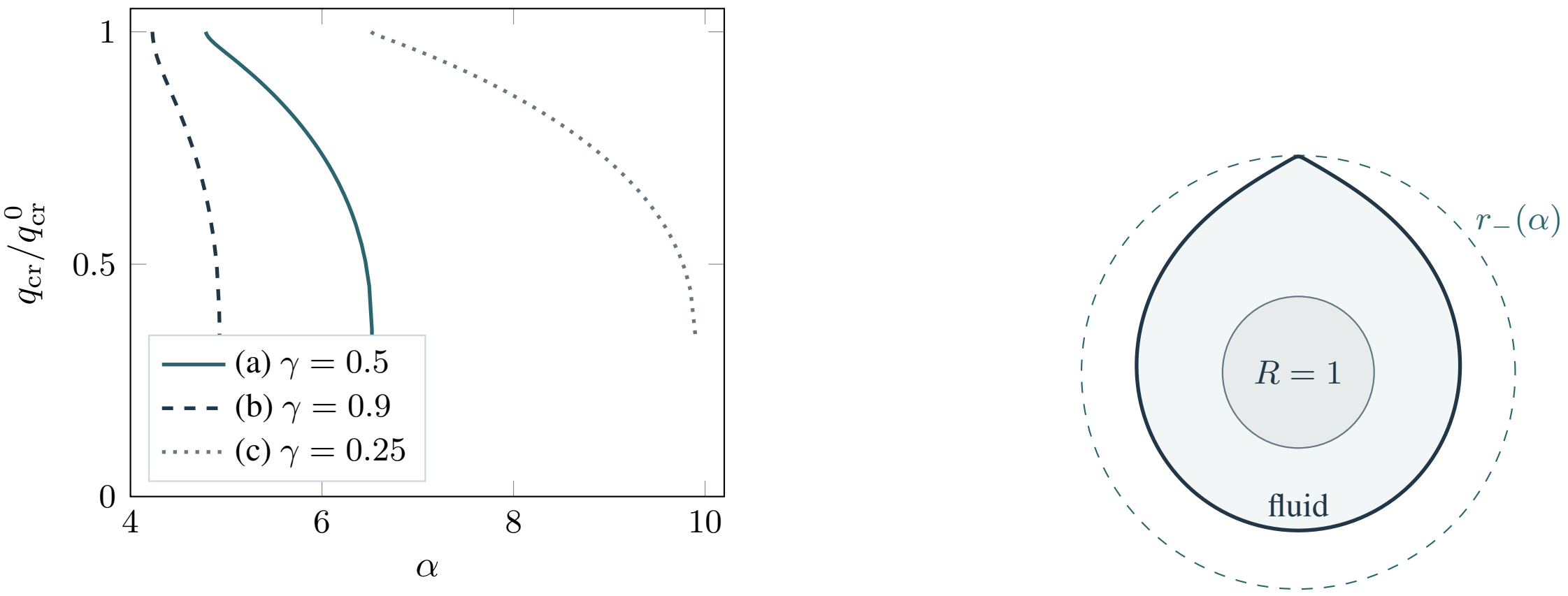

Along all three branches we observe the following.
\begin{enumerate}
\item The parameter stays in $[\alpha_c,\max\alpha]$. This interval lies strictly above $\max\{\alpha_+,\alpha_\circ\}$, as \autoref{thm_crossing} requires, and alternative (III) of \autoref{thm_global} is not seen. Each computed wave crosses the laminar circle; for instance, the last wave of the high-resolution run of \autoref{subsec_hp} has $S_t=2.091<R_0=2.460<S_c=2.859$. The values of $\max\alpha$ lie far below the bound $A(\gamma,p_0)$ of \autoref{thm_alpha_bound}, which shows how far from sharp that bound is.
\item In runs (a) and (b), $\alpha$ reaches a maximum and then decreases while the wave keeps steepening (\autoref{fig_numerics}, left). This is the annular counterpart of the speed maximum of near-extreme flat Stokes waves \cite{LonguetHigginsFox1977}. It is the reason why continuation in the arclength, not in $\alpha$, is necessary.
\item The waves stay strictly monotone from crest to trough. Their slope increases steadily and remains below $28^\circ$, well away from the vertical tangents of alternative (II).
\item The crest approaches the circle $r_-(\alpha)$ and the crest speed decreases, as predicted by \autoref{thm_main_stokes} and the localization of stagnation points.
\end{enumerate}
These observations support alternative (I) for these parameters. They do not reach the extreme wave itself. The singularity forms at the crest, and the uniform grid in the conformal variable cannot resolve it beyond $q_{\rm cr}/q_{\rm cr}^{\,0}\approx0.2$. Following the branch closer to the $120^\circ$ corner would require a crest-adapted variable, in the spirit of the high-precision computations for flat Stokes waves \cite{DLK2016,Lushnikov2016}, and that is left for future work.

\subsection{High-resolution continuation}\label{subsec_hp}
The dense collocation scheme above costs $O(N^3)$ per Newton step, which limits it to $N\le4096$. To go further we use a second scheme.
\begin{itemize}
\item The Bernoulli residual is evaluated on $2N+1$ points by fast cosine and sine transforms. The equations are its first $N+1$ cosine coefficients (a Galerkin projection), with the same unknowns $(b_1,\dots,b_N,\varrho,\alpha)$.
\item Newton steps are solved by GMRES with Jacobian-vector products by finite differences.
\item The preconditioner has two levels: an LU factorization of the dense Jacobian on the first $2048$ modes, together with $\varrho$, $\alpha$ and the arclength row, and a diagonal approximation on the remaining modes.
\item $N$ is doubled whenever the top eighth of the spectrum exceeds $10^{-9}$ of the largest coefficient.
\end{itemize}
Starting from the last wave of the refined run~(b), the branch was continued with $N=16384$ and then $N=32768$, that is, $65537$ grid points. At the last computed wave:
\begin{itemize}
\item $q_{\rm cr}/q_{\rm cr}^{\,0}=0.091$ and $r_--S_c=3.6\times10^{-3}$, with $\alpha=4.9126>\alpha_c=4.2276$;
\item the spectral tail is $1.6\times10^{-10}$, and the Bernoulli residual on a grid twice as fine is $1.9\times10^{-6}$, three orders of magnitude below that of the dense runs;
\item the wave is strictly monotone, and its maximal slope is $29.22^\circ$, attained at $\phi/\pi=7.6\times10^{-4}$, very close to the crest, where $h_\alpha'(S)<0$, as \autoref{prop_max_slope} requires;
\item the total-weight identity of \autoref{prop_total_weight} holds with $2.4622$ against $2.4625$, and the dispersion relation of \autoref{rem_nonlinear_dispersion} holds to ten digits.
\end{itemize}
Along this continuation the maximal slope increases monotonically: $28.3^\circ$, $28.7^\circ$, $29.0^\circ$, $29.2^\circ$ at $q_{\rm cr}/q_{\rm cr}^{\,0}=0.13$, $0.11$, $0.10$, $0.091$. This is consistent with convergence to a $120^\circ$ corner, whose limiting slope at the crest is $30^\circ$. The data do not decide whether the slope overshoots $30^\circ$ before the limit, as it does for near-extreme flat waves \cite{LonguetHigginsFox1977,StraussWheeler2016}. Below $q_{\rm cr}/q_{\rm cr}^{\,0}\approx0.09$ the GMRES iterations stagnate. Further progress needs a preconditioner adapted to the crest, or a reformulation of the singularity in the spirit of \cite{DLK2016,Lushnikov2016}. The scripts \texttt{hp\_core.py} and \texttt{hp\_continue.py} accompany the paper.

\subsection{A parameter scan}\label{subsec_scan}
To see whether the behaviour of runs (a)--(c) is typical, we repeated the dense continuation of \autoref{sec_numerics} with $N\le2048$ on a grid of $20$ parameter pairs. We took $\gamma\in\{0.2,0.4,0.6,0.8,0.95\}$ and $c_0=f\,e^{2\gamma}$ with $f\in\{0.2,0.5,0.8,0.95\}$. All $20$ pairs satisfy the standing assumptions \eqref{eq_standing}, since $f<1$ and $\gamma<1$, and also $c_0>\sinh^2\gamma$. Two changes were made to the stopping rules of the dense runs.
\begin{itemize}
\item If the corrector fails after the grid is doubled, the arclength step is halved instead of the run continuing on the coarse grid. Without this change, one run continued on an under-resolved grid and showed a spurious Gibbs oscillation next to the crest.
\item Monotonicity is tested relative to the amplitude: a run stops only if $S$ increases somewhere by more than $10^{-7}\,(S_c-S_t)$. For $\gamma=0.2$ the trough is long and almost flat. Increases of order $10^{-11}$ there had stopped the first version of the scan, but they are rounding effects.
\end{itemize}

In the table, $f=c_0e^{-2\gamma}$. The column $q_{\rm fold}$ gives $q_{\rm cr}/q_{\rm cr}^{\,0}$ at the maximum of $\alpha$ along the branch, and $q_{\rm end}$ gives it at the last computed wave. The last two columns give $\alpha-\alpha_c$ at that wave and the largest slope along the branch.
\begin{center}
\small
\begin{tabular}{@{}ccccccccc@{}}
\toprule
$\gamma$ & $f$ & $c_0$ & $\alpha_c$ & $\max\alpha-\alpha_c$ & $q_{\rm fold}$ & $q_{\rm end}$ & $\alpha_{\rm end}-\alpha_c$ & max slope\\
\midrule
$0.20$ & $0.20$ & $0.298$ & $2.703$ & $0.980$ & $0.329$ & $0.311^{*}$ & $0.980$ & $25.4^\circ$\\
$0.20$ & $0.50$ & $0.746$ & $4.926$ & $2.468$ & $0.330$ & $0.312^{*}$ & $2.467$ & $25.3^\circ$\\
$0.20$ & $0.80$ & $1.193$ & $7.149$ & $3.956$ & $0.330$ & $0.312^{*}$ & $3.955$ & $25.3^\circ$\\
$0.20$ & $0.95$ & $1.417$ & $8.261$ & $4.699$ & $0.330$ & $0.313^{*}$ & $4.698$ & $25.3^\circ$\\
$0.40$ & $0.20$ & $0.445$ & $2.575$ & $0.629$ & $0.332$ & $0.235$ & $0.624$ & $26.7^\circ$\\
$0.40$ & $0.50$ & $1.113$ & $4.201$ & $1.606$ & $0.339$ & $0.239$ & $1.595$ & $26.5^\circ$\\
$0.40$ & $0.80$ & $1.780$ & $5.826$ & $2.584$ & $0.338$ & $0.239$ & $2.565$ & $26.4^\circ$\\
$0.40$ & $0.95$ & $2.114$ & $6.639$ & $3.073$ & $0.340$ & $0.239$ & $3.051$ & $26.4^\circ$\\
$0.60$ & $0.20$ & $0.664$ & $2.865$ & $0.518$ & $0.333$ & $0.203$ & $0.512$ & $27.1^\circ$\\
$0.60$ & $0.50$ & $1.660$ & $4.430$ & $1.349$ & $0.344$ & $0.209$ & $1.332$ & $26.8^\circ$\\
$0.60$ & $0.80$ & $2.656$ & $5.994$ & $2.180$ & $0.348$ & $0.211$ & $2.154$ & $26.7^\circ$\\
$0.60$ & $0.95$ & $3.154$ & $6.776$ & $2.596$ & $0.348$ & $0.211$ & $2.564$ & $26.7^\circ$\\
$0.80$ & $0.20$ & $0.991$ & $3.341$ & $0.467$ & $0.322$ & $0.182$ & $0.459$ & $27.5^\circ$\\
$0.80$ & $0.50$ & $2.477$ & $5.014$ & $1.240$ & $0.344$ & $0.192$ & $1.219$ & $27.0^\circ$\\
$0.80$ & $0.80$ & $3.962$ & $6.687$ & $2.015$ & $0.353$ & $0.194$ & $1.981$ & $26.9^\circ$\\
$0.80$ & $0.95$ & $4.705$ & $7.524$ & $2.403$ & $0.354$ & $0.194$ & $2.363$ & $26.9^\circ$\\
$0.95$ & $0.20$ & $1.337$ & $3.802$ & $0.444$ & $0.324$ & $0.169$ & $0.435$ & $27.7^\circ$\\
$0.95$ & $0.50$ & $3.343$ & $5.626$ & $1.198$ & $0.354$ & $0.182$ & $1.174$ & $27.1^\circ$\\
$0.95$ & $0.80$ & $5.349$ & $7.450$ & $1.955$ & $0.356$ & $0.185$ & $1.916$ & $27.0^\circ$\\
$0.95$ & $0.95$ & $6.352$ & $8.363$ & $2.334$ & $0.360$ & $0.185$ & $2.287$ & $27.0^\circ$\\
\bottomrule
\end{tabular}
\end{center}
A star marks the four runs with $\gamma=0.2$. They stopped on the relative monotonicity test, but at that point the spectral tail was $1.2\times10^{-6}$, at the resolution limit. The increase of $S$ was $1.1\times10^{-7}$ of the amplitude, located in the flat trough region. We therefore regard these stops as resolution stops as well. The pair $\gamma=0.95$, $f=0.2$ is the only one in the band $\sinh^2\gamma<c_0\le\gamma^2e^{2\gamma}/(1+\gamma)^2$, since $(0.95/1.95)^2>0.2$; on its computed part $\alpha-\alpha_c$ stays positive, as for all other pairs. All other runs stopped because the spectral tail exceeded $2\times10^{-6}$ at $N=2048$. For $\gamma\ge0.4$, the value of $\max\alpha$ agrees to within $1.1\times10^{-6}$ with a second scan that used $N\le1024$. The fold is therefore resolved.

On all $20$ branches we observe the following.
\begin{enumerate}
\item \emph{Supercriticality.} The inequality $\alpha>\alpha_c$ holds at every computed wave. The smallest value of $\alpha_{\rm end}-\alpha_c$ is $0.435$, at $\gamma=0.95$, $f=0.2$. Together with runs (a)--(c), \autoref{conj_supercritical} is thus confirmed on the computed part of $23$ branches.
\item \emph{A universal fold.} Along every branch $\alpha$ has exactly one interior maximum, after which it decreases slowly while the wave keeps steepening. The fold occurs at $q_{\rm cr}/q_{\rm cr}^{\,0}\in[0.322,0.360]$ in every case, while $\max\alpha-\alpha_c$ varies by a factor of more than $10$. The fold position is therefore almost independent of the parameters. For near-extreme flat Stokes waves the phase speed also has its first maximum before the highest wave \cite{LonguetHigginsFox1977}. A consequence is that near the extreme wave the continuum is not a graph over $\alpha$, so any proof of (A) must control $\alpha$ along a curve that turns back.
\item \emph{Slope.} The maximal slope increases monotonically along every branch. Its largest value at the end of a run is $27.7^\circ<30^\circ$. Consistently with \autoref{subsec_hp}, these runs end well before the regime in which an overshoot of $30^\circ$ could be detected.
\end{enumerate}
The script \texttt{scan\_branches.py} reproduces the table.

\section{Open problems}\label{sec_open}

By \autoref{thm_global} and \autoref{thm_band}, when $c_0>\sinh^2\gamma$ the monotone continuum ends in an extreme wave unless its parameter approaches the degenerate value $\alpha_+$, and under \eqref{eq_floor_condition} it ends in an extreme wave (\autoref{cor_thin}). The Stokes existence statement for $c_0>\sinh^2\gamma$ therefore reduces to the single parameter floor (A) in the band $\sinh^2\gamma<c_0\le\gamma^2e^{2\gamma}/(1+\gamma)^2$. \autoref{tab_status} records the status of the main statements; all of them assume $0<\gamma<1$.

\begin{enumerate}
\item \emph{(A) Parameter floor:} $\alpha\ge\alpha_++\eta$. By \autoref{thm_crossing}, $\alpha>\alpha_\circ$ on the continuum in the lower region, so (A) is proved under \eqref{eq_floor_condition}, that is, when $\alpha_\circ>\alpha_+$. It remains open for $\sinh^2\gamma<c_0\le\gamma^2e^{2\gamma}/(1+\gamma)^2$, where $\alpha_\circ\le\alpha_+$. There the floor must keep the continuum away from $\alpha_+$, where the forbidden band collapses. A natural sufficient condition is the following annular analogue of supercriticality, the counterpart of the bound $F>1$ on the Froude number of solitary waves \cite{Wheeler2015,CWW2018}.

\begin{conjecture}[Global supercriticality]\label{conj_supercritical}
If $c_0>\sinh^2\gamma$, then $\alpha>\alpha_c$ at every non-laminar element of $\mathcal C^+_\delta$.
\end{conjecture}

\autoref{thm_supercritical} proves this near the bifurcation point.

\begin{remark}[A nonlinear dispersion relation]\label{rem_nonlinear_dispersion}
Pairing \eqref{eq_nekrasov} with $\sin\phi$ and using the symmetry of $\mathcal K_L$ gives, for every $(\Psi,S)\in\mathcal N$,
\[
\coth L\int_0^\pi\vartheta\sin\phi\,d\phi=\int_0^\pi\mathfrak g_\alpha\big(S(\phi)\big)\sin\vartheta\,\sin\phi\,d\phi .
\]
On the computed wave of run~(b) with $N=16384$ the two sides agree to $13$ digits. At the laminar flow this relation is the bifurcation condition $\coth\gamma=\mathfrak g_\alpha(R_0)$ of \autoref{rem_nekrasov_checks}. \autoref{conj_supercritical} asks whether this relation forces $\alpha>\alpha_c$ along the whole continuum. We have no global argument. A direct comparison with the laminar flow is obstructed by two facts. First, the weight $\mathfrak g_\alpha(S)$ is large near the crest. Second, the conformal modulus $L$ is not $\gamma$: in run~(b), $L=0.785$ while $\gamma=0.9$.
\end{remark}

Since $\alpha_c>\alpha_+$, the conjecture implies (A) with $\eta=\alpha_c-\alpha_+$, and then \autoref{thm_band} gives the extreme wave. The band is thin: for instance, at $\gamma=0.3$ and $c_0$ in the band, $\alpha_c-\alpha_+<0.008$. Simple integral constraints do not suffice. We checked numerically that the crest--trough comparison, the energy sandwich and the two-sided angular flow-force estimates of \autoref{lem_trough_bound} remain compatible as $\alpha\downarrow\alpha_+$ in the band. On all $23$ computed branches of \autoref{sec_numerics} (runs (a)--(c) and the scan of \autoref{subsec_scan}), $\alpha-\alpha_c$ stays positive, including after the turning point in $\alpha$.
\item \emph{(B) Slope bound.} (B) is proved whenever (A) holds, for every $c_0>\sinh^2\gamma$ (\autoref{thm_band}); the bound is $\vartheta<60^\circ$ when moreover $f\ge f_1(\gamma)$ (\autoref{rem_C1_region}). As $\alpha\downarrow\alpha_+$ the constant $\kappa_\alpha$ tends to $\frac12$, but by \autoref{rem_kappa_crest} the slope bound stays uniform as long as the crests stay away from stagnation. If the crests approach stagnation as $\alpha\to\alpha_\pm$, \autoref{thm_degenerate_limit} and \autoref{cor_IV} produce a right-angle extreme wave. What is missing is control of the continuum when it approaches or enters the window $\Delta\le0$.
\item \emph{Degenerate parameters $\alpha=\alpha_\pm$.} The crest shape is settled by \autoref{thm_degenerate}: every stagnation point is a symmetric right-angle corner. Moreover, by \autoref{thm_degenerate_limit} and \autoref{cor_IV}, if alternative (III) occurs with crests approaching stagnation, it produces such a right-angle extreme wave. What remains open is whether the continuum actually reaches it.
\item \emph{Full density for general variational solutions.} Every stagnation point of an overhanging classical wave is a corner (\autoref{cor_overhanging}). The same holds for variational solutions with regular interfaces and no air bubbles, or with finitely many air components and an $\mathcal H^1$-null singular set (\autoref{thm_variational_crest}(c)). Full density remains possible, a priori, only when infinitely many air components accumulate at the crest and the interface is not regular (\autoref{rem_full_density}).
\item \emph{The case $c_0\le\sinh^2\gamma$.} For $c_0<\sinh^2\gamma$, \autoref{thm_window} shows that inside the window $\Delta<0$ the continuum cannot stagnate. It must overturn or reach $\alpha_\pm$; secondary bifurcation is excluded by \autoref{thm_crossing}. The critical case $c_0=\sinh^2\gamma$ is reduced to the same alternatives as $c_0>\sinh^2\gamma$ (\autoref{thm_critical}). A global theory beyond the window remains open.
\end{enumerate}

\begin{remark}[Scope of the methods]\label{rem_scope}
The tools of this paper use the specific model to different degrees. The weighted subharmonicity of \autoref{lem_subharmonic}, the radial phase monotonicity and the blow-up analysis use only that $\Psi$ is harmonic and that the Bernoulli function is a smooth function $Q(|x|)$ of the radius with a simple zero at the stagnation radius; the rotational Jacobi field and the angular flow force use only rotation invariance, and the flow force holds with $H'=\sigma Q(\sigma)$ for any such $Q$. We therefore expect the crest theory of \autoref{sec_corner} and the monotonicity of \autoref{sec_continuum} to extend to general radial Bernoulli functions, but we have not checked the details. In contrast, the crossing invariant uses the parameter path \eqref{eq_parameter_path}, on which the laminar radius solves the circular relation for every $\alpha$, and the structural inequality, the monotonicity of the weight and \autoref{cor_thin} use that $Q_\alpha$ is quadratic. Non-constant vorticity would destroy the harmonicity of $\Psi$, on which the subharmonicity, the Nekrasov equation and the conformal parametrization all rest.
\end{remark}

\section*{Acknowledgements}
The author was supported by the Natural Science Foundation of Sichuan Province (Grant No.~2025ZNSFSC0072).
The author gratefully acknowledges the use of Claude (Anthropic) for auxiliary calculations and proof derivations; the core proof idea and strategy were provided by the author, who takes full responsibility for the final proof.
\section*{Conflict of interest}
The author declares that there is no conflict of interest.

\section*{Data availability}
No data were generated other than by the scripts described in the paper.\sloppy\ The scripts \texttt{annular\_continuation\_core.py}, \texttt{run\_branch.py}, \texttt{run\_branch\_fine.py} (the refined run with $N\le4096$), \texttt{verify\_identities.py}, \texttt{verify\_tables.py}, \texttt{verify\_A.py}, \texttt{hp\_core.py}, \texttt{hp\_continue.py}, \texttt{nekrasov\_check.py}, \texttt{scan\_branches.py}, \texttt{cone\_check.py} and \texttt{region\_check.py} reproduce every numerical value reported here; they use only Python with NumPy and SciPy. The script \texttt{convergence\_test.py} reproduces the convergence table of \autoref{subsec_convergence}. The scripts \texttt{cone\_check.py} and \texttt{region\_check.py} are numerical illustrations; no proof uses them, and in particular the proof of \autoref{cor_thin} uses no interval arithmetic.

\end{document}